\pdfoutput=1
\documentclass[pdf, 10pt]{article}
\usepackage[square,authoryear]{natbib}
\usepackage{marsden_article}
\usepackage{amscd}

\input xy
\xyoption{all}

\newtheorem{Theorem}{Theorem}[section]

\theoremstyle{remark}
\newtheorem{Remark}{Remark}[section]

\theoremstyle{definition}

\def\R{\mathbb{R}}
\def\E{\mathbb{E}}
\def\T{\mathbb{T}}
\def\P{\mathbb{P}}

\def\Cov{\operatorname{Cov}}
\def\Var{\operatorname{Var}}
\def\eref#1{(\ref{#1})}

\begin{document}
\title{Ballistic Transport at Uniform Temperature}

\author{Nawaf Bou-Rabee$^*$  \and Houman Owhadi\thanks{Applied \& Computational Mathematics (ACM), Caltech, Pasadena, CA 91125 ({\tt nawaf@acm.caltech.edu}, \tt{owhadi@acm.caltech.edu}).} }

\maketitle

\begin{abstract}
A paradigm for isothermal, mechanical rectification of stochastic
fluctuations is introduced in this paper.  The central idea is to
transform energy injected by random perturbations into
rigid-body rotational kinetic energy. The prototype considered in
this paper is a mechanical system consisting of a set of rigid
bodies in interaction through magnetic fields. The system is
stochastically forced by white noise and dissipative through
mechanical friction. The Gibbs-Boltzmann distribution at a specific
temperature defines the unique invariant measure under the 
flow of this stochastic process and allows us to define 
``the temperature'' of the system.  This measure is also ergodic 
and strongly mixing.  Although the system does not exhibit global 
directed motion, it  is shown that global ballistic motion is possible (the 
mean-squared displacement grows like $t^2$).  More precisely, although work cannot be
extracted from thermal energy by the second law of thermodynamics,
it is shown that ballistic transport from thermal energy is possible.
In particular, the dynamics is characterized by a meta-stable state 
in which the system exhibits directed motion over random time scales.  
This phenomenon is caused by interaction of three attributes of the
system: a non flat (yet bounded) potential energy landscape, a rigid body
effect (coupling translational momentum and angular momentum through
friction) and the degeneracy of the noise/friction tensor on the momentums
(the fact that noise is not applied to all degrees of freedom).
\end{abstract}

\section{Introduction}

Mechanical rectification was introduced to describe devices that
convert small-amplitude mechanical vibrations into directed rotary
or rectilinear mechanical motion \citep{Br1989}.  The purpose of
this paper is to extend this concept to mechanical systems in
isothermal environments subjected to stochastic fluctuations.
Although the second law of thermodynamics prevents directed
mechanical motion from thermal fluctuations in an isothermal
environment, violations of the second law of thermodynamics on
certain time-scales in microscopic systems have been theoretically
predicted \citep{Gallavotti1995a,MR1705590} and experimentally
observed \citep{Ciliberto1998, WaSeMiSeEv2002}.

In this paper we exhibit a system characterized by ballistic non-directed 
motion at uniform temperature.  Set $\mu$ to be the unique
invariant measure of the system and $x(t)$ the displacement of the system 
at time $t$ then $\mu$ a.s.~$\lim_{t \to \infty} (x(t)-x(0))/t \to 0$  
but $\E_{\mu}[\{ x(t)-\E(x(t)) \}^2]\sim t^2$.  Moreover the system is characterized by ``dynamic'' meta-stable states in which it experiences ballistic directed motion over random time scales (flights) in addition to ``static'' meta-stable states where the system is 
``captured'' in potential wells.

The mechanism introduced in this paper is different from the one
associated to the Gallavotti-Cohen fluctuation theorem \citep{Gallavotti1995a}. 
It is based on the fact that one can obtain anomalous diffusion by introducing
degenerate noise and friction in the momentums.  This anomalous diffusion
manifests itself in the just mentioned flights over time-scales which are random and the probability distributions of the flight-durations have a heavy tail.  
Moreover we expect that it will be of relevance to Brownian motors \citep{Ast97,Pet02} 
since it shows how thermal noise can be used to get ballistic
transport, in other words how to obtain fast and efficient
intracellular transport without using chemical energy.  Although
global ballistic transport is achieved it is not directed and hence
work is not produced from thermal fluctuations (this mechanism is
not in violation of the second law of thermodynamics). However this
mechanism would be sufficient for intracellular transport
of a large fraction of the material being carried to target areas with
ballistic speed without draining cellular energy reserves.

\begin{figure}[htbp]
\begin{center}
\includegraphics[scale=0.35,angle=0]{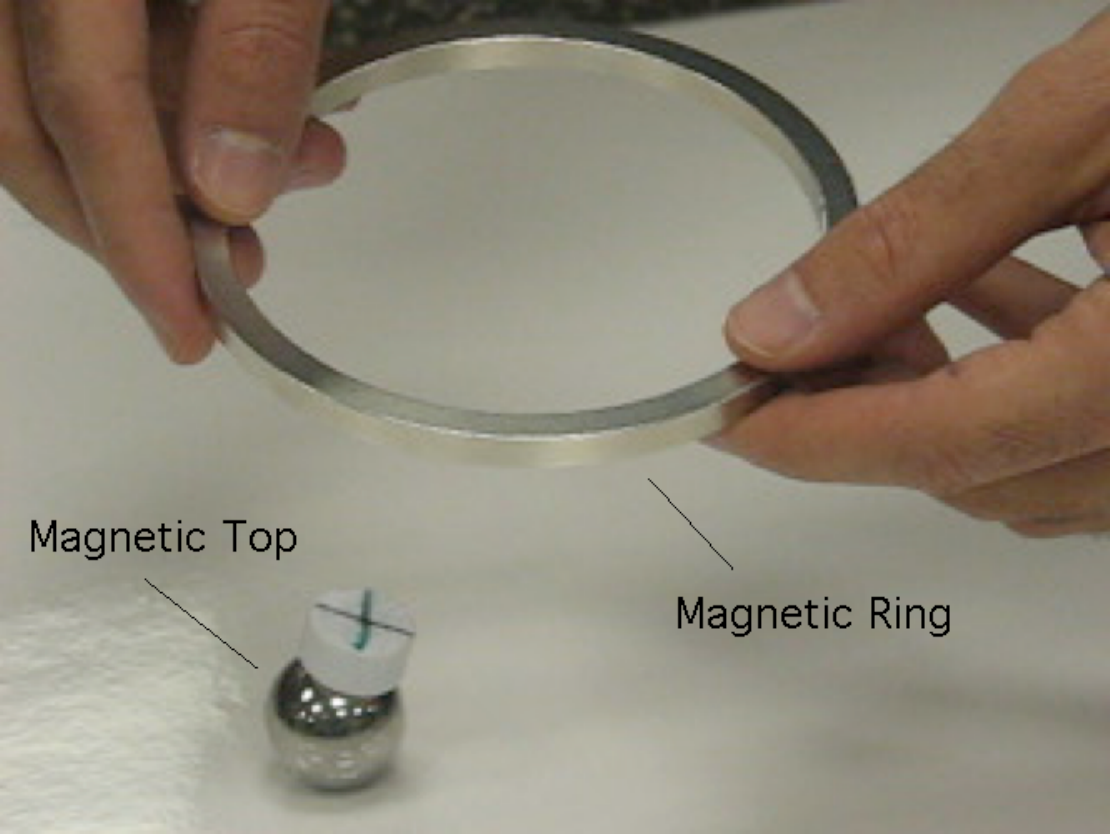}
\caption{\footnotesize  {\bf Picture of mechanical system.}   The
physical system consists of a magnetized axisymmetric top on a flat
surface and a magnetized ring as shown above.   The dipole moments
in the ring are oriented radially. } \label{fig:moviesnapshot}
\end{center}
\end{figure}

The prototype analyzed in this paper consists of a magnetized top
(ball) on a surface interacting with a suspended magnetized ring as
shown in Fig.~\ref{fig:moviesnapshot}.   The following dynamics is
observed: when one lowers the magnetic ring to within a certain range of 
heights and then tilts the ring, one observes the top transition from a 
state of no spin about its axis of symmetry to a state of nonzero spin.   
The reader is referred to the following urls for movies and simulations of the phenomenon: 
\begin{quote}
\url{http://www.acm.caltech.edu/~nawaf/BallisticTransport/} 

\url{http://www.acm.caltech.edu/~owhadi/BallisticTransport/}
\end{quote}
This phenomenon seems counterintuitive and non-Hamiltonian, since
one would expect the angular momentum of the top to be conserved
(since without friction the Hamiltonian of the system is invariant
under an $S^1$ rotation of the ball). The origin of this mechanical
device can be traced back to the work of David Hamel on magnetic
motors in the unofficial sub-scene  of physics \citep{MaSi1995}.  The 
mechanism presented in this paper when the device is not at uniform
temperature could in principle be used as a method of extracting energy
from macroscopic fluctuations.

\begin{figure}[htbp]
\begin{center}
\includegraphics[scale=0.2,angle=0]{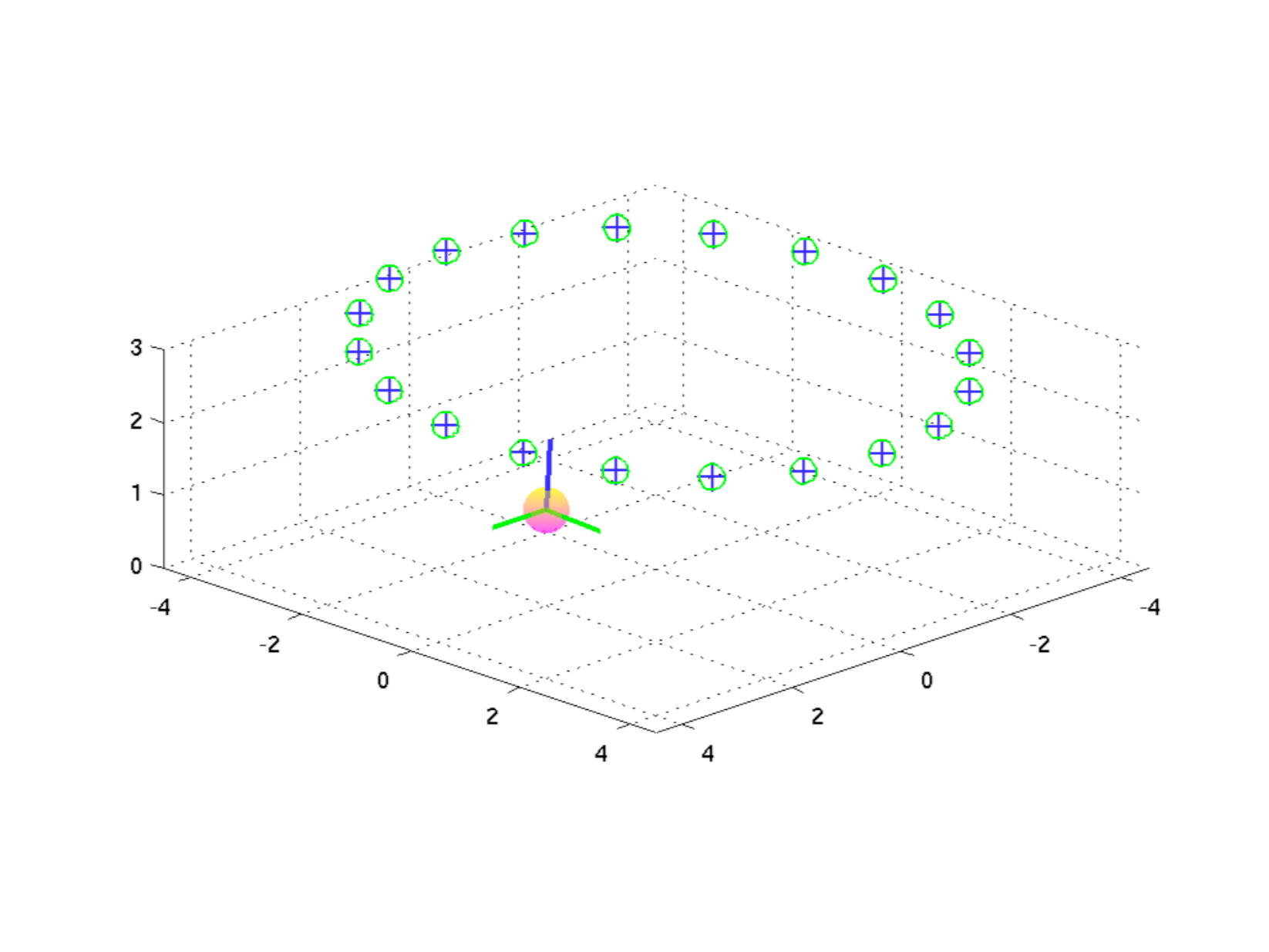}
\includegraphics[scale=0.2,angle=0]{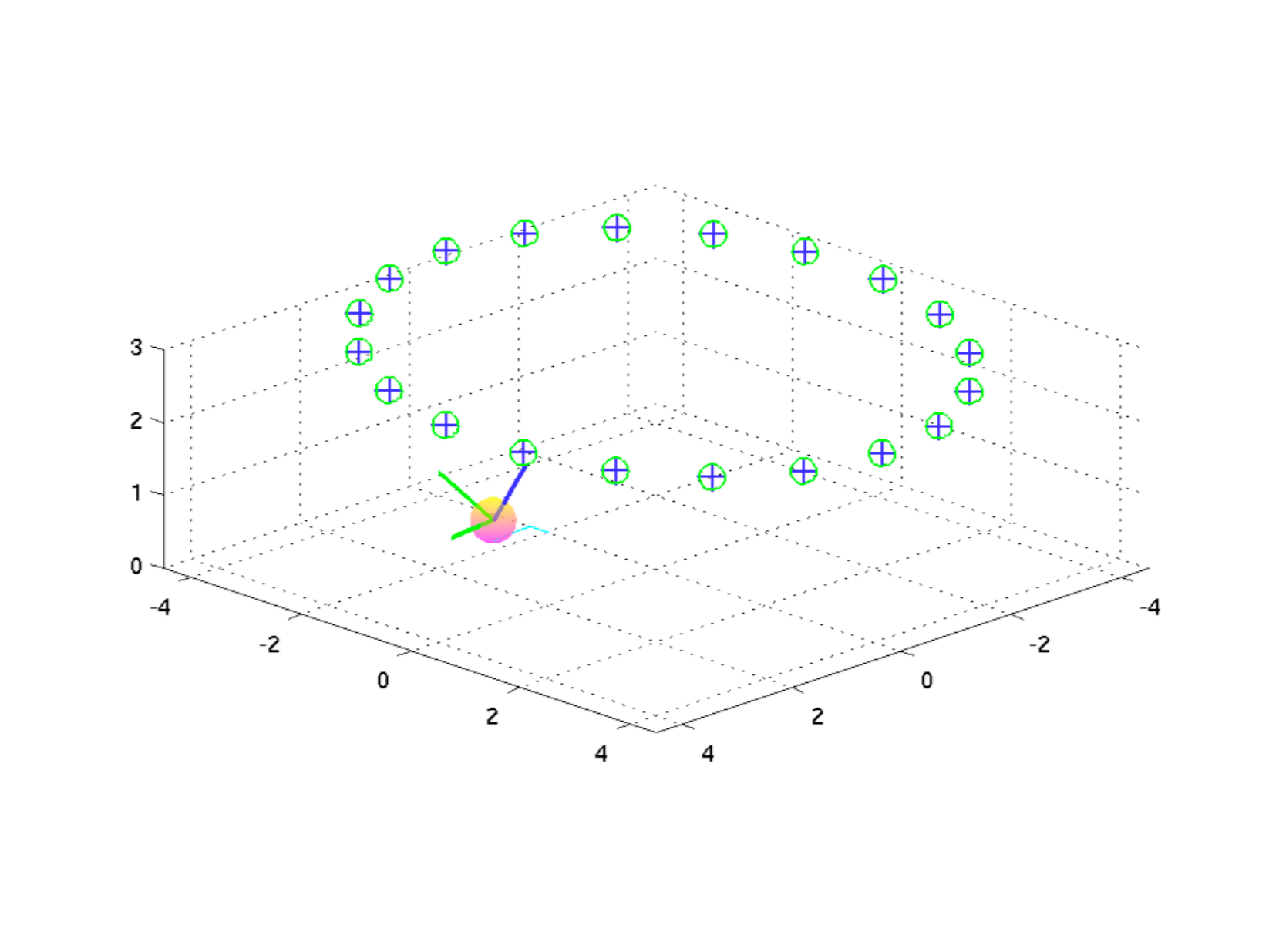}
\includegraphics[scale=0.2,angle=0]{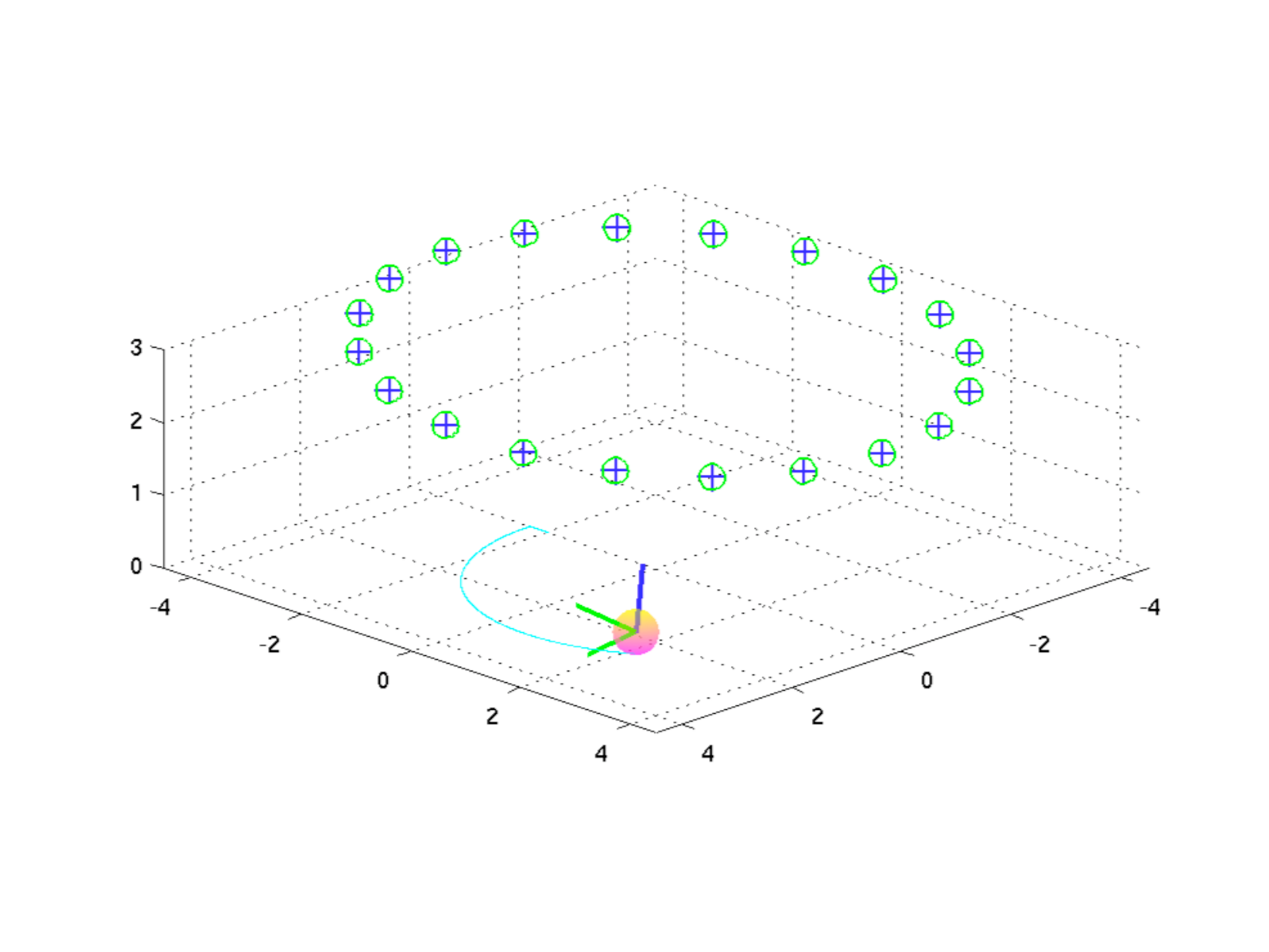}
\includegraphics[scale=0.2,angle=0]{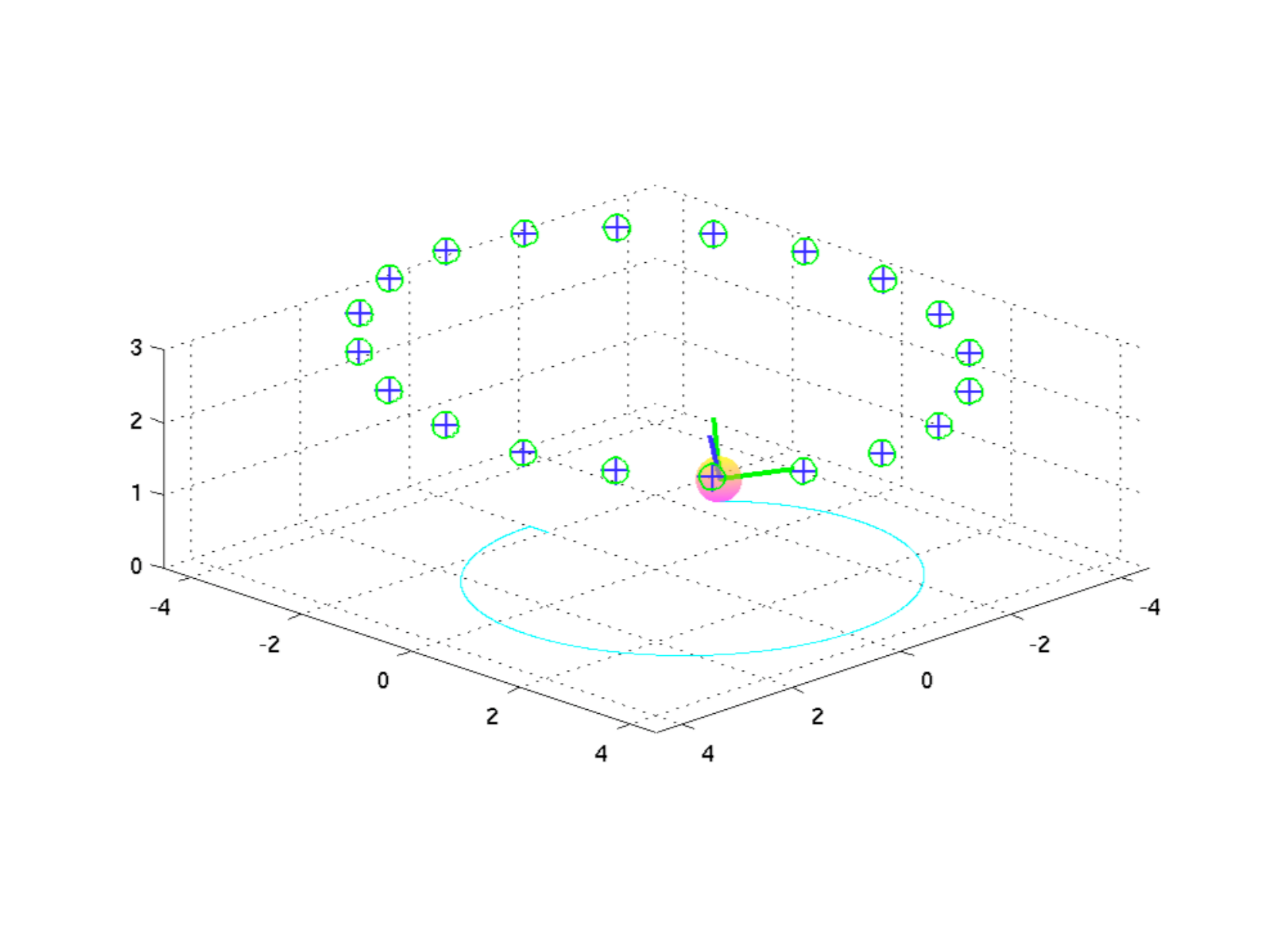}
\includegraphics[scale=0.2,angle=0]{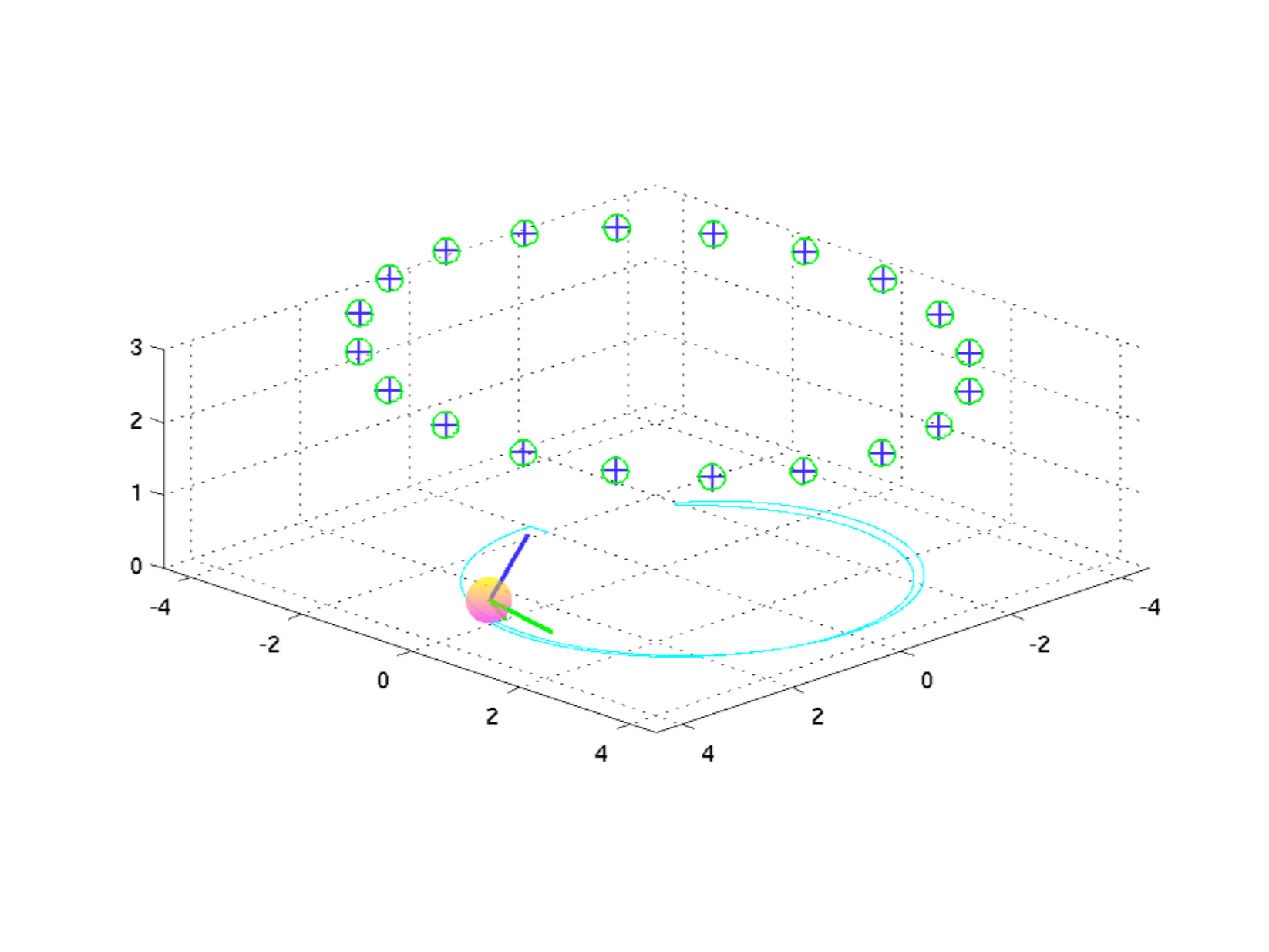}
\includegraphics[scale=0.2,angle=0]{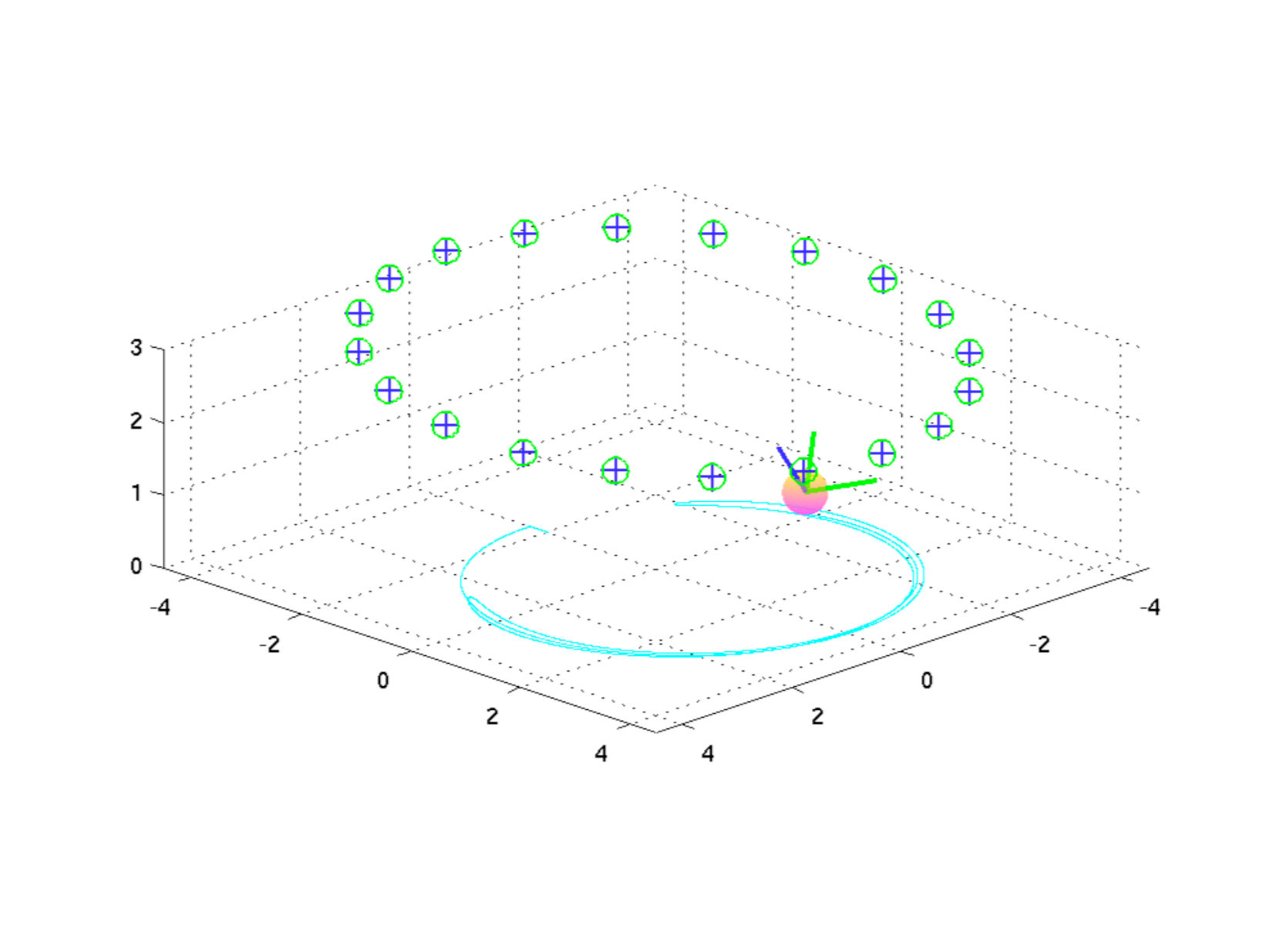}
\includegraphics[scale=0.2,angle=0]{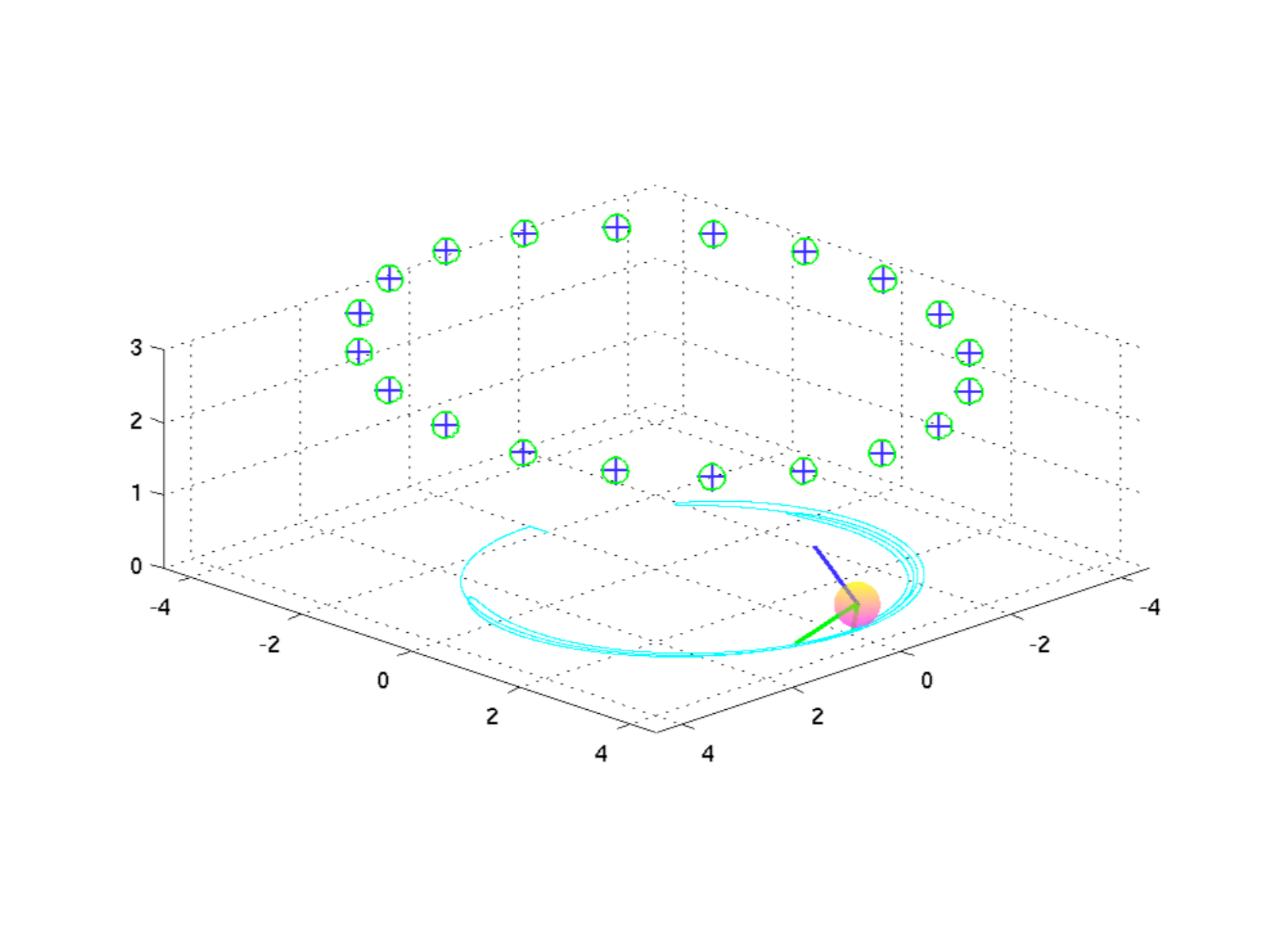}
\includegraphics[scale=0.2,angle=0]{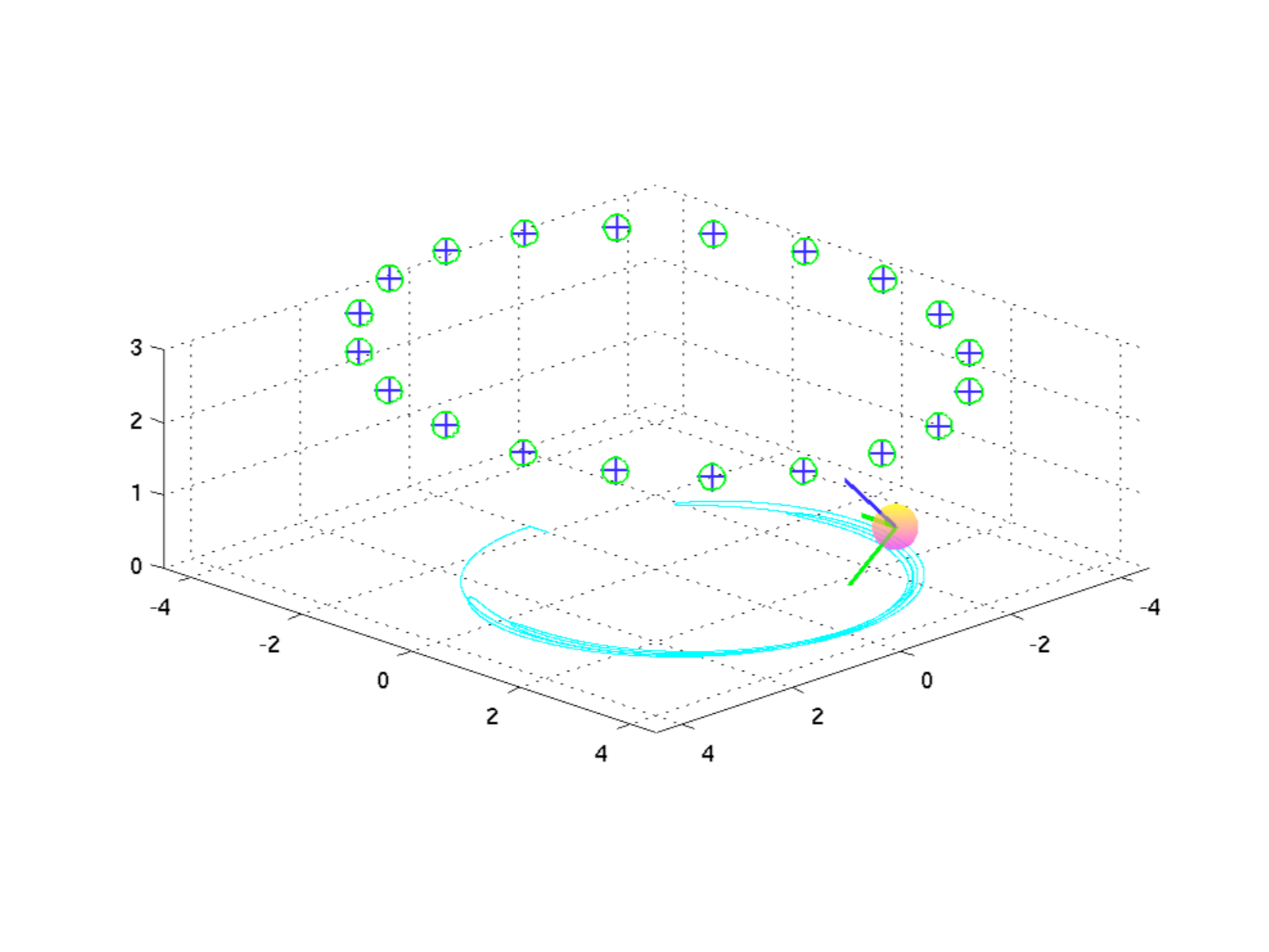}
\caption{\footnotesize  {\bf Snapshots of Nonconservative Tilted
Simulation.} The above are snapshots of the simulation described in
Fig.~\ref{fig:tilt}.  The top is initially set with its axis of
symmetry pointing nearly vertical.  The axis of symmetry then aligns
with the ambient magnetic field.  However, the state is not a
minimum of $V_e$, and hence, the state is unstable.  The top moves
towards a state that minimizes its magnetic potential energy, and
acquires spin in this process.
   }
\label{fig:snapshots}
\end{center}
\end{figure}

\section{Preview of the Paper}

In  section \ref{HamelsDevice}  Hamel's device is analyzed through
an idealized model based on magnetostatics and the spinning of the
top is caused by the introduction of surface frictional forces.
Simulations (figure \ref{fig:snapshots}) are done using variational
integrators and concur with experimental observations.

In section \ref{kkshshkjhe} the magnetic ring is allowed to be
dynamic, and a fixed outer ring of a finite number of magnetic
dipoles is introduced to stabilize it (see figure
\ref{fig:magmotor}). The inner magnetic ring is excited through
white noise applied as a torque. The steel ball and inner ring are
coupled through a magnetostatic potential and the motion of the
steel ball is dissipative through slip friction.

The mathematical description of the fluctuation driven motor is
obtained by generalizing Langevin processes from a system
of particles on a linear configuration space to rigidified particles
whose configuration space is the Lie group $\operatorname{SE}(3)$.
The question of random perturbations of a rigid body was treated by
previous investigators who added perturbations to the Lie-Poisson equations
without potential or dissipative torques \citep{Li1997,LiWa2005}.   Hence, 
they do not consider generalizing Langevin processes to the Lie-Poisson setting.

\begin{figure}[htbp]
\begin{center}
\includegraphics[scale=0.2,angle=0]{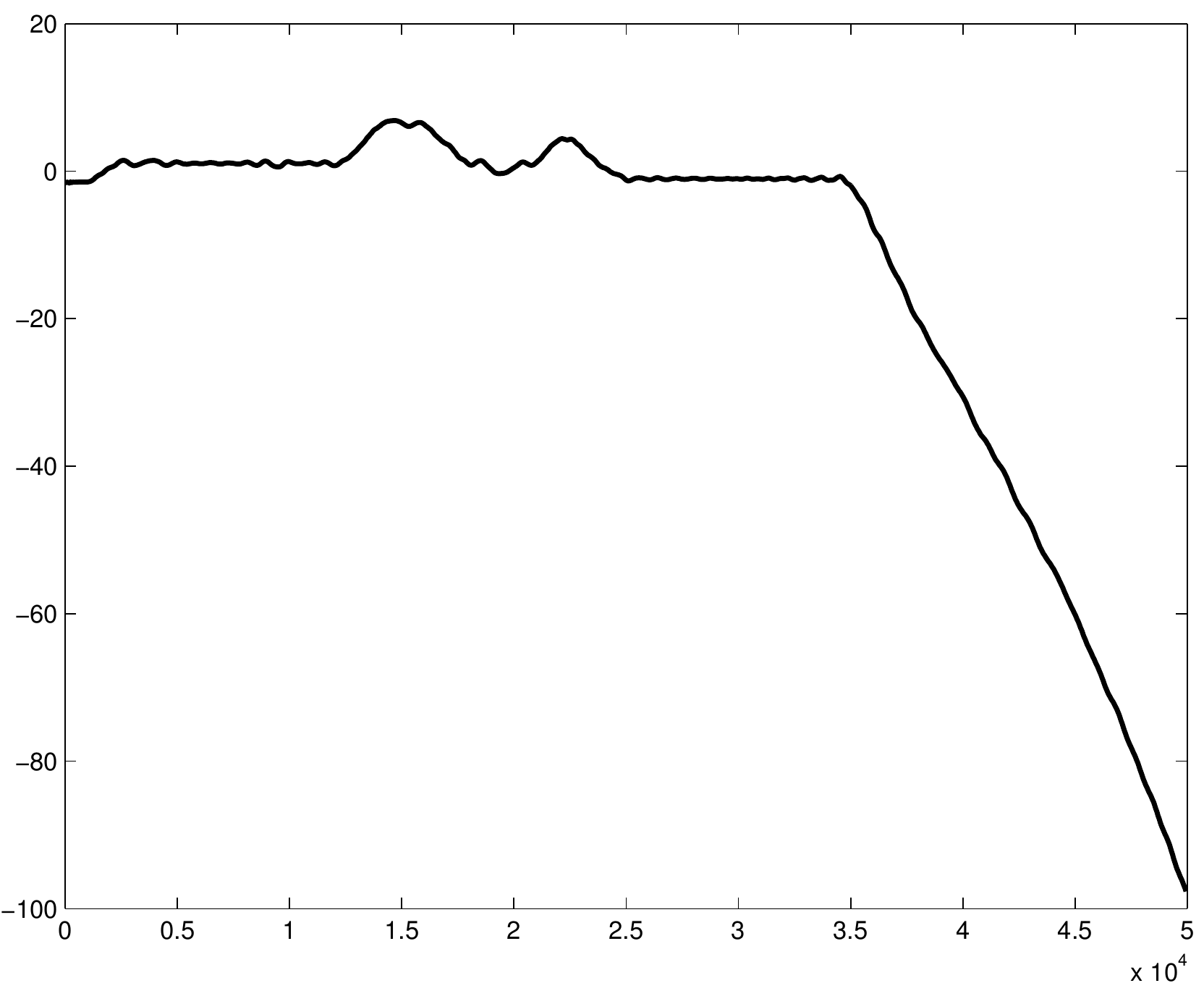}
\includegraphics[scale=0.2,angle=0]{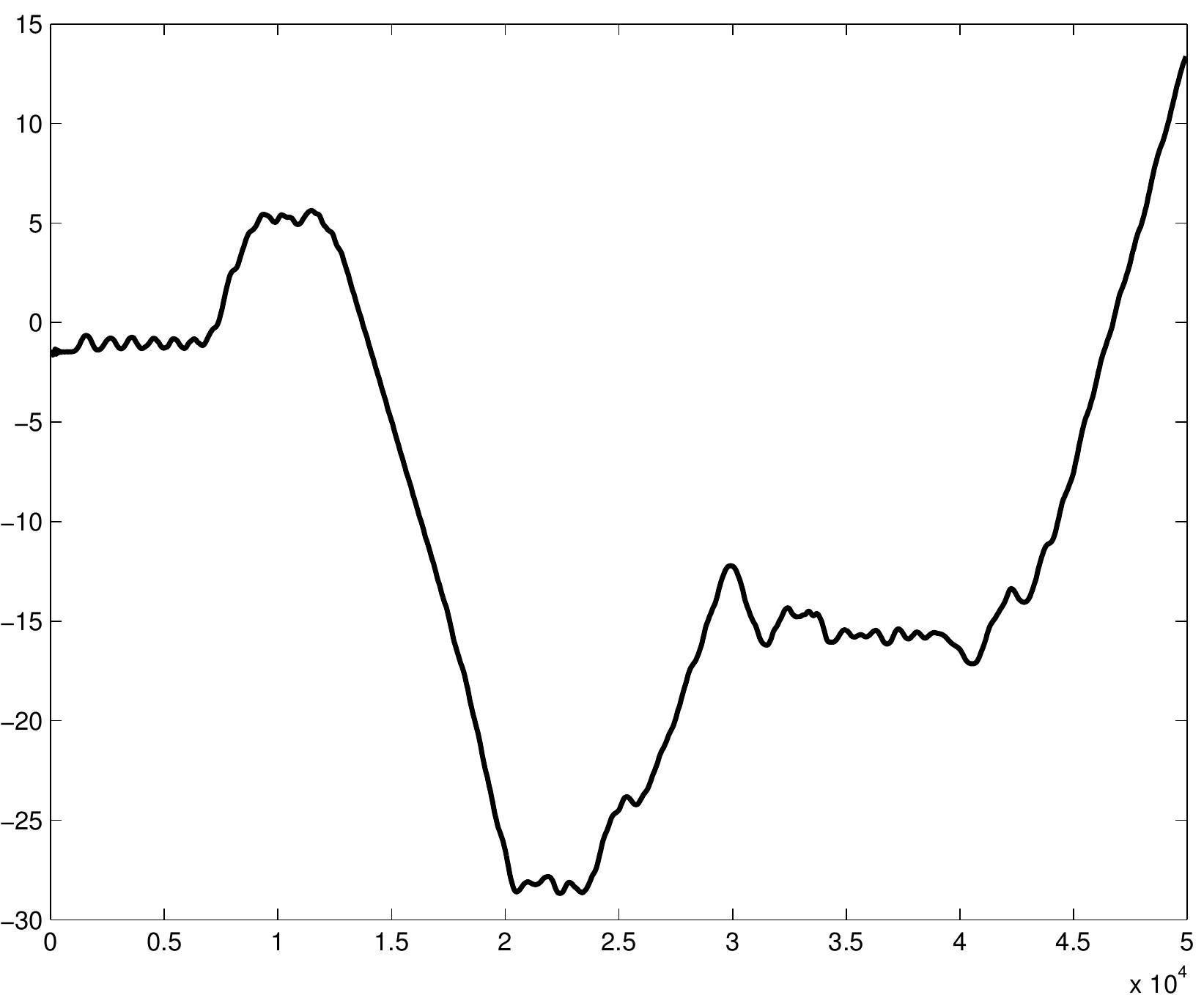}
\includegraphics[scale=0.2,angle=0]{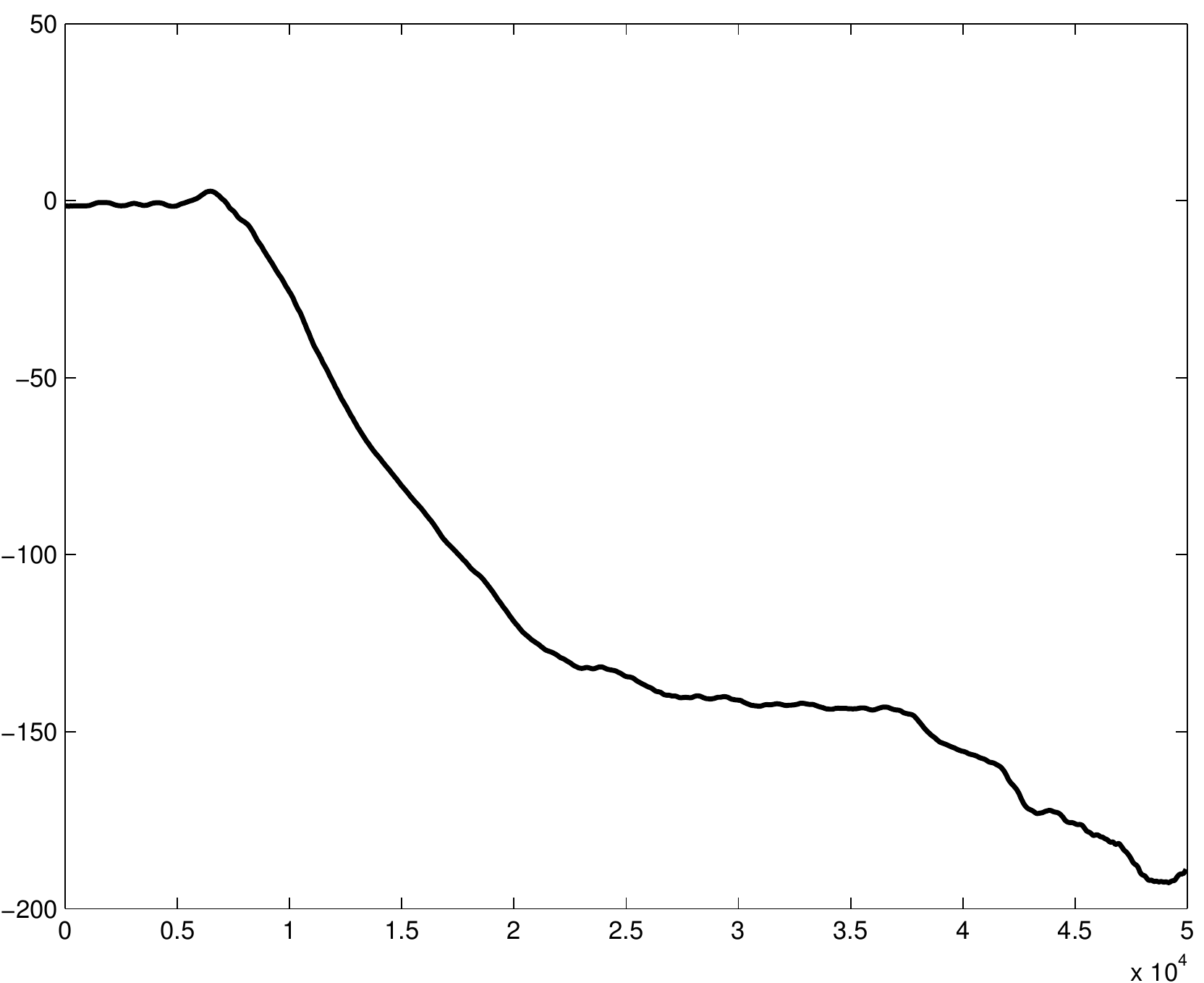} \\
\hbox{\hspace{0.8in}   (a)  $\alpha=0.0002$ \hspace{0.3in} (b)
$\alpha=0.0002$ \hspace{0.3in} (c)  $\alpha=0.00025$  \hspace{0.3in}
} \caption{\footnotesize  {\bf Angular-position of magnetic ball
(non-uniform temperature).}  The angular component of the center of
mass is plotted for three different realizations. The plots show the
magnetic ball  transition between two meta-stable states:
noise-driven and inertia-driven motion. }
\label{fig:nctmagmotortheta}
\end{center}
\end{figure}

We refer to figure \ref{fig:nctmagmotortheta} for a plot of the
angular position of the ball versus time for different values of
noise amplitude $\alpha$. The simulation of the system exhibits two
distinct kinds of metastable states. In the first kind the ball is
``stuck'' in a magnetic potential well whose depth depends on the
position of the inner ring, and moves into another potential well
when the energy barrier between the two wells is close to minimum
(a phenomenon known as stochastic resonance) or transitions to a
metastable state of the second kind.  In the second kind the ball
spins in circles clockwise or counter-clockwise in a directed way
for a random amount of time (that numerically and heuristically
observe to be exponential in law) until it gets stuck in a potential
well.

Since the system associated to figure \ref{fig:nctmagmotortheta} is
described by a generalized Langevin process it is possible
to introduce a notion of temperature defined as $\frac{(\text{noise
amplitude})^2}{\text{friction constant}}$.  However, since noise but
no friction is applied at the level of the inner ring and friction
but no noise is applied at the level of the ball, the system is not
at uniform temperature, i.e., the inner ring is at infinite
temperature and the ball is at zero temperature.  Nevertheless, one could
in principle use such a mechanism as a way to extract energy from 
macroscopic fluctuations.

\begin{figure}[htbp]
\begin{center}
\includegraphics[scale=0.2,angle=0]{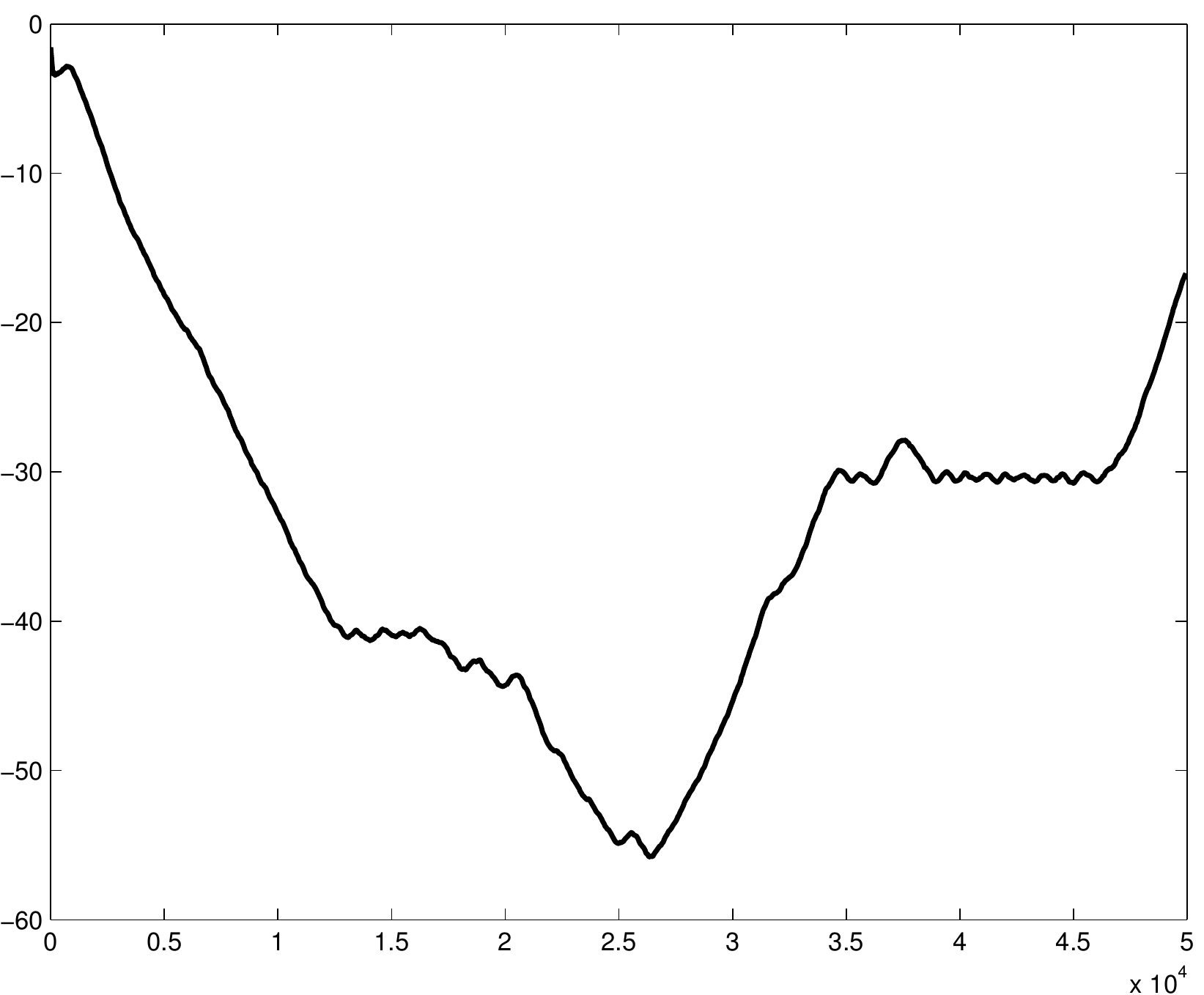}
\includegraphics[scale=0.2,angle=0]{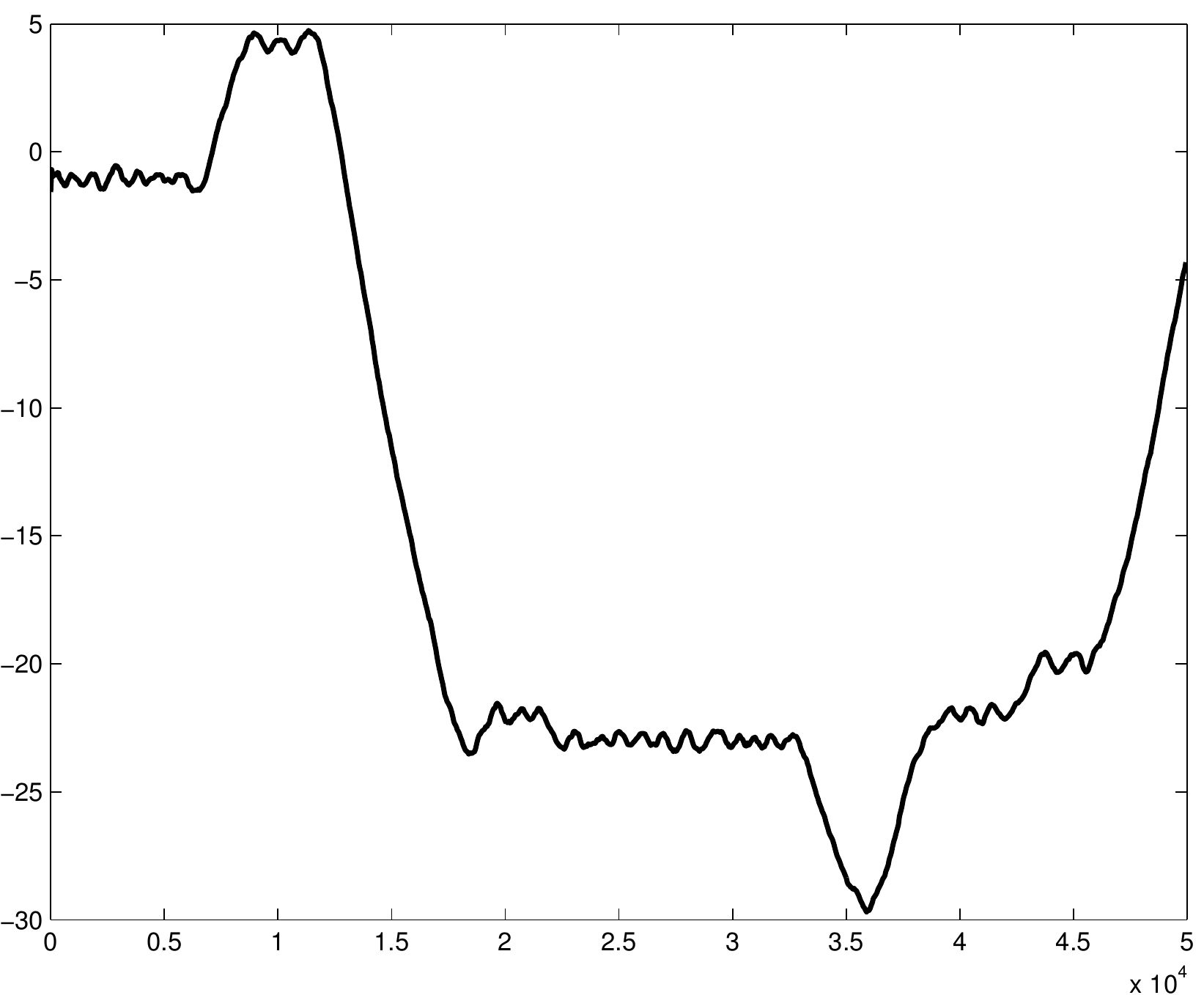}
\includegraphics[scale=0.2,angle=0]{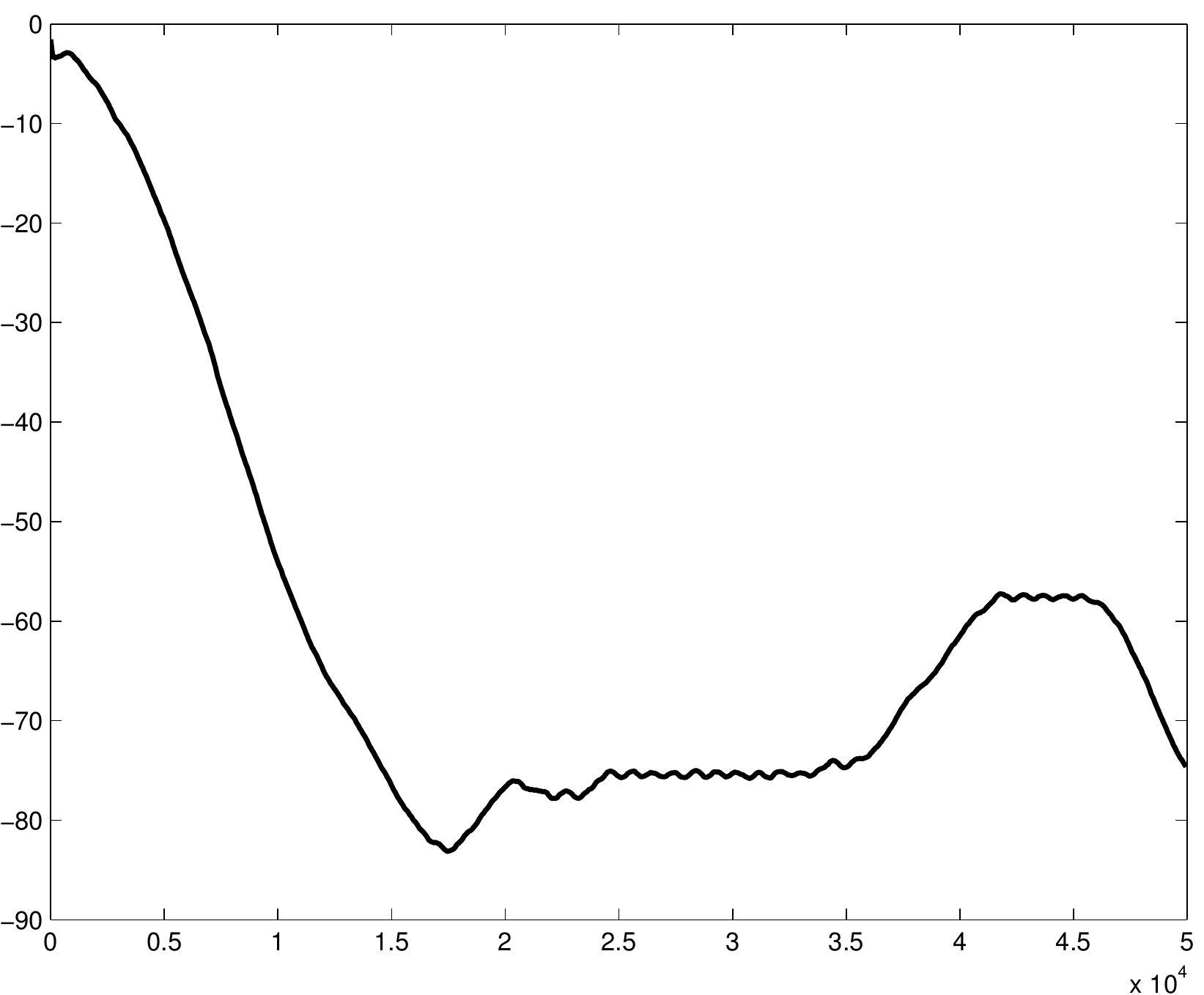}
\includegraphics[scale=0.2,angle=0]{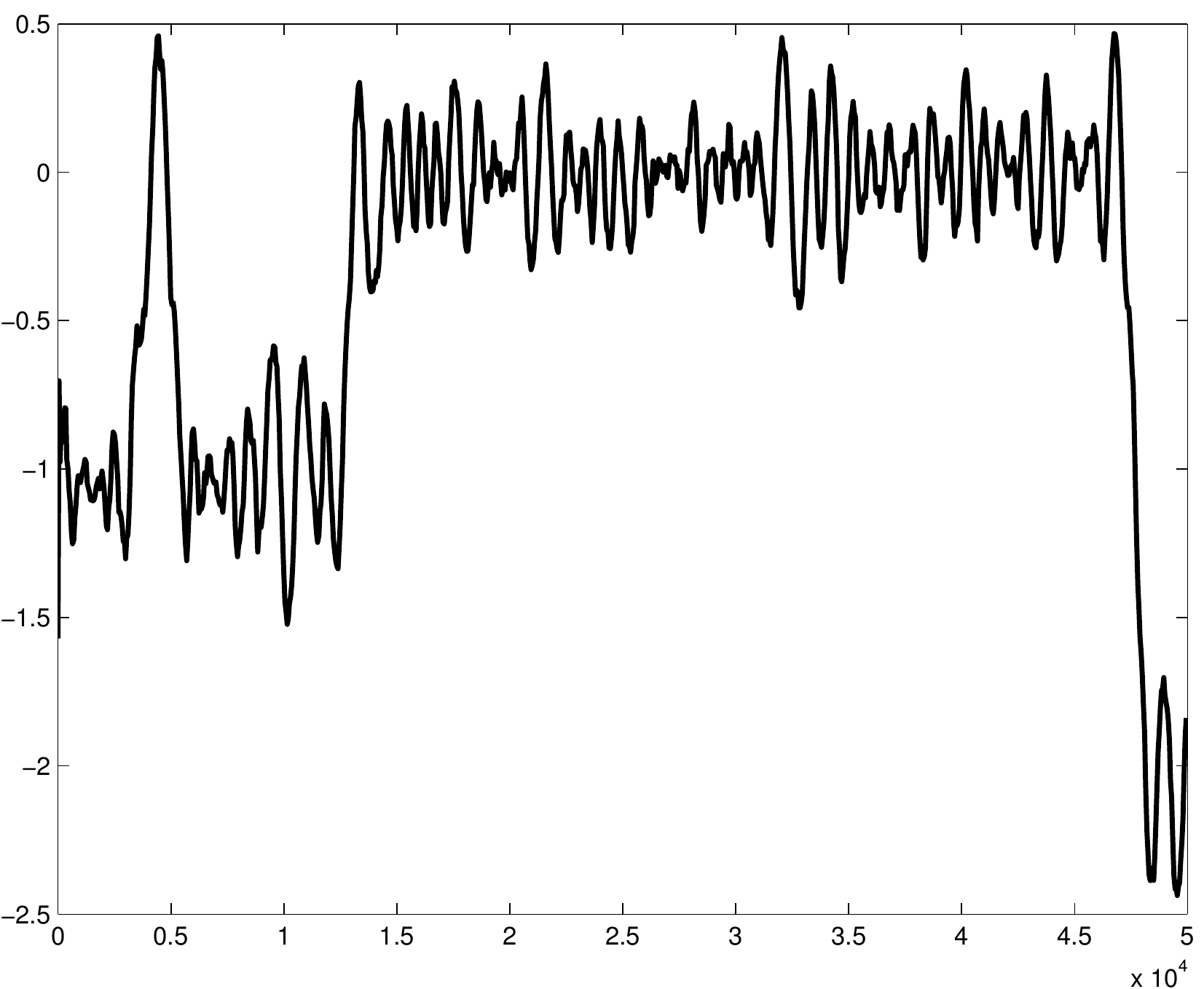} \\
\hbox{\hspace{0.15in}   (a)  $\alpha=0.001$ \hspace{0.3in} (b)
$\alpha=0.001$ \hspace{0.3in} (c)  $\alpha=0.00075$ \hspace{0.3in}
(d) $\alpha=0.00075$  \hspace{0.3in}    } \caption{\footnotesize
{\bf Angular-position of magnetic ball (uniform temperature).} Four
different realizations of the angular component of the center of
mass are plotted.    The plots for $\alpha=0.0002$ show the magnetic
ball transition between two metastable states: noise-driven and
inertia-driven  motion.  If $\alpha$ is increased to
$\alpha=0.00025$, realization (c) shows a similar transition.
Realization (d), however, does not show such a transition, i.e., the
ball never transitions to a metastable state of the second kind.
 }
\label{fig:ctmagmotortheta}
\end{center}
\end{figure}

In a second step frictional torque is introduced to the inner ring
and thermal torque (white noise) to the ball, so that the generator
of the process is characterized by a unique Gibbs-Boltzmann
invariant distribution.  Throughout this paper we interpret this
property as placing the mechanical system at uniform temperature.
We refer to figure \ref{fig:ctmagmotortheta} for a plot of the
angular position of the ball versus time for different values of
noise amplitude $\alpha$ (and temperature). The system is still
characterized by ballistic directed motion meta-stable states.

\begin{figure}[htbp]
\begin{center}
\includegraphics[scale=0.5,angle=0]{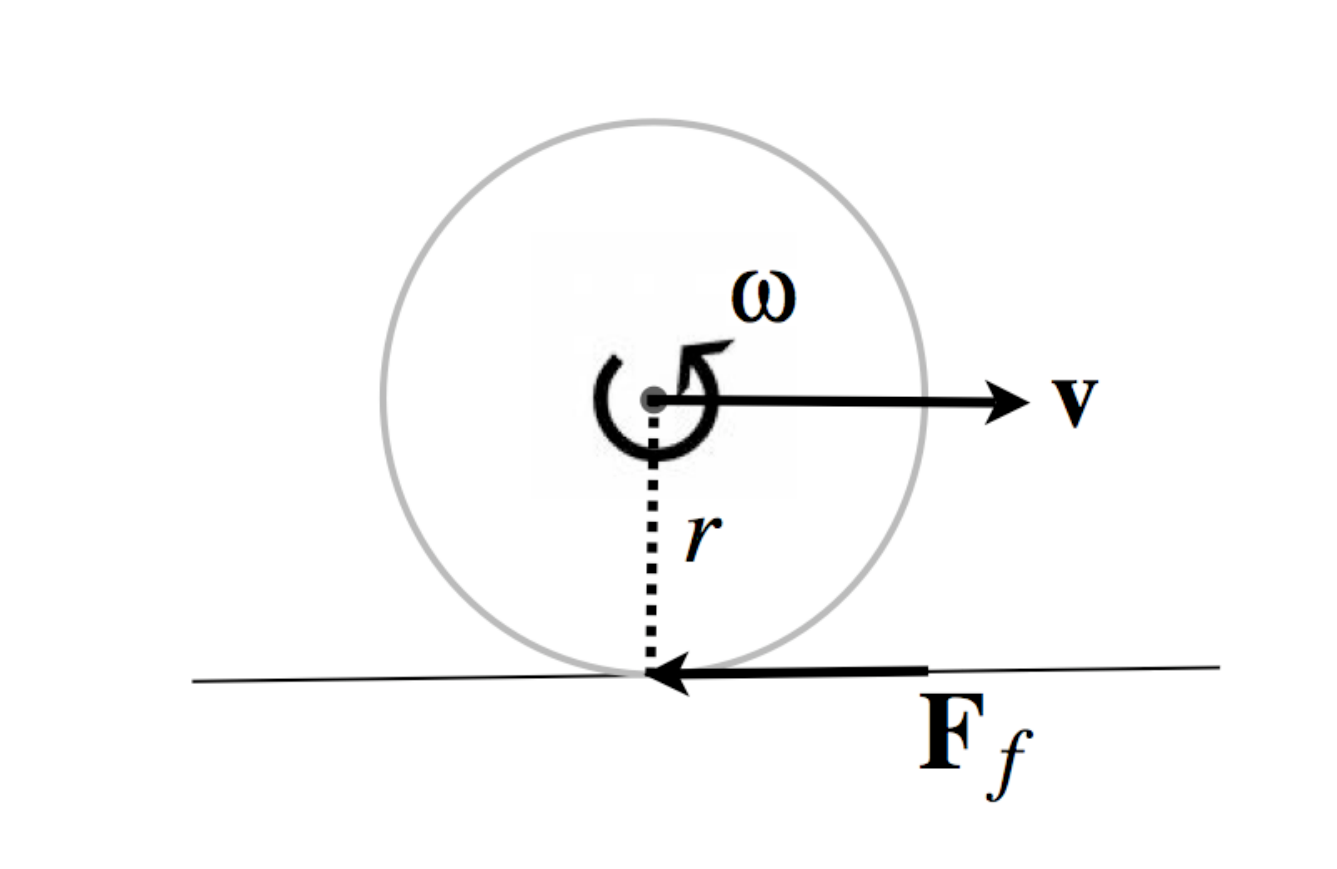}
\caption{\footnotesize  {\bf Sliding Disk.}   Consider a sliding
disk of radius $r$ that is free to translate and rotate on a
surface.   We assume the disk is in sliding frictional contact with
the surface.  The configuration space of the system is
$\operatorname{SE}(2)$, but with the surface constraint the
configuration space is just $\mathbb{R} \times
\operatorname{SO}(2)$. } \label{fig:slidingdisk}
\end{center}
\end{figure}

\begin{figure}[htbp]
\begin{center}
\includegraphics[scale=0.25,angle=0]{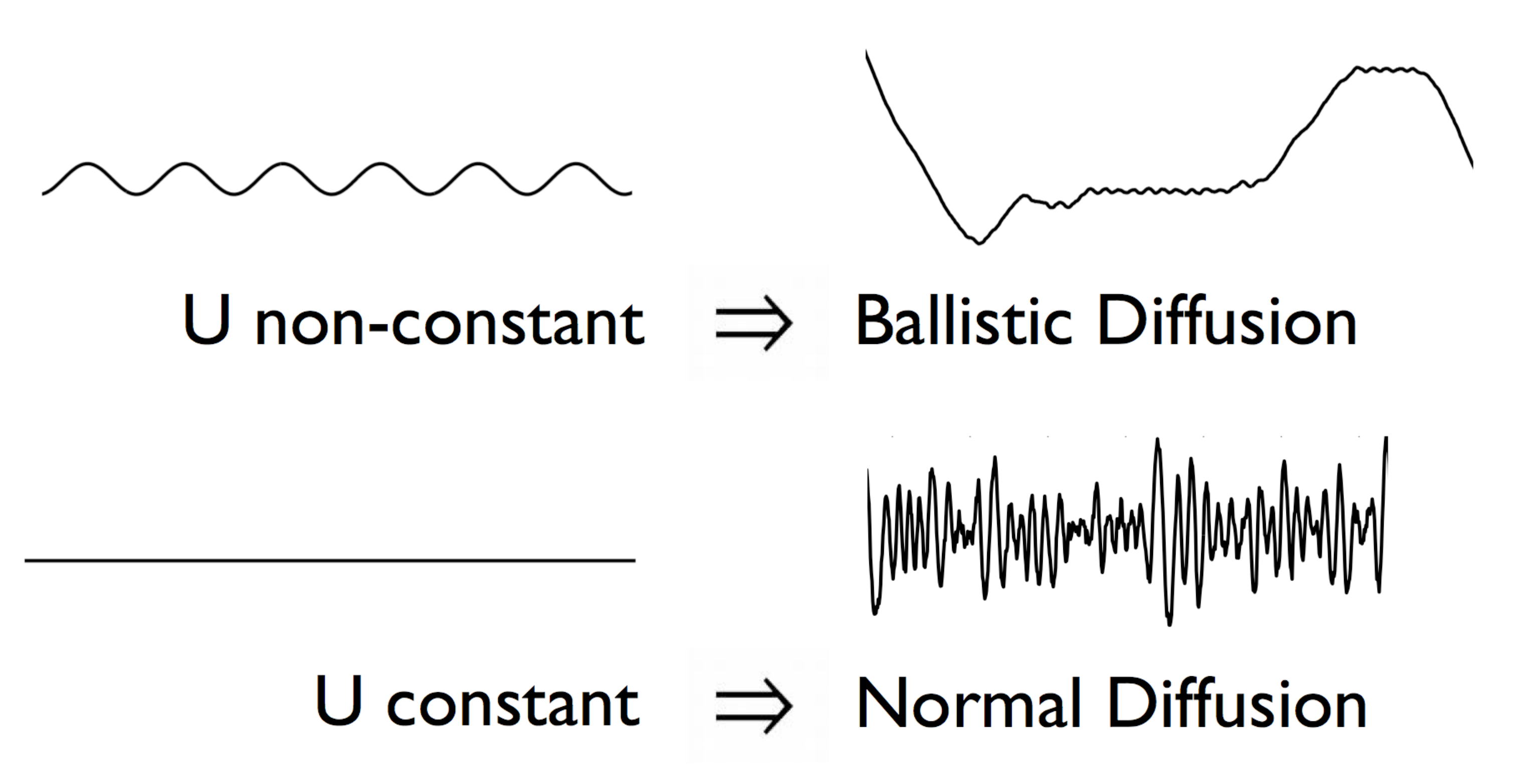}
\caption{\footnotesize  {\bf Ballistic vs.~Normal Diffusion.}   If $U$ is non-constant theorem~\ref{thm:ballistic} implies that the mean squared displacement with respect to the invariant law is ballistic.  Numerically, flights are observed in the $x$-displacement as shown in the diagram.  However, if $U$ is constant this diffusion is normal, i.e., the $x$-displacement behaves like Brownian motion. } 
\label{fig:ballisticvsnormal}
\end{center}
\end{figure}

\begin{figure}[htbp]
\begin{center}
\includegraphics[scale=0.4,angle=0]{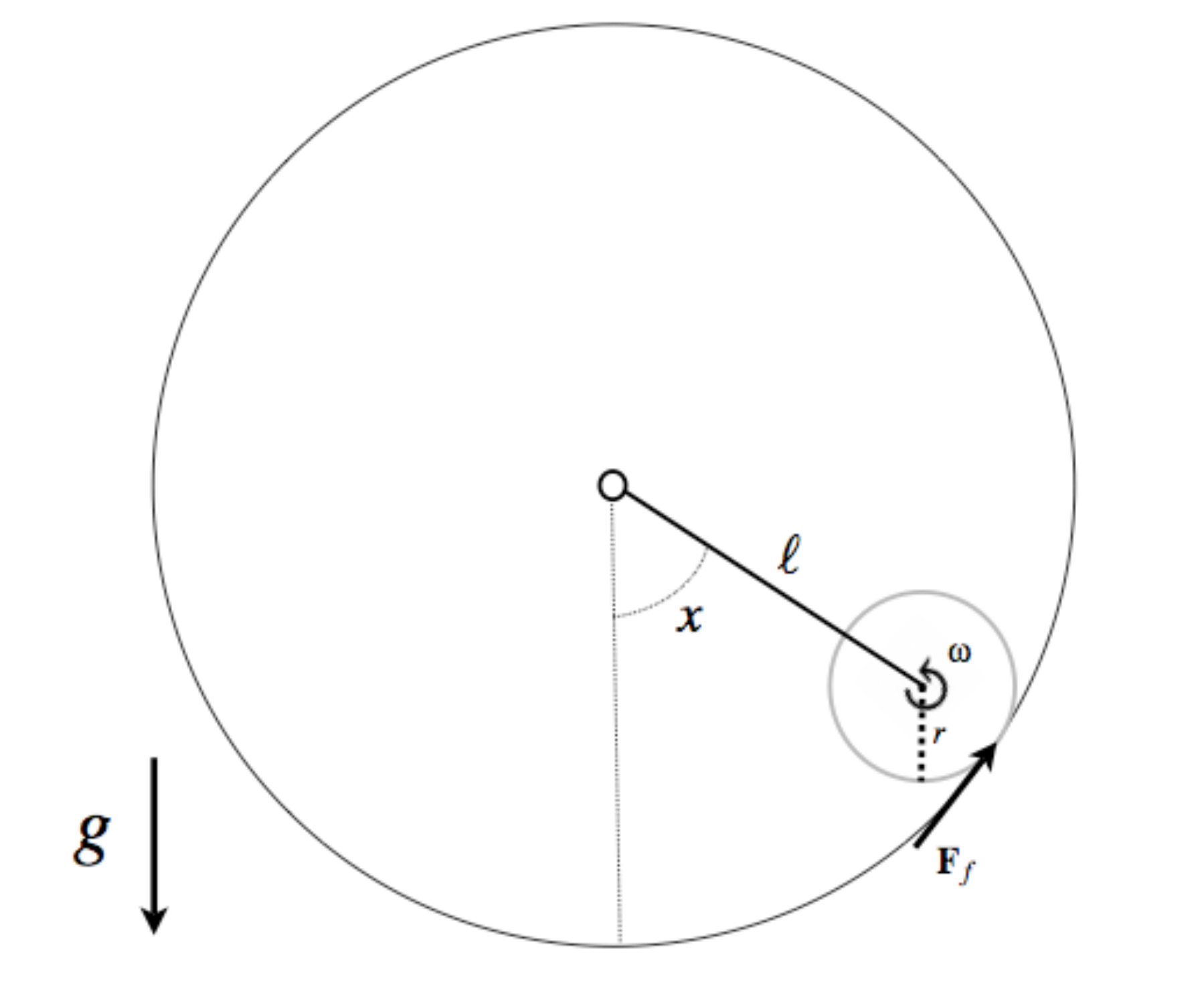}
\caption{\footnotesize  {\bf Ballistic Pendulum.}   If the dimensionless
potential is $U =  \cos(x)$, then the sliding disk is simply a pendulum in which
the bob in the pendulum is replaced by a disk and the pendulum is placed 
within a cylinder as shown.     } \label{fig:ballisticpendulum}
\end{center}
\end{figure}

To understand the behavior of the fluctuation driven motor
prototype, the paper considers in section \ref{kjskjsss8332} a
sliding disk (figure \ref{fig:slidingdisk}) that has the same
essential behavior, but whose configuration space is
$\operatorname{SE}(2)$.  The solution of this simplified system is
a Langevin process with degenerate noise and friction
matrices in the momentums.  The disk is free to slide and rotate.   
Assume that one rescales position by $r$ and time by some characteristic frequency of rotation or other time-scale.  The dimensionless Lagrangian is given by
the difference of kinetic and potential energy
\[
L(x, v, \theta, \omega) = \frac{1}{2} v^2 + \frac{\sigma}{2} \omega^2 - U(x)
\]
where $U: \mathbb{R} \to \mathbb{R}$ is assumed to be smooth and
periodic.   If $U = \cos(x)$, figure~\ref{fig:ballisticpendulum} shows that
the sliding disk is a modified one-dimensional pendulum.  The contact with the 
surface is modelled using a sliding friction law. For this purpose we introduce a symmetric matrix
$\mathbf{C}$ defined as,
\[
\mathbf{C} = \begin{bmatrix} 1 &  1/\sigma \\
                             1/\sigma & 1/\sigma^2  \end{bmatrix} \text{.}
\]
$\mathbf{C}$ is degenerate since the frictional force is actually
applied to only a single degree of freedom, and hence, one of its
eigenvalues is zero.  In addition to friction the system is excited
by white noise so that the governing equations become
 \begin{equation}   \label{eq:slidingdisk}
 \begin{cases}
 \begin{array}{rcl}
 d x &=& v dt  \\
 d \theta &=& \omega dt  \\
\begin{bmatrix} d v \\   d  \omega \end{bmatrix} &=&
\begin{bmatrix} - \partial_x U  \\ 0 \end{bmatrix} dt
- c \mathbf{C} \begin{bmatrix}   v \\    \sigma \omega  \end{bmatrix} dt
+ \alpha \mathbf{C}^{1/2} \begin{bmatrix} d B_v \\ d B_{\omega}
\end{bmatrix}  
 \end{array}
 \end{cases}
 \end{equation}
 where $\mathbf{C}^{1/2}$ is the matrix square root of $\mathbf{C}$. 
 Let $E$ denote the energy of the mechanical system given by,
\begin{equation} \label{eq:slidingdiskE}
E = \frac{1}{2} v^2 + \frac{1}{2} \sigma \omega^2  + U(x) \text{,}
\end{equation}
and let  $\beta = 2 c / \alpha^2$.

In \citep{BoOw2007c}, we show that the Gibbs-Boltzmann distribution
\[
\mu = \exp\left( - \beta E \right) \text{,}
\]
defines the unique invariant measure under the flow of the sliding disk.  This measure
is also ergodic and strongly mixing.  Using this result we prove in \S~\ref{kjskjsss8332}, that if $U$ is non-constant then the $x$-displacement of the sliding disk is $\mu$ a.s.~not ballistic (cf. proposition~\ref{thm:nondirected}).   However, the mean-squared displacement with respect to the invariant law is ballistic (cf. theorem~\ref{thm:ballistic}).  More precisely, we show that the squared standard deviation of the $x$-displacement with respect to its noise-average grows like $t^2$. This implies that the process exhibits not only ballistic transport but also ballistic diffusion.   If $U$ is constant
then the squared standard deviation of the $x$-displacement is diffusive (grows like $t$).

\begin{figure}[htbp]
\begin{center}
\includegraphics[scale=0.25,angle=0]{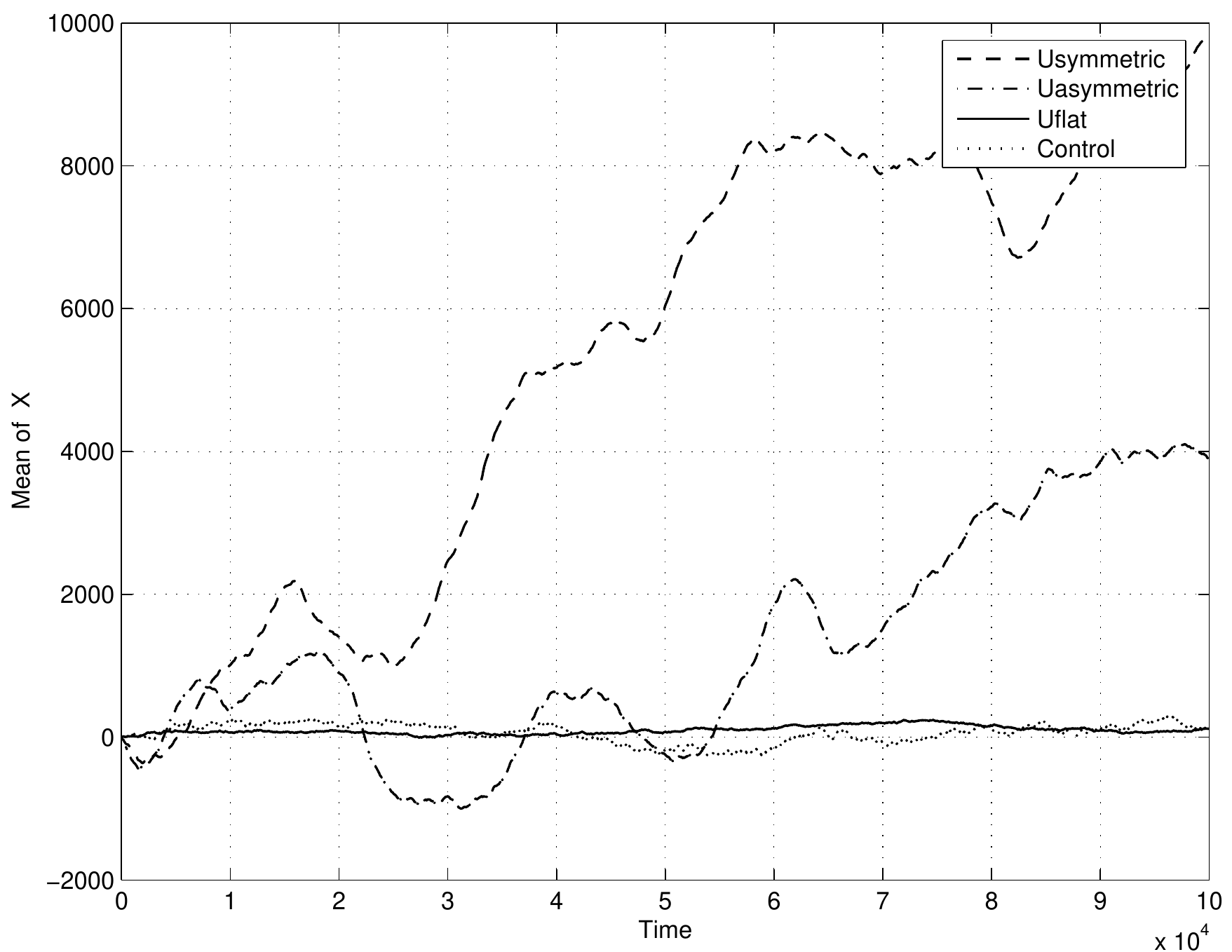} \\
\includegraphics[scale=0.25,angle=0]{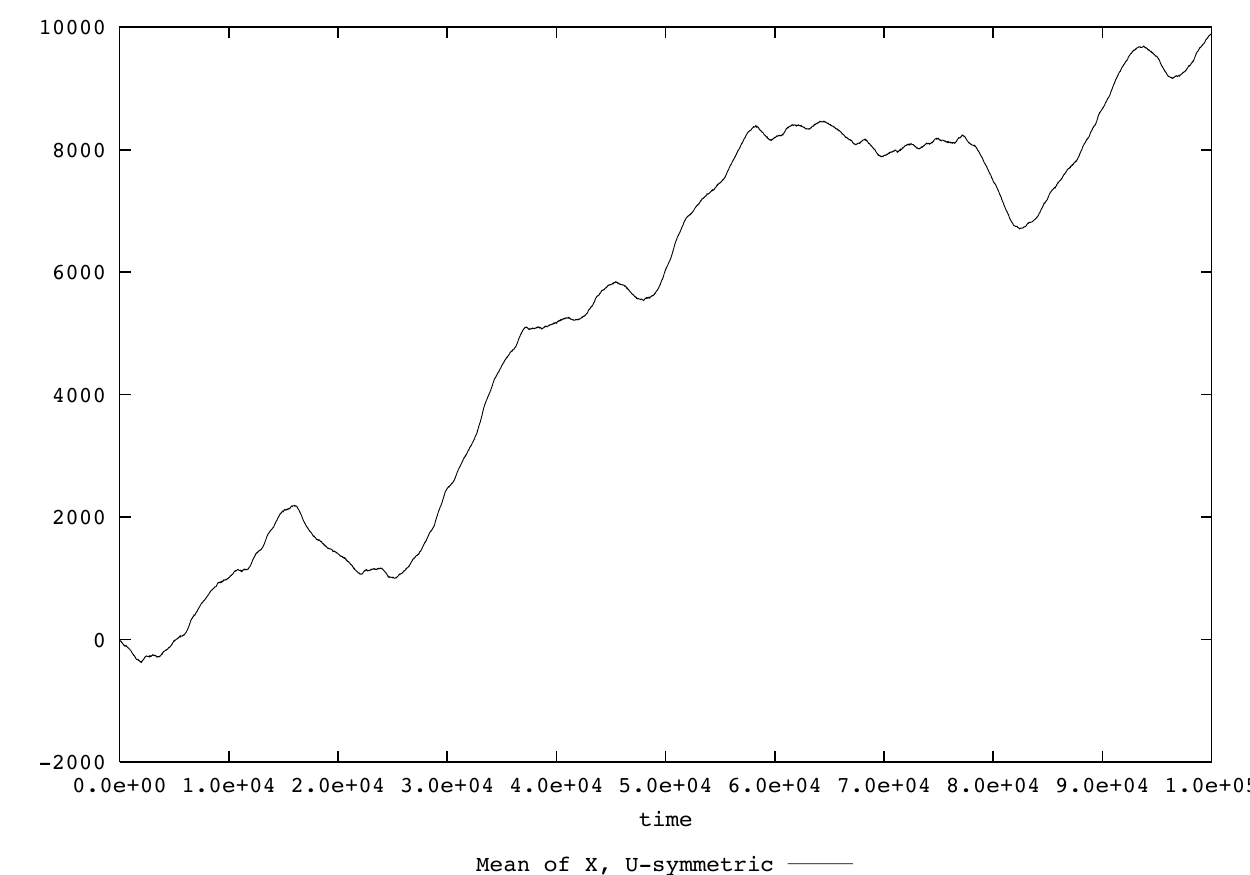}
\includegraphics[scale=0.25,angle=0]{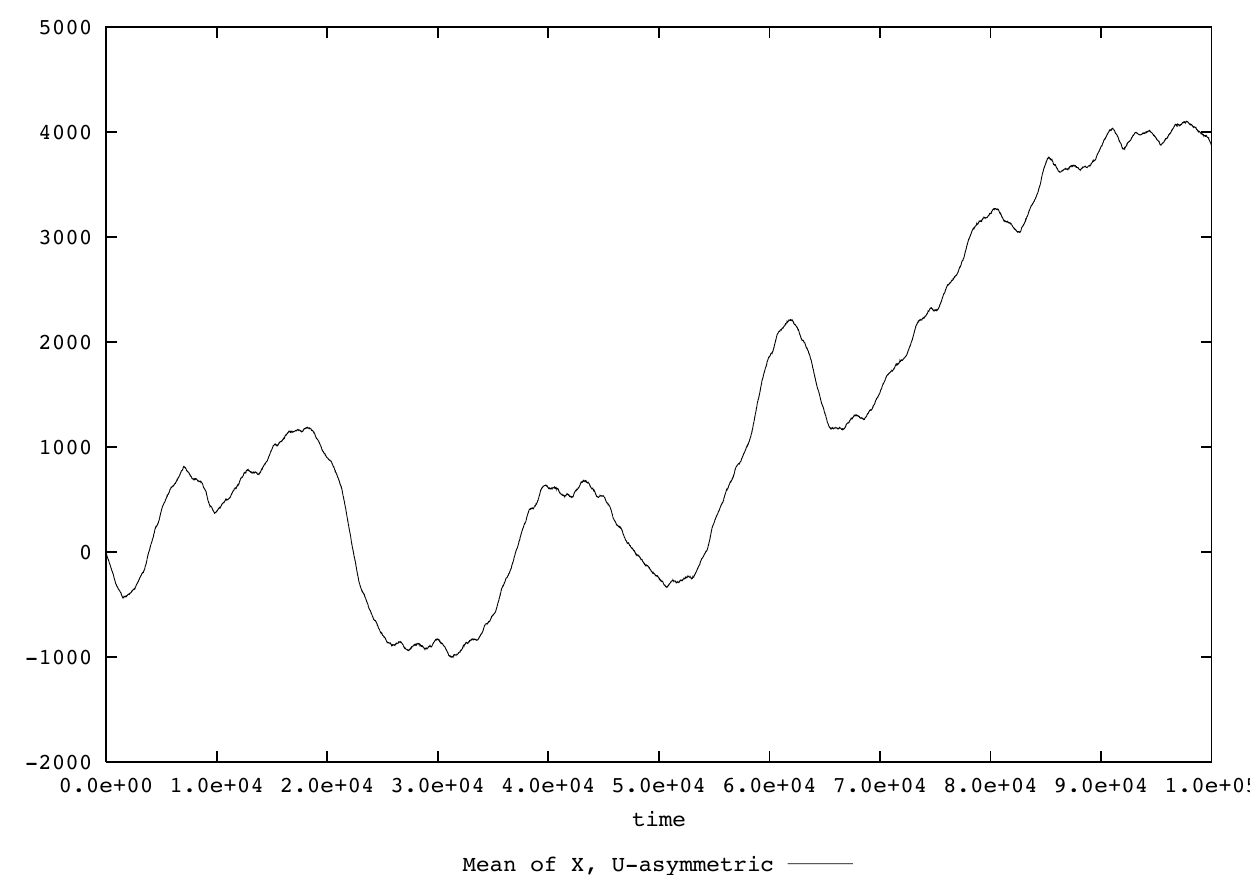}
\includegraphics[scale=0.25,angle=0]{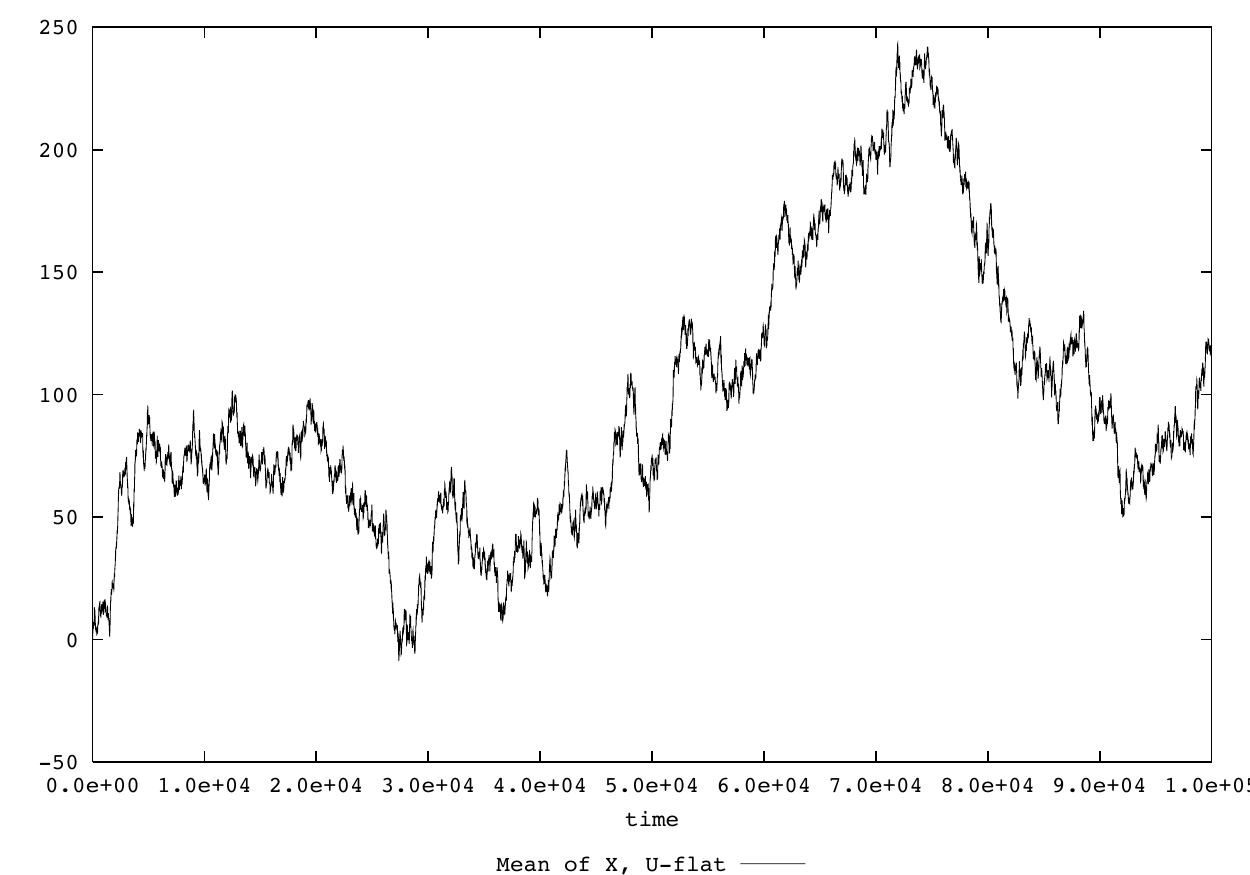}
\includegraphics[scale=0.25,angle=0]{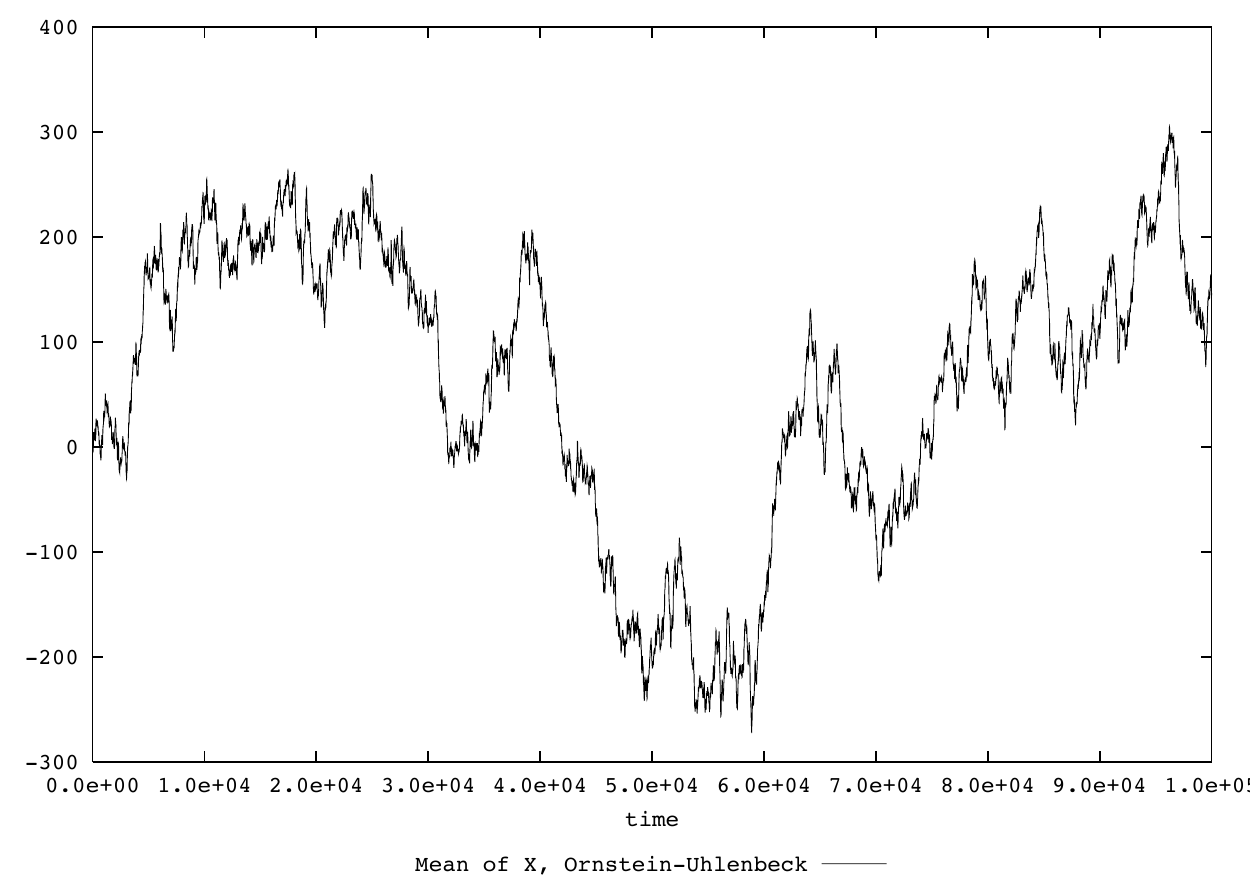} \\
\hbox{\hspace{0.15in}   (a)  U-symmetric \hspace{0.3in} (b)
U-asymmetric \hspace{0.3in} (c)  U-flat \hspace{0.3in}    (d)
Control  \hspace{0.3in}    } \caption{\footnotesize {\bf Sliding Disk at Uniform
Temperature, $h=0.01$, $\alpha=5.0$, $c=0.1$.} The mean of the
$x$-displacement of the disk for $U$ symmetric, asymmetric, flat and
control case.   Figures (a) and (b) show that when $U$ is
non-constant, directed motion as a meta-stable state is possible.
On the other hand, (c) and (d) do not show such behavior. The figure on the top 
superposes these graphs in a single plot for comparison. }
\label{fig:meanX}
\end{center}
\end{figure}

\begin{figure}[htbp]
\begin{center}
\includegraphics[scale=0.3,angle=0]{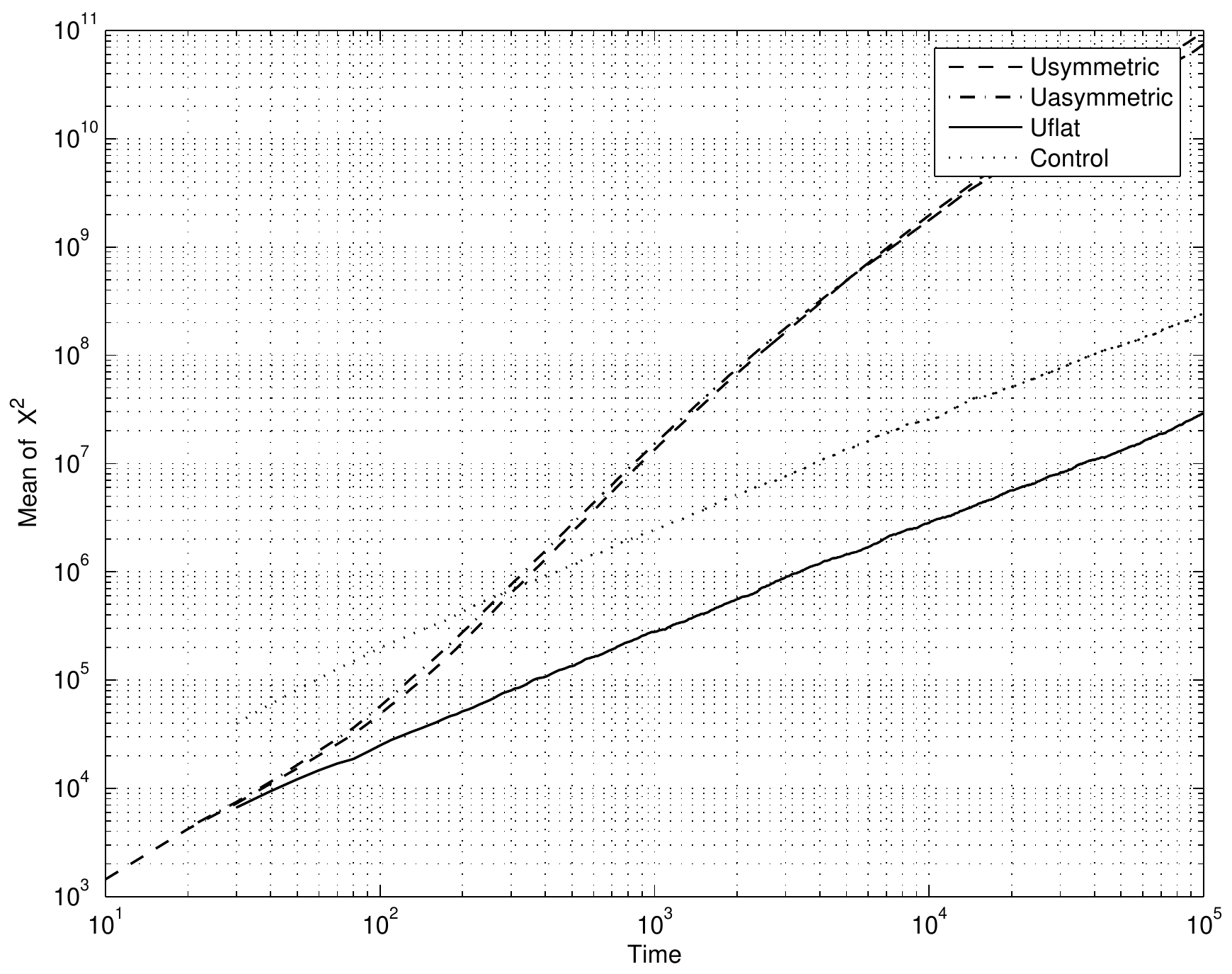}
\caption{\footnotesize {\bf Sliding Disk at Uniform temperature, $h=0.01$,
$\alpha=5.0$, $c=0.1$.} A log-log plot of the mean squared
displacement of the ball.  It clearly shows that the x-position
exhibits anomalous diffusion when $U$ is symmetric or asymmetric. In
the control and flat $U$ cases the diffusion is normal. }
\label{fig:loglogX2}
\end{center}
\end{figure}

\begin{figure}[htbp]
\begin{center}
\includegraphics[scale=0.25,angle=0]{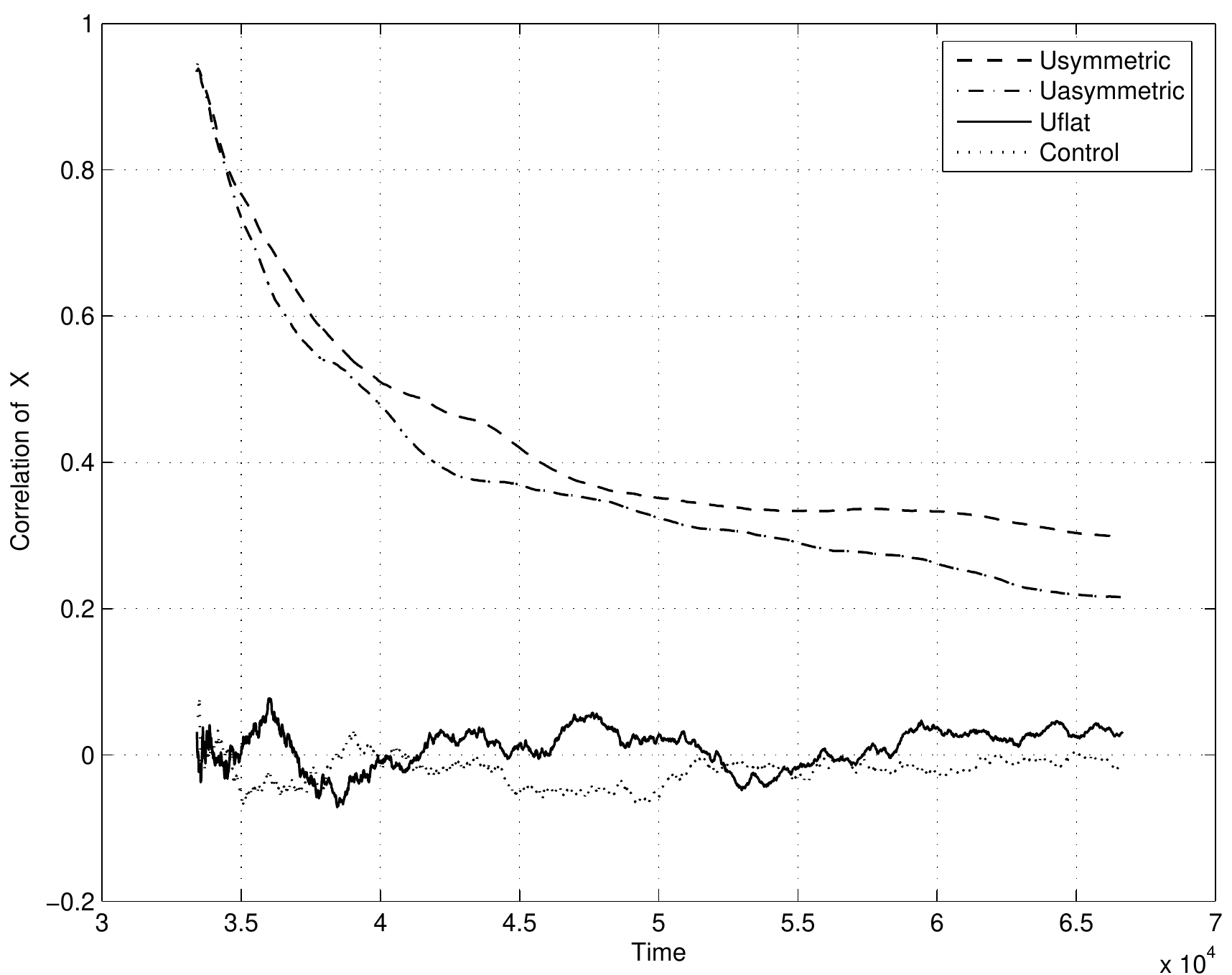} \\
\includegraphics[scale=0.25,angle=0]{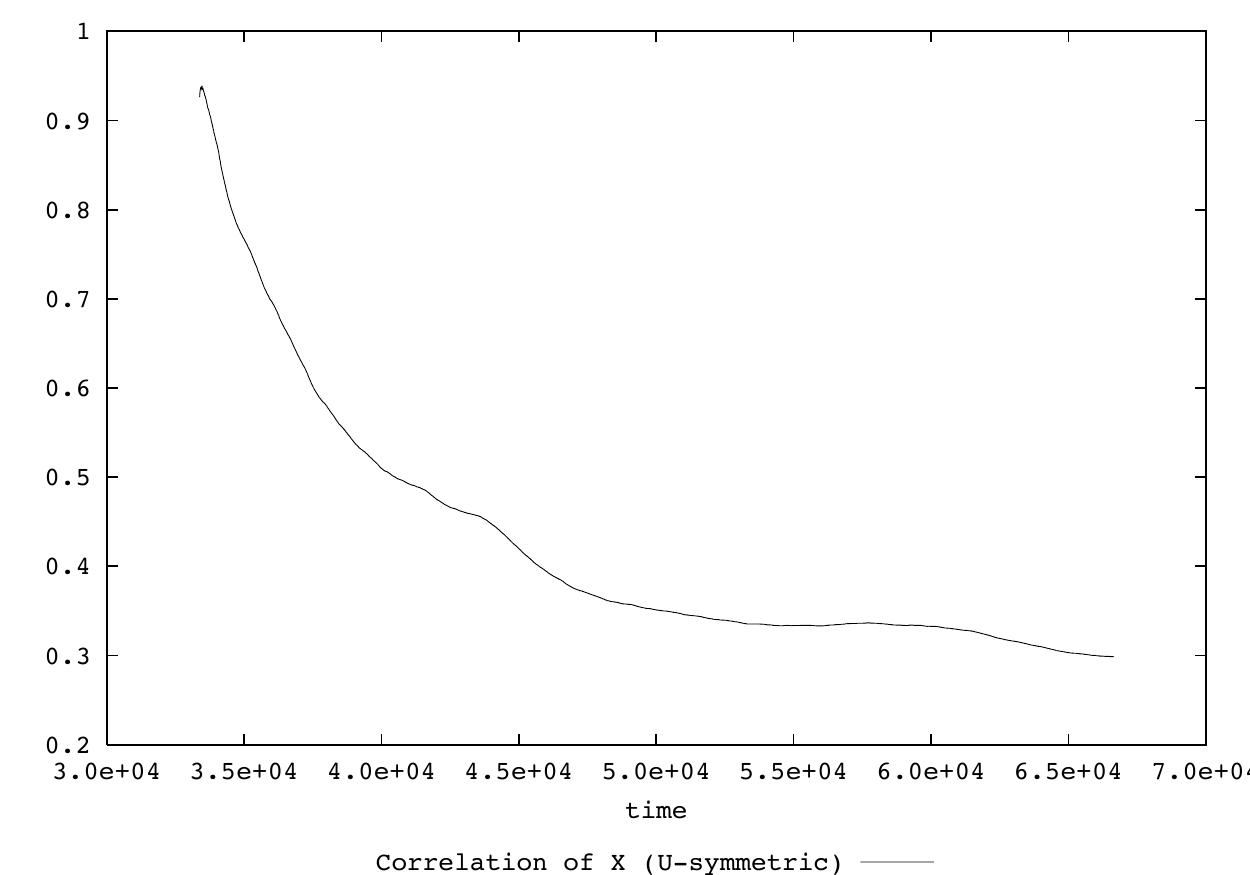}
\includegraphics[scale=0.25,angle=0]{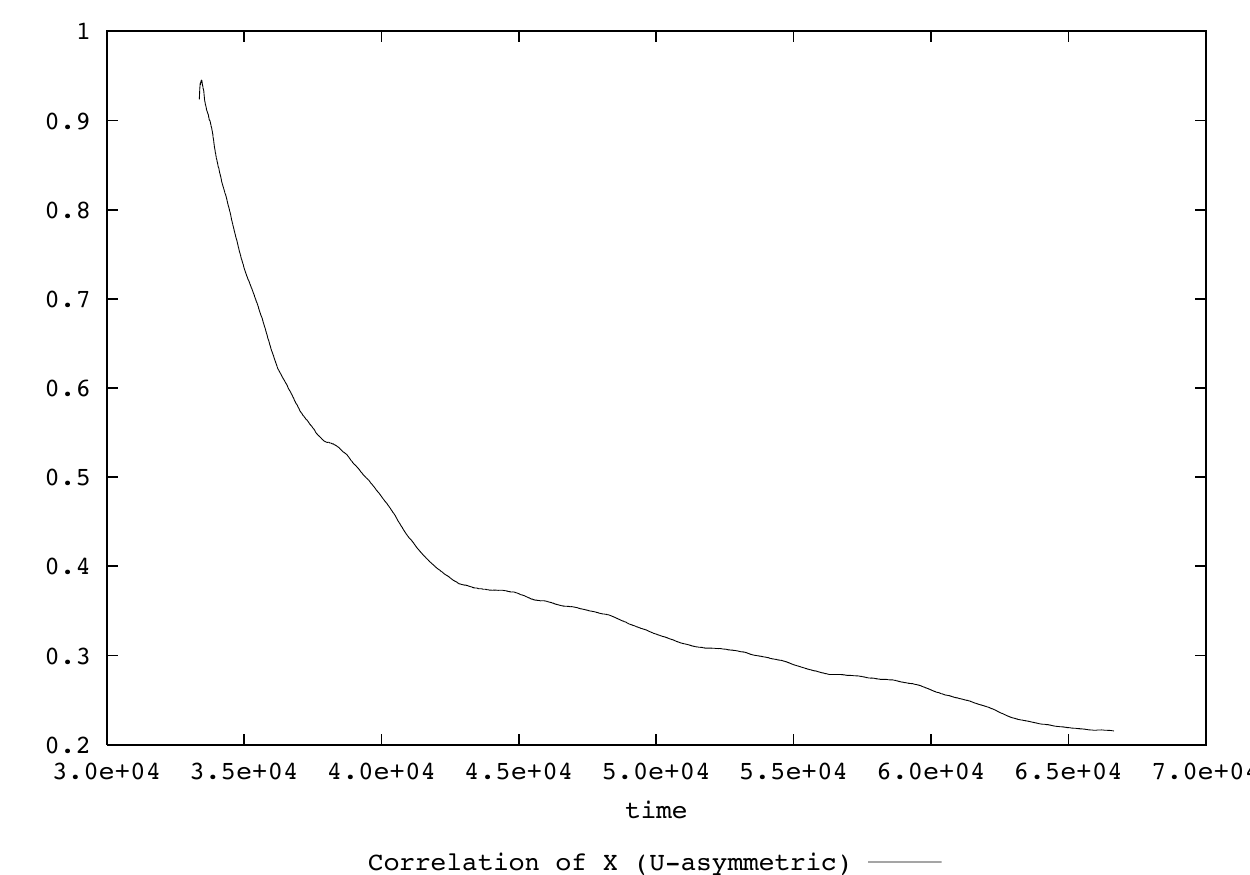}
\includegraphics[scale=0.25,angle=0]{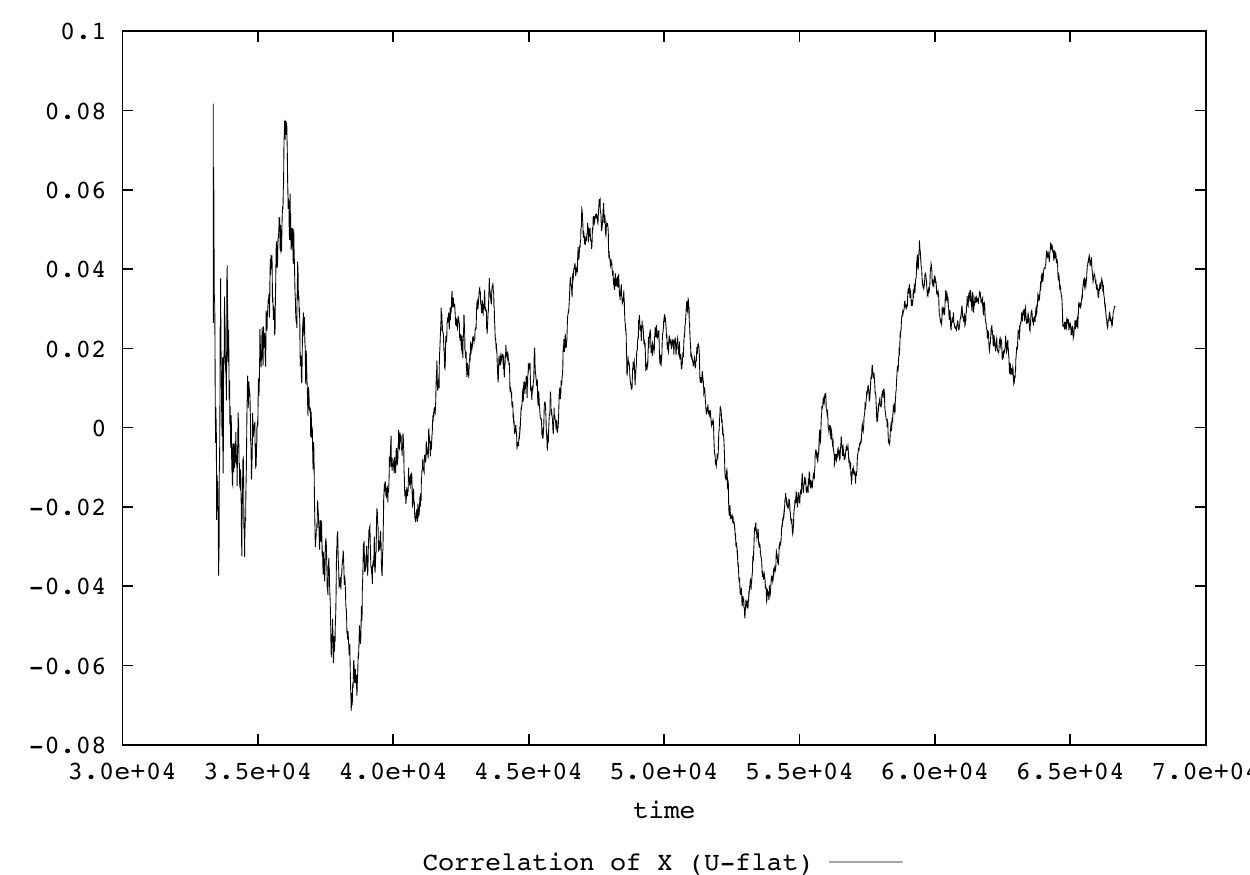}
\includegraphics[scale=0.25,angle=0]{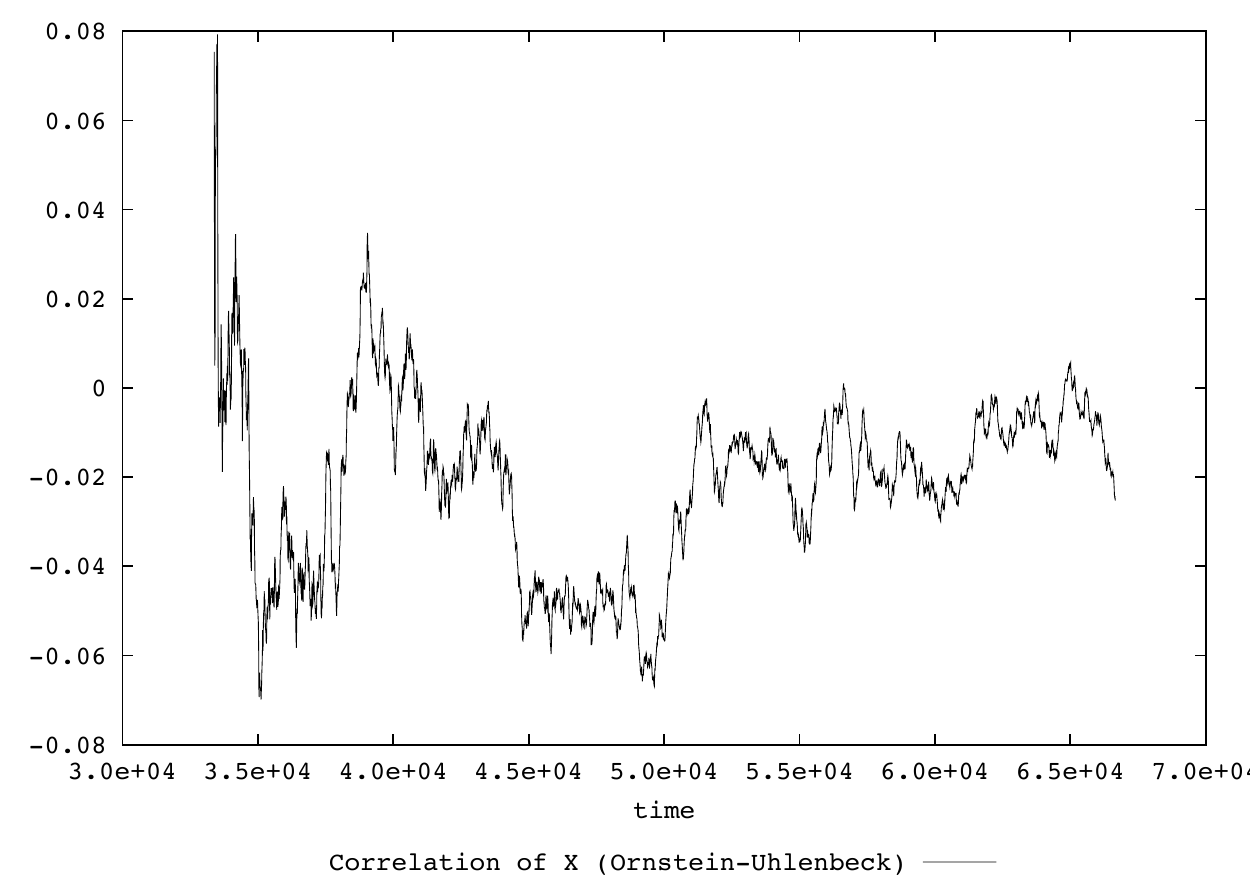} \\
\hbox{\hspace{0.15in}   (a)  U-symmetric \hspace{0.3in} (b)
U-asymmetric \hspace{0.3in} (c)  U-flat \hspace{0.3in}    (d)
Control  \hspace{0.3in}    } \caption{\footnotesize {\bf Sliding Disk at Uniform
Temperature, $h=0.01$, $\alpha=5.0$, $c=0.1$.} The correlation of
the $x$ displacement of the disk for $U$ symmetric, asymmetric,
flat, and the control case.  When $U$ is symmetric or asymmetric the
correlation in the $x$ displacement is nonzero for a certain
time-scale which is larger than the characteristic time-scales
associated with $U$ and with the friction factor. However, in the
other cases the correlation is negligible.  The figure on the top superposes
these graphs in a single plot for comparison.  } \label{fig:corX}
\end{center}
\end{figure}

In the numerics the {\em sliding disk is initially at rest} and averages are computed with respect to realizations.  Numerically one observes the following consequences of this ballistic behavior.  Figure~\ref{fig:meanX} shows that when $U$ is non-constant then the motion is characterized by meta-stable directed motion states. Figure~\ref{fig:loglogX2} shows that the mean square displacement $\E[x(t)^2]$ grows like $t$ (with respect to time) when $U$ is constant or in the case of a one dimensional standard Langevin process (control) whereas 
it grows like $t^2$ as soon as $U$ is non-constant (the motion becomes ballistic).  A diagram comparing the solution behavior in the $U$ constant and non-constant cases is provided in Figure~\ref{fig:ballisticvsnormal}.  Figure~\ref{fig:corX} corresponds to the plot of
\begin{equation}
\frac{\Cov(x(t+2s)-x(t+s),x(t+s)-x(t))}{(\Var(x(t+2s)-x(t+s)))^\frac{1}{2}(\Var(x(t+s)-x(t)))^\frac{1}{2}}
\end{equation}
for $t$ large as a function of $s$. It clearly shows that the system is characterized by long time memory/correlation when $U$ is non-constant whereas it has almost no memory when $U$ is constant or in in the case of a one dimensional standard Langevin process.

\newpage

\section{Sliding Disk: Simplified Model}\label{kjskjsss8332}

To understand the behavior of the prototype stochastic mechanical 
rectifier which will be treated in detail in subsequent sections, 
we designed a simplified model whose configuration space is
$\operatorname{SE}(2)$.  The system consists of a disk sliding on a
surface as shown in Figure \ref{fig:slidingdisk}.  The effect of the
outer ring is modelled as a periodic, one-dimensional potential;
while the effect of the inner ring is incorporated into white noise.
The SDE for the isothermal, sliding disk is a Langevin
process, but the noise and friction matrices are degenerate.  A
statistical numerical analysis discussed below shows that this
process has two interesting statistical properties which are both
linked to rigid-body inertial effects: 1) a non-trivial correlation
on  certain time-scales and 2) ballistic diffusion.

\paragraph{Ballistic Transport at Uniform Temperature}

The proof that the sliding disk is at uniform temperature is based on
finding the generator for (\ref{eq:slidingdisk}) and showing that the 
Gibbs measure is the unique, invariant measure under the flow of this generator.

\begin{theorem} \label{thm:ergodic}
Let $\chi:=\T\times \R\times \T\times \R$ denote the phase space of the sliding disk ($\T$ standing for the torus of dimension one).  Set $\xi:=(x,v,\theta,\omega)$ to be the solution of (\ref{eq:slidingdisk}).  Let $E$ denote the energy of the sliding disk given by,  $E = \frac{1}{2} v^2 + \frac{1}{2} \sigma \omega^2  + U(x)$.  Set  $\beta = 2 c / \alpha^2$ and let $\mu$ be the Gibbs probability
measure defined by
\begin{equation}\label{kjhss}
\mu(d\xi):=\frac{e^{-\beta E}}{Z} d\xi
\end{equation}
where $Z:=\int_{\chi} e^{-\beta E}d\xi$.  If $U$ is non-constant, then the Gibbs measure $\mu$ is ergodic and strongly mixing with respect to the stochastic process $\xi$.  Furthermore, it is the unique, invariant probability measure for the stochastic process $\xi$.
\end{theorem}

This result is a special case of the proof provided in \citep{BoOw2007c}.  To determine the measure is invariant, one computes the infinitesimal generator $\mathcal{L}$ associated to the stochastic process $\xi$ and proves that if $f$ is $\mu$ - measureable then,
\[
\int_{\chi} \mathcal{L} f \mu(d \xi) = 0
\]  
Using the fact that the measure is ergodic, one can readily prove that the $x$-displacement itself is not ballistic.

\begin{proposition} \label{thm:nondirected}
Provided that $U$ is non-constant, then $\mu$ a.s.~
\[
\lim_{t \to \infty} \frac{x(t) - x(0)}{t} \to 0 \text{.}
\]
\end{proposition}

\begin{proof}
Since $\mu$ is ergodic with respect to $\xi$, 
\[
\lim_{t \to \infty} \frac{x(t) - x(0)}{t} = \lim_{t \to \infty} \frac{1}{t} \int_0^t v dt = 
\int_{\chi} v \mu(d \xi) = 0 \text{.}
\]
\end{proof}

Moreover, as stated in the following proposition, the squared standard deviation of the $x+\theta$-degree of freedom of the sliding disk grows like $t$.   Let $\E$ denote the expectation with respect to the Brownian noise and $\E_\mu$ the expectation with respect to the Brownian noise and the initial configuration (sampled from the invariant measure $\mu$). Set $x_t = x(t)$ and $\theta_t = \theta(t)$.

\begin{proposition} 
The squared standard deviation of the $x_t + \theta_t$-degree of freedom is diffusive, i.e.,
\begin{equation}\label{kskjhs3}
\lim_{t\rightarrow
\infty}\frac{\E_\mu[(x_t+\theta_t-\E[x_t+\theta_t])^2]}{t}=\frac{2
\alpha^2 \sigma^2}{c^2 (\sigma^2+1) } \text{.}
\end{equation}

\end{proposition}
\begin{proof}
First, diagonalize the diffusion and friction matrices in
(\ref{eq:slidingdisk}) using the following invertible matrix:
\[
\mathbf{V} = \frac{1}{1+\sigma}  \begin{bmatrix} - 1   & \sigma \\
1 & 1 \end{bmatrix}
\]
as follows,
\begin{align*}
\mathbf{V}  \begin{bmatrix} d v \\   d  \omega \end{bmatrix} =&
\mathbf{V}  \begin{bmatrix} - \partial_x U  \\ 0 \end{bmatrix} dt -
c
\begin{bmatrix}
0 & 0 \\
0 & \frac{\sigma + 1}{\sigma}   \end{bmatrix} \mathbf{V}
\begin{bmatrix} v \\    \omega  \end{bmatrix} dt  + \frac{\alpha
(\sigma+1)}{\sqrt{\sigma^2+1}}
\begin{bmatrix}
0 & 0 \\
0 &  1  \end{bmatrix} \mathbf{V}
 \begin{bmatrix} d B_v \\ d B_{\omega} \end{bmatrix}
\text{.}
 \end{align*}
 Simplifying this expression yields the following pair of equations:
 \begin{align}
d \left( - v + \sigma \omega \right) &=   \partial_x U dt \text{,}\label{eq:1} \\
d \left(  v +  \omega \right) &= -  \partial_x U dt - c \gamma
\left( v + \omega \right) + \bar{\alpha} ( d B_v + d B_{\omega} )
\text{.}\label{eq:2}
 \end{align}
where
\[
\gamma = \frac{ (\sigma+1)}{ \sigma},~~~\bar{\alpha}= \frac{ \alpha
(\sigma +1) }{ \sqrt{\sigma^2 +  1} } \text{.}
\]

Set $B_s:=(B_v +  B_{\omega})/\sqrt{2}$. Integrating \eref{eq:1} and
\eref{eq:2} gives
 \begin{align}
 - v + \sigma \omega  &= -v_0+\sigma \omega_0 +\int_0^t  \partial_x U ds \text{,}\label{eq:1w} \\
 v +  \omega  &=(v_0+\omega_0)e^{-c\gamma t}- \int_0^t e^{-c\gamma (t-s)}\partial_x U ds + \sqrt{2}\bar{\alpha}  \int_0^t e^{-c\gamma (t-s)} dB_s \text{.}\label{eq:2w}
 \end{align}
Integrating~\eref{eq:2w} gives
\begin{align*}
(x_t+\theta_t)-(x_0+\theta_0)=&\frac{v_0+\omega_0}{c
\gamma}(1-e^{-c\gamma t})- \int_0^t \frac{1- e^{-c\gamma (t-s)}}{c
\gamma}\partial_x U ds \\&+ \sqrt{2}\bar{\alpha} \int_0^t \frac{1-
e^{-c\gamma (t-s)}}{c \gamma} dB_s \text{.}
\end{align*}
The result follows by observing that
\begin{equation}
\int_0^t \partial_x U ds=(-v_t+\sigma \omega_t)-(-v_0+\sigma
\omega_0)
\end{equation}

\end{proof}

However, one can prove that the squared standard deviation of the $-x_t + \sigma \theta_t$-degree of freedom grows like $t^2$.  This implies that the process exhibits not only ballistic transport but also ballistic diffusion along this degree of freedom.

\begin{proposition}
Assume that $U$ is non constant, then
\begin{equation}\label{ksddssdkjhs3}
\lim\sup_{t\rightarrow \infty}\frac{\E_\mu\Big[\big(-x_t+\sigma
\theta_t-\E[-x_t+\sigma \theta_t]\big)^2\Big]}{t^2}\leq
4\frac{1+\sigma}{\beta}
\end{equation}
and
\begin{equation}\label{ksddssdkjhs3}
\lim\inf_{t\rightarrow \infty}\frac{\E_\mu\Big[\big(-x_t+\sigma
\theta_t-\E[-x_t+\sigma \theta_t]\big)^2\Big]}{t^2}\geq
\frac{1}{4}\frac{1+\sigma}{\beta}
\end{equation}
\end{proposition}

\begin{proof}
From Cauchy-Schwartz inequality one obtains that
\[
\left(\frac{1}{t}\int_0^t (-v_s+\sigma \omega_s) \,ds\right)^2\leq
\frac{1}{t}\int_0^t (-v_s+\sigma \omega_s)^2\,ds.
\]
Hence
\begin{equation}
\lim\sup_{t\rightarrow \infty}\frac{\E_\mu\Big[\big(-x_t+\sigma
\theta_t-(-x_0+\sigma \theta_0)\big)^2\Big]}{t^2} \leq
\mu\big[(-v_0+\sigma \omega_0)^2\big]\text{.}
\end{equation}
We obtain the first inequality of the proposition by observing that
\[
\mu\big[(-v_0+\sigma \omega_0)^2\big]=\frac{1+\sigma}{\beta}
\text{.}
\]
and
\begin{equation}
\E_\mu\Big[\big(\E[-x_t+\sigma \theta_t]-(-x_0+\sigma
\theta_0)\big)^2\Big] \leq \E_\mu\Big[\big(-x_t+\sigma
\theta_t-(-x_0+\sigma \theta_0)\big)^2\Big]
\end{equation}
Let us now prove the lower bound. Integrating equation \eref{eq:1w}
gives
 \begin{equation}\label{jhsjhge2}
 - x_t + \sigma \theta_t-\E[-x_t+\sigma \theta_t]  = t \int_0^t  (1-\frac{s}{t})(\partial_x U-\E[\partial_x U]) ds
 \end{equation}
Write
 \begin{equation}
A_t =  \int_0^t  (1-\frac{s}{t})(\partial_x U-\E[\partial_x U]) ds
 \end{equation}
 and
  \begin{equation}
B_t =  \int_0^t  \frac{s}{t}(\partial_x U-\E[\partial_x U]) ds
 \end{equation}
 Observe that
  \begin{equation}
A_t+B_t = (-v_t+\sigma \omega_t)-\E[-v_t+\sigma \omega_t]\text{.}
 \end{equation}
Since $\mu$ is strongly mixing when $U$ is non constant, it follows that \citep{DaZa1996}
  \begin{equation}
\lim\inf_{t\rightarrow \infty}\E_\mu[(A_t+B_t)^2] = \mu[(-v_t+\sigma
\omega_t)^2]
 \end{equation}
Furthermore
  \begin{equation}
\E_\mu[(A_t+B_t)^2] \leq 2 \big(\E_\mu[A_t^2]+\E_\mu[B_t^2]\big)
 \end{equation}
and since the law of the process $(x_t,\theta_t,v_t,\omega_t)$
remains invariant under $\P_\mu$ by reversing time and flipping the
velocities $v_t,\omega_t$ we deduce that
$\E_\mu[A_t^2]=\E_\mu[B_t^2]$ and
  \begin{equation}
\lim\inf_{t\rightarrow \infty}\E_\mu[A_t^2] =
\frac{1}{4}\mu[(-v_t+\sigma \omega_t)^2]
 \end{equation}
We conclude by the taking the expectation of square of
\eref{jhsjhge2} with respect to $\E_\mu$.
\end{proof}

The following theorem is a straightforward consequence of the
previous propositions.

\begin{Theorem}  \label{thm:ballistic}
We have
\begin{itemize}
\item If $U$ is constant then
\begin{equation}\label{kskjhs3}
\lim_{t\rightarrow
\infty}\frac{\E_\mu\big[(x_t-\E[x_t])^2\big]}{t}=\frac{2 \alpha^2
\sigma^2}{c^2 (\sigma^2+1) (\sigma+1)^2}
\end{equation}
\item If $U$ is non constant then
\begin{equation}\label{ksddssdkjhs3}
\lim\sup_{t\rightarrow
\infty}\frac{\E_\mu\big[(x_t-\E[x_t])^2\big]}{t^2}\leq
\frac{4}{\beta(1+\sigma)}
\end{equation}
and
\begin{equation}\label{ksddssdkjhs3}
\lim\inf_{t\rightarrow
\infty}\frac{\E_\mu\big[(x_t-\E[x_t])^2\big]}{t^2}\geq
\frac{1}{4\beta(1+\sigma)}
\end{equation}
\end{itemize}
\end{Theorem}

\begin{Remark}
Using Ito's formula and (\ref{eq:slidingdiskE}) one obtains that,
\[
d E = - c (v+\omega)^2 dt + \frac{\alpha^2}{2} \left(1 +
\frac{1}{\sigma} \right) dt + \text{martingales} \text{.}
\]
Integrating this expression gives
\[
E(t) - E(0) = -c \int_0^t  (v+ \omega)^2 ds  +  \frac{\alpha^2}{2}
\left(1 + \frac{1}{\sigma} \right) t + \text{martingales.}
\]
The first and second terms represent the energy loss due to friction
and the energy injected due to the noise respectively.
\end{Remark}

\begin{Remark}
Setting $v_t = v(t)$, it follows from (\ref{eq:1w}) and (\ref{eq:2w}) that
\begin{align} 
v_t=&\frac{\sigma (v_0+\omega_0)e^{-c\gamma t}+v_0-\sigma
\omega_0}{\sigma+1}- \int_0^t \frac{\sigma e^{-c\gamma
(t-s)}+1}{\sigma+1}\partial_x U ds \nonumber \\
&+ \frac{\sqrt{2}\bar{\alpha}}{\sigma+1}  \int_0^t
e^{-c\gamma (t-s)} dB_s \text{.} \label{eq:v}
\end{align}
The long term memory effect exhibited in Figure~\ref{fig:corX}
has its origin in the term $\int_0^t \partial_x U(x(s))  ds$ in (\ref{eq:v}). That term is equal to zero when $U$ is constant, and itself has its origin in the rigid
body interaction between rotation and translation and the fact the
friction matrix is singular.
\end{Remark}

\paragraph{Stochastic Variational Integrator}

To simulate the dynamics of the sliding disk at constant
temperature, a stochastic variational Euler method is applied \citep{BoOw2007a}.  The
discrete scheme for the isothermal case is given explicitly by:
\begin{equation} \label{eq:svi} \begin{cases} \begin{array}{rcl}
  x_{n+1} &=& x_n + h v_{n+1}  \text{,} \\
  \theta_{n+1} &=& \theta_n + h \omega_{n+1}  \text{,} \\
\begin{bmatrix} v_{n+1} \\     \omega_{n+1} \end{bmatrix} &=&
\begin{bmatrix} v_n \\     \omega_n \end{bmatrix} +
h \begin{bmatrix} - \partial_x U(x_n)  \\ 0 \end{bmatrix} - h c
\mathbf{C} \begin{bmatrix}   v_n \\    \sigma \omega_n  \end{bmatrix}  +
 \alpha \mathbf{C}^{1/2} \begin{bmatrix} d B_v \\ d B_{\omega} \end{bmatrix} \text{.} 
 \end{array} \end{cases}
 \end{equation}
 It is an explicit, first-order strongly convergent method.

\paragraph{Simulation}

\begin{figure}[htbp]
\begin{center}
\includegraphics[scale=0.3,angle=0]{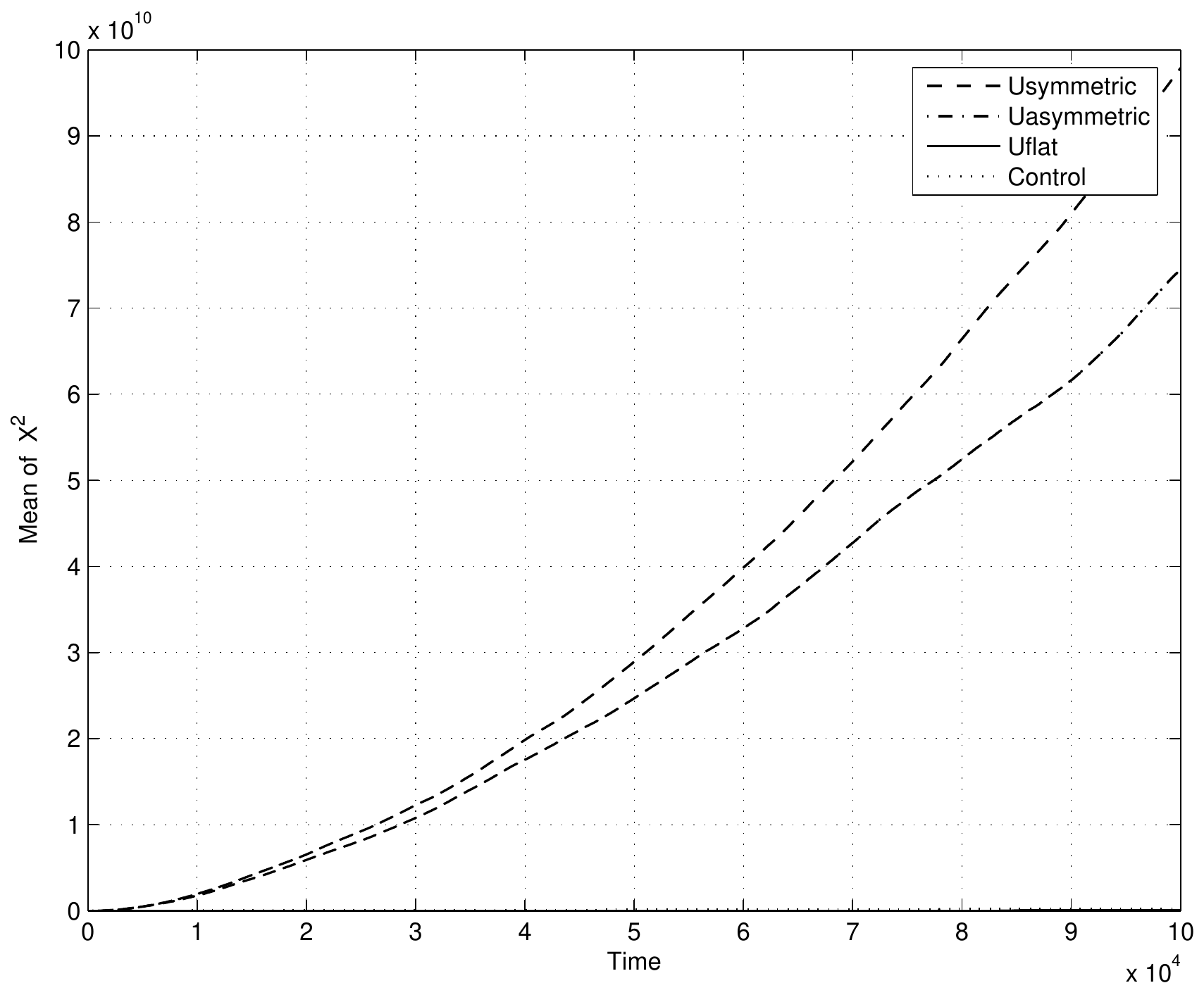} \\
\includegraphics[scale=0.25,angle=0]{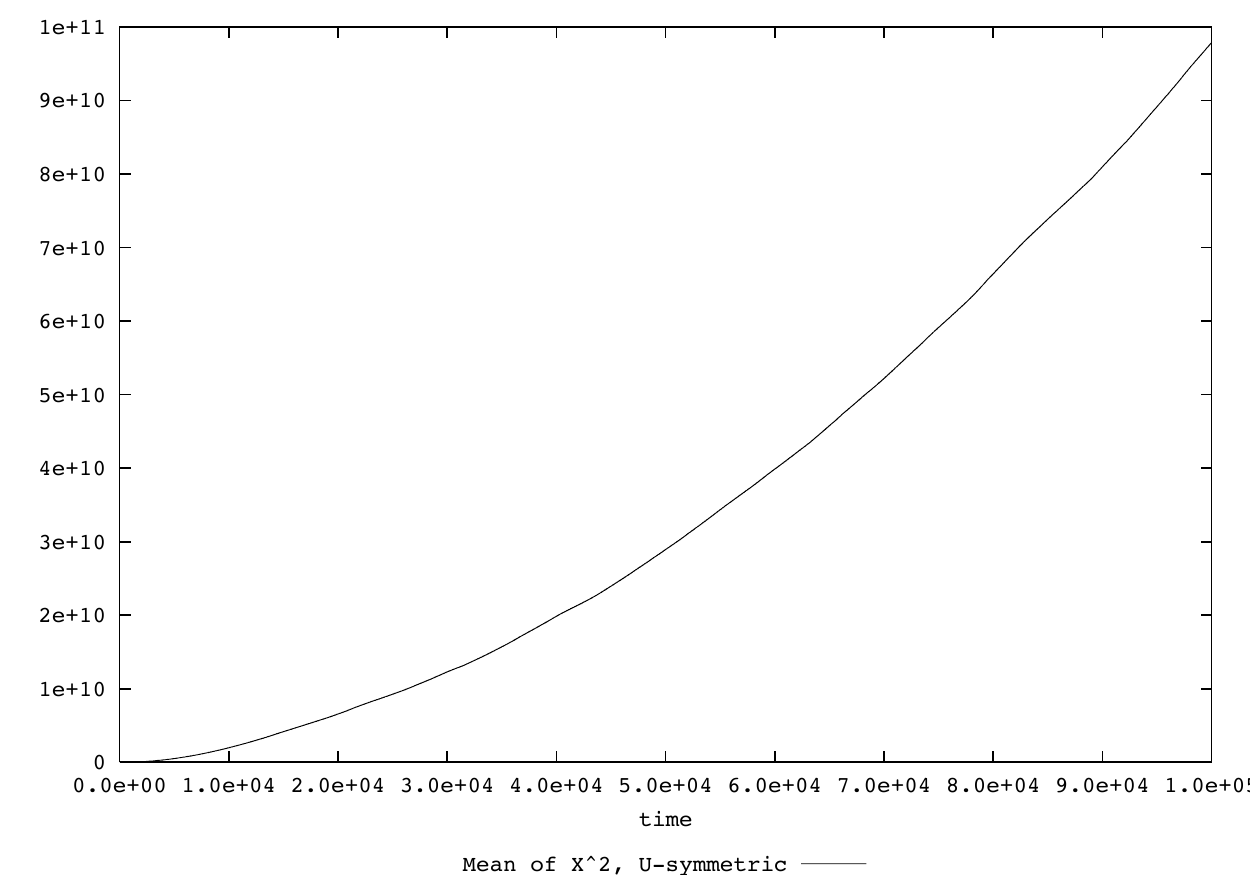}
\includegraphics[scale=0.25,angle=0]{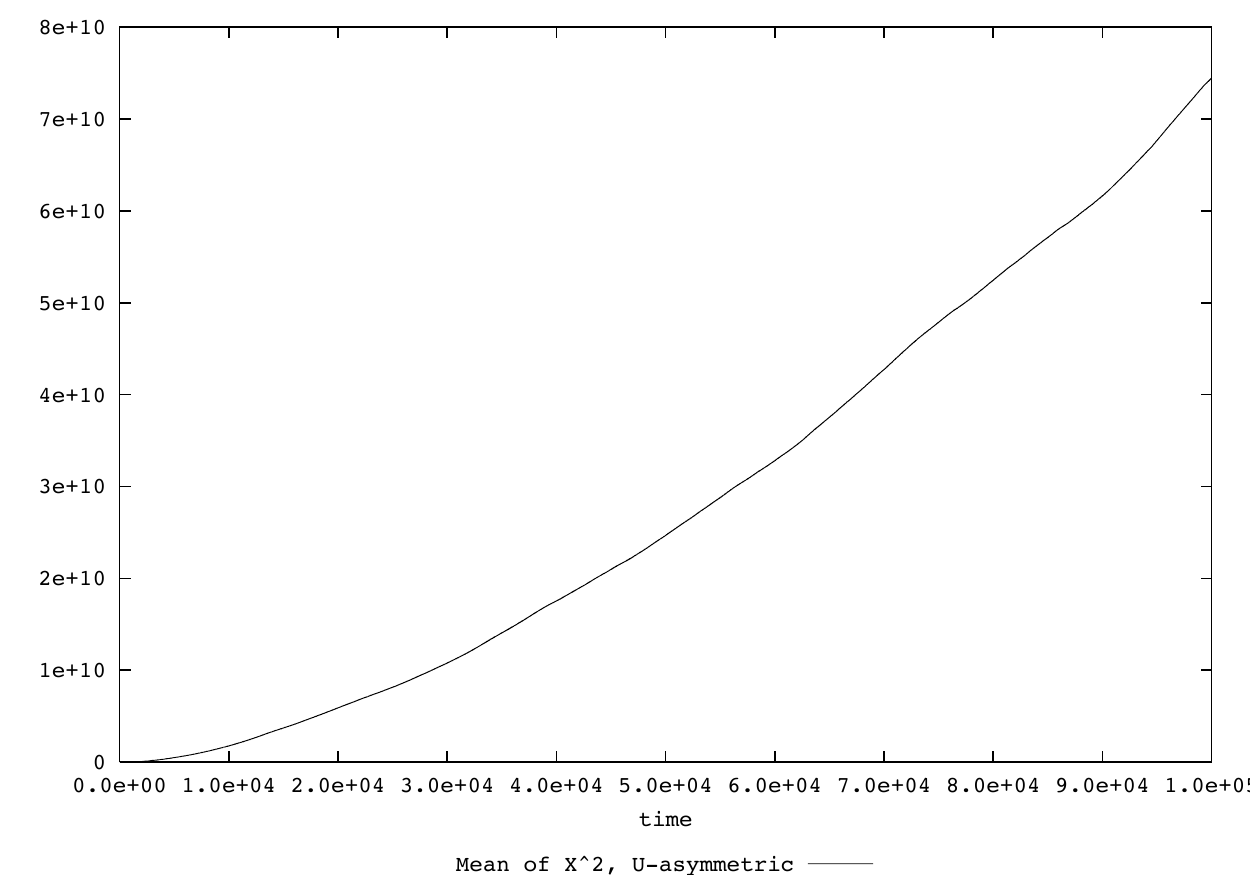}
\includegraphics[scale=0.25,angle=0]{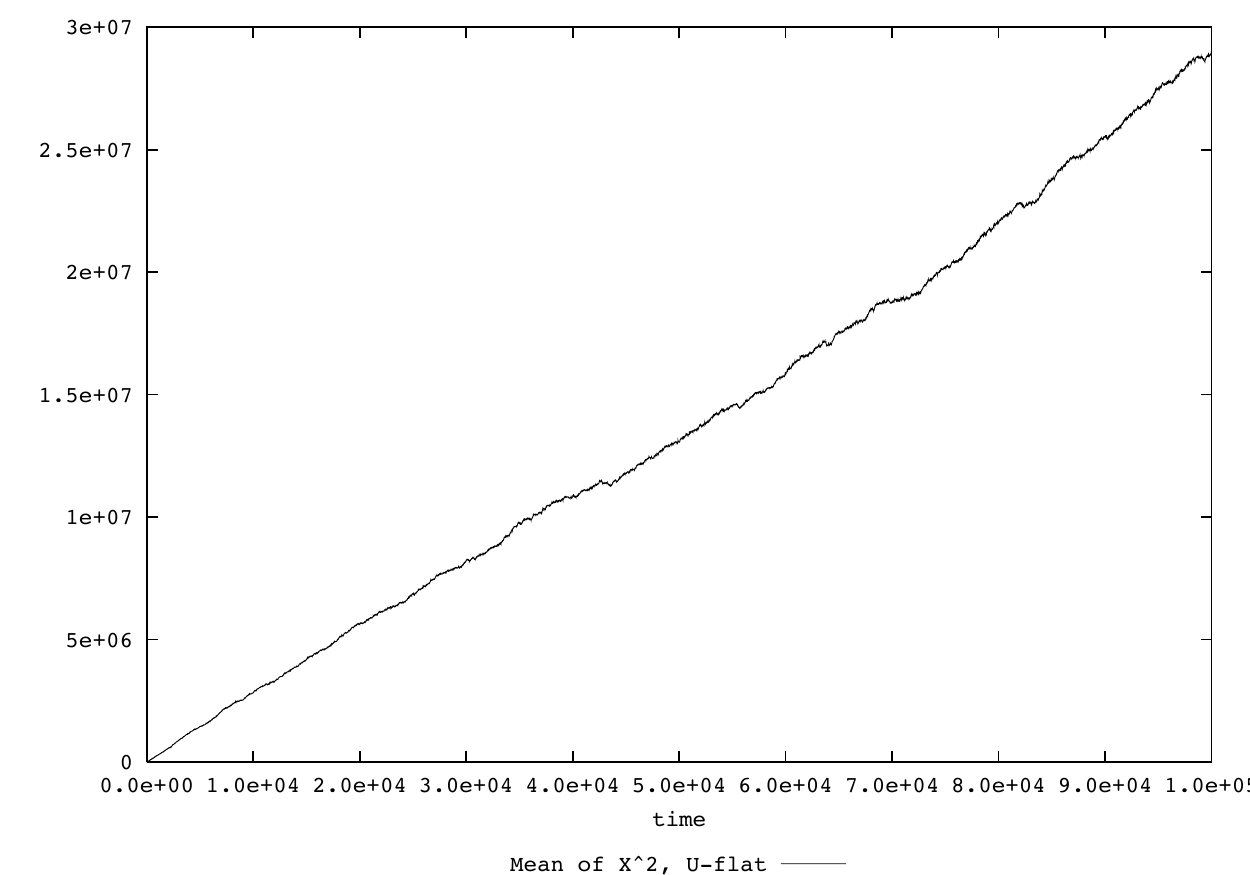}
\includegraphics[scale=0.25,angle=0]{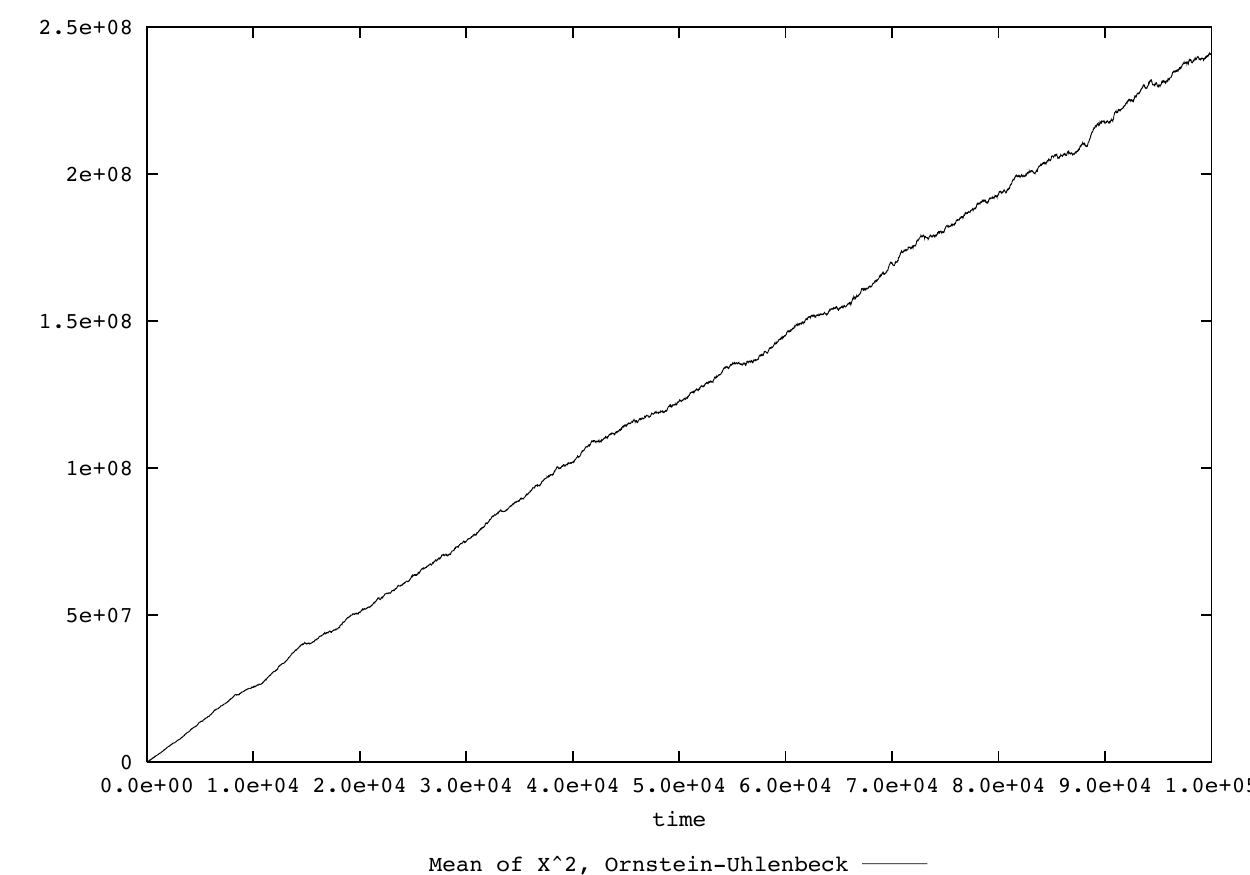}  \\
\hbox{\hspace{0.15in}   (a)  U-symmetric \hspace{0.3in} (b)
U-asymmetric \hspace{0.3in} (c)  U-flat \hspace{0.3in}    (d)
Control  \hspace{0.3in}    } \caption{\footnotesize {\bf Sliding Disk at Uniform
Temperature, $h=0.01$, $\alpha=5.0$, $c=0.1$.} From left: the mean
squared position of the disk for $U$ symmetric, asymmetric, flat,
and the control. The figure on the top superposes
these graphs in a single plot for comparison.  } \label{fig:meanX2}
\end{center}
\end{figure}

\begin{figure}[htbp]
\begin{center}
\includegraphics[scale=0.2,angle=0]{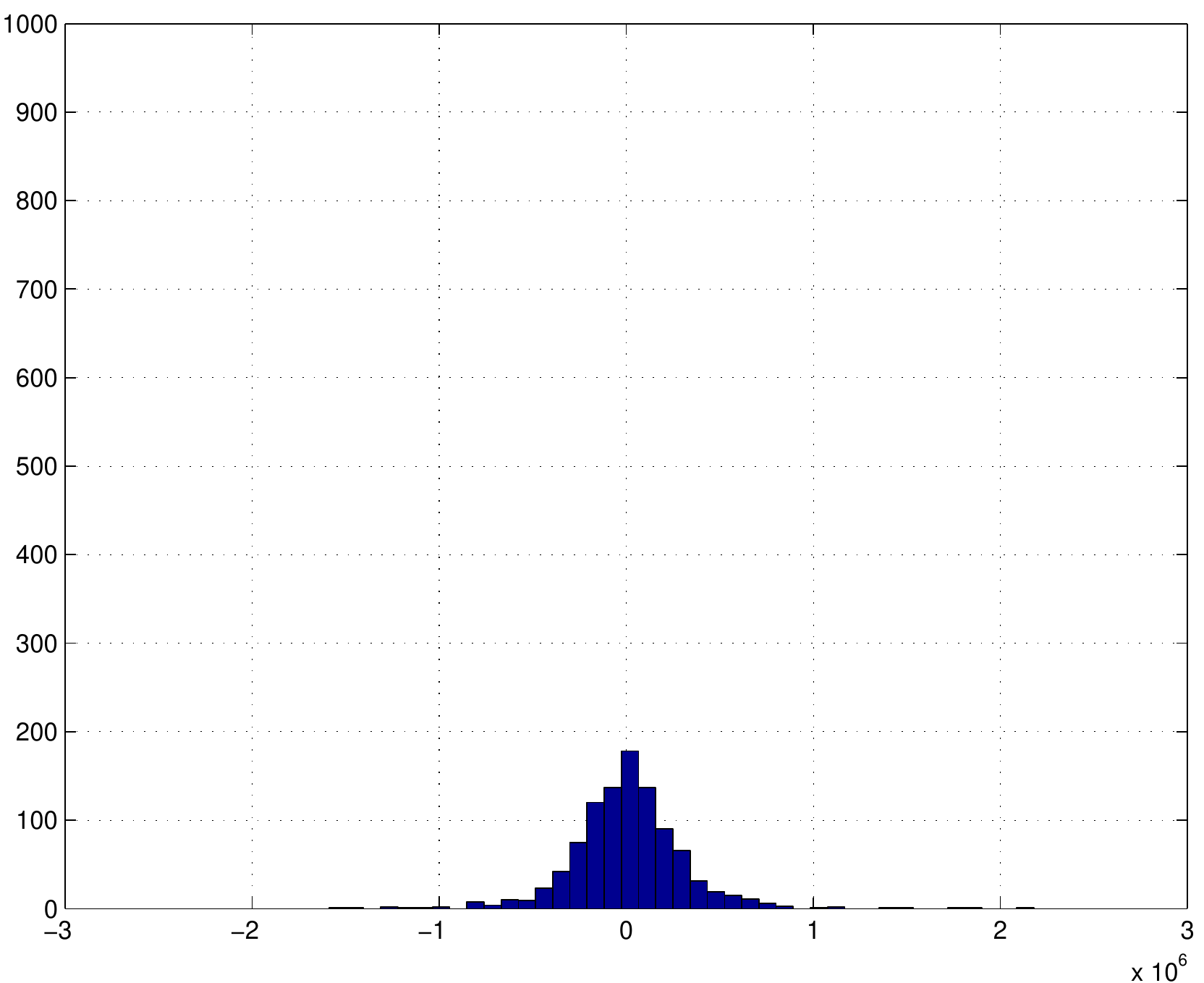}
\includegraphics[scale=0.2,angle=0]{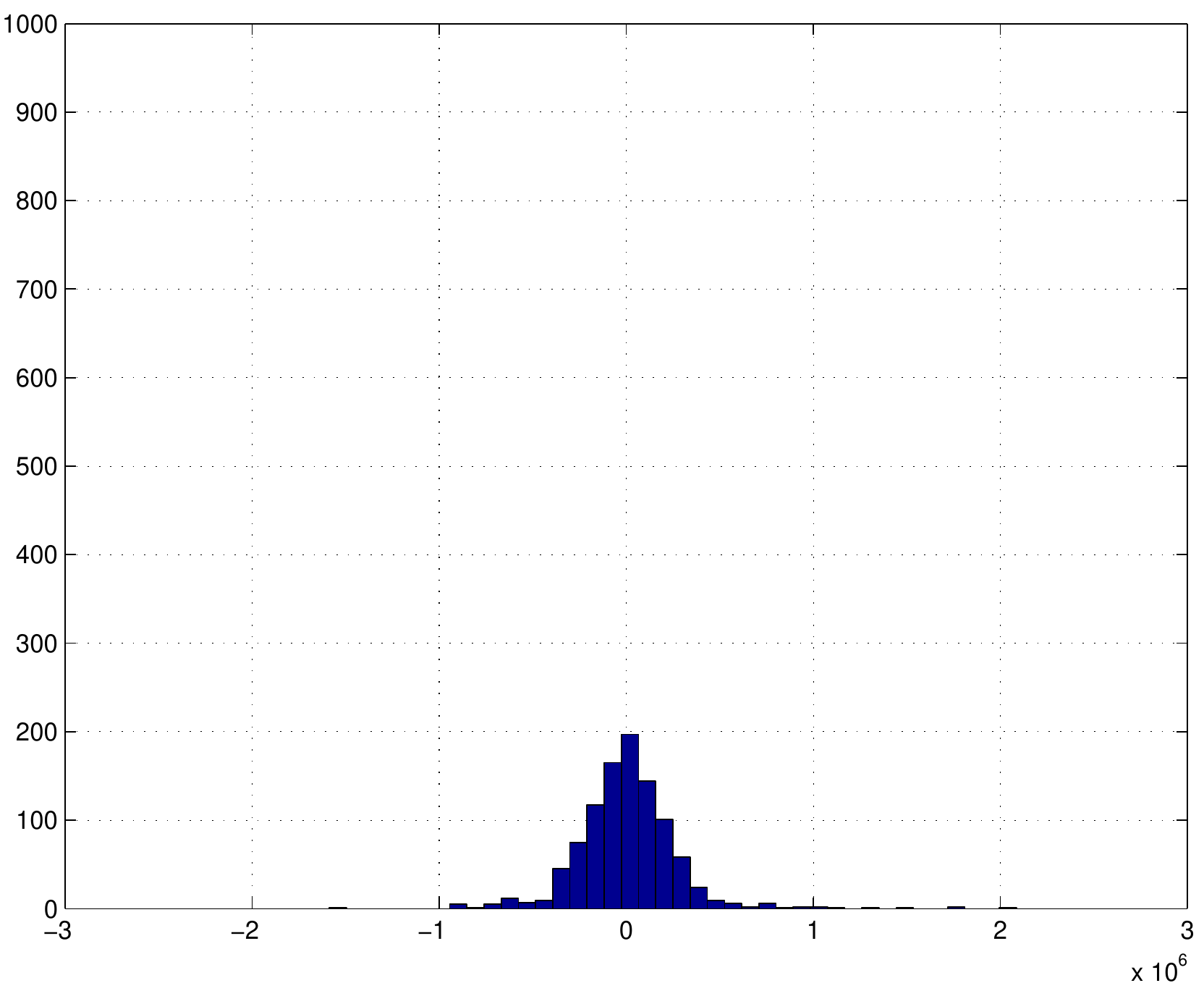}
\includegraphics[scale=0.2,angle=0]{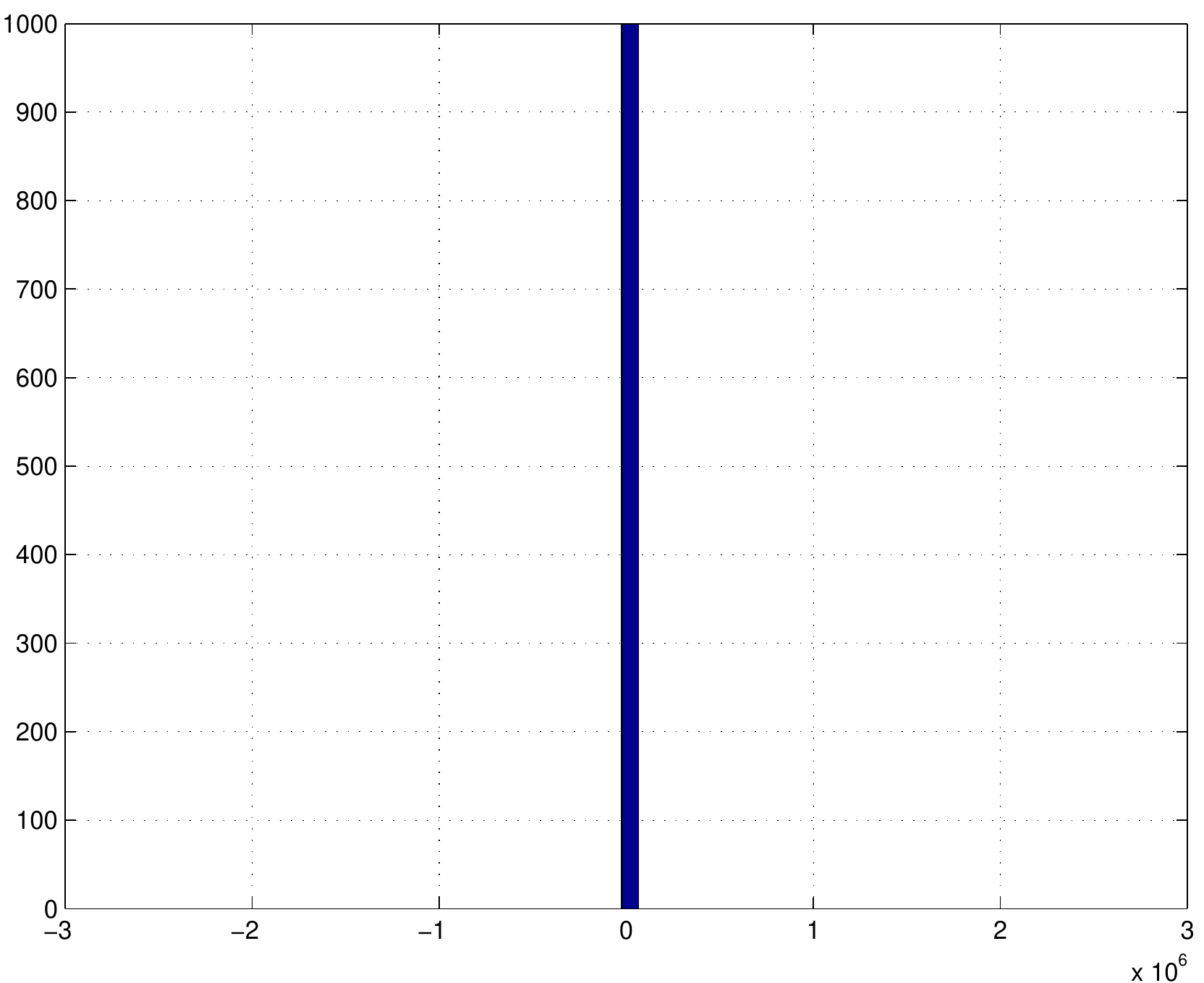}
\includegraphics[scale=0.2,angle=0]{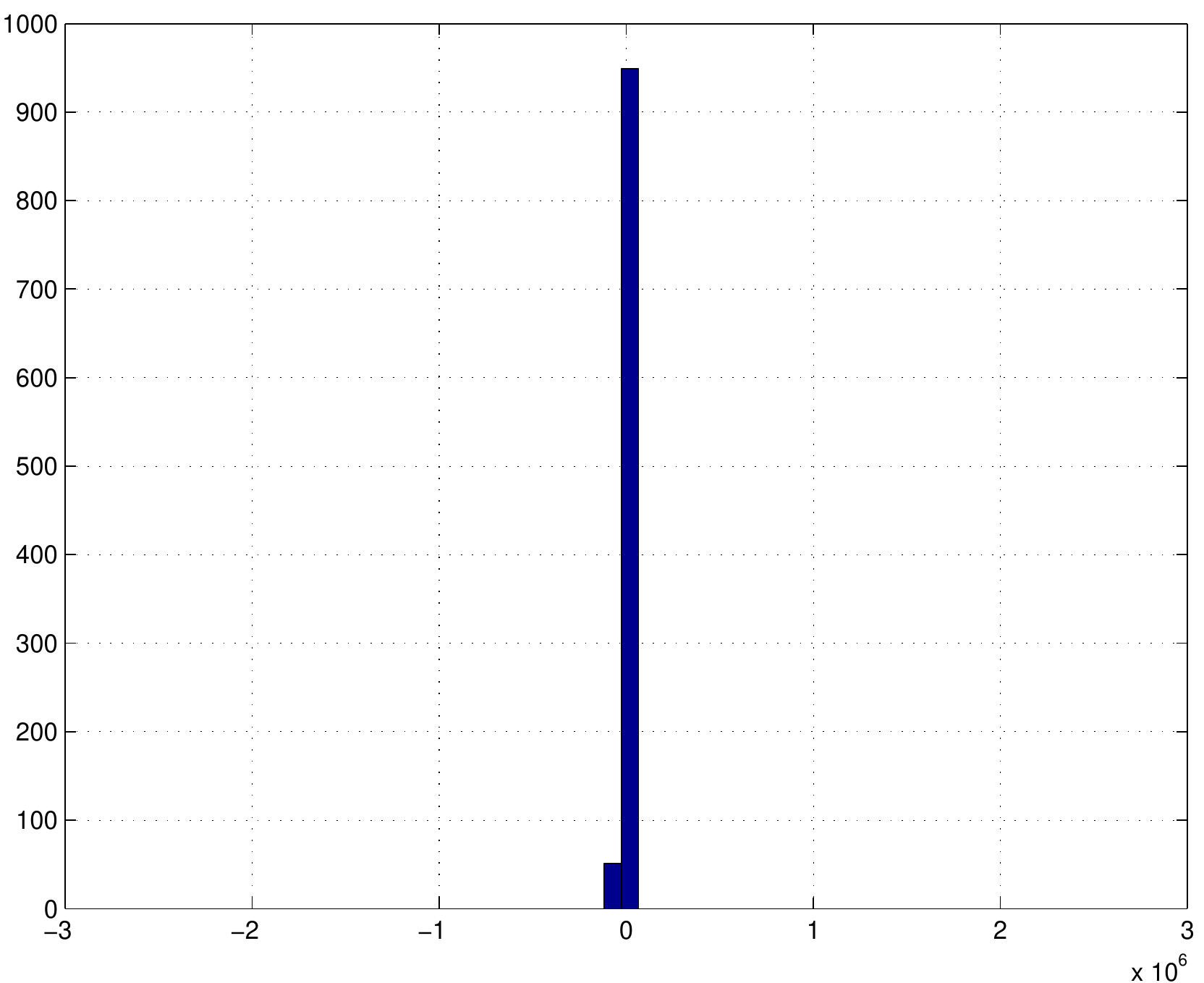}  \\
\hbox{\hspace{0.15in}   (a)  U-symmetric \hspace{0.3in} (b)
U-asymmetric \hspace{0.3in} (c)  U-flat \hspace{0.3in}    (d)
Control  \hspace{0.3in}    } \caption{\footnotesize {\bf Sliding Disk at Uniform
Temperature, $h=0.01$, $\alpha=5.0$, $c=0.1$.} From left: histogram
of the $x$-displacement of the disk at $T=50000$ for $U$ symmetric,
asymmetric, and flat.  We observe a wider spread in the cases when
$U$ is symmetric or asymmetric.    } \label{fig:histX}
\end{center}
\end{figure}

We will consider four different systems to simulate.  The first
three are sliding disks with a symmetric, asymmetric, and flat
potentials:
\begin{center}
\begin{tabular}{c|c}
$U(x)$  &  \\
\hline
 $\sin(x)$ & symmetric \\
$\sin(x) + 0.4 \sin(2 x)$  & asymmetric \\
$0$ & flat  
\end{tabular}
\end{center}
The fourth case is a control consisting of a simulation of a 1-D
Langevin process at the same temperature and with a
symmetric potential:
\begin{align*}
dX &= V dt \text{,} \\
d V &= -c V dt - \cos(X) dt+ \alpha dB_V \text{.}
\end{align*}
For all of the simulations the {\em sliding disk is initially at rest},
$r=0.25$, $m=1.0$, $J=m r^2/2$ and the friction and noise factors
are as indicated in the figures.   The means of the $x$ displacement
of the disk as shown in Fig.~\ref{fig:meanX} are very small compared
with the spread as, e.g., shown in the histogram of the final
position of the ball as shown in Fig.~\ref{fig:histX}.   However,
the mean squared displacement shows ballistic diffusion in the cases
when $U$ is symmetric or asymmetric and normal diffusion otherwise
(see Figures \ref{fig:meanX2}-\ref{fig:loglogX2}).   The time-scale
associated with this ballistic diffusion is plotted in
Fig.~\ref{fig:corX} which shows the correlation in the
x-displacement when $U$ is symmetric or asymmetric.   This
time-scale is much greater than the characteristic time-scale
associated with the friction or the potential. Recall from the integral expression of the velocity, that when $U$ is zero a rigid-body term is neglected.  This demonstrates the
important role of the rigid-body effect in the ballistic diffusion
of the $x$-displacement when $U$ is symmetric and asymmetric.

Finally we also consider adding a non-degenerate, but anisotropic
dissipation matrix to the sliding disk.  It is numerically observed that if the 
anisotropy is large enough, $\E( x_t^2)$ is ballistic for a long period
of time.

\begin{figure}[htbp]
\begin{center}
\includegraphics[scale=0.3,angle=0]{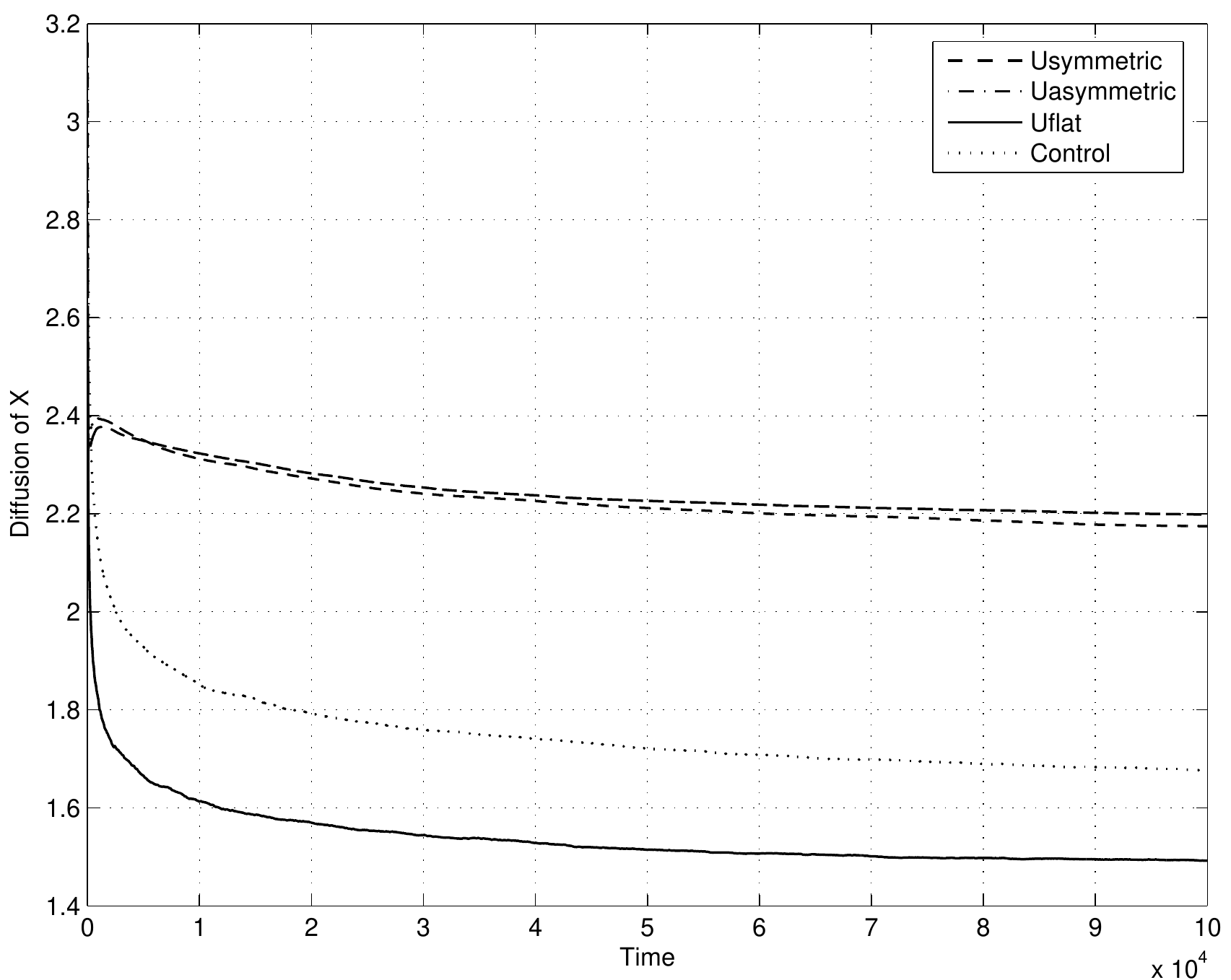} \\
\includegraphics[scale=0.2,angle=0]{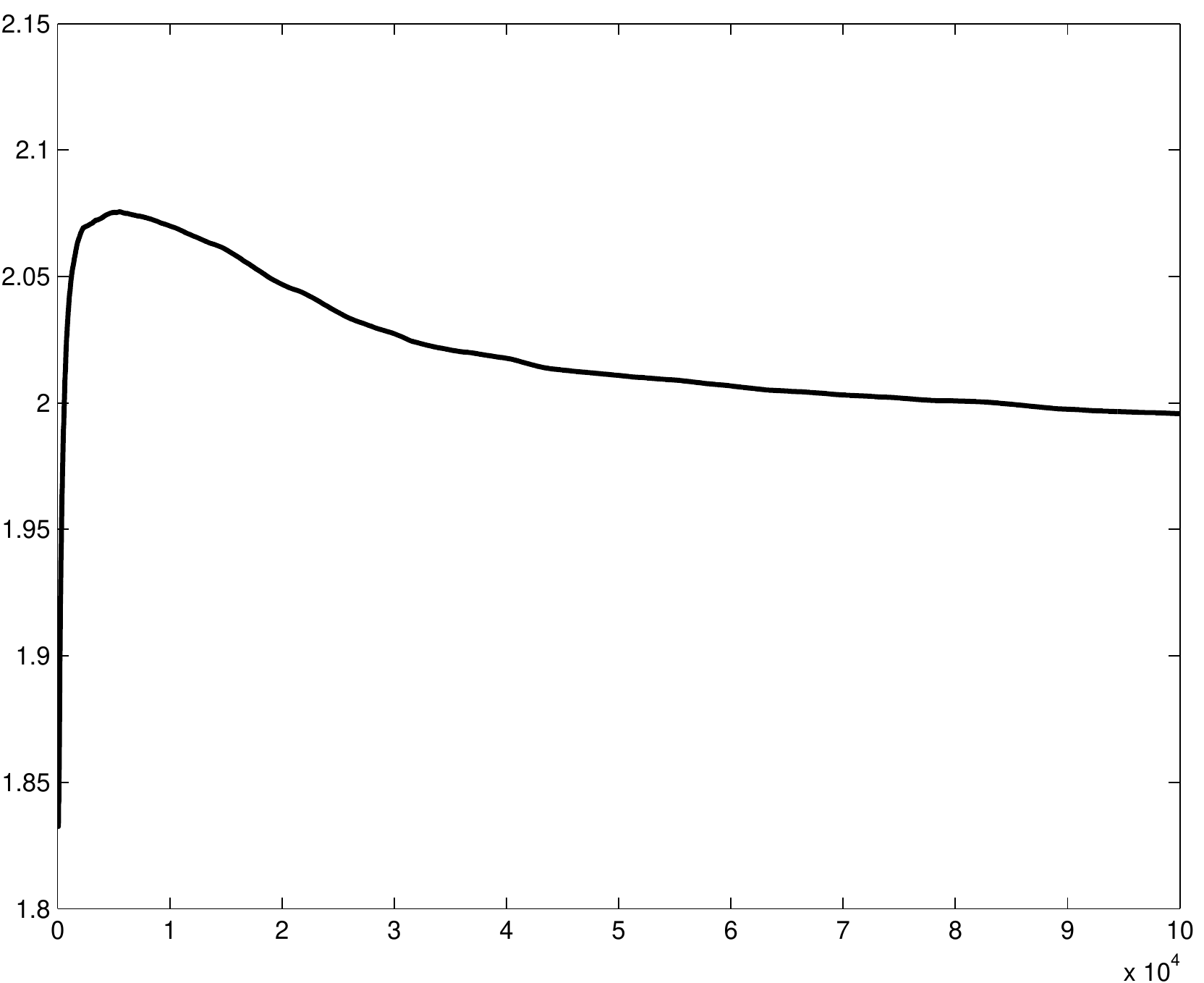}
\includegraphics[scale=0.2,angle=0]{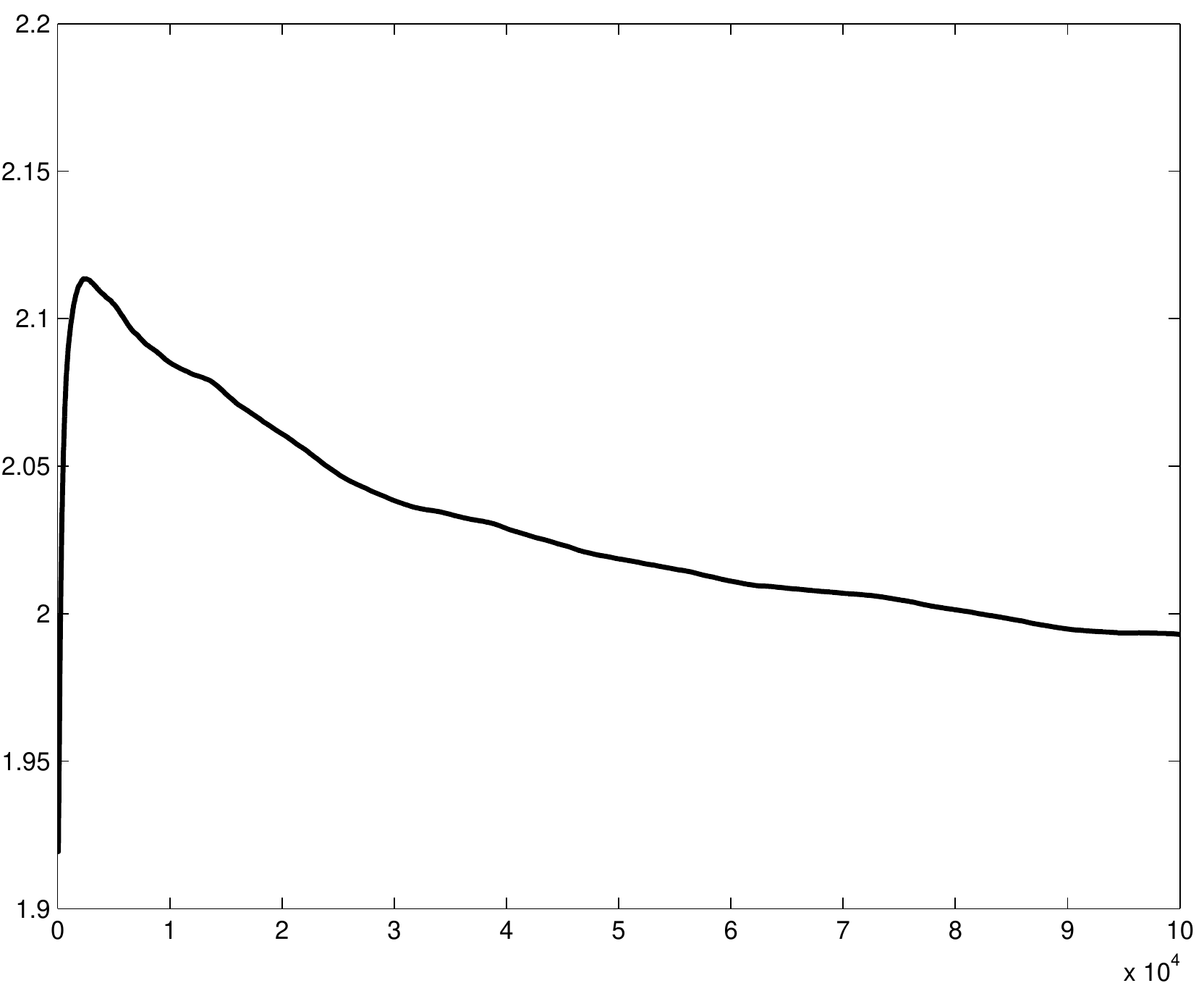}
\includegraphics[scale=0.2,angle=0]{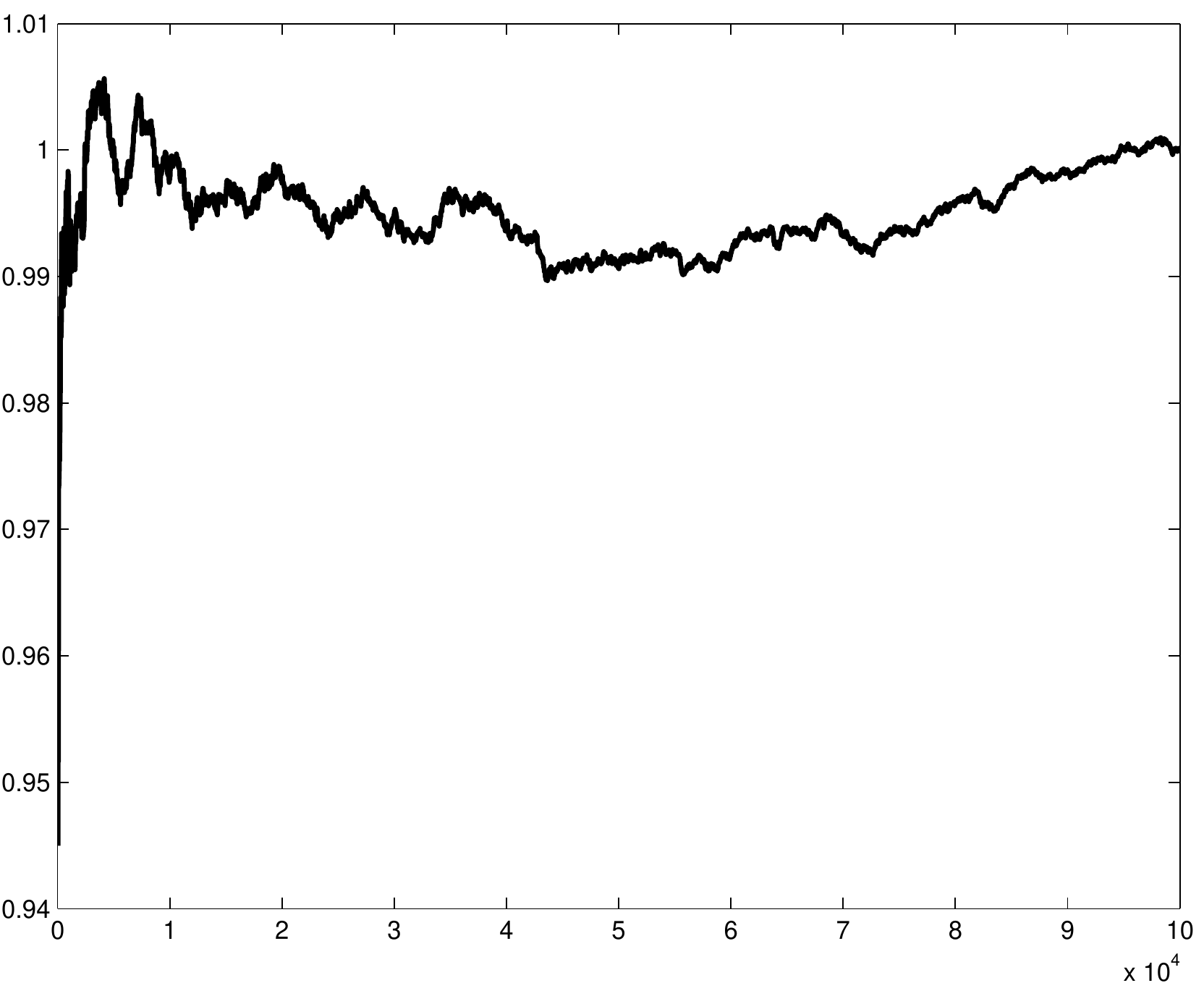}
\includegraphics[scale=0.2,angle=0]{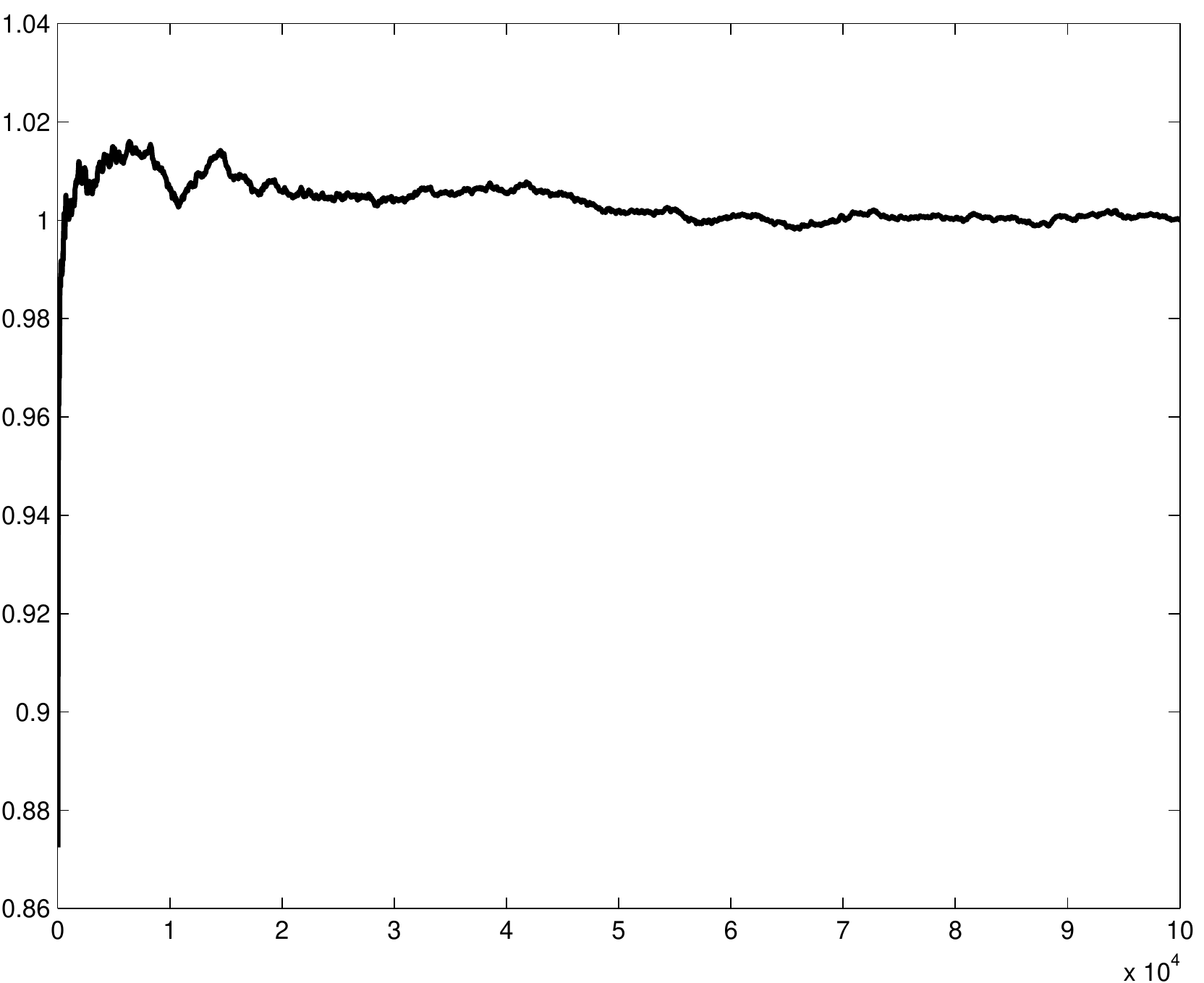} \\
\hbox{\hspace{0.15in}   (a)  U-symmetric \hspace{0.3in} (b)
U-asymmetric \hspace{0.3in} (c)  U-flat \hspace{0.3in}    (d)
Control  \hspace{0.3in}    } \caption{\footnotesize {\bf Sliding Disk at Uniform
Temperature, $h=0.01$, $\alpha=5.0$, $c=0.1$.} From left:
``diffusion'' of $x$-displacement of the disk for $U$ symmetric,
asymmetric, flat, and the control.  For the cases when $U$ is flat
and the control, the diffusion is normal.  Whereas in the other
cases the diffusion is ballistic.  The figure on the top superposes
these graphs in a single plot for comparison. } \label{fig:diffusionX}
\end{center}
\end{figure}

\newpage

\section{Hamel's magnetic top}\label{HamelsDevice}

In this section an observed magnetism induced spinning phenomenon is analyzed.   
Up to our knowledge, the first magnetism induced spinning device was  ``Tesla's egg of colombus'' exhibited in 1892 \citep{Ge1919}.  Using two phase AC energizing coils in quadrature, Tesla placed a copper plated ellipsoid on a wooden plate above a
rotating magnetic field. The egg stood on its pointed end without
cracking its shell and began to spin at high speed to the amazement
of the scientists who witnessed the experiment. This effect was
caused by induced eddy currents on the surface of the ellipsoid. A
related magnetism induced spinning phenomenon is the Einstein-de Hass 
effect in which the rotation of an object is caused
by a change in magnetization \citep{Fr1979}. In this effect the
magnetic field causes an alignment of electronic spins and their
angular momenta is transferred to the atomistic lattice.
We will show below  through an idealized model based on
magnetostatics that  the spinning of the top in Hamel's device is 
due to surface friction.

\paragraph{Mechanism behind curious rotation}

The following observation was made on a simple mechanical system consisting of a magnetized top and ring as shown in Fig.~\ref{fig:moviesnapshot}.   When the ring is held above the top within a certain range of heights, and then tilted, one observes the 
top transitions from a state of no spin to a state of nonzero spin about its axis of symmetry.

\begin{figure}[htbp]
\begin{center}
\includegraphics[scale=0.25,angle=0]{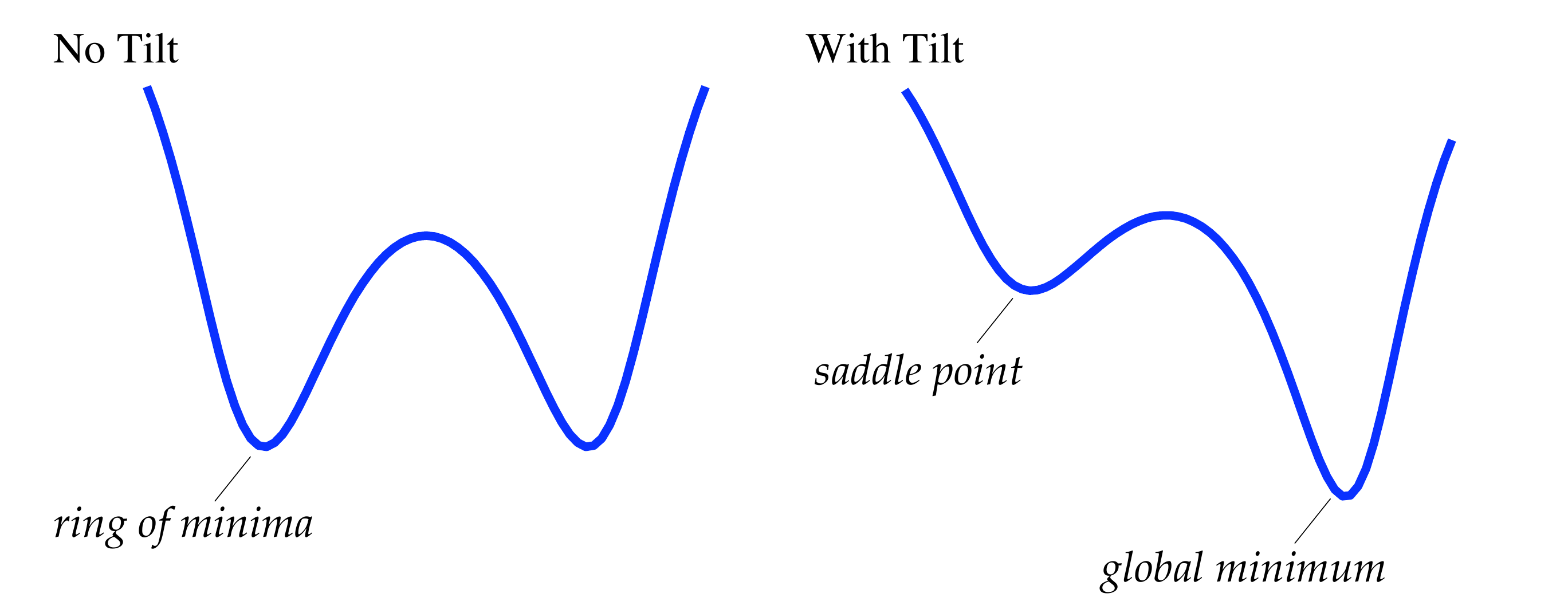}
\caption{\footnotesize  {\bf Cross-Section of Potential Energy
Surface.} Cross-sectional sketches of the magnetic potential energy
with and without tilt.  With no tilt in the ring, the potential
energy surface is axisymmetric and has a ring of minima.  If the
ring is tilted slightly, a unique local minimum exists across
from a saddle point. } \label{fig:potentialcrossection}
\end{center}
\end{figure}

The system is modeled as two rigid bodies in magnetostatic interaction: a 
ball with a magnetic dipole aligned to one of its axes, and a fixed magnetized ring of radially aligned magnetic dipoles.  The main tool used to analyze the observed curious rotation is a Lagrange-Dirichlet stability criterion \citep{MaRa1999}.  It is
shown that the fixed points of the system's governing equations
correspond to the magnetic top being at rest with its axis of
symmetry aligned with the local magnetic field and its translational position at a critical point of the magnetic potential energy.  Stability of this point is determined by analyzing the nature of this critical point.  If the attitude of the ring is normal to the surface, the magnetic potential energy is very nearly axisymmetric.   A cross-sectional
sketch of this potential energy is shown in Fig.~\ref{fig:potentialcrossection}.   In this case there exists a ring of minima that are not individually stable.  However, if the ring is tilted
slightly, the potential has a unique local minimum opposite a saddle point as
shown in the same figure.

\begin{figure}[htbp]
\begin{center}
\includegraphics[scale=0.25,angle=0]{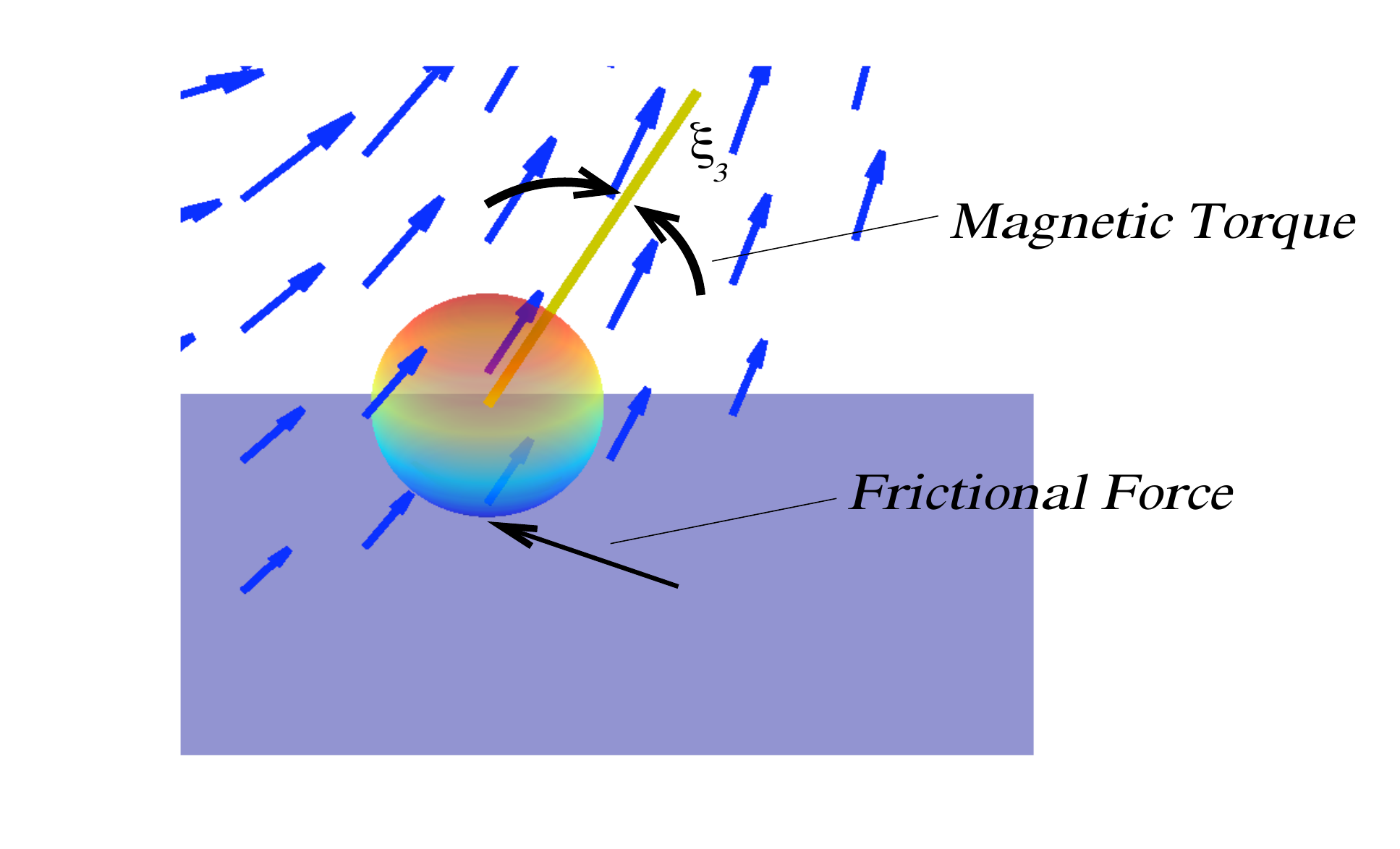}
\caption{\footnotesize  {\bf Mechanism behind spin.} A sketch of the
ball, its axis of symmetry  $\boldsymbol{\xi}_3$, the ambient
magnetic field (blue arrows), the restoring magnetic torques, and 
frictional force.  The ball tends to a position that minimizes its
potential energy. If the ball is initially unstable, then as it
moves towards the minimum a frictional force resists this motion.
The frictional force causes torques about axes orthogonal to the
moment arm $\mathbf{q}$.  However, the restoring magnetic torque
counters the torques about the axes orthogonal to
$\boldsymbol{\xi}_3$ as shown in the sketch.   This argument
clarifies why the spin is primarily along $\boldsymbol{\xi}_3$. }
\label{fig:whyspin}
\end{center}
\end{figure}

If the ring is tilted and the ball's initial position is unstable,
the ball moves towards the stable fixed point which induces a
resisting frictional force (see Fig.~\ref{fig:whyspin}). The
torque due to the sliding friction is in directions orthogonal
to the moment arm $\mathbf{q}$.  However, the torque due to the
magnetic field counters the torques about the axes perpendicular
to $\boldsymbol{\xi}_3$.  This magnetic torque keeps
$\boldsymbol{\xi}_3$ aligned with the local magnetic  field. Thus,
the torque due to friction mainly  causes a spin about the
$\boldsymbol{\xi}_3$ axis.

Thus far, the ring has been kept fixed.  If the ring is perturbed
slightly the position of the local minimum will change.  Hence this
system is unstable with respect to perturbations of the ring.  This
instability will be utilized to design a prototype stochastic
mechanicial rectifier.

\begin{figure}[htbp]
\begin{center}
\includegraphics[scale=0.35,angle=0]{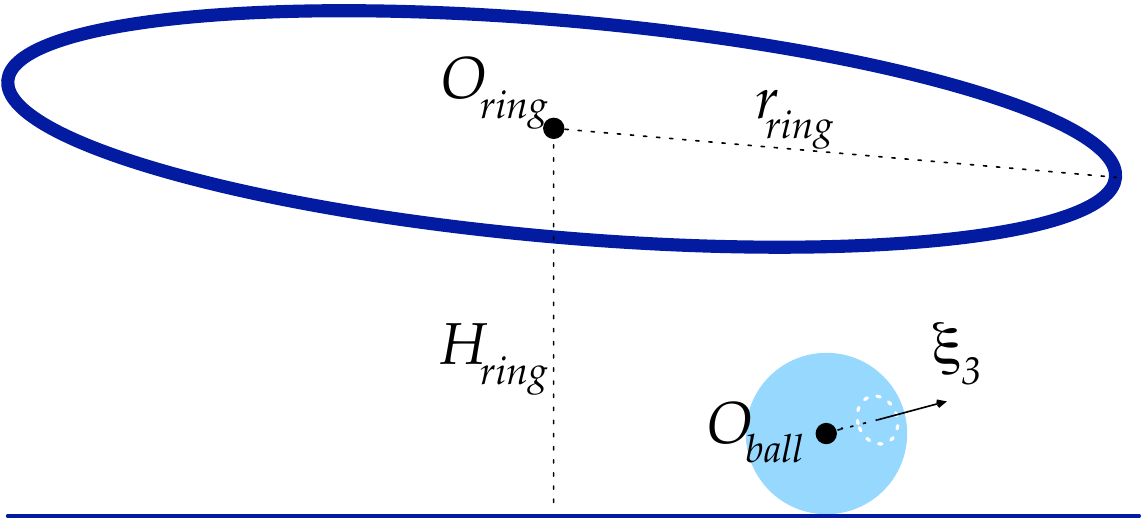}
\caption{\footnotesize  {\bf Illustration of Magnetized Rigid Ball
and Ring.} The figure depicts a rigid ball on a flat surface with a
dipole fixed at the ball's centroid $O_{ball}$ and in the direction
$\boldsymbol{\xi}_3$.   We will assume that the center of mass of
the ball is coincident with $O_{ball}$.  The ring is located at a
height $H_{ring}$ above a reference point $O$ on the surface, and
its radius is $r_{ring}$.   The magnetic dipoles are placed
symmetrically around this ring with dipole moments in the radial
direction. } \label{fig:magball}
\end{center}
\end{figure}

Hamel's magnetic top is modeled as an axisymmetric rigid ball of radius $r$ and mass $m$.   A magnetic dipole is attached to the
center of mass, and in the  direction of the axis of symmetry
$\boldsymbol{\xi}_3 \in \mathbb{R}^3$ as shown in
Fig.~\ref{fig:magball}. For simplicity, the surface friction is
modelled using a sliding friction law.   In what follows the
continuous and discrete model are derived using the 
Hamilton-Pontryagin (HP) variational principle \citep{BoMa2007a}.

\paragraph{Magnetostatic Field} Let $\mathbf{B} : \mathbb{R}^3 \to \mathbb{R}^3$ be the ambient magnetostatic field.  For simplicity, we assume the magnetostatic
field is due to a magnetic ring consisting of $N$ magnetic dipoles
equally spaced around a ring of radius $r_{ring}$ centered at the
point $O_{ring}$ and in the plane defined by the vector
$\boldsymbol{\zeta}_3 \in \mathbb{R}^3$.   The point $O_{ring}$ is
at a height $H_{ring}$ above a reference point $O$ on the surface of
the ball.  To the reference point $O$, an orthonormal frame
$(\mathbf{e}_x, \mathbf{e}_y, \mathbf{e}_z)$ is attached. For
$i=1,...,N$, let $\mathbf{d}_i \in \mathbb{R}^3$ denote the location
of the ith-dipole with respect to $O_{ring}$ and $\mathbf{m}_i \in
\mathbb{R}^3$ denote its dipole moment.   The dipole moments are
assumed to be aligned in the radial direction of the ring.  Let
$\mu$ be the permeability of free space. The magnetic field of the
ith-dipole at a  field point $\mathbf{x} \in \mathbb{R}^3$  can be
computed from  the vector $\mathbf{r}_i$ connecting the field point
to the ith-dipole given by
\[
\mathbf{r}_{i} = \mathbf{x} - r_{ring} \mathbf{d}_i - H_{ring}
\mathbf{e}_z
\]
using the standard formula for the magnetostatic field due to a
dipole:
\begin{equation}
\mathbf{B}(\mathbf{r}_i, \mathbf{m}_i) =  \frac{\mu}{4 \pi \|
\mathbf{r}_i \|^5} \left( 3 (\mathbf{r}_i^\mathrm{T} \mathbf{m}_i)
\mathbf{r}_i -  \| \mathbf{r}_i \|^2  \mathbf{m}_i \right) \text{.}
\label{eq:magneticfield}
\end{equation}
Using the principle of superposition, the magnetic field at a field
point $\mathbf{x} \in \mathbb{R}^3$ due to the $N$ dipoles is
determined as the vector sum of the magnetic fields due to each
dipole:
\[
\mathbf{B}( \mathbf{x} ) = \sum_{i=1}^{N} \mathbf{B}(
\mathbf{r}_i, \mathbf{m}_i) \text{.}
\]

\begin{figure}[htbp]
\begin{center}
\includegraphics[scale=0.2,angle=0]{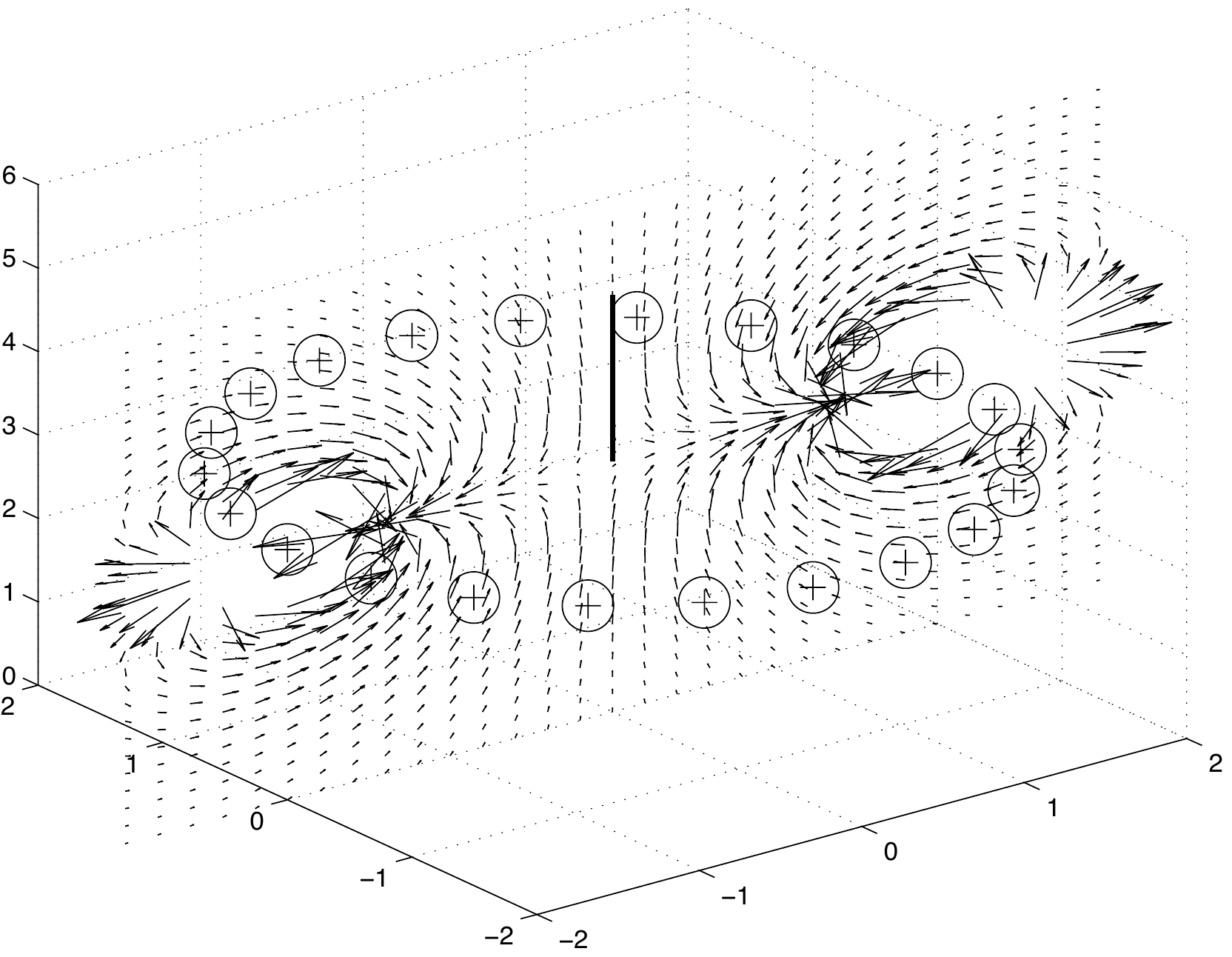}
\includegraphics[scale=0.2,angle=0]{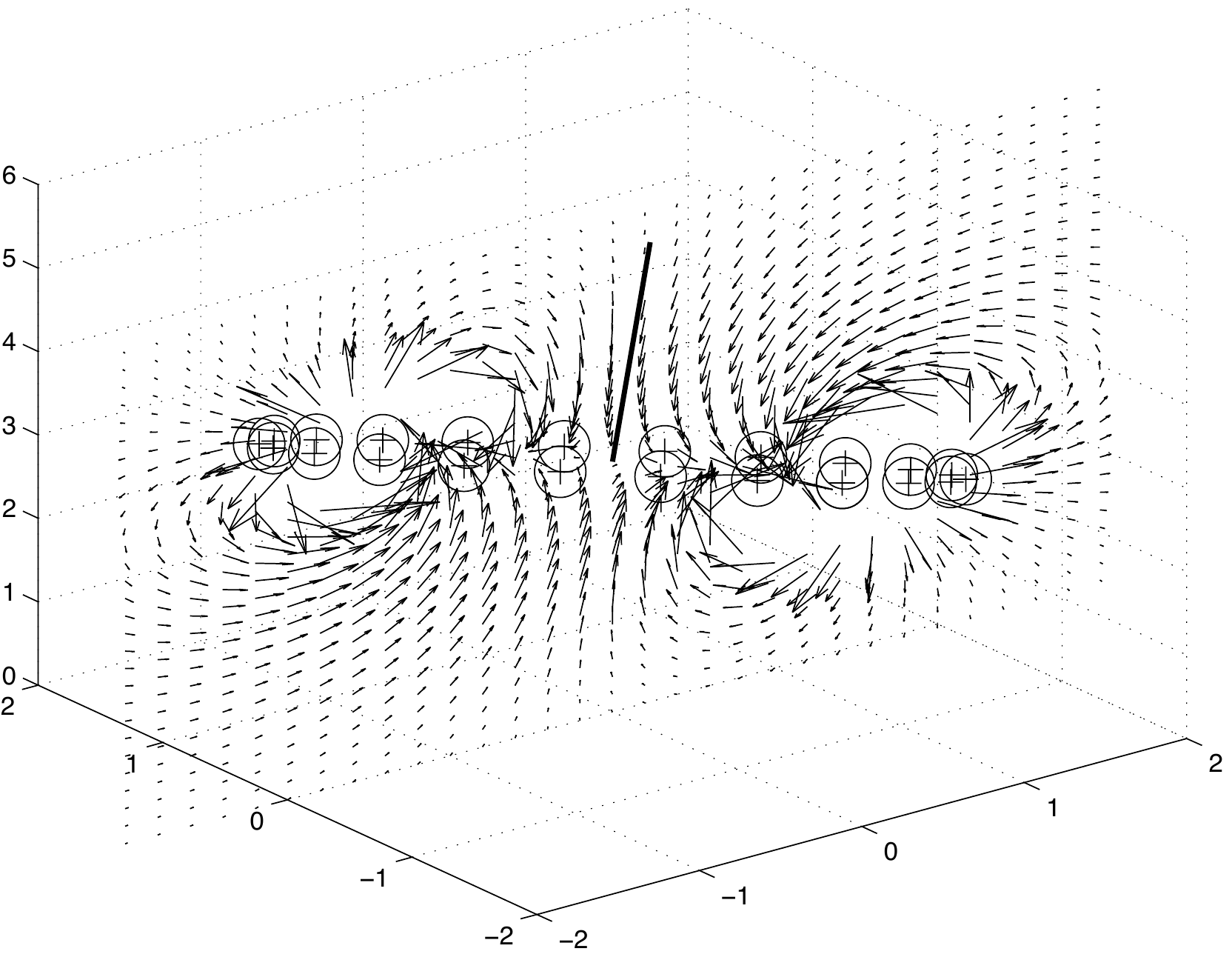}
\includegraphics[scale=0.2,angle=0]{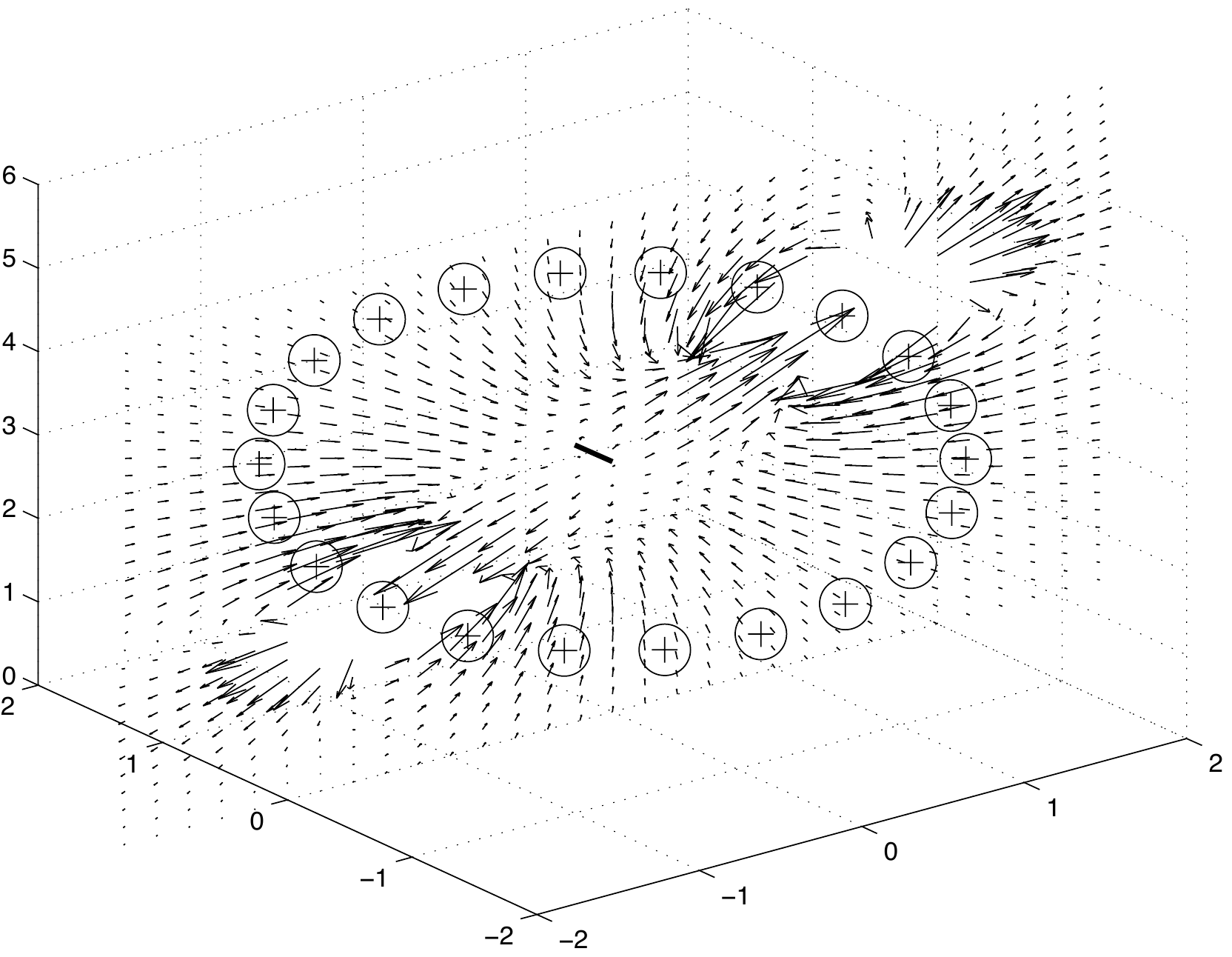}
\hbox{\hspace{0.1in}   $\begin{array}{c} \text{(a) axisymmetric} \\
\phi_R=0, \theta_R=0 \end{array}$ \hspace{0.3in}  $\begin{array}{c}
\text{(b) tilted} \\ \phi_R=\pi/4, \theta_R=0 \end{array}$
\hspace{0.3in}   $\begin{array}{c} \text{(c) tilted} \\
\phi_R=\pi/4, \theta_R=\pi \end{array}$  } \caption{\footnotesize
{\bf Illustration of Magnetized Ring and its Magnetic Field.} The
figures depict various orientations of a magnetized ring consisting
of a  discrete number of dipoles equally spaced around a circle in
the plane defined by $\boldsymbol{\zeta}_3$ which is also depicted.
The magnetic fields lines in the plane $y=0$ are plotted.   The
angles $\phi_R$ and $\theta_R$ shown are the altitude and azimuth of
the normal to the plane in which the ring is in.  The height of the
center of the ring, the radius of the ring, the number of dipoles,
and the orientation of the ring are all parameters in this
magnetostatic field model.    } \label{fig:magfield}
\end{center}
\end{figure}

\paragraph{Lagrangian of Magnetic Ball}

The configuration space of the system is $Q = \mathbb{R}^3 \times
\operatorname{SO}(3)$ and its Lagrangian is denoted $L: TQ \to
\mathbb{R}$. Let $(\mathbf{x}(t), \mathbf{\dot{x}}(t)) \in
\mathbb{R}^3 \times \mathbb{R}^3$ denote the translational position
and velocity of the ball measured with respect to the reference
frame attached to $O$.  Let $(\mathbf{e}_1, \mathbf{e}_2,
\mathbf{e}_3)$ denote an inertial orthonormal frame attached to
$O_{ball}$ and related to a body-fixed frame  $(\boldsymbol{\xi}_1,
\boldsymbol{\xi}_2, \boldsymbol{\xi}_3)$ via the rotation matrix
$R(t) \in \operatorname{SO}(3) $.   Let $J =
\operatorname{diag}(I_1, I_2, I_3)$ be the standard diagonal inertia
matrix of the body and $ \boldsymbol{\xi}_3$ the axis of symmetry of
the ball.  We assume the mass distribution of the body is  symmetric
with respect to the axis of symmetry $ \boldsymbol{\xi}_3$.  This
assumption implies that the principal moments of inertia about
$\boldsymbol{\xi}_1$ and $\boldsymbol{\xi}_2$ are equal, i.e.,
$I=I_1=I_2$.  Let $(\mathbf{x}(t), \mathbf{\dot{x}}(t)) \in
\mathbb{R}^3 \times \mathbb{R}^3$ denote the translational position
and velocity of the ball.

We will use the isomorphism between $\mathbb{R}^3$ and the Lie
algebra of $\operatorname{SO}(3)$, $\mathfrak{so}(3)$, given by the
hat map (See appendix.).  In terms of this identification, we define
the reduced Lagrangian $\ell: \mathbb{R}^3 \times \mathbb{R}^3
\times  \operatorname{SO}(3) \times \mathbb{R}^3 \to \mathbb{R}$ as
\[
\ell(\mathbf{x}, \mathbf{\dot{x}}, R, \boldsymbol{\omega})  =
L(\mathbf{x}, \mathbf{\dot{x}}, R, \widehat{\boldsymbol{\omega}} R)
\text{.}
\]
For the free magnetic ball, i.e., magnetic ball without dissipation,
it is given explicitly by,
\begin{equation}
\ell(\mathbf{x}, \mathbf{\dot{x}}, R, \boldsymbol{\omega}) =
\frac{m}{2} \mathbf{\dot{x}}^\mathrm{T}  \mathbf{\dot{x}} +
\frac{1}{2}  \boldsymbol{\omega}^\mathrm{T} R J R^\mathrm{T}
\boldsymbol{\omega}  - mg \mathbf{x}^\mathrm{T} \mathbf{e}_3 +
\boldsymbol{\xi}_3^\mathrm{T} \mathbf{B}(\mathbf{x}) \text{.}
\label{eq:ell}
\end{equation}
From left the terms represent the translational and rotational
kinetic energy of the ball, the gravitational potential energy and
the dipole potential energy. The ball is also subject to the surface
constraint  $\varphi: \mathbb{R}^3 \to \mathbb{R}$ given by:
\begin{align}
\varphi( \mathbf{x}) &= -r  +  \mathbf{e}_3^\mathrm{T} \mathbf{x}
\text{.} \label{eq:constraint}
\end{align}
This holonomic constraint restricts the translational motion of the 
ball to a plane.

\paragraph{Governing Conservative Equations}

The equations of motion will be determined using a HP description \citep{BoMa2007a}.   The constrained HP action integral is given by,
\[
s = \int_a^b \left[ \ell(\mathbf{x}, \mathbf{v}, R,
\boldsymbol{\omega}) +\left\langle  \mathbf{p}, \mathbf{\dot{x}} -
\mathbf{v} \right\rangle +\left\langle  \widehat{\boldsymbol{\pi}},
\dot{R} R^\mathrm{T}  - \widehat{\boldsymbol{\omega}} \right\rangle
+ \lambda \varphi( \mathbf{x} )  \right] dt \text{.}
\]
The HP principle states that
 \[
\delta s = 0
\]
where the variations are arbitrary except that the endpoints
$(\mathbf{x}(a),R(a))$ and  $(\mathbf{x}(b),R(b))$ are held fixed.
The equations are given by,
 \begin{align}
\mathbf{\dot{x}}  = \mathbf{v}~~~ & \text{(reconstruction equation),} \\
\mathbf{p} = \frac{ \partial \ell }{ \partial \mathbf{v} }  ~~~& \text{(Legendre transform),}  \\
 \frac{d}{dt} \mathbf{p} = \frac{\partial \ell}{\partial \mathbf{x}} + \lambda \frac{ \partial \varphi }{ \partial \mathbf{x} }  \label{eq:eulerlagrange}  ~~~&\text{(Euler-Lagrange equations),} \\
 \varphi( \mathbf{x} ) = 0 \label{eq:surfaceconstraint}  ~~~&\text{(constraint equation),} \\
  \frac{d}{dt}    R =  \widehat{\boldsymbol{\omega}}  R \label{eq:reconstruction} ~~~&\text{(reconstruction equation),} \\
      \boldsymbol{\pi} = \frac{\partial \ell}{\partial  \boldsymbol{\omega}} \label{eq:legendre} ~~~&\text{(reduced Legendre transform),} \\
    \frac{d}{dt} \widehat{\boldsymbol{\pi}}  =  \frac{\partial \ell}{\partial R} R^{\mathrm{T}} -  \widehat{ \widehat{\boldsymbol{\pi}} \boldsymbol{\omega} } \label{eq:lp} ~~~&\text{(Lie-Poisson equations).}
 \end{align}
 Evaluating these equations at $\ell$ as defined in (\ref{eq:ell}) yields
 \begin{equation} \label{eq:magtopc}
 \begin{cases}
 \begin{array}{rcl}
 \mathbf{\dot{x}} &=& \mathbf{v} \\
m \mathbf{\dot{v}} &=&   \left( \mathbf{D} \mathbf{B} (\mathbf{x}) \right)^{\mathrm{T}}  \boldsymbol{\xi}_3  + (\lambda - mg) \mathbf{e}_3  \\
 \mathbf{x}^{\mathrm{T}} \mathbf{e}_3 &=&  r  \\
      \dot{R} &= &  \widehat{\boldsymbol{\omega}} R \\
      \boldsymbol{\pi} &= & R J R^{\mathrm{T}} \boldsymbol{\omega} \\
\dot{\boldsymbol{\pi}}  &= &  \widehat{\boldsymbol{\xi}_3}
\mathbf{B}(\mathbf{x})   \text{.}
\end{array}
\end{cases}
 \end{equation}

 \begin{Remark}
 Since the system is axisymmetric the Legendre transform in (\ref{eq:magtopc}) simplifies:
 \begin{equation}
\boldsymbol{\pi} =  I \boldsymbol{\omega} + (I_3 - I)  (
\boldsymbol{\omega}^{\mathrm{T}}  \boldsymbol{\xi}_3 )
\boldsymbol{\xi}_3  \implies \boldsymbol{\omega}  = \frac{1}{I}
\left( \boldsymbol{\pi} + \frac{I - I_3}{I_3} (
\boldsymbol{\pi}^{\mathrm{T}}  \boldsymbol{\xi}_3 )
\boldsymbol{\xi}_3 \right) \label{eq:spatiallt} \text{.}
 \end{equation}
As a consequence one does not need to solve for the
evolution of all three columns of $R(t)$ to integrate the ODE in $\boldsymbol{\pi}$.  Instead one just needs to solve for the evolution of the third column,
$\boldsymbol{\xi}_3$, using:
 \begin{equation}
\dot{\boldsymbol{\xi}}_3 = \widehat{ \boldsymbol{\omega} }
\boldsymbol{\xi}_3 \label{eq:xi3} \text{.}
 \end{equation} 
 \end{Remark}

The following conservation law follows from axisymmetry.

 \begin{proposition} The following momentum map is conserved under the flow of (\ref{eq:magtopc}),
 \[
 J = \boldsymbol{\pi}^{\mathrm{T}} \boldsymbol{\xi}_3 \text{.}
 \]
\end{proposition}

\begin{proof}
This momentum map is due to an $S^1$ symmetry of the Lagrangian
about the axis $\boldsymbol{\xi}_3$ which can be
computed by formula (12.2.1) of  \citep{MaRa1999}.   The group
acts on $Q$ by:
\[
\Phi^{Q}_s(\mathbf{x},R) = ( \mathbf{x},  \exp(s
\widehat{\boldsymbol{\xi}_3}) R)  \text{.}
\]
The corresponding infinitesimal generator is given by:
\[
\psi^{Q}(\mathbf{x}, R)  = \frac{d}{ds} \left.
\Phi^{Q}_s(\mathbf{x}, R)  \right|_{s=0} \text{.}
\]
Thus, one can invoke Noether's theorem to conclude the following
momentum map $J: TQ \to (S^1)^*$ is preserved
\[
J =\left\langle \frac{\partial L}{\partial \dot{R}},
\psi^{Q}(\mathbf{x}, R)  \right\rangle= \frac{\partial
\ell}{\partial \boldsymbol{\omega}}^{\mathrm{T}} \boldsymbol{\xi}_3
=   \boldsymbol{\pi}^{\mathrm{T}} \boldsymbol{\xi}_3 \text{.}
\]
\end{proof}

This conservation law indicates that if initially the top is not
spinning about its axis of symmetry, then it will never spin about
this axis.  Thus, one cannot obtain the curious rotation using
magnetic and gravitational effects alone.  This result suggests
surface friction will also play a role in producing this phenomenon.

\paragraph{Governing Nonconservative Equations}

Let $\mathbf{q} = - r \mathbf{e}_3$ denote the vector connecting the center of mass $C$ to the contact point $Q$ as shown in Fig.~\ref{fig:magball}. We model the surface frictional force using
a sliding friction law proportional to the slip velocity, i.e.,  the
velocity of the contact  point on the rigid body relative to the
center of mass $\mathbf{V}_Q$:
\[
\mathbf{F}_f = - c \mathbf{V}_Q
\]
This law assumes zero static friction but, is nevertheless
reasonable on very slippery surfaces where static friction is
negligible. Moreover, it is an experimental fact that the magnetic
top exhibits this curious rotation even on oily surfaces. In reality
one must keep in mind that the ball is not in point-contact  with
the surface; rather, a finite area of the ball is in contact with
the surface and is moving relative to the surface to make spinning
possible. In this case this sliding model of friction is quite
reasonable.  A more refined model of friction would include
rotational torque and dry frictional effects (Coulomb friction).

The slip velocity is given by:
\begin{equation}
\mathbf{V}_Q = \dot{\mathbf{x}} + \widehat{ \boldsymbol{\omega} }
\mathbf{q} \text{.} \label{eq:slip}
\end{equation}
The force of friction is therefore,
\[
\mathbf{F}_f = - c \mathbf{V}_Q \text{,}
\]
and the torque due to friction is,
\[
\boldsymbol{\tau}_f = \widehat{\mathbf{q}} \mathbf{F}_f  \text{.}
\]
The governing dynamical equations of the magnetic ball with friction
are given by:
\begin{equation} \label{eq:magtopnc}
\begin{cases}
\begin{array}{rcl}
\mathbf{\dot{x}} &=& \mathbf{v} \text{,} \\
m \mathbf{\dot{v}} &=&   \left( \mathbf{D} \mathbf{B}(\mathbf{x}) \right)^{\mathrm{T}}  \boldsymbol{\xi}_3   + (\lambda - mg) \mathbf{e}_3 - c ( \mathbf{v} + r \widehat{ \mathbf{e}_3 }\boldsymbol{\omega} )    \text{,} \\
 \mathbf{x}^{\mathrm{T}} \mathbf{e}_3 &=& r \text{,} \\
\dot{\boldsymbol{\xi}}_3 &=&  \widehat{\boldsymbol{\omega}} \boldsymbol{\xi}_3  \text{,} \\
 \boldsymbol{\omega}  &=& \frac{1}{I} \left( \boldsymbol{\pi} + \frac{I - I_3}{I_3} ( \boldsymbol{\pi}^{\mathrm{T}}  \boldsymbol{\xi}_3 )  \boldsymbol{\xi}_3 \right)  \text{,} \\
\dot{\boldsymbol{\pi}} &=&   \widehat{\boldsymbol{\xi}_3}
\mathbf{B}(\mathbf{x})  + c r \widehat{\mathbf{e}_3} (
\mathbf{v} + r \widehat{ \mathbf{e}_3 }\boldsymbol{\omega}  )
 \text{.}
\end{array}
\end{cases}
 \end{equation}

\begin{figure}[htbp]
\begin{center}
\includegraphics[scale=0.25,angle=0]{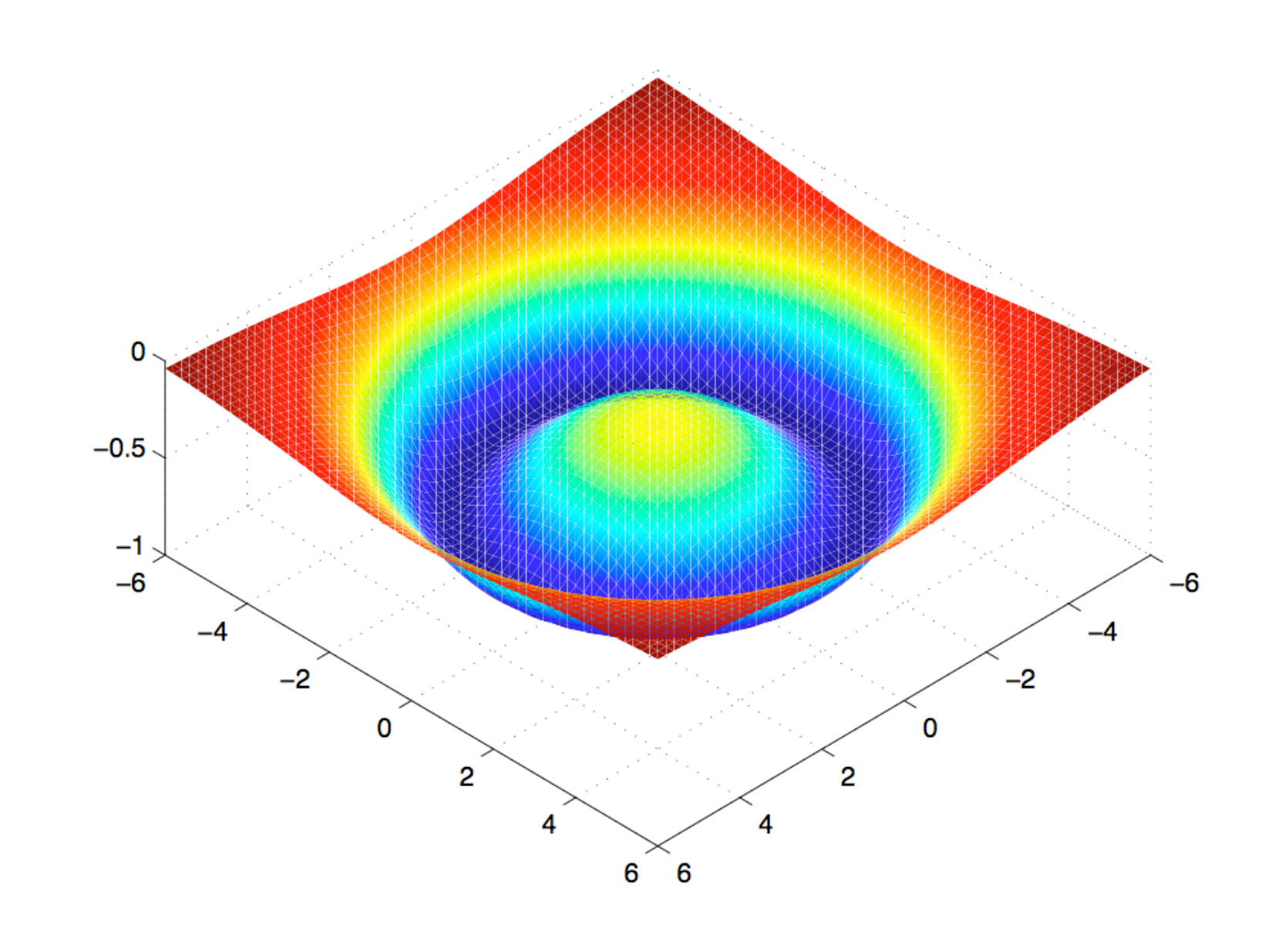}
\includegraphics[scale=0.25,angle=0]{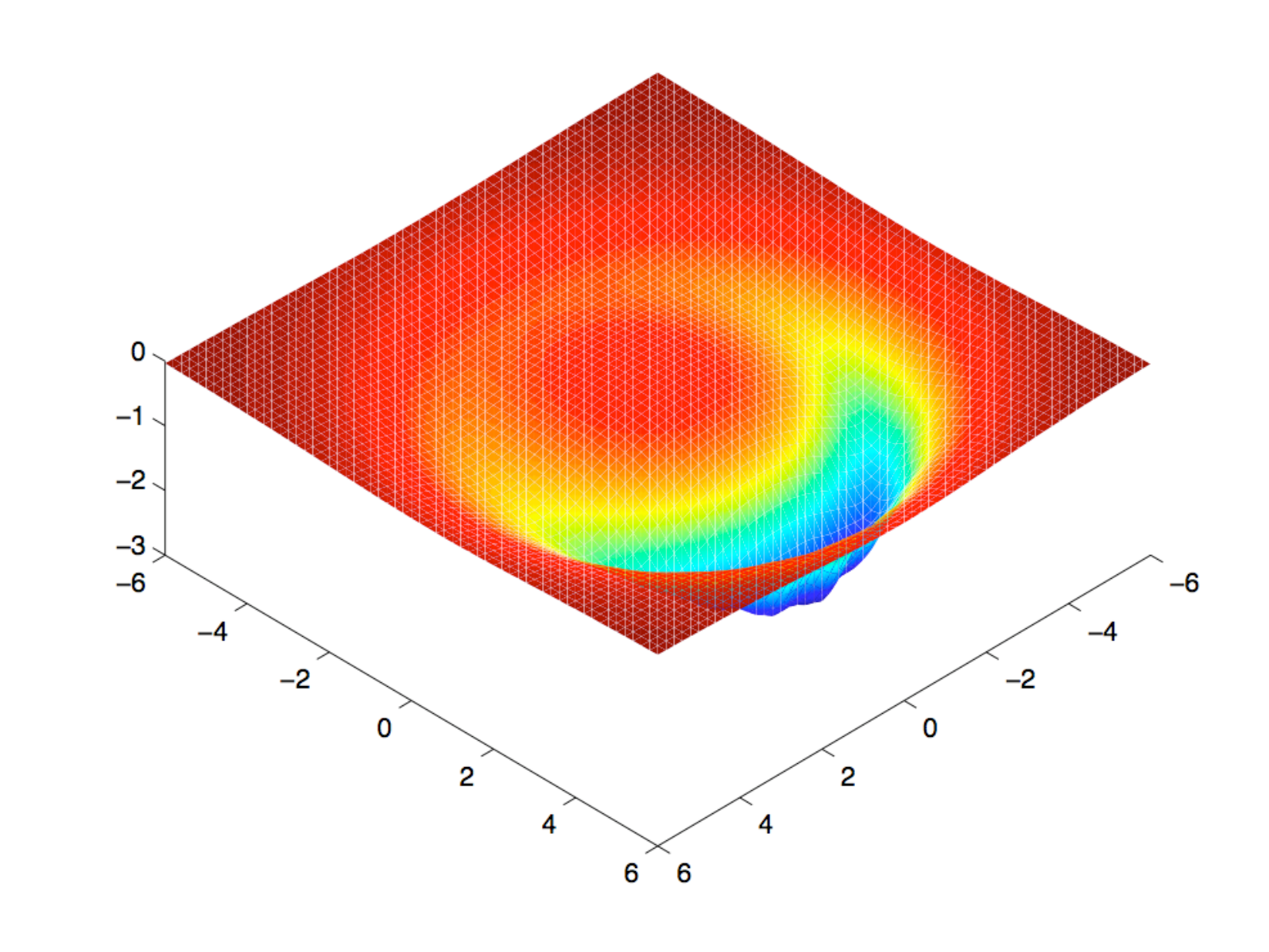}
\includegraphics[scale=0.25,angle=0]{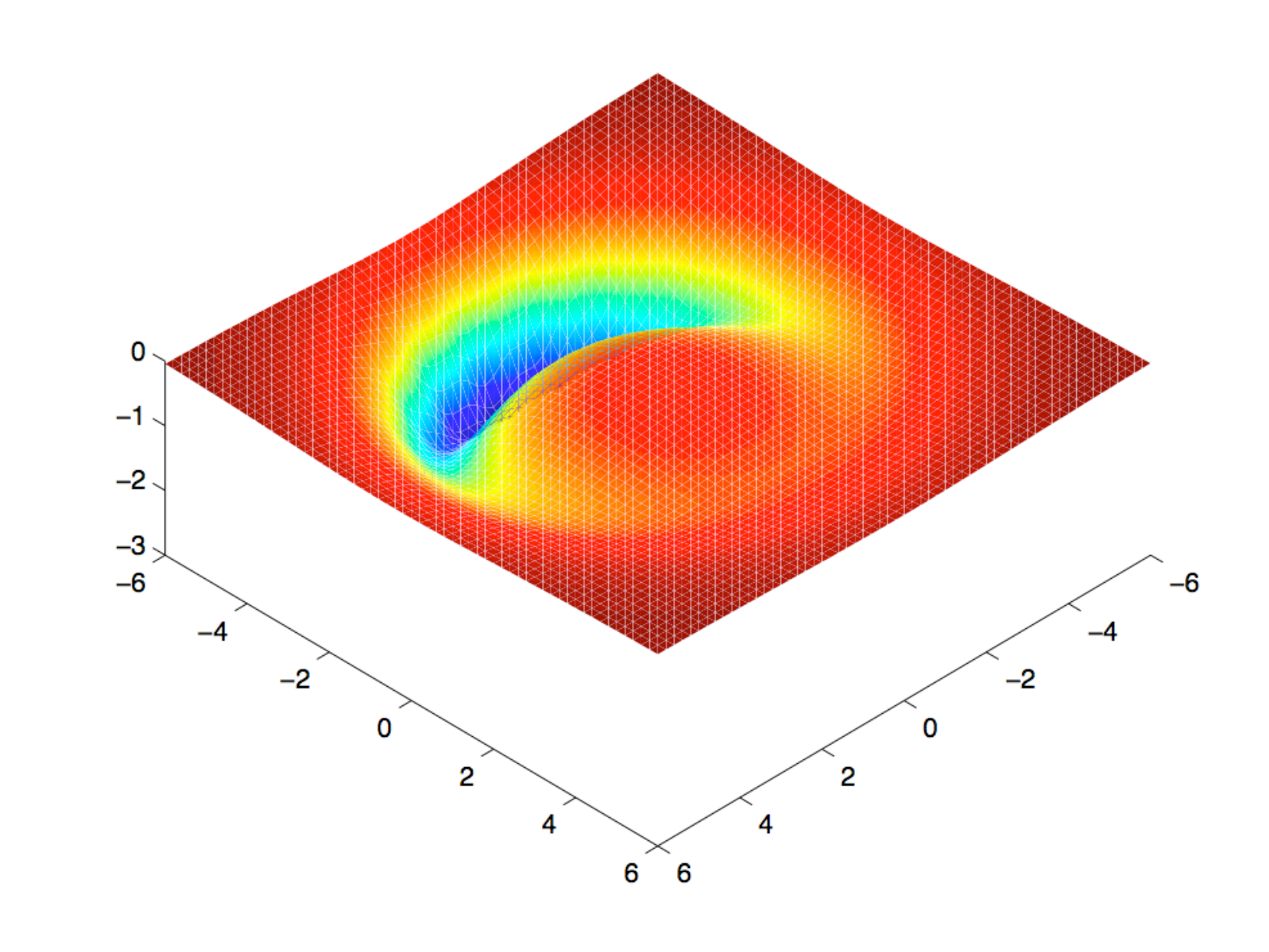}
\hbox{\hspace{0.1in}   $\begin{array}{c} \text{(a) axisymmetric} \\
\phi_R=0, \theta_R=0 \end{array}$ \hspace{0.3in}  $\begin{array}{c}
\text{(b) tilted} \\ \phi_R=\pi/16, \theta_R=0 \end{array}$
\hspace{0.3in}   $\begin{array}{c} \text{(c) tilted} \\
\phi_R=\pi/16, \theta_R=\pi \end{array}$  } \caption{\footnotesize
{\bf Surface plot of $V_e$.}  The potential energy surface is
plotted as a function of $x$ and $y$ for ball radius $r = 0.3$, ring
height $H=2$, $25$ dipoles on the ring, and ring radius
$r_{ring}=4$.   The angles $\phi$ and $\theta$ shown are the
altitude and azimuth of the normal to the plane in which the ring is
in.  By tilting the ring one can obtain a critical point which is a
local minimum of $V_e$.
   }
\label{fig:potentialsurface}
\end{center}
\end{figure}

\begin{Remark} The fixed points of (\ref{eq:magtopnc}) satisfy:
\begin{align*}
\dot{\mathbf{x}} &= 0,~~~\boldsymbol{\omega}=0 \\
\dot{\boldsymbol{\pi}} &= 0 \implies \boldsymbol{\xi}_3 =  \kappa_2 \mathbf{B}(\mathbf{x}),~~~\kappa_2 \in \mathbb{R},~ \text{constant}\text{,} \\
\dot{\mathbf{v}} &= 0 \implies \left( \mathbf{D} \mathbf{B}(\mathbf{x}) \right)^{\mathrm{T}}   \mathbf{B}(\mathbf{x}) =  0 \text{.}
\end{align*}   
The Lagrange-Dirichlet criterion can be used to analyze the stability of these states since the kinetic energy vanishes at these fixed points and since the dissipation is proportional to velocity \citep{MaRa1999}.  In particular, one can invoke a classical theorem due to Thomson-Tait-Chetaev, to conclude that if the fixed point is potentially stable or unstable, then it remains stable or unstable after introducing arbitrary dissipative forces proportional to velocities \citep{Me2002}.  By this criterion, the fixed points are stable provided that the  equilibrium is a strict local minimum of the potential energy.  The potential energy evaluated at this equilibrium point is given by $V_e: \mathbb{R}^3 \to \mathbb{R}$,
 \[
V_e( \mathbf{x} ) = - \kappa_2  \mathbf{B}(\mathbf{x} )^{\mathrm{T}}
\mathbf{B}( \mathbf{x}  )  \text{.}
\]
In the limit as the number of dipoles along the ring is infinite,
$\mathbf{B}$ is axisymmetric with respect to the surface provided
the attitude of the ring is normal to the surface on which the ball
is on.  In this case $V_e$ can be written as a function of the
distance to the origin: $f(\rho) = V_e(\mathbf{x})$ where $\rho=\|
\mathbf{x} \| $ with $\mathbf{x}^{\mathrm{T}} \mathbf{e}_3 = r$.
Assume that the ring is above the ball and the potential energy is
$V_e(\mathbf{x}) \le 0$ for all $\mathbf{x} \in \mathbb{R}^3$. The
origin is an unstable critical point, since $f$ has a local maximum
at that point.   Moreover, $f(\rho)$ tends to zero as $\rho$ becomes
large.  One can pick a distance $M \gg r_{ring}$ (the radius of the
ring) sufficiently large  to make $f(\rho)$ arbitrarily close to
zero.  In the interval $\rho \in (0, M)$, $V_e(\rho)$ is continuous,
and therefore, the function has a local minimum in this interval.
However, since $f$ is axisymmetric this critical point is not a
local minimum in the plane, but rather a circle of critical points.
As shown in Fig.~\ref{fig:potentialsurface}, the ball will be
unstable with respect to circumferential perturbations.  The figure
also shows that when the axisymmetry of the magnetic field is broken
by changing the orientation of the ring, there is a critical point
that is a local minimum of $V_e( \mathbf{x} )$.
\end{Remark}

\paragraph{Simulations}

\begin{figure}[htbp]
\begin{center}
\includegraphics[scale=0.25,angle=0]{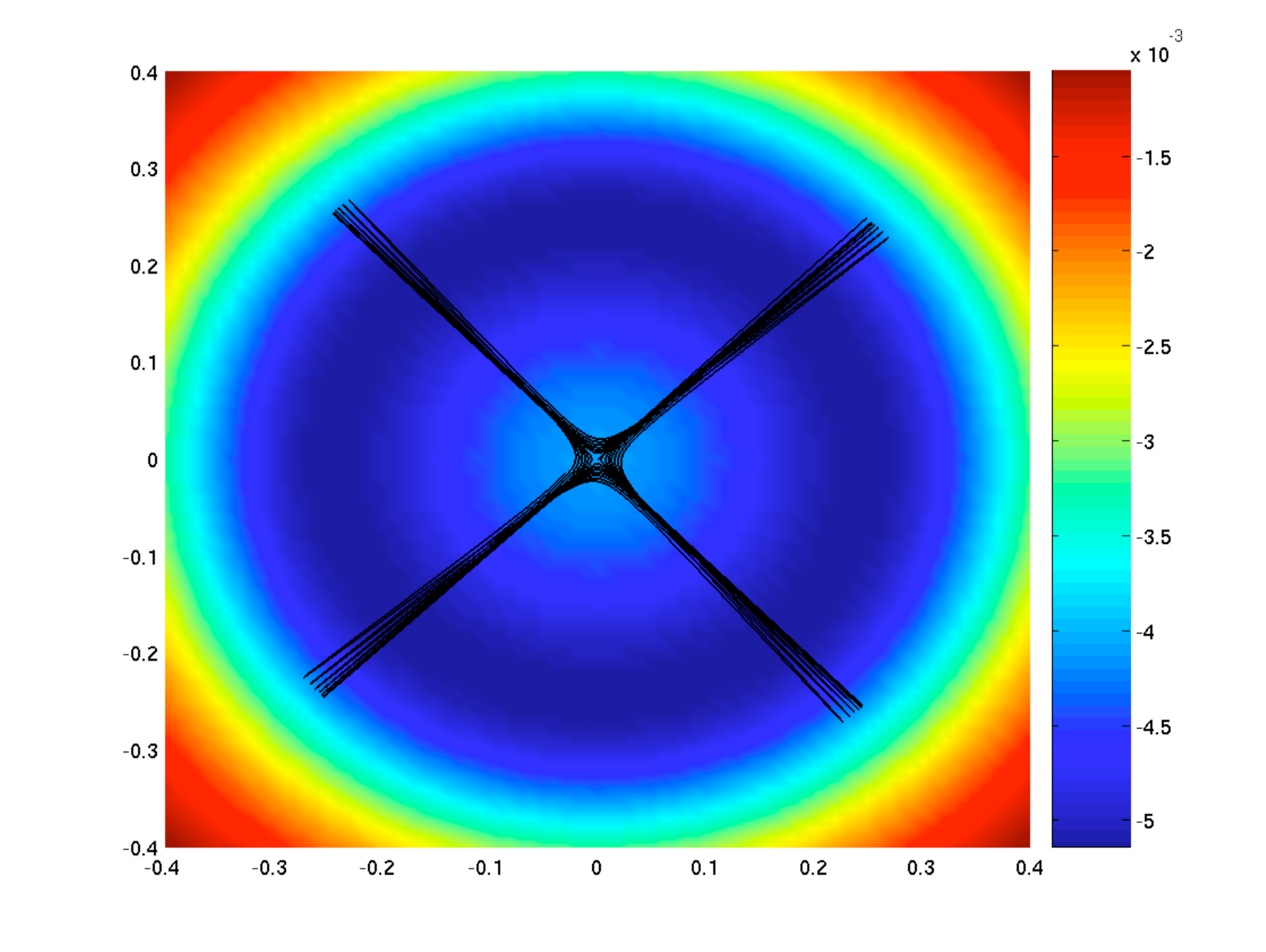}
\includegraphics[scale=0.25,angle=0]{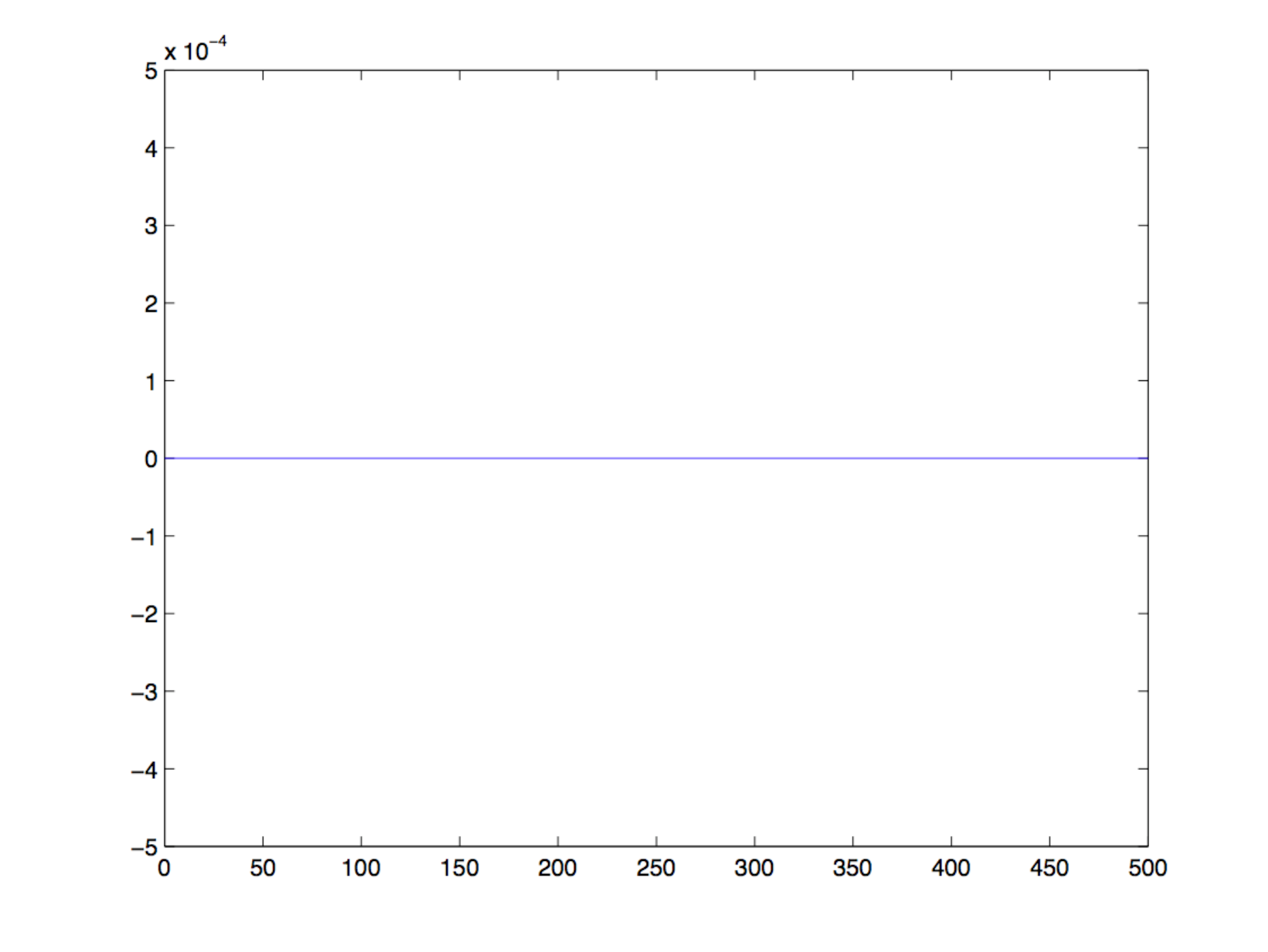}
\hbox{\hspace{0.5in}   (a) potential surface $V_e$   \hspace{0.6in}
(b) momentum map $J(t)$ \hspace{0.3in}    } \caption{\footnotesize
{\bf Conservative Axisymmetric.}  The figures plot $x(t), y(t)$ and
$J(t)$ of the top in the absence of friction, and for $\phi_R=0$ and
$\theta_R=0$.    Superimposed on (a) is the potential energy level
set which appears axisymmetric.  (a) shows that the top bounces
between the local maximum at the center and the circle of minima
(dark blue ring).  (b) shows that $J(t)$ is preserved as predicted
by the theory. The accompanying simulation vividly illustrates this
motion and demonstrates that one does not get curious rotation in
this case. } \label{fig:conaxi}
\end{center}
\end{figure}

\begin{figure}[htbp]
\begin{center}
\includegraphics[scale=0.25,angle=0]{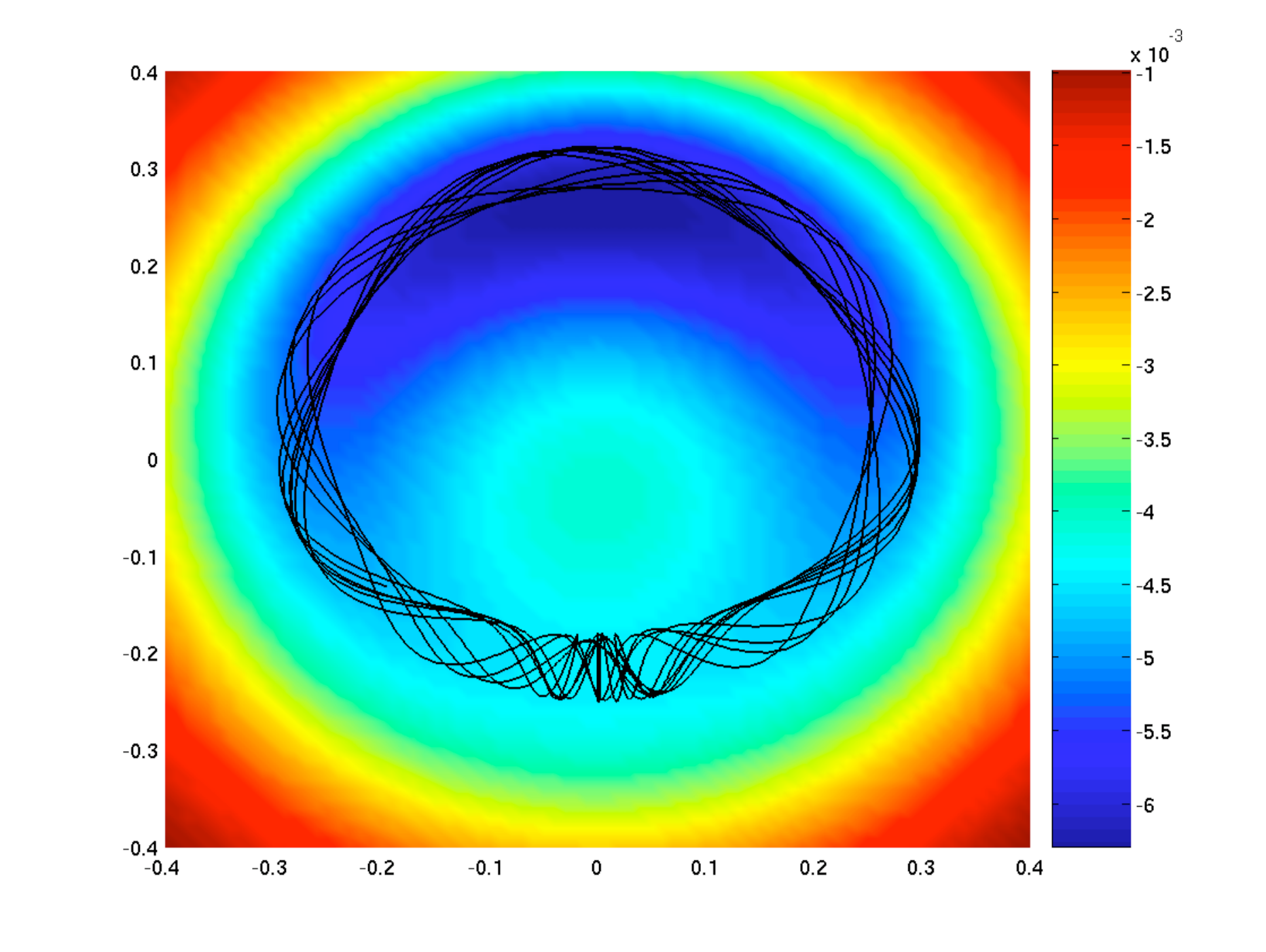}
\includegraphics[scale=0.25,angle=0]{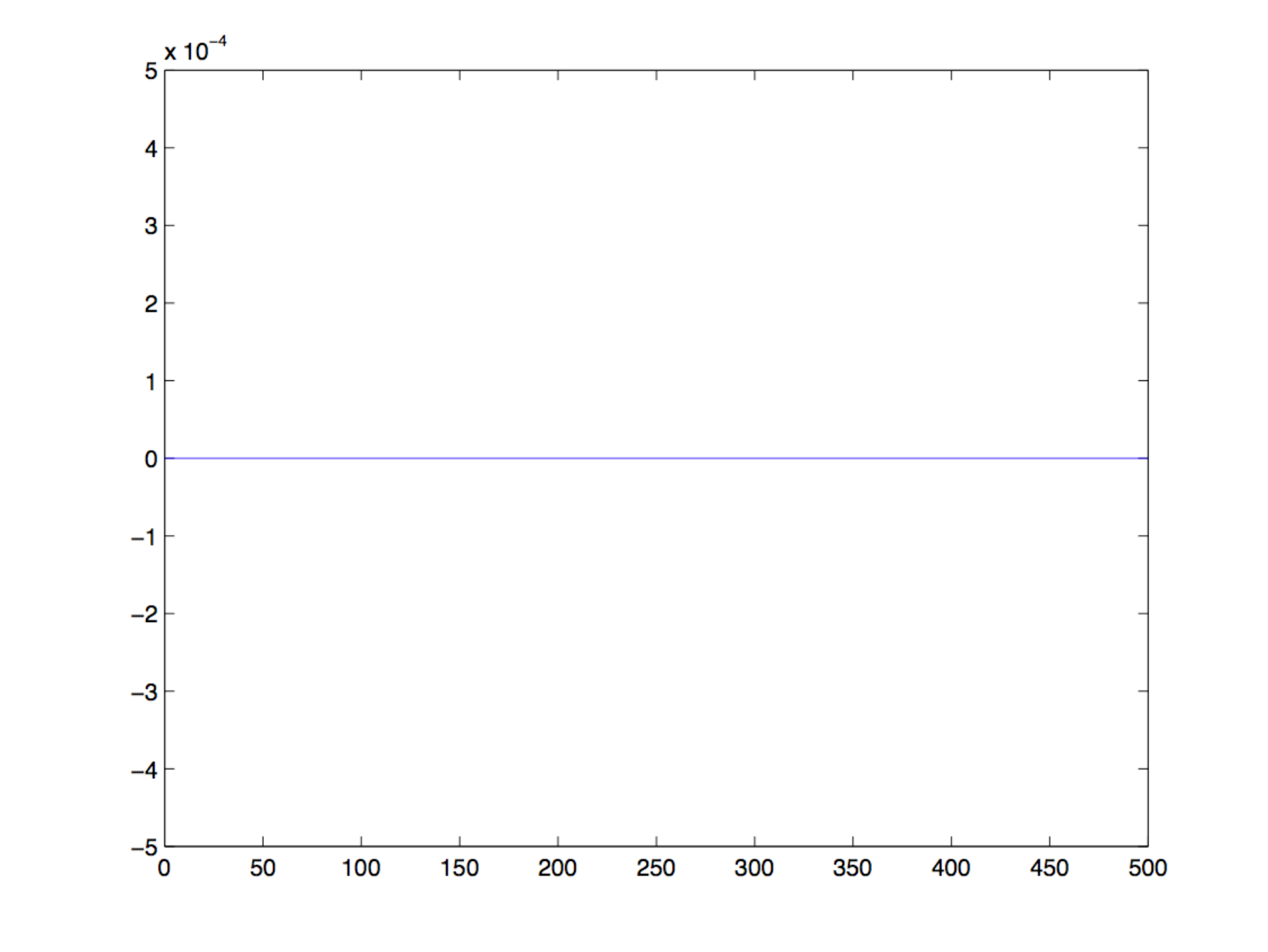}
\hbox{\hspace{0.5in}   (a) potential surface $V_e$   \hspace{0.6in}
(b) momentum map $J(t)$ \hspace{0.3in}    } \caption{\footnotesize
{\bf Conservative Tilted.}   The figures plot $x(t), y(t)$ and
$J(t)$ of the top in the absence of friction, and for $\phi_R=0$ and
$\theta_R=\pi/64$.    Superimposed on (a) is the potential energy
level set which is no longer axisymmetric.  (a) shows that the top
moves erratically within an annulus.  (b) confirms that $J(t)$ is
preserved. The accompanying simulation vividly illustrates this
motion and demonstrates that one does not get curious rotation in
this case. } \label{fig:contilt}
\end{center}
\end{figure}

Four simulations are performed on the top to confirm the theory and
explain the curious rotation of the physical system.   We use the HP 
integrator introduced in \citep{BoMa2007a}.  Let $\phi_B$
and $\theta_B$ denote the azimuth and altitude of the ball's
attitude $\boldsymbol{\xi}_3$. In all of the simulations the top is
initialized as follows:
\begin{align*}
&c= 0.1~\text{kg} \cdot \text{m}/\text{s},
r = 1.8~\text{cm}, m=200 ~\text{g}, \phi_B = 0, \theta_B = 0,  I_3=I_1=2/5 m r^2, \\
 &\mu \| \mathbf{m}_i \| /(4 \pi)  = \mu \| \boldsymbol{\xi}_3 \| /(4 \pi)   = 10^{-5} ~\text{Tesla}~\text{meter}^3,  \\
 & H_{ring} =20 ~\text{cm}, N=20~ \text{dipoles},  r_{ring}=34 ~\text{cm}
\end{align*}
The timestep size and number of timesteps are $h=0.025$ and
$N=20000$ respectively. In the first four simulations the
orientation of the ring is kept fixed.   The first case simulated is
when the magnetostatic field is axisymmetric and friction is absent.
In this case the top undergoes what appears to be chaotic motion
involving a balance between the magnetic potential and translational
kinetic energies as shown in the simulation and
Fig.~\ref{fig:conaxi}.  In the second simulation the altitude of the
attitude of the ring is slightly perturbed causing the magnetostatic
field to be noticeably asymmetric, see Fig.~\ref{fig:contilt}.
However, in both of these cases the top does not acquire spin about
the symmetry axis $\boldsymbol{\xi}_3$ as predicted by the theory.

\begin{figure}[htbp]
\begin{center}
\includegraphics[scale=0.25,angle=0]{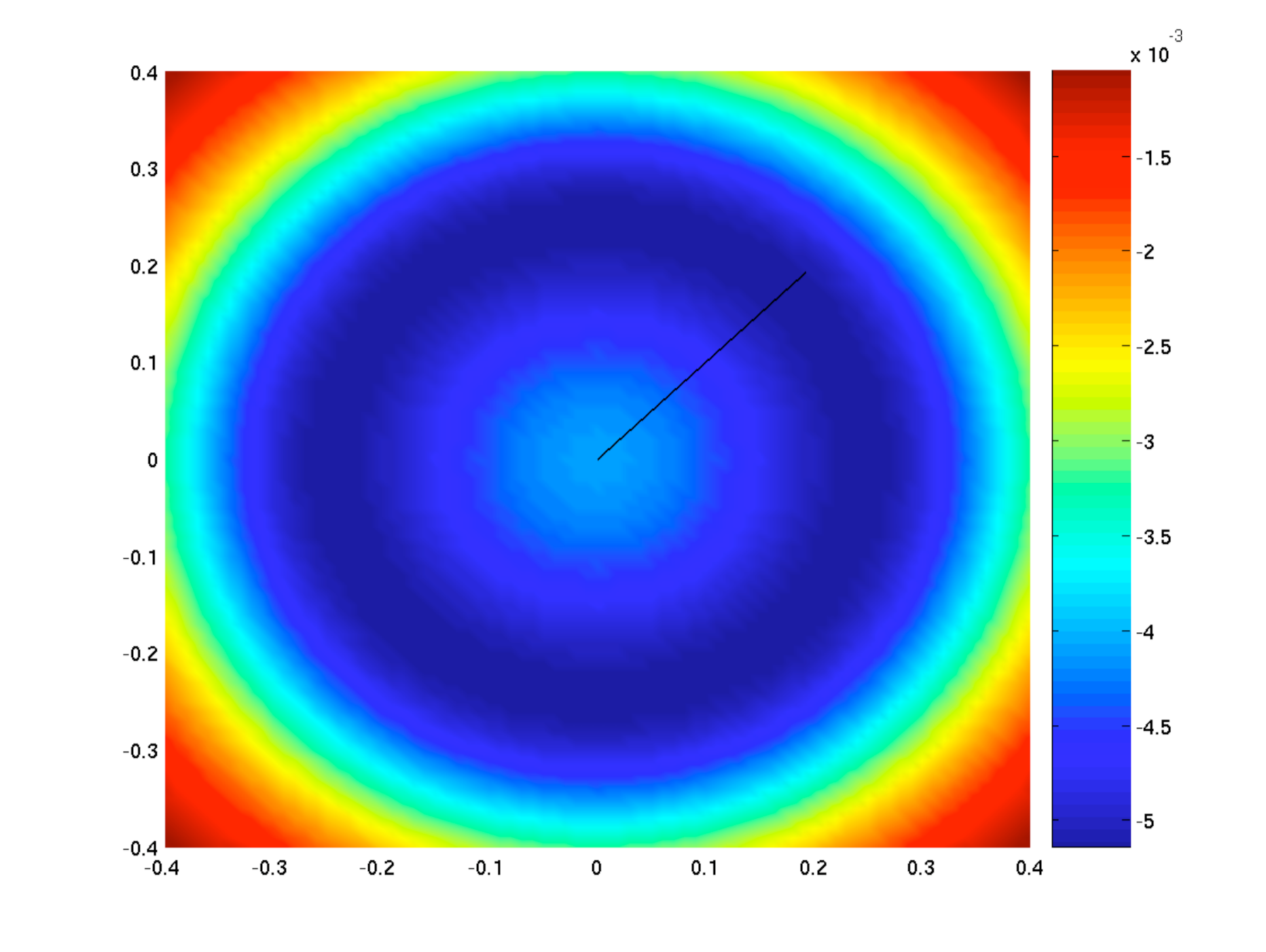}
\includegraphics[scale=0.25,angle=0]{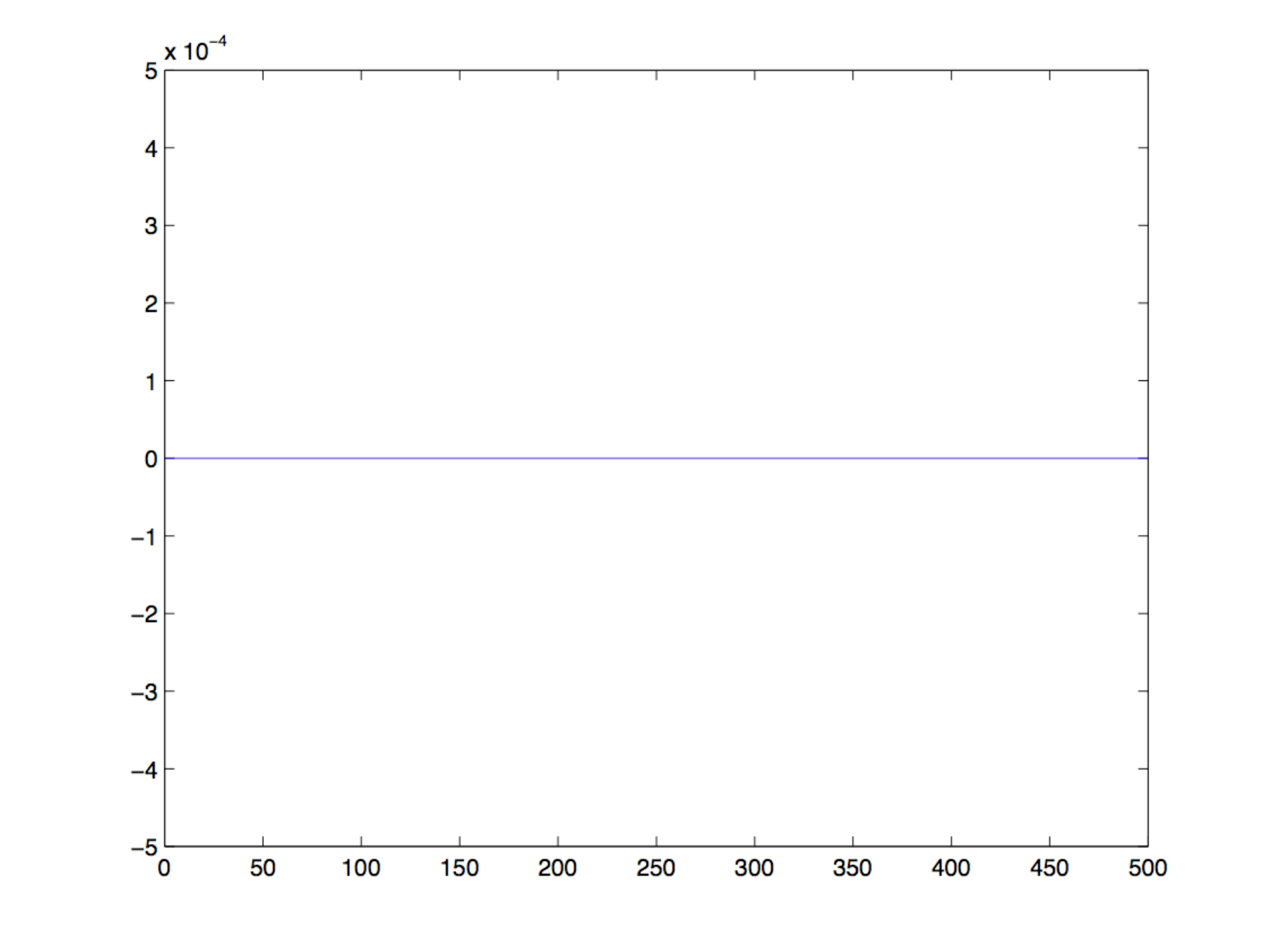}
\hbox{\hspace{0.5in}   (a) potential surface $V_e$   \hspace{0.6in}
(b) momentum map $J(t)$ \hspace{0.3in}    } \caption{\footnotesize
{\bf Nonconservative Axisymmetric.}  The figures plot $x(t), y(t)$
and $J(t)$ of the top in the presence of friction, and for
$\phi_R=0$ and $\theta_R=0$.    Superimposed on (a) is the potential
energy level set which appears axisymmetric.  (a) shows that the top
moves from the local maximum at the center to a point that minimizes
the potential energy (dark blue ring). However, (b) shows that no
spin develops. The accompanying simulation vividly illustrates this
motion and demonstrates that one does not get curious rotation in
this case.
   }
\label{fig:axi}
\end{center}
\end{figure}

\begin{figure}[htbp]
\begin{center}
\includegraphics[scale=0.25,angle=0]{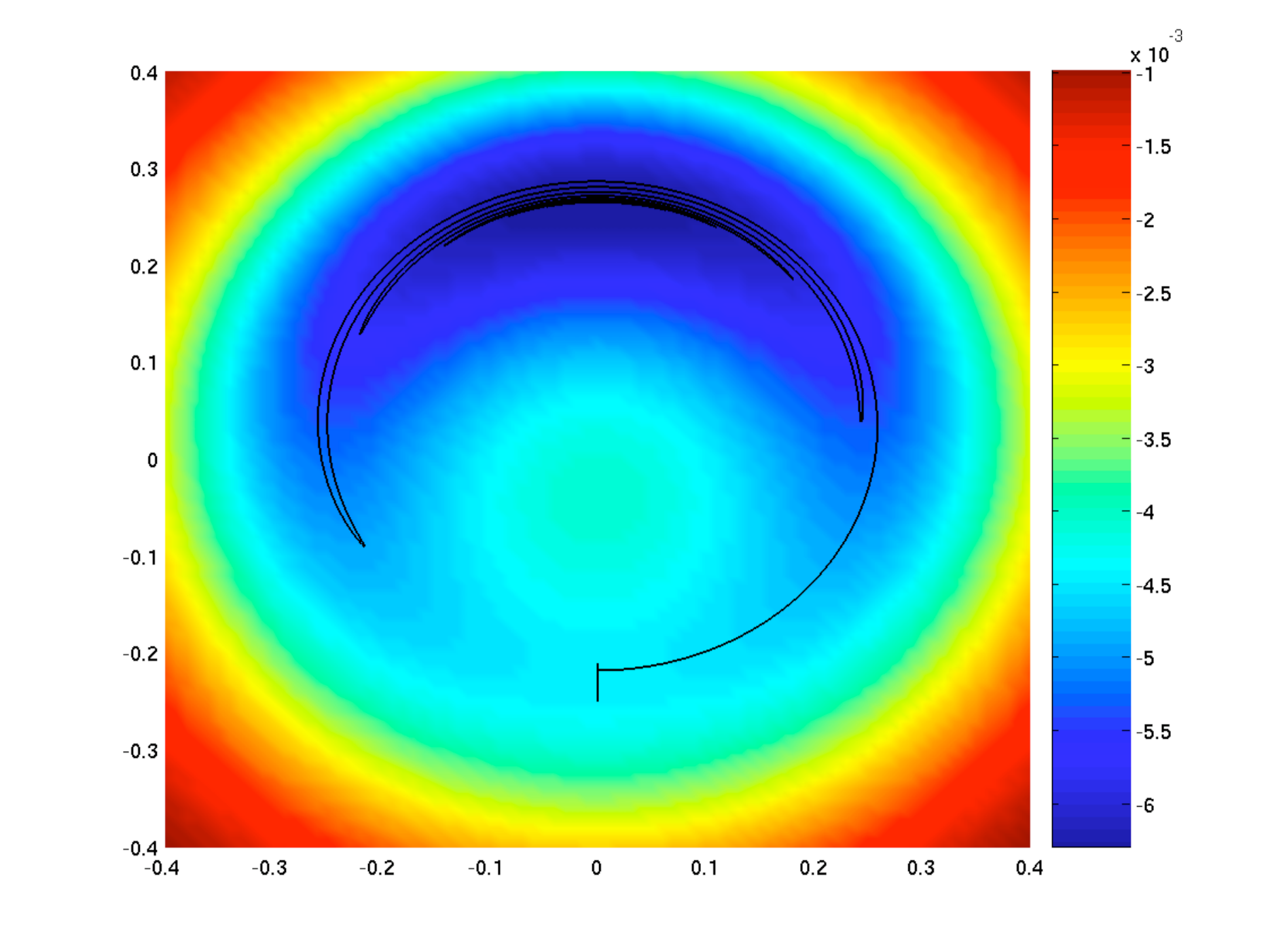}
\includegraphics[scale=0.25,angle=0]{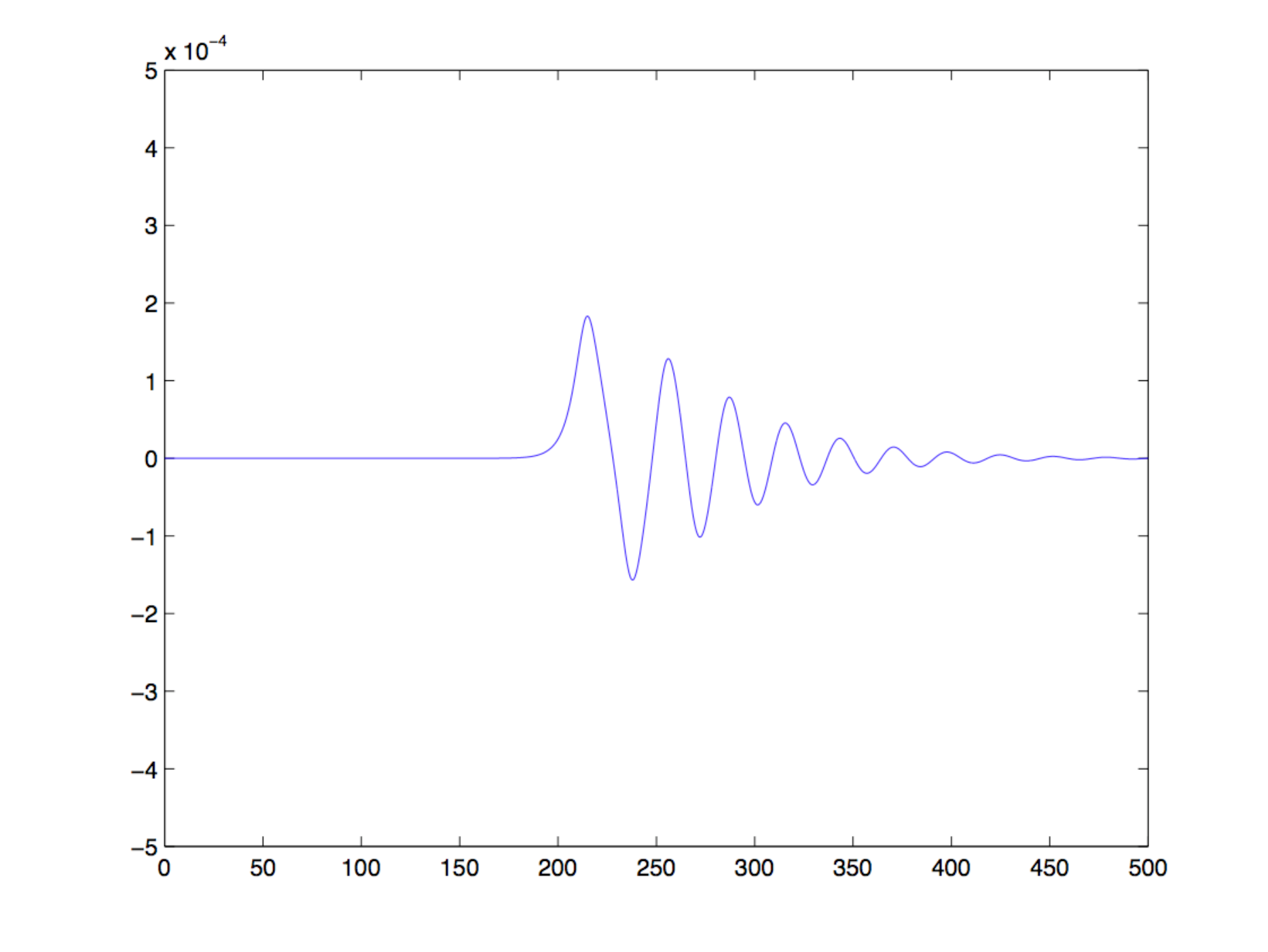}
\hbox{\hspace{0.5in}   (a) potential surface $V_e$   \hspace{0.6in}
(b) momentum map $J(t)$ \hspace{0.3in}    } \caption{\footnotesize
{\bf Nonconservative Tilted.} The figures plot $x(t), y(t)$ and
$J(t)$ of the top in the presence of friction, and for $\phi_R=0$
and $\theta_R=\pi/64$.    Superimposed on (a) is the potential
energy level set which is no longer axisymmetric.  (a) shows that
the top moves along arcs.  (b) confirms that one does get spin in
this case. The accompanying simulation vividly illustrates this
motion and clearly demonstrates that the ball goes from a state of
rest to a state of nonzero spin about its attitude.   }
\label{fig:tilt}
\end{center}
\end{figure}

In the presence of friction $c=0.3$, the motion changes.  The third
case considered is the same configuration as the first case, but
with friction.  In this case the kinetic energy of the top is
dissipated and the top moves towards a point where the magnetic
potential energy is minimum, see Fig.~\ref{fig:axi}.  However,  no
spin develops.   If the attitude of the ring is slightly perturbed,
spin does develop as shown in Fig.~\ref{fig:tilt}.  The origin of
this spin in case 4 is explained here.

It is very clear from (\ref{eq:magtopnc}) that one can get spin about the axis
of symmetry in the presence of friction.  However, the frictional
forces may produce torques about the other axes as well.   So what
is not clear is why the ball spins primarily about
$\boldsymbol{\xi}_3$.  This phenomenon is clarified in the following remark.

\begin{Remark}
If the initial position of the ball is not a local minimum of the
magnetic potential energy (as in the simulation), the position will
be unstable as predicted by a Lagrange-Dirichlet criterion.   The
ball's position then adjusts to minimize its magnetic potential
energy which causes sliding friction (see Fig.~\ref{fig:whyspin}).
The torque due to the sliding friction will introduce a torque in
directions orthogonal to the moment arm $\mathbf{q}$. However, the
torque due to the magnetic field will counter the torques about the
axes perpendicular to $\boldsymbol{\xi}_3$.  Keep in mind that the
magnetic torque keeps $\boldsymbol{\xi}_3$ aligned with the local
magnetic  field.   Thus, the torque due to friction mainly causes a
spin about the $\boldsymbol{\xi}_3$ axis as shown in the simulation.
Snapshots of the simulation are provided in
Fig.~\ref{fig:snapshots}.
\end{Remark}

\newpage

\section{Fluctuation Driven Magnetic Motor}\label{kkshshkjhe}

The instability of the top with respect to perturbations of the ring
is the key idea behind the fluctuation driven magnetic motor.   
These random fluctuations are modelled as a white
noise torque on the ring.  However, if the ring is allowed to move
freely, it could possibly turn on its side.  To stabilize the ring,
a fixed magnetized outer ring is installed.

By adjusting the radii and heights of the rings and the inertia of
the inner ring, one can obtain a configuration in which the attitude
of the inner ring can be randomly torqued without undergoing large
excursions from the vertical position.   In this case numerical
experiments reveal that one can adjust the amplitude of the white
noise so that the ball undergoes directed motion on certain
time-scales. To be precise the numerics indicates that starting from
a position of rest initially the top's motion is dominated by the effect 
of white noise.   If the outer ring is close enough one can see the top oscillate between the wells in the magnetic potential caused by the dipoles in the outer
ring.  This behavior is reminiscent of stochastic resonance.

After some time, the top accumulates enough kinetic energy that it
displays directed motion along a circle of certain radius.  This
motion is inertia-driven and the energy injected into the system by
the white noise mainly adds to the speed of the top.  Provided that
the inner ring does not turn on its side, one of two things can
happen: 1) the top reaches a critical velocity in which the amount
of energy dissipated by the surface friction is on average equal to
the amount of energy  injected by the thermal noise or 2) the top
gathers enough kinetic energy to escape from the potential well
created by the inner and outer rings.  Conducting this same
experiment at uniform temperature, reveals that this phenomenon
persists. The isothermal, magnetic ball-ring system is a prototype
fluctuation-driven motor.

The fluctuation-driven motor considered in this paper is
related to the granular, magnetic balls in ferrofluidic thermal
ratchets.  In such ratchets a time-varying magnetostatic field
transfers angular momentum to magnetic spherical grains  in a
ferrofluid \citep{EnMuReJu2003}. The ferromagnetic grains are
modelled in the same way as the rigid ball in the prototype.
However, our work is different in that the magnetic potential in the
prototype is autonomous.  In fact, the use of a non-autonomous
potential to design a fluctuation driven motor is well understood
\citep{AjPr1992}.

In this section the equations of motion for the dynamics of a
magnetic ball interacting with a dynamic, inner magnetized
ring are derived.   The magnetic effects of a fixed outer ring are
also taken into account.  As before the magnetic ball is assumed to
be constrained to a flat surface that resists the motion of the top
via surface sliding friction.  The center of mass of the inner ring is kept at a fixed
height, but otherwise it is free to rotate and subjected to white noise torques. 
The governing stochastic differential equations are analyzed using energy
arguments.  Simulations validate this analysis and show that one can
get sustained directed motion of the top from random perturbations of
the inner ring.  Moreover, it is shown that this phenomenon persists even 
when the system is at uniform temperature.  This system is called a 
fluctuation driven magnetic motor.

\begin{figure}[htbp]
\begin{center}
\includegraphics[scale=0.35,angle=0]{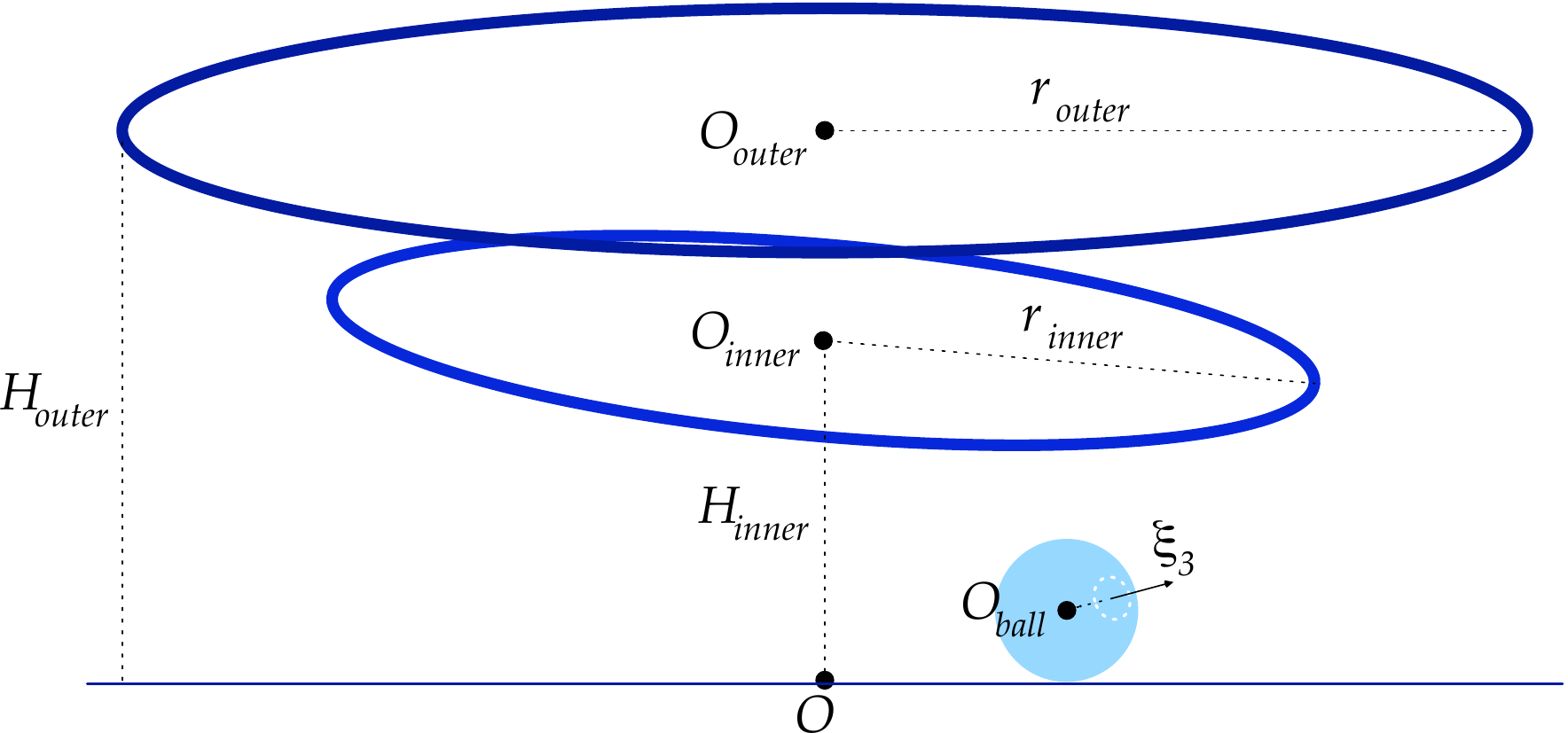}
\caption{\footnotesize  {\bf Illustration of Magnetized Rings and
Magnetic Top.}   The figure shows the magnetic top on the surface
with two rings above it.    The inner and outer rings have radii
$r_{\text{outer}}$ and $r_{\text{inner}}$ and heights
$H_{\text{outer}}$ and $H_{\text{inner}}$ respectively.  The
centroids of the ball, inner and outer rings are also labeled.   The
purpose of the outer ring is to prevent the trivial equilibrium in
which the inner ring turns on its side.
 }
\label{fig:magmotor}
\end{center}
\end{figure}

The mass and radius of the magnetic ball are $m$ and $r$
respectively.  The mass of the inner ring is $m_{\text{inner}}$. The
heights and radii of the inner and outer rings are plotted and
defined in Fig.~\ref{fig:magmotor}.

\paragraph{Lagrangian of Magnetic Ball \& Inner Ring}

The configuration space of the system is $Q = \mathbb{R}^3 \times
\operatorname{SO}(3) \times \operatorname{SO}(3)$.  Let $O_{ball}$,
$O_{inner}$ and $O_{outer}$ denote the centroids of the ball, inner
and outer rings respectively. Let $(\mathbf{e}_1, \mathbf{e}_2,
\mathbf{e}_3)$, $(\mathbf{f}_1, \mathbf{f}_2, \mathbf{f}_3)$, and
$(\mathbf{g}_1, \mathbf{g}_2, \mathbf{g}_3)$ denote inertial
orthonormal frames attached to $O_{ball}$, $O_{inner}$ and
$O_{outer}$, and related to body-fixed frames $(\boldsymbol{\xi}_1,
\boldsymbol{\xi}_2, \boldsymbol{\xi}_3)$, $(\boldsymbol{\zeta}_1,
\boldsymbol{\zeta}_2, \boldsymbol{\zeta}_3)$, and
$(\boldsymbol{\eta}_1, \boldsymbol{\eta}_2, \boldsymbol{\eta}_3)$
via the rotation matrices  $R_B(t)$, $R_R(t)$, and $R_O \in
\operatorname{SO}(3)$.   Let $\mathbb{I}_B =
\operatorname{diag}(I_1, I_2, I_3)$ and $\mathbb{I}_R=
\operatorname{diag}(J_1,J_2, J_3)$ be the standard diagonal inertia
matrices of the body and inner ring, and let $ \boldsymbol{\xi}_3$
and $\boldsymbol{\zeta}_3$ be the attitude of the ball and inner
ring.  The following assumptions are made: the outer ring is fixed,
the mass distribution of the inner ring is  symmetric with respect
to its attitude (or axis of symmetry), and the ball's mass
distribution is spherically symmetric.  This  assumption implies
that $J_B=I_1=I_2=I_3$ and $J=J_1=J_2$.

The magnetostatic field at any field point is due to the dipole in
the ball, and an inner and outer ring consisting of $N$ and $M$
magnetic dipoles respectively.   On each ring the dipoles are
equally spaced around a circle of  radius $r_{inner}$ and
$r_{outer}$ in the planes defined by the vectors
$\boldsymbol{\eta}_3$ and $\boldsymbol{\zeta}_3 \in \mathbb{R}^3$.
We assume that there is no self-interaction between the dipoles
within each body. The heights of the rings above the surface are
denoted by $H_{\text{inner}}$ and $H_{\text{outer}}$.

For $i=1, \cdots, N$ and $j=1, \cdots, M$, let
$\mathbf{d}_i^{\text{inner}}(t)$ and  $\mathbf{d}_j^{\text{outer}}
\in \mathbb{R}^3$ denote the location of the ith and jth-dipoles on
the inner and outer rings with respect to the points
$O_{\text{inner}}$ and $O_{\text{outer}}$ and let
$\mathbf{m}_i^{\text{inner}}(t)$ and $\mathbf{m}_j^{\text{outer}}
\in \mathbb{R}^3$ denote the orientation of their respective dipole
moments.   The dipole moments are assumed to be in the radial
direction of the ring.  The magnetic field of each dipole at a field
point can be determined from  the vector $\mathbf{r}$ connecting the
field point to the dipole using (\ref{eq:magneticfield})
\begin{align*}
\mathbf{B}_i^{\text{inner}} (\mathbf{r}) =   \mathbf{B}(\mathbf{r}, \mathbf{m}_i^{\text{inner}})  ~~~&\text{(field of inner ring dipole),}  \\
\mathbf{B}_j^{\text{outer}} (\mathbf{r}) =   \mathbf{B}(\mathbf{r}, \mathbf{m}_i^{\text{outer}})~~~&\text{(field of outer ring dipole),} \\
\mathbf{B}_0 (\mathbf{r}) = \mathbf{B}(\mathbf{r},
\boldsymbol{\xi}_3) ~~~&\text{(field of ball).}
\end{align*}
Define the following vectors $\mathbf{r}_{i j} \in \mathbb{R}^3$
$i=0,\cdots,N$ and $j=0,\cdots,M$ which joins the inner and outer
dipoles, the inner dipoles to the ball, and outer dipoles to the
ball as follows,
\[
\mathbf{r}_{ij} = \left\{   \begin{array}{lll}
-\mathbf{d}_j^{\text{outer}} - H_{\text{outer}} \mathbf{e}_3
+ \mathbf{d}_i^{\text{inner}} + H_{\text{inner}} \mathbf{e}_3 & \text{if $i,j>0$} & \text{(ith outer to jth inner dipole),} \\
 \mathbf{x} - \mathbf{d}_i^{\text{inner}} - H_{\text{inner}} \mathbf{e}_3 & \text{if $j=0$} & \text{(jth outer to ball),} \\
  \mathbf{x} -  \mathbf{d}_j^{\text{outer}} - H_{\text{outer}} \mathbf{e}_3 & \text{if $i=0$}  & \text{(ith inner to ball).} \end{array} \right.
\]

The spatial representation of the reduced Lagrangian  $\ell: T
\mathbb{R}^3 \times  \operatorname{SO}(3) \times \mathbb{R}^3 \times
\operatorname{SO}(3) \times \mathbb{R}^3  \to \mathbb{R}$ of the
free magnetic ball, i.e., magnetic ball without dissipation, is
given by,
\begin{align}
& \ell(  \mathbf{x}, \mathbf{\dot{x}}, R_B, \boldsymbol{\omega}_B,
R_R, \boldsymbol{\omega}_R) = \frac{m}{2}
\mathbf{\dot{x}}^\mathrm{T}  \mathbf{\dot{x}} +  \frac{1}{2}  J_B
\boldsymbol{\omega}_B^\mathrm{T} \boldsymbol{\omega}_B
+  \frac{1}{2}  \boldsymbol{\omega}_R^\mathrm{T} R_R \mathbb{I}_R R_R^\mathrm{T} \boldsymbol{\omega}_R  \nonumber \\
&+ 2 \sum_{i=1}^N  \boldsymbol{\xi}_3^\mathrm{T}
\mathbf{B}_i^{\text{inner}}(\mathbf{r}_{i0}) + \sum_{j=1}^M
\boldsymbol{\xi}_3^\mathrm{T}  \mathbf{B}_j^{\text{outer}}(
\mathbf{r}_{0j}) + \sum_{i=1}^N  \sum_{j=1}^M
(\mathbf{m}_i^{\text{inner}})^\mathrm{T} \mathbf{B}_j^{\text{outer}}
\left(  \mathbf{r}_{ij} \right) \label{eq:ellringnball}
\end{align}
From left the terms represent the translational and rotational
kinetic energy of the ball, the rotational kinetic energy of the
ring, the gravitational potential energy of the ball and the
magnetic potential energy of the inner ring and ball dipoles.

\paragraph{Governing Equations for Magnetic Ball \& Ring}

The equations of motion are determined from a Hamilton-Pontryagin (HP) principle \citep{BoMa2007a}.   The HP action integral is given by,
\[
s = \int_a^b \left[ \ell(\mathbf{x}, \mathbf{\dot{x}}, R_B,
\boldsymbol{\omega}_B, R_R, \boldsymbol{\omega}_R) + \left\langle
\widehat{\boldsymbol{\pi}}_B, \dot{R}_B R_B^\mathrm{T}  -
\widehat{\boldsymbol{\omega}_B} \right\rangle + \left\langle
\widehat{\boldsymbol{\pi}}_R, \dot{R}_R R_R^\mathrm{T}  -
\widehat{\boldsymbol{\omega}_R} \right\rangle  + \lambda \varphi(
\mathbf{x} )  \right] dt \text{.}
\]
The HP principle states that
 \[
\delta s = 0
\]
where the variations are arbitrary except that the endpoints
$(\mathbf{x}(a),R_B(a),R_R(a))$ and $(\mathbf{x}(b),R_B(b),R_R(b))$
are held fixed.   The equations are given by,
 \begin{align}
 \frac{d}{dt} \frac{\partial \ell}{\partial \mathbf{\dot{x}}} = \frac{\partial \ell}{\partial \mathbf{x}} + \lambda \frac{\partial \varphi}{\partial \mathbf{x}}   \label{eq:balleulerlagrange}  ~~~&\text{(Euler-Lagrange equations for ball),} \\
 \varphi( \mathbf{x} )  = 0 \label{eq:constraintball} ~~~&\text{(constraint equation),} \\
  \frac{d}{dt}    R_B =  \widehat{\boldsymbol{\omega}_B} R_B \label{eq:reconstructionball} ~~~&\text{(reconstruction equation for ball),} \\
    \frac{d}{dt}    R_R =  \widehat{\boldsymbol{\omega}_R}  R_R \label{eq:reconstructionring} ~~~&\text{(reconstruction equation for ring),} \\
      \boldsymbol{\pi}_B = \frac{\partial \ell}{\partial  \boldsymbol{\omega}_B} \label{eq:legendreball} ~~~&\text{(reduced Legendre transform for ball),} \\
            \boldsymbol{\pi}_R = \frac{\partial \ell}{\partial  \boldsymbol{\omega}_R} \label{eq:legendrering} ~~~&\text{(reduced Legendre transform for ring),} \\
    \frac{d}{dt} \widehat{\boldsymbol{\pi}_B}  = \frac{\partial \ell}{\partial R_B} R_B^{\mathrm{T}}   -   \widehat{\widehat{\boldsymbol{\pi}_B} \boldsymbol{\omega}_B} \label{eq:lpball} ~~~&\text{(Lie-Poisson equations for ball),}  \\
        \frac{d}{dt} \widehat{\boldsymbol{\pi}_R}  = \frac{\partial \ell}{\partial R_R} R_R^{\mathrm{T}}  -   \widehat{\widehat{\boldsymbol{\pi}_R} \boldsymbol{\omega}_R} \label{eq:lpring} ~~~&\text{(Lie-Poisson equations for ring).}
 \end{align}
The nonconservative system is obtained by adding the frictional
force and torque derived earlier to the translational and rotational
equations of the ball: (\ref{eq:balleulerlagrange}) and
(\ref{eq:lpball}).

\paragraph{Random Perturbations \& Uniform Temperature}

Consider driving the inner ring by the following Wiener process. 
 Let $\mathbf{W} \in \mathbb{R}^3$ denote Brownian
motion in $\mathbb{R}^3$, and append the following random torque:
\[
\mathbf{T}_R = \alpha d \mathbf{W}  \in \mathbb{R}^3  
\]
to the Lie-Poisson equation for the ring (\ref{eq:lpring}).

The work performed by $\mathbf{T}_R$
is equal to the change in total energy of the system.  On the other
hand, the work done by the kicks over a time interval $[a, b]$ is
given by the following formula:
\begin{equation}
\text{Work of kicks} = \alpha \int_a^b  \mathbf{T}_R^{\mathrm{T}}
\boldsymbol{\omega}_R \text{.}
\end{equation}
The work transferred from the inner ring to the top is given by:
\begin{equation}
\text{Work transferred to top} = \Delta \text{Total Energy} - \Delta
\text{Ring Energy} \text{.}
\end{equation}
Thus a measure of the efficiency of the magnetic motor is given by
the following ratio:
\begin{equation}
\text{Efficiency} = \frac{\text{Work transferred to top} }{
\text{Work of kicks}}  \text{.}
\end{equation}

As a next step frictional torque is introduced to the inner ring
and thermal torque (white noise) to the ball, so that the generator
of the process is characterized by a unique Gibbs-Boltzmann
invariant distribution.   The reader is referred to the Appendix for 
the governing equations of the magnetic motor at uniform and
non-uniform temperatures.

\paragraph{Simulations}

\begin{figure}[htbp]
\begin{center}
\includegraphics[scale=0.2,angle=0]{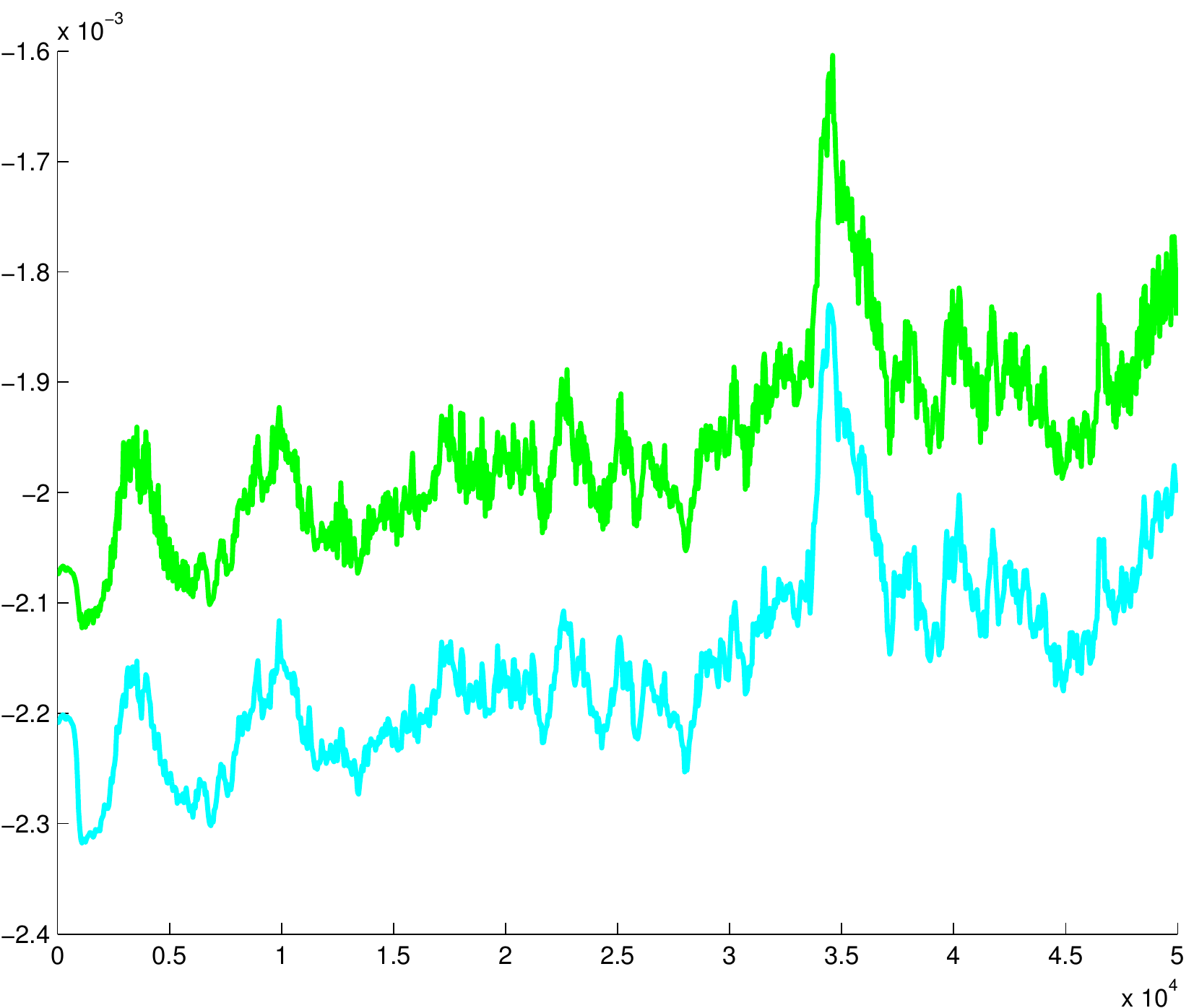}
\includegraphics[scale=0.2,angle=0]{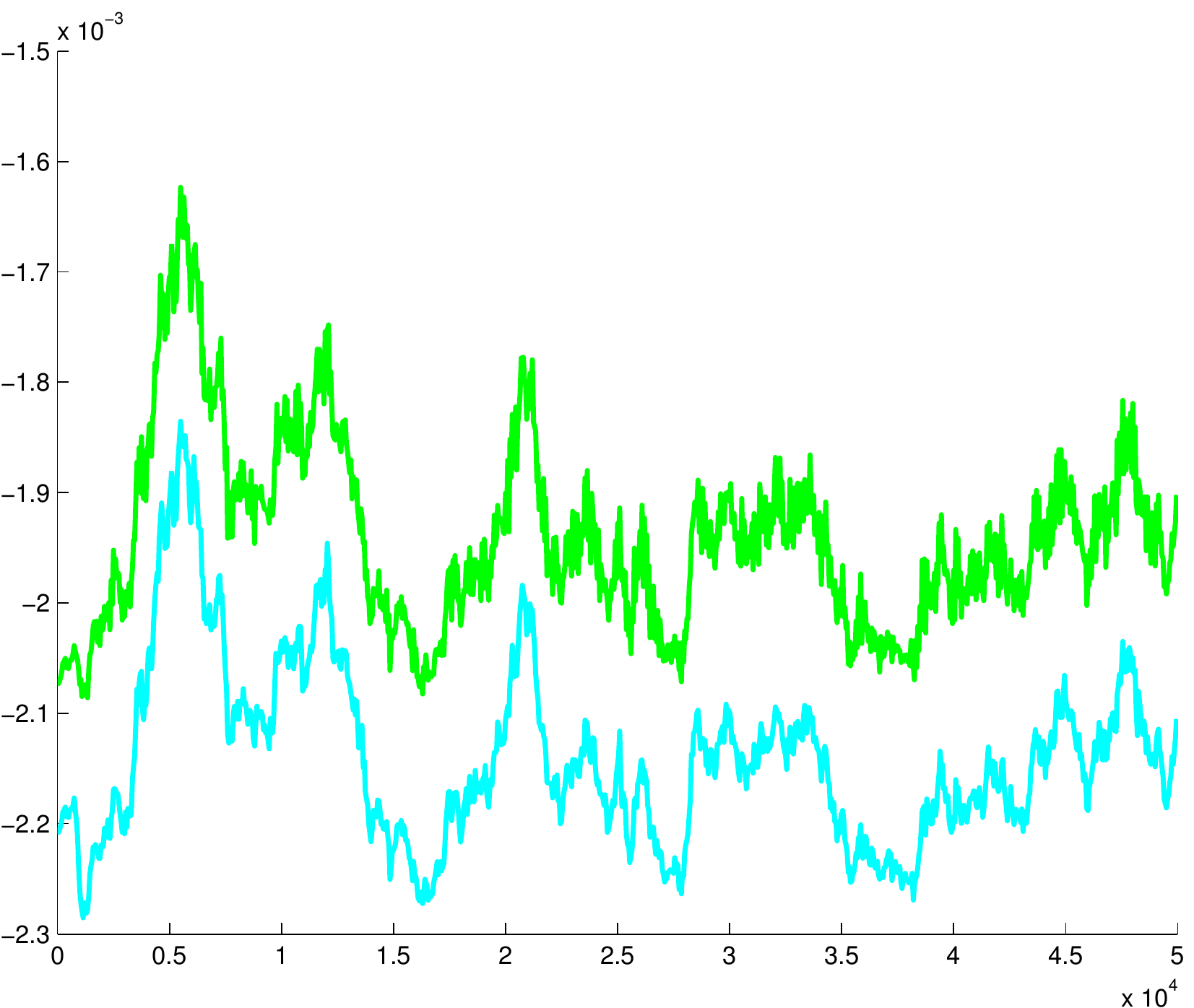}
\includegraphics[scale=0.2,angle=0]{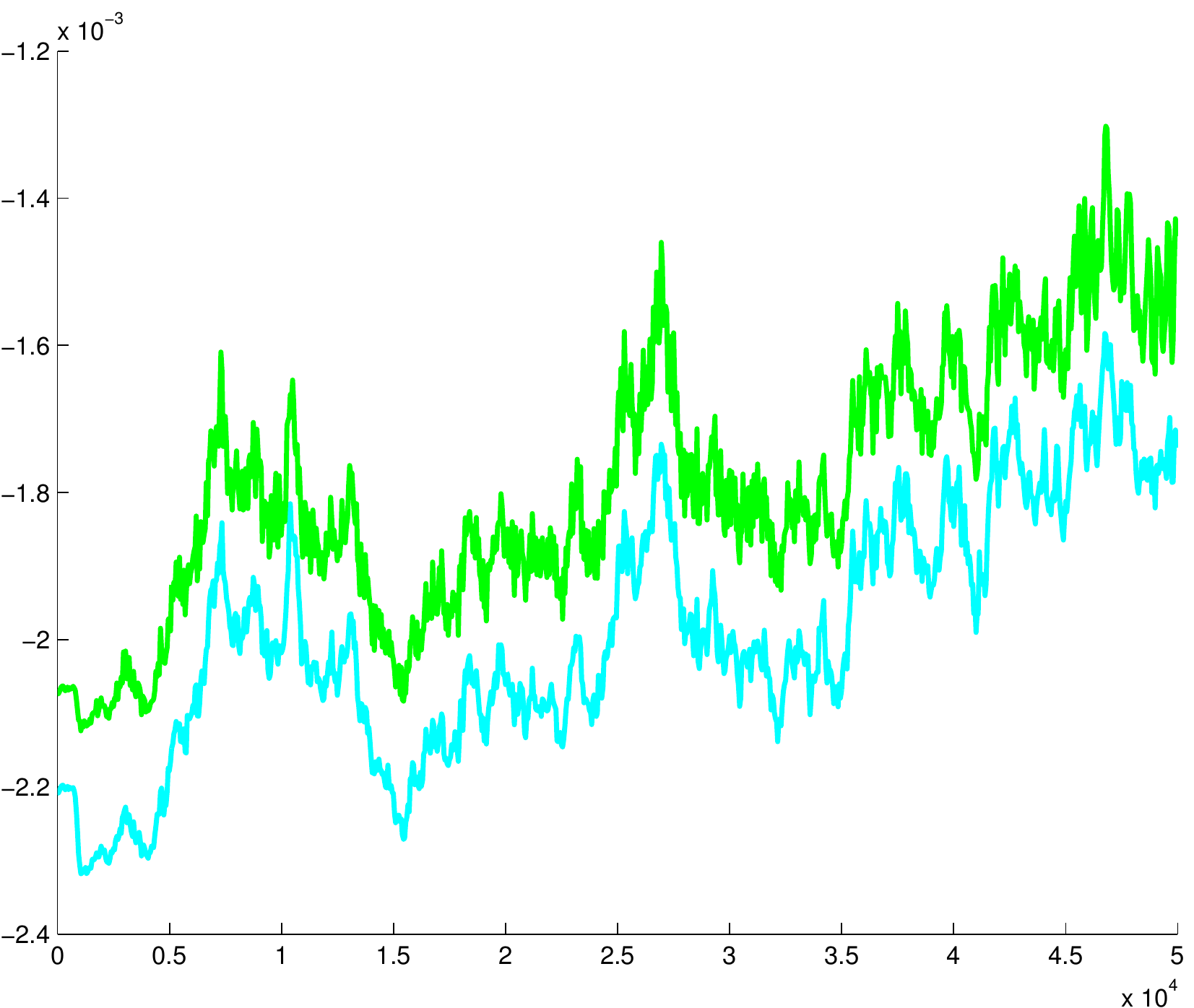}  \\
\hbox{\hspace{0.8in}   (a)  $\alpha=0.0002$ \hspace{0.3in} (b)
$\alpha=0.0002$ \hspace{0.3in} (c)  $\alpha=0.00025$  \hspace{0.3in}
} \caption{\footnotesize  {\bf Total and inner ring energy
(non-uniform temperature).}   The total energy of the system is
shown in green and the ring energy in cyan for different amplitudes
of the noise and for two different samples.   Although there is
dissipation, the thermal fluctuations inject energy into the system.
} \label{fig:nctmagmotorenergy}
\end{center}
\end{figure}

\begin{figure}[htbp]
\begin{center}
\includegraphics[scale=0.2,angle=0]{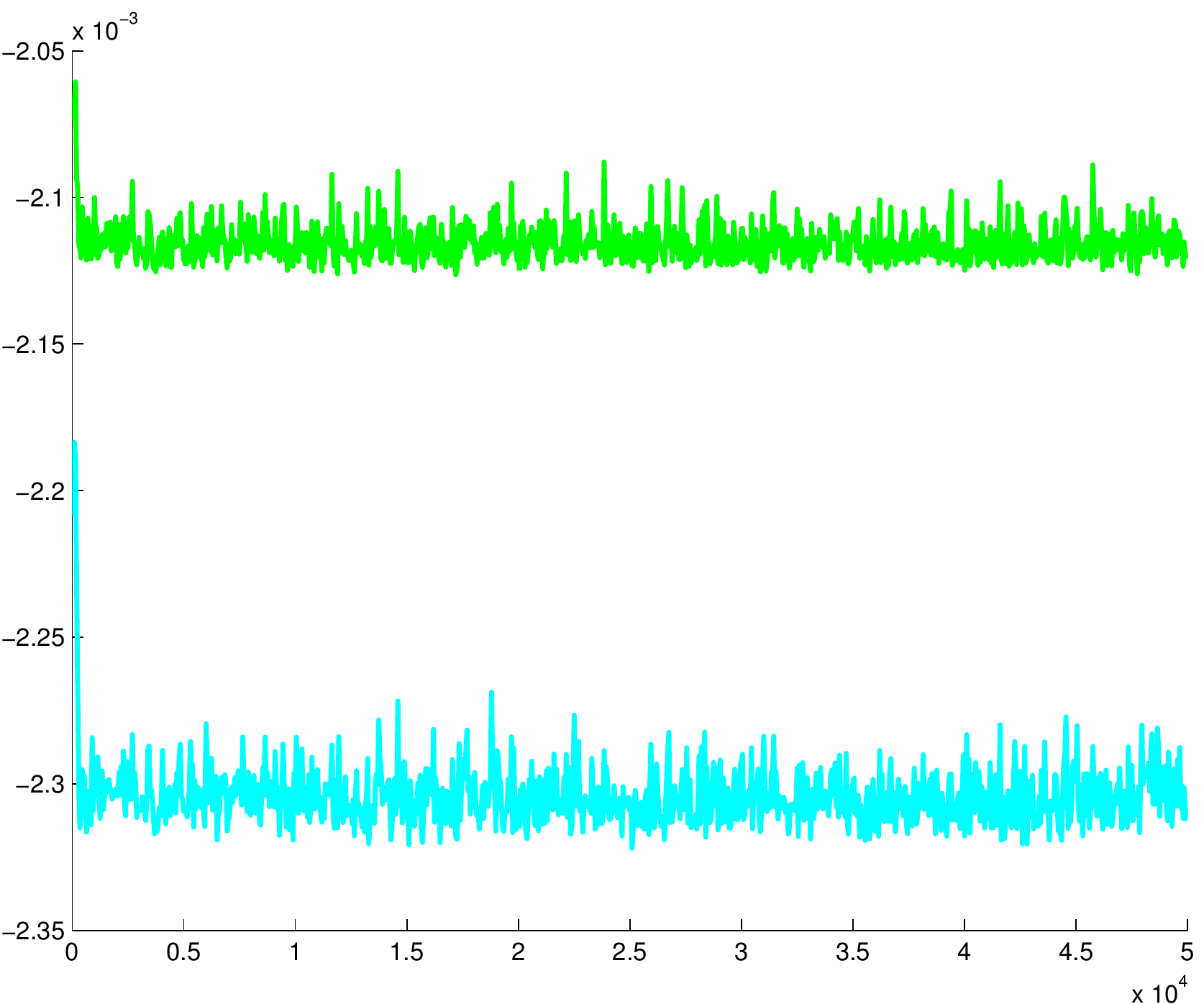}
\includegraphics[scale=0.2,angle=0]{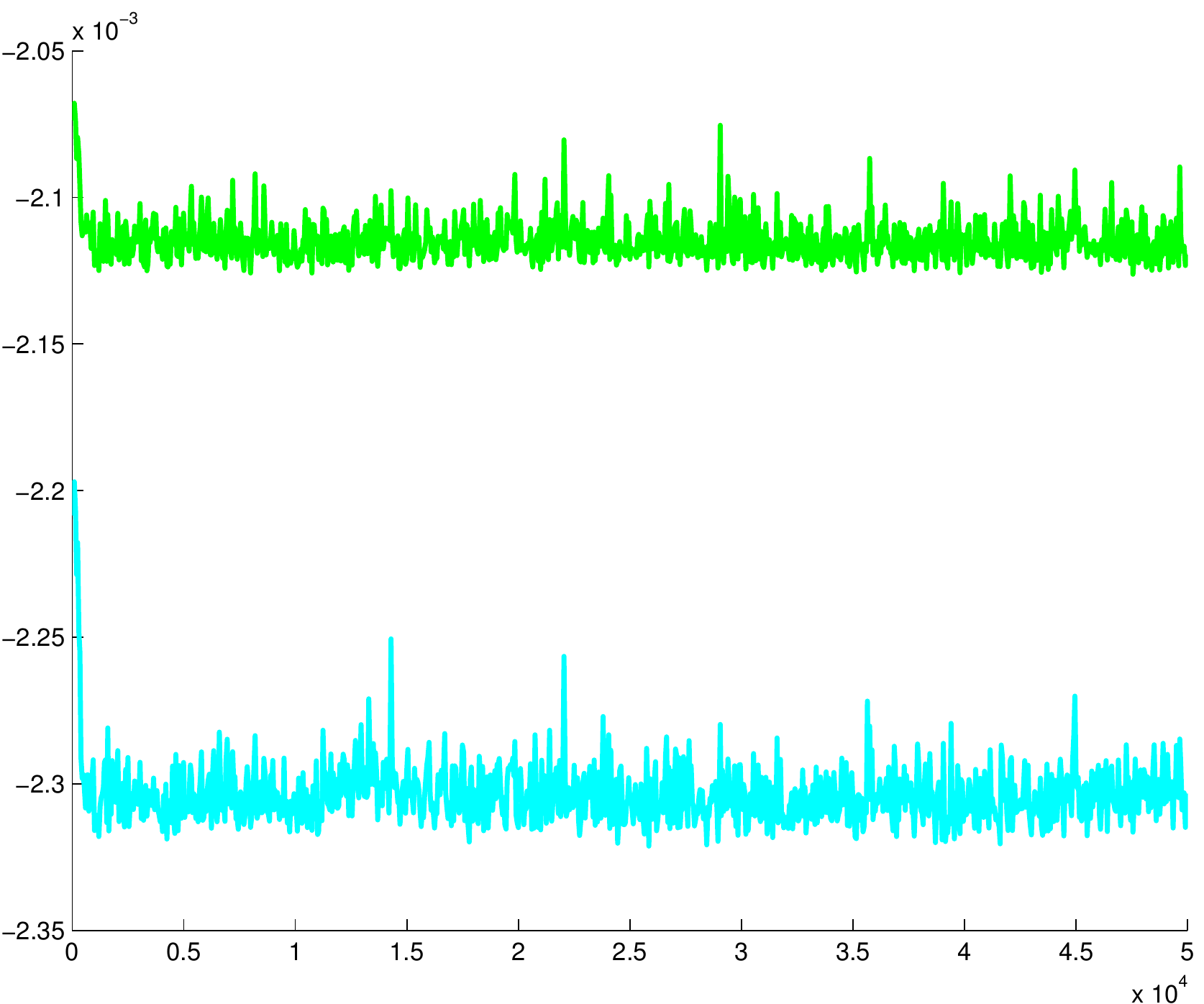}
\includegraphics[scale=0.2,angle=0]{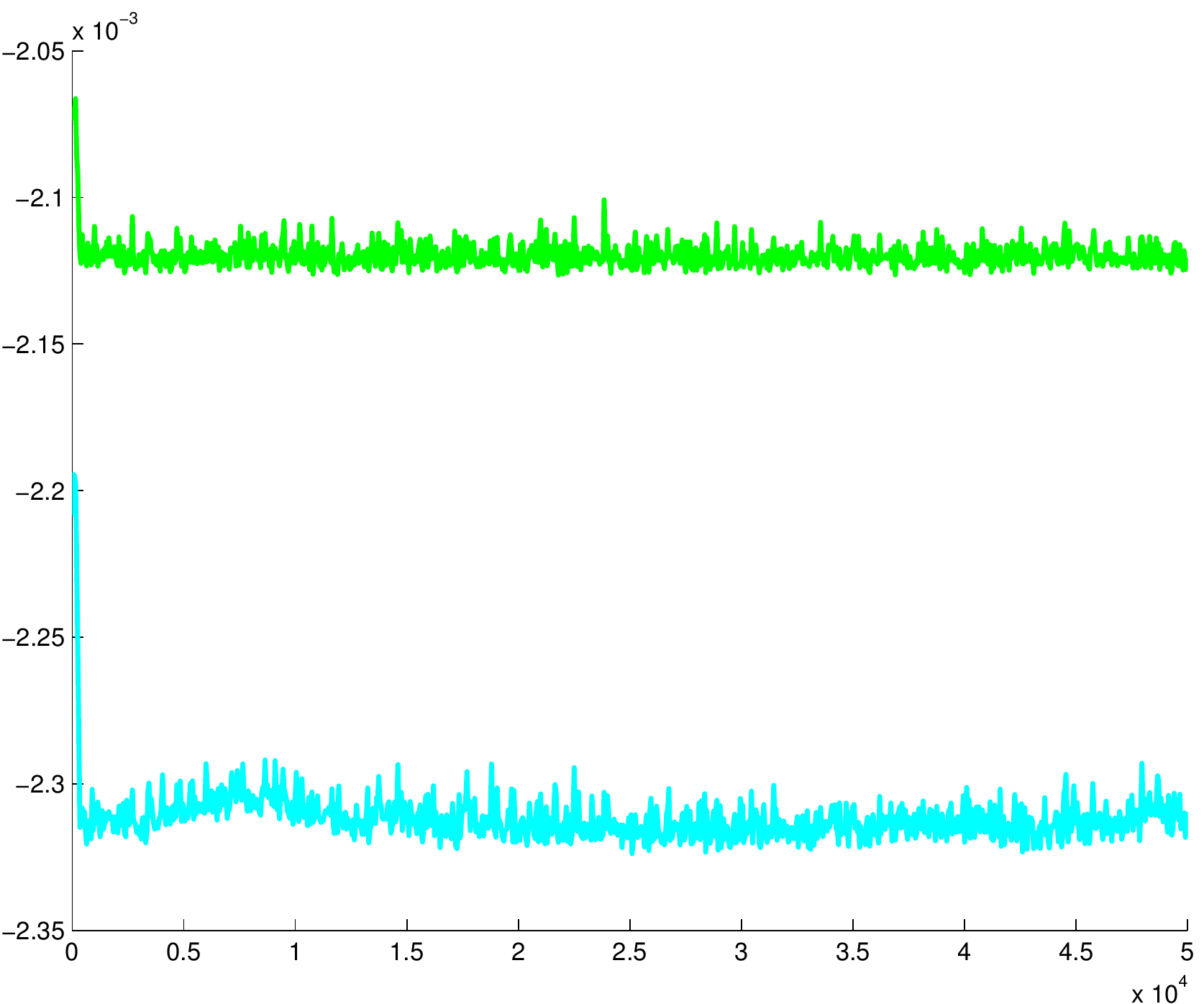}
\includegraphics[scale=0.2,angle=0]{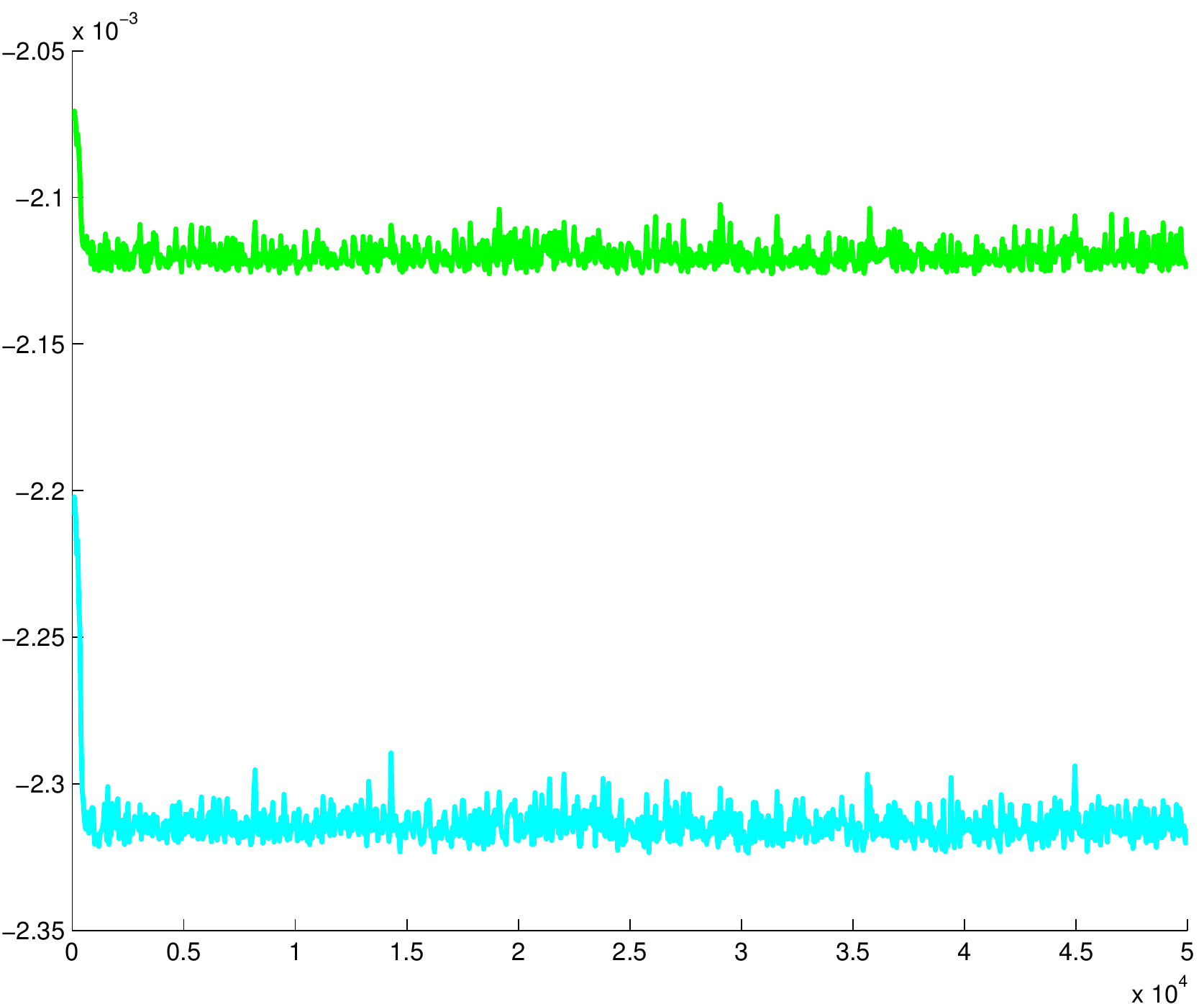} \\
\hbox{\hspace{0.15in}   (a)  $\alpha=0.001$ \hspace{0.3in} (b)
$\alpha=0.001$ \hspace{0.3in} (c)  $\alpha=0.00075$ \hspace{0.3in}
(d) $\alpha=0.00075$  \hspace{0.3in}    } \caption{\footnotesize
{\bf Total and inner ring energy (uniform temperature).}   The total
energy of the system is shown in green and the ring energy in cyan
for different amplitudes of the noise.   Although there is
dissipation, the thermal fluctuations inject energy into the system.
} \label{fig:ctmagmotorenergy}
\end{center}
\end{figure}

\begin{figure}[htbp]
\begin{center}
\includegraphics[scale=0.2,angle=0]{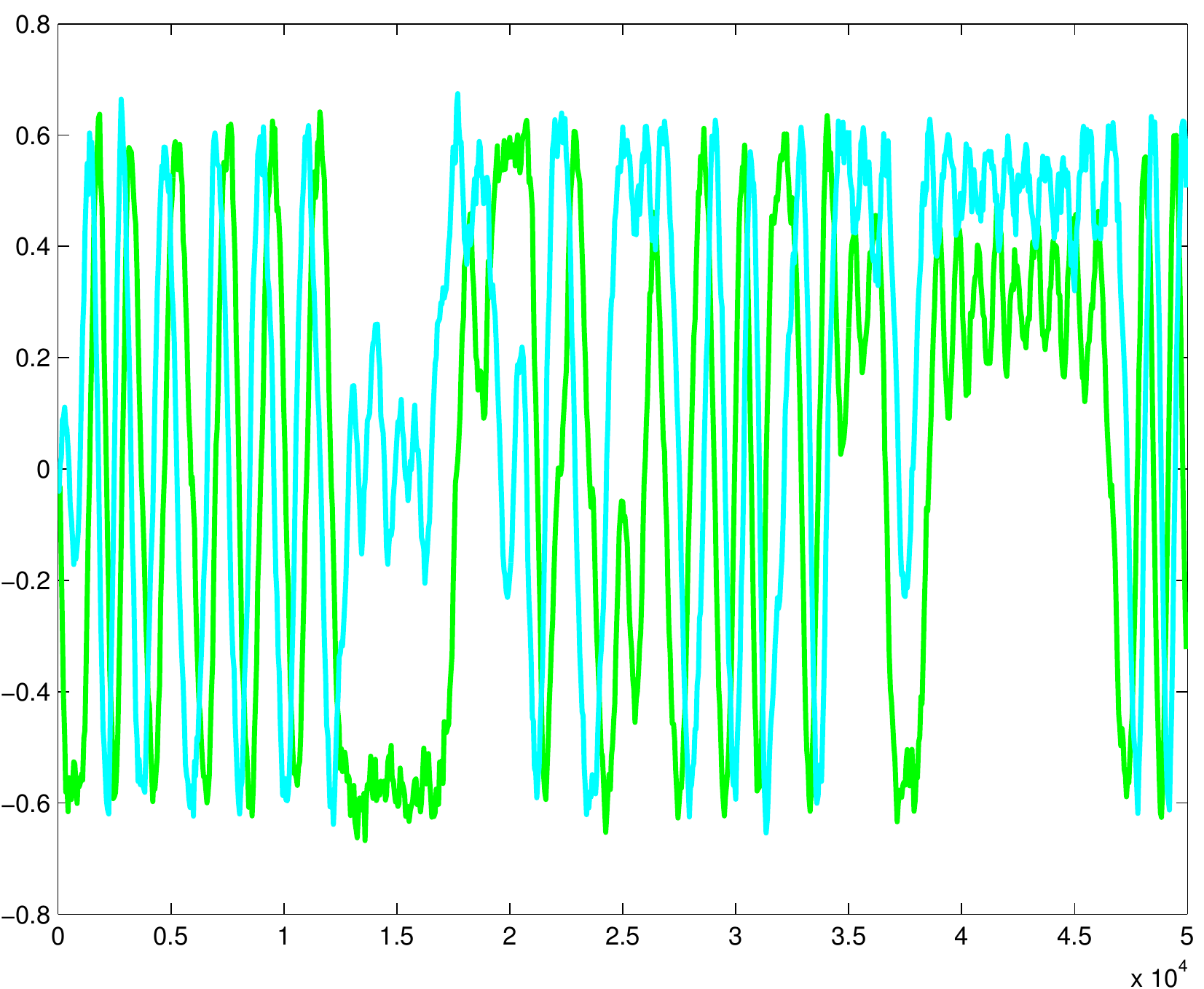}
\includegraphics[scale=0.2,angle=0]{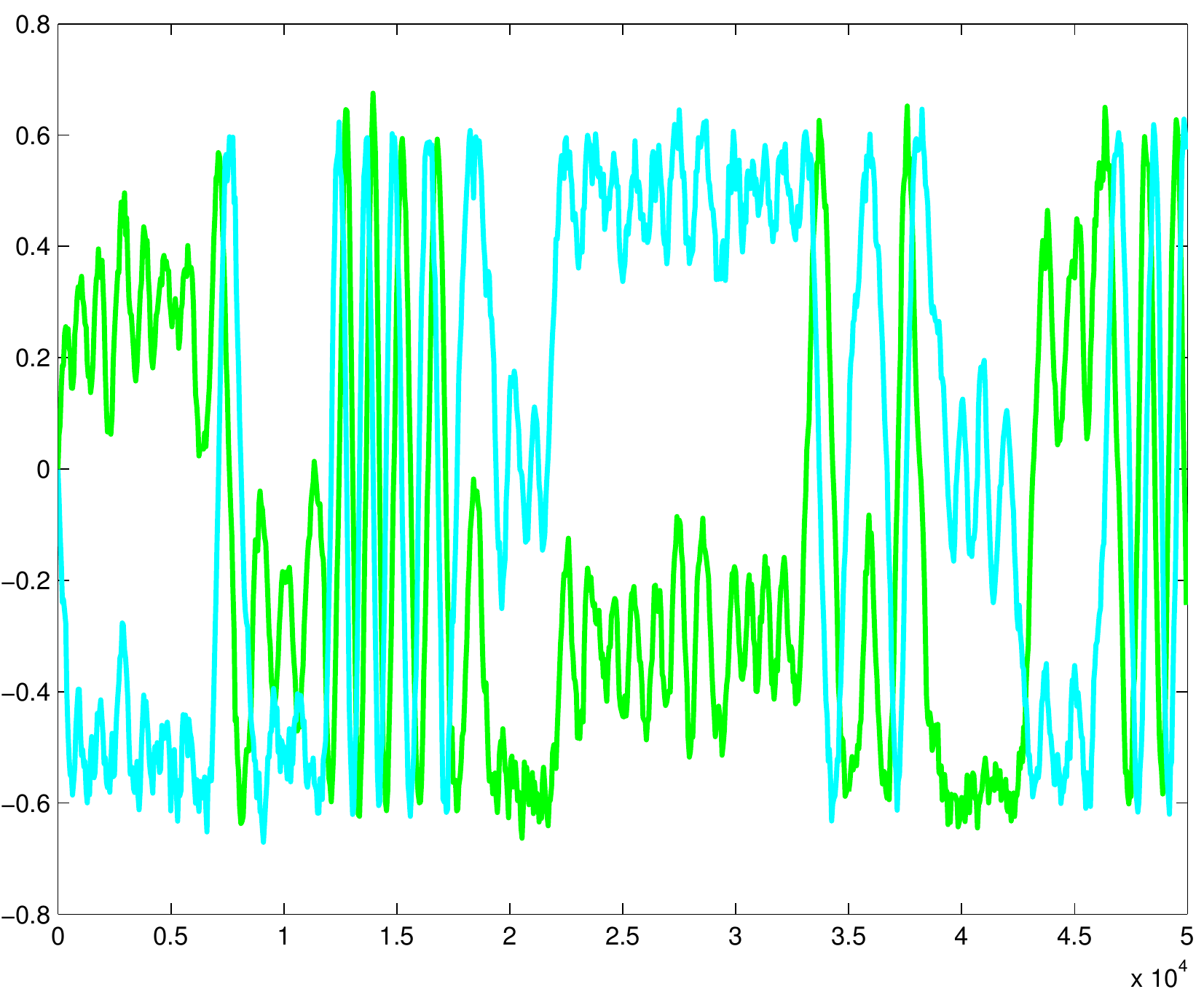}
\includegraphics[scale=0.2,angle=0]{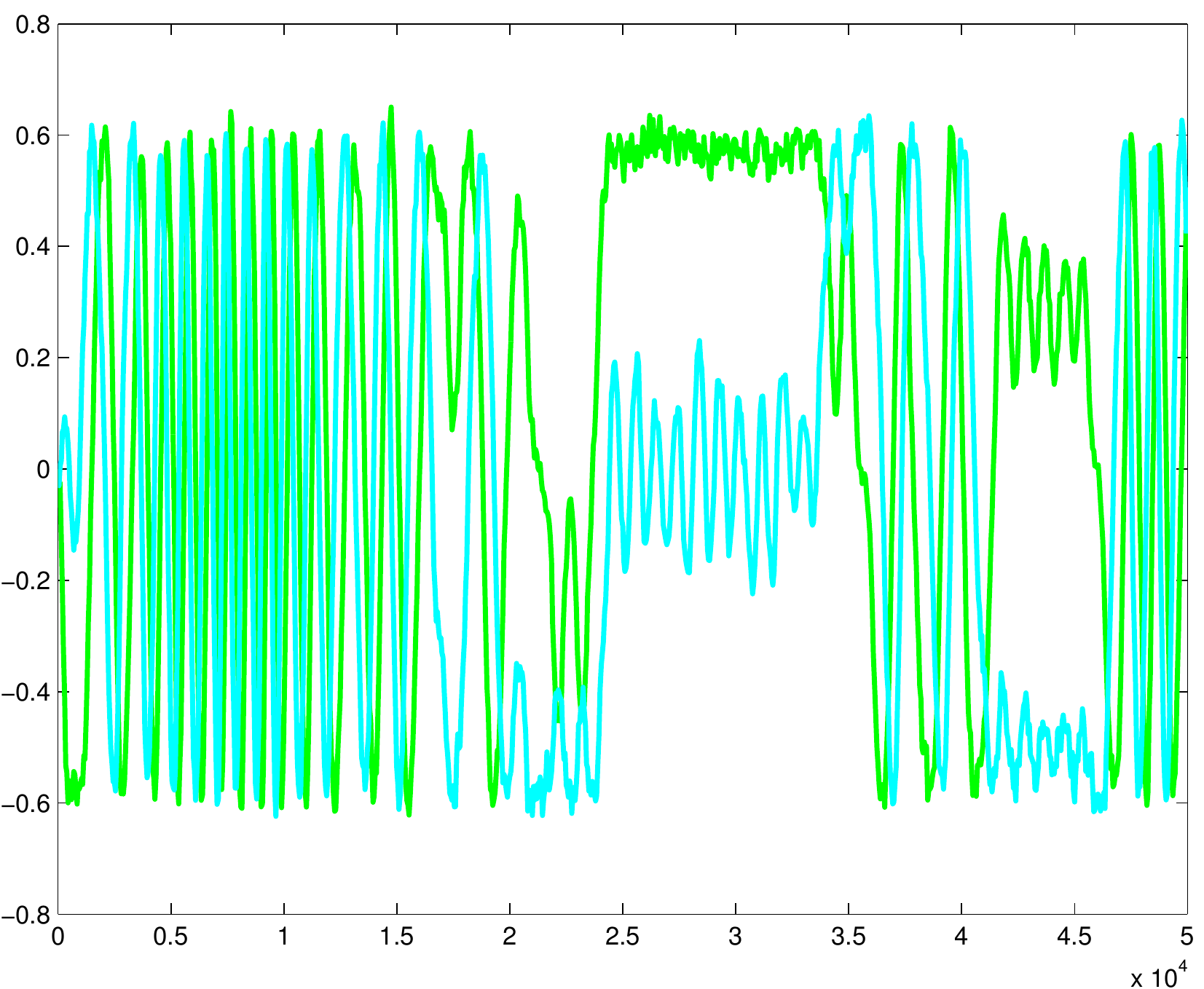}
\includegraphics[scale=0.2,angle=0]{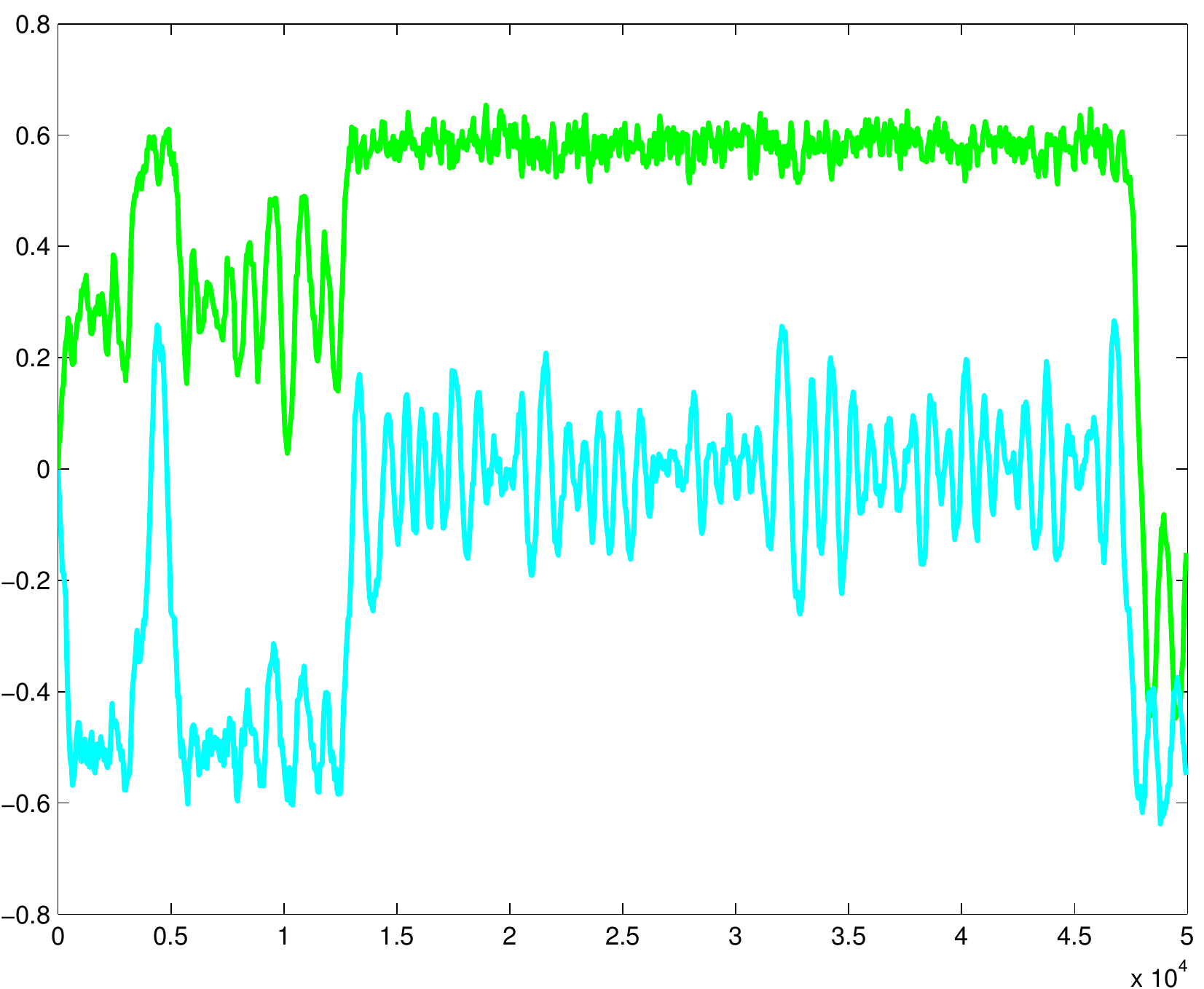} \\
\hbox{\hspace{0.15in}   (a)  $\alpha=0.001$ \hspace{0.3in} (b)
$\alpha=0.001$ \hspace{0.3in} (c)  $\alpha=0.00075$ \hspace{0.3in}
(d) $\alpha=0.00075$  \hspace{0.3in}    } \caption{\footnotesize
{\bf xy-position of magnetic ball planar view (uniform
temperature).} The  $x$ and $y$ components of the center of mass are
plotted in green and blue respectively.   It is very clear from this
plot that the motion transitions from noise driven to being inertia
driven. } \label{fig:ctmagmotorxyplots}
\end{center}
\end{figure}

A stochastic variational integrator is used to carry out these simulations
\citep{BoOw2007a}.  Simulations are conducted at uniform and non-uniform temperature as described below. The initial conditions and parameters used are given by
\begin{align*}
c&= 0.15~\text{kg} \cdot \text{m}/\text{s}, r = 4.0~\text{cm}, m=500
~\text{g}, \phi_B = 0, \theta_B = 0,
 I_3=I_1=2/5 m r^2, \\
 \mu & \| \mathbf{m}_i^{\text{inner, outer}} \| /(4 \pi)  = \mu \| \boldsymbol{\xi}_3 \| /(4 \pi)   = 2 \times 10^{-6} ~\text{T} \cdot \text{m}^3,  \\
 H_{\text{inner}} &=H_{\text{outer}} =48 ~\text{cm},
 N=20~ \text{dipoles},
 M=5~\text{dipoles}, \\
  r_{\text{inner}}&= 75 ~\text{cm},
 r_{\text{outer}}= 98 ~\text{cm},
 m_{\text{inner}} = 6~\text{kg}, J_3=m_{\text{inner}} r_{\text{inner}}^2, J = 1/2 J_3.
\end{align*}
In all of the simulations the magnetic top is initially at rest with
its attitude aligned to the vertical.   The difference in the
simulations is the amplitude of the oscillations and the discrete
sample from the normal distribution.  Figures
\ref{fig:ctmagmotorenergy}-\ref{fig:ctmagmotortheta} present data
from the uniform temperature simulations.  The key point about the
uniform temperature simulations is as seen in Figures
\ref{fig:ctmagmotorxyplots} and \ref{fig:ctmagmotortheta}  the
magnetic top transitions from noise to ballistic motion.  As
expected the non-uniform temperature case also exhibits this
behavior as shown in Figures
\ref{fig:nctmagmotorenergy}-\ref{fig:nctmagmotortheta}.

\begin{figure}[htbp]
\begin{center}
\includegraphics[scale=0.2,angle=0]{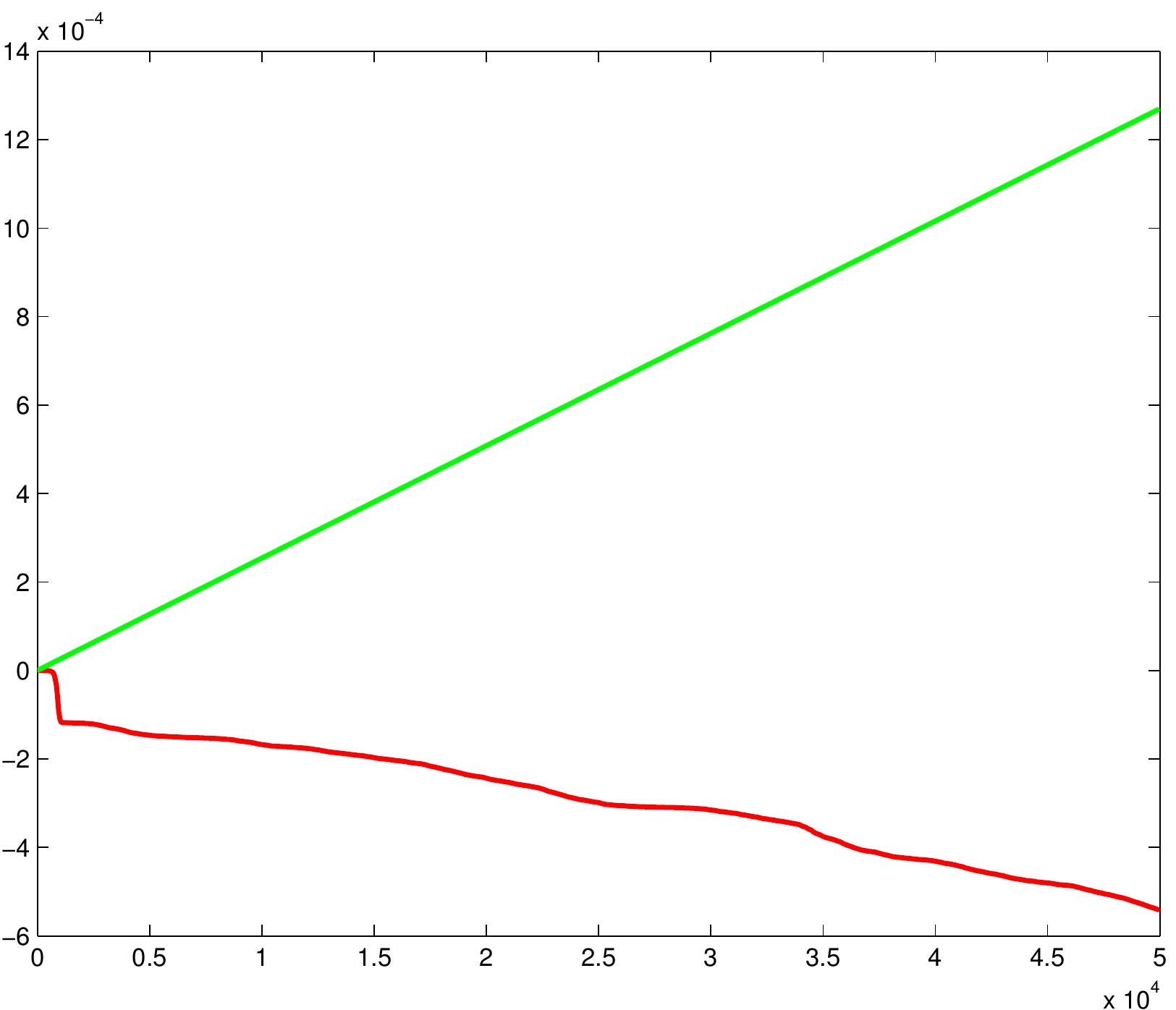}
\includegraphics[scale=0.2,angle=0]{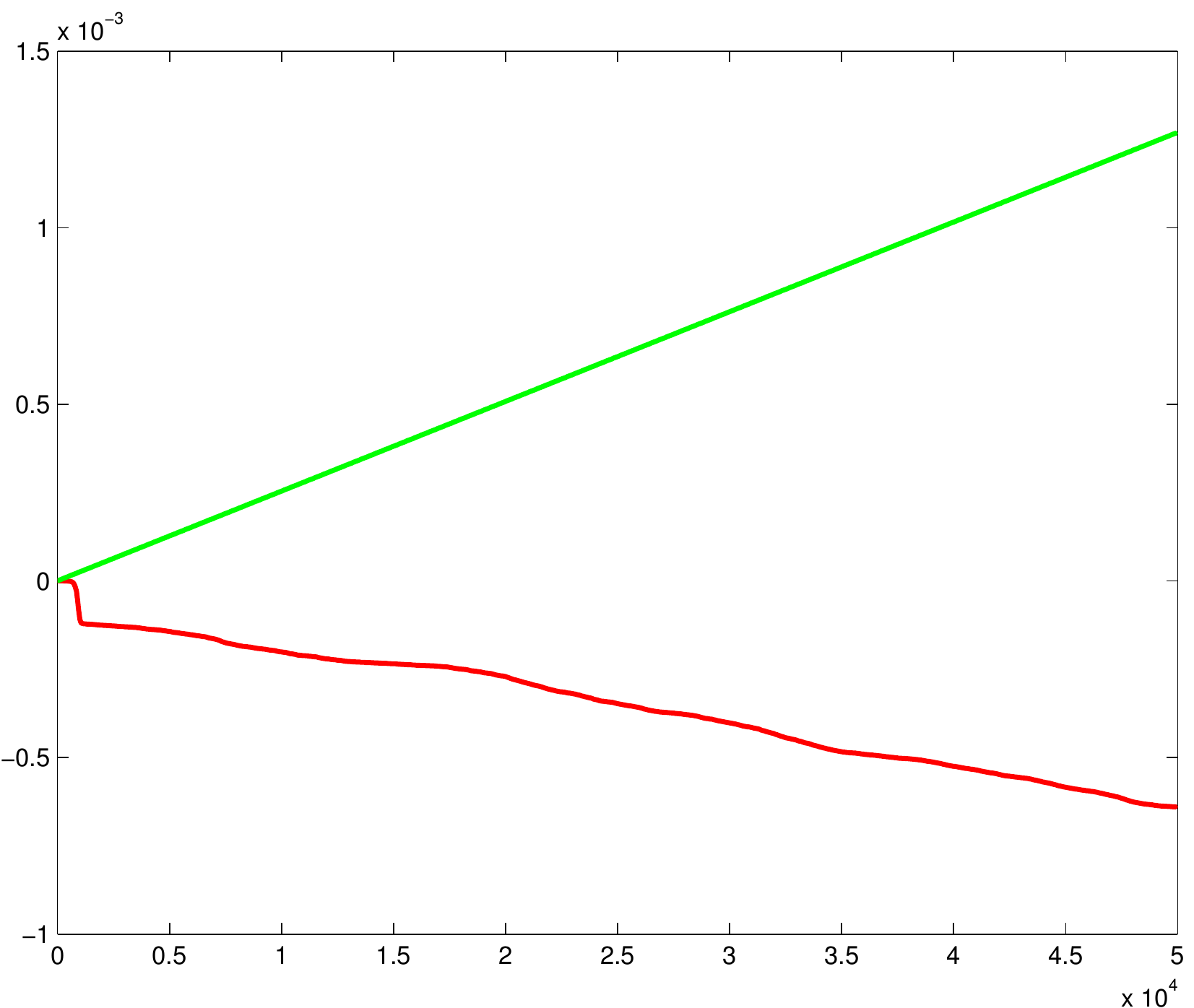}
\includegraphics[scale=0.2,angle=0]{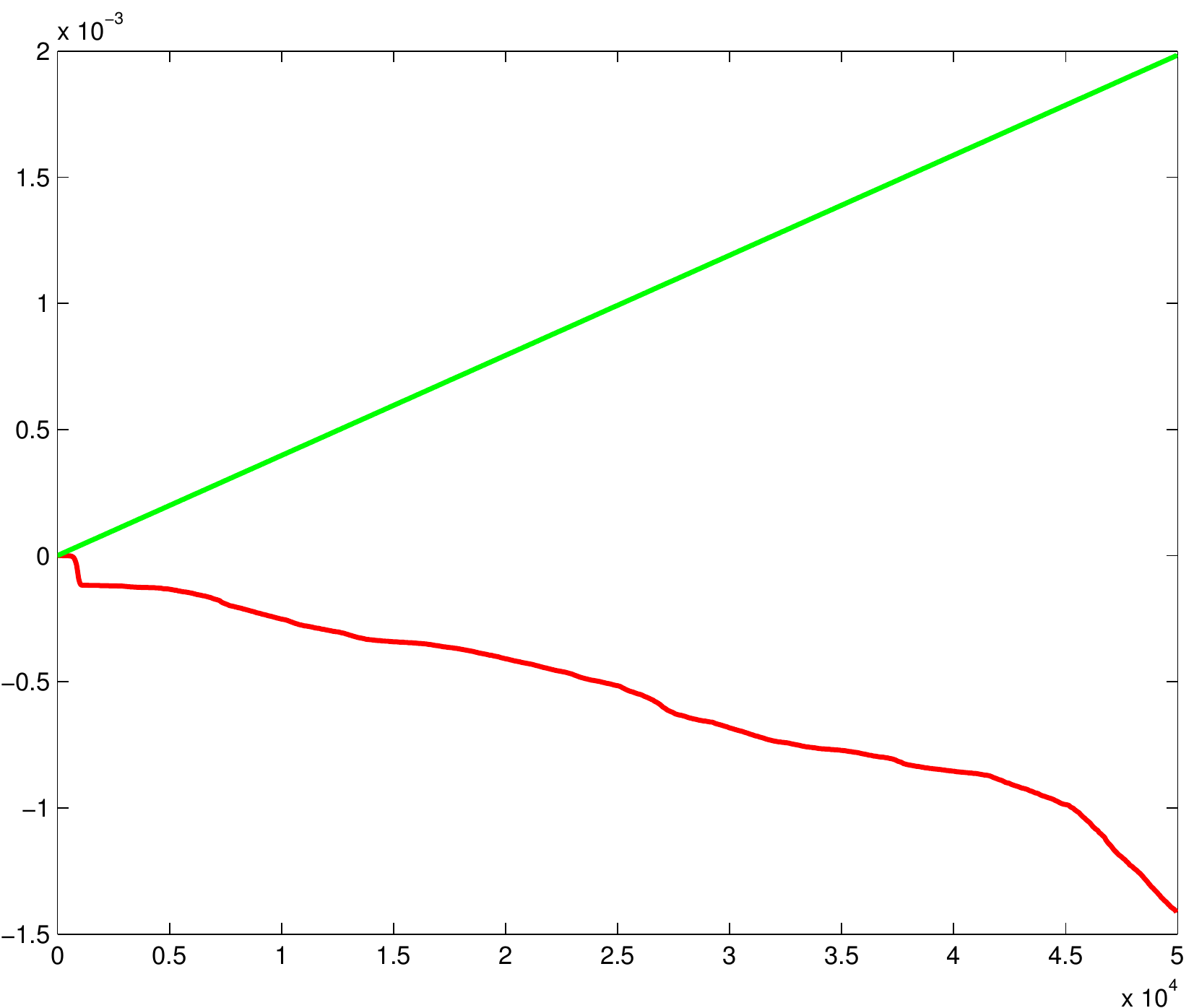}  \\
\hbox{\hspace{0.8in}   (a)  $\alpha=0.0002$ \hspace{0.3in} (b)
$\alpha=0.0002$ \hspace{0.3in} (c)  $\alpha=0.00025$  \hspace{0.3in}
} \caption{\footnotesize  {\bf Energy injected and dissipated
(non-uniform temperature).}   A plot of the energy injected by the
white noise as computed analytically and dissipated by friction
computed numerically. } \label{fig:nctmagmotorenergyloss}
\end{center}
\end{figure}

\begin{figure}[htbp]
\begin{center}
\includegraphics[scale=0.2,angle=0]{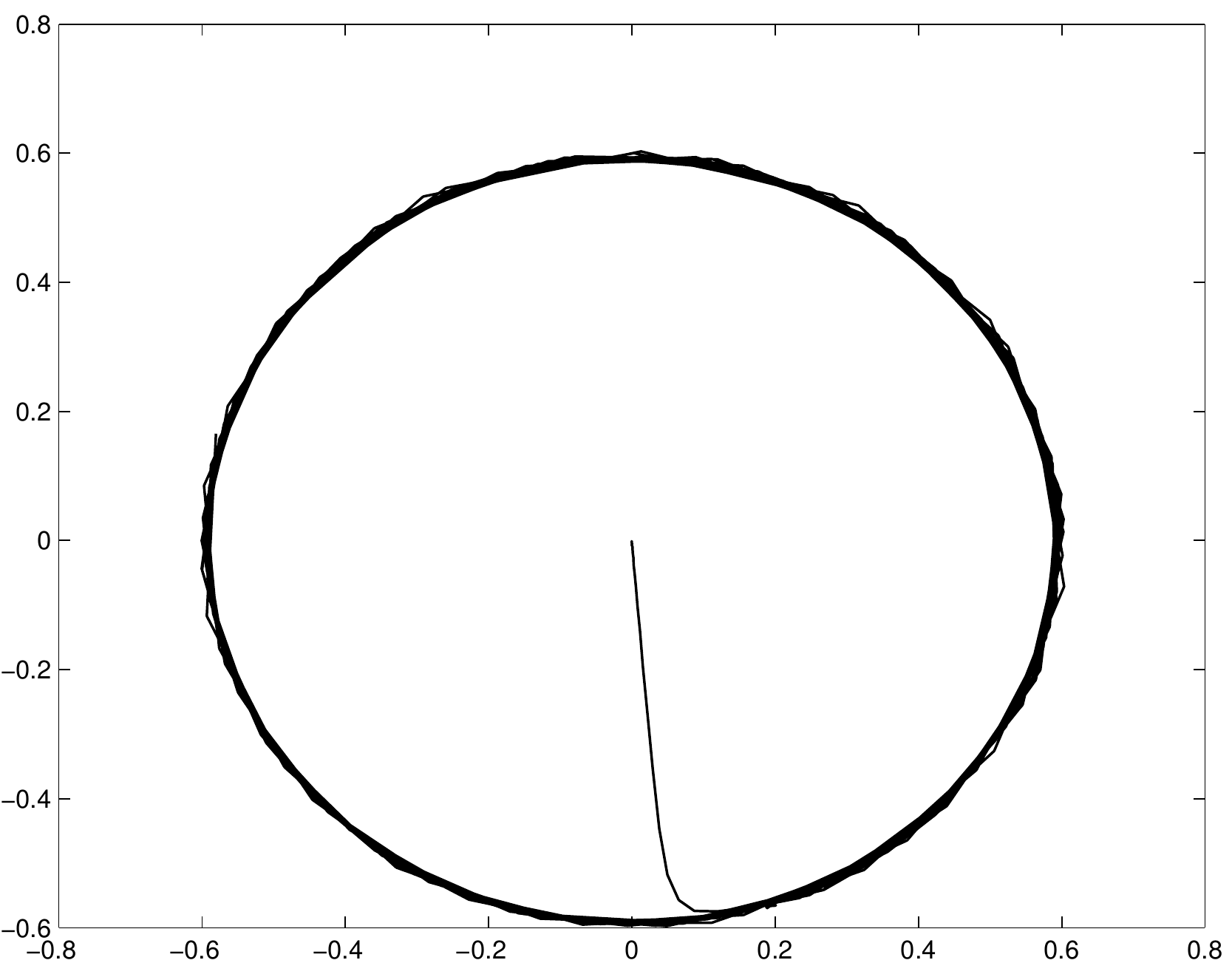}
\includegraphics[scale=0.2,angle=0]{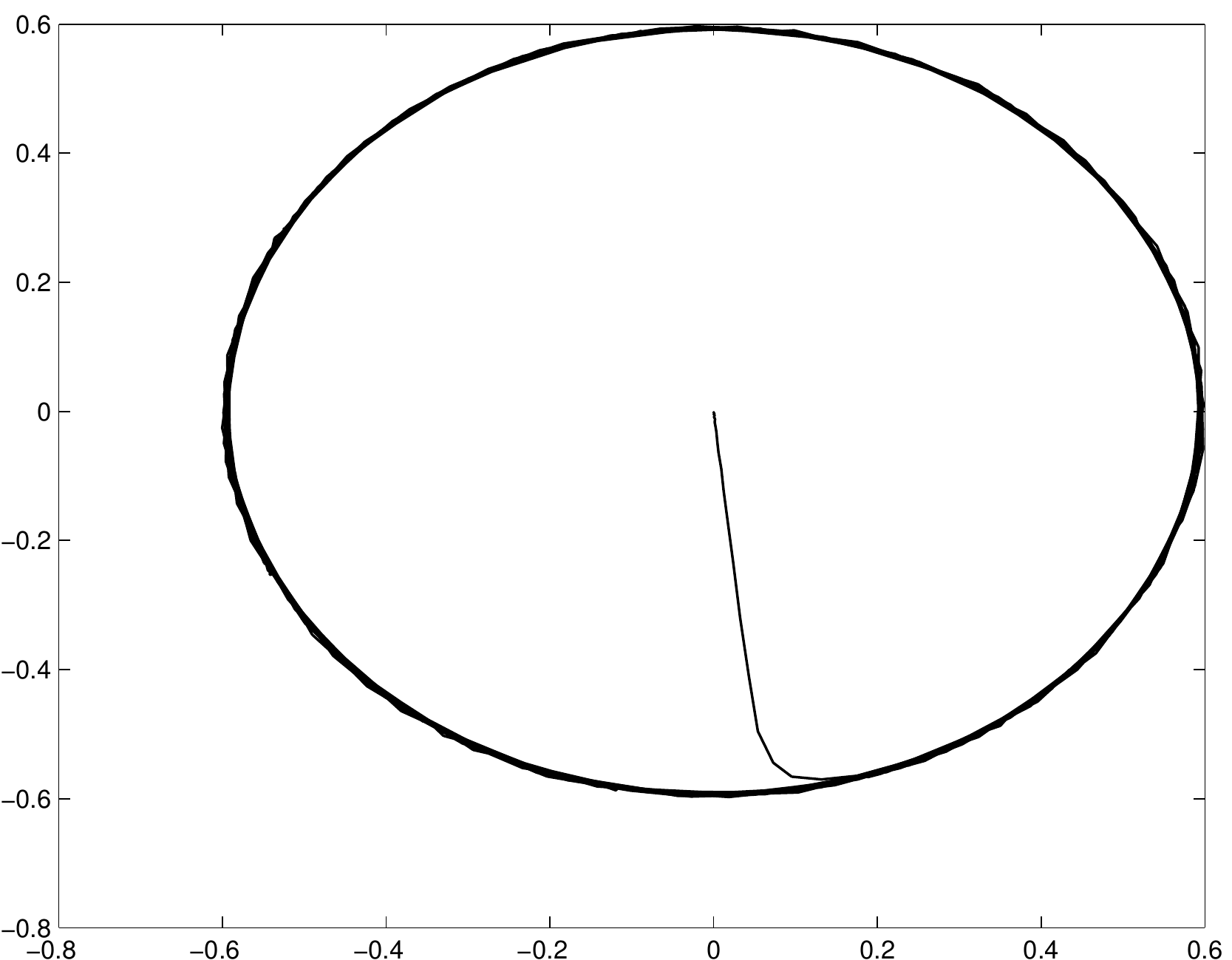}
\includegraphics[scale=0.2,angle=0]{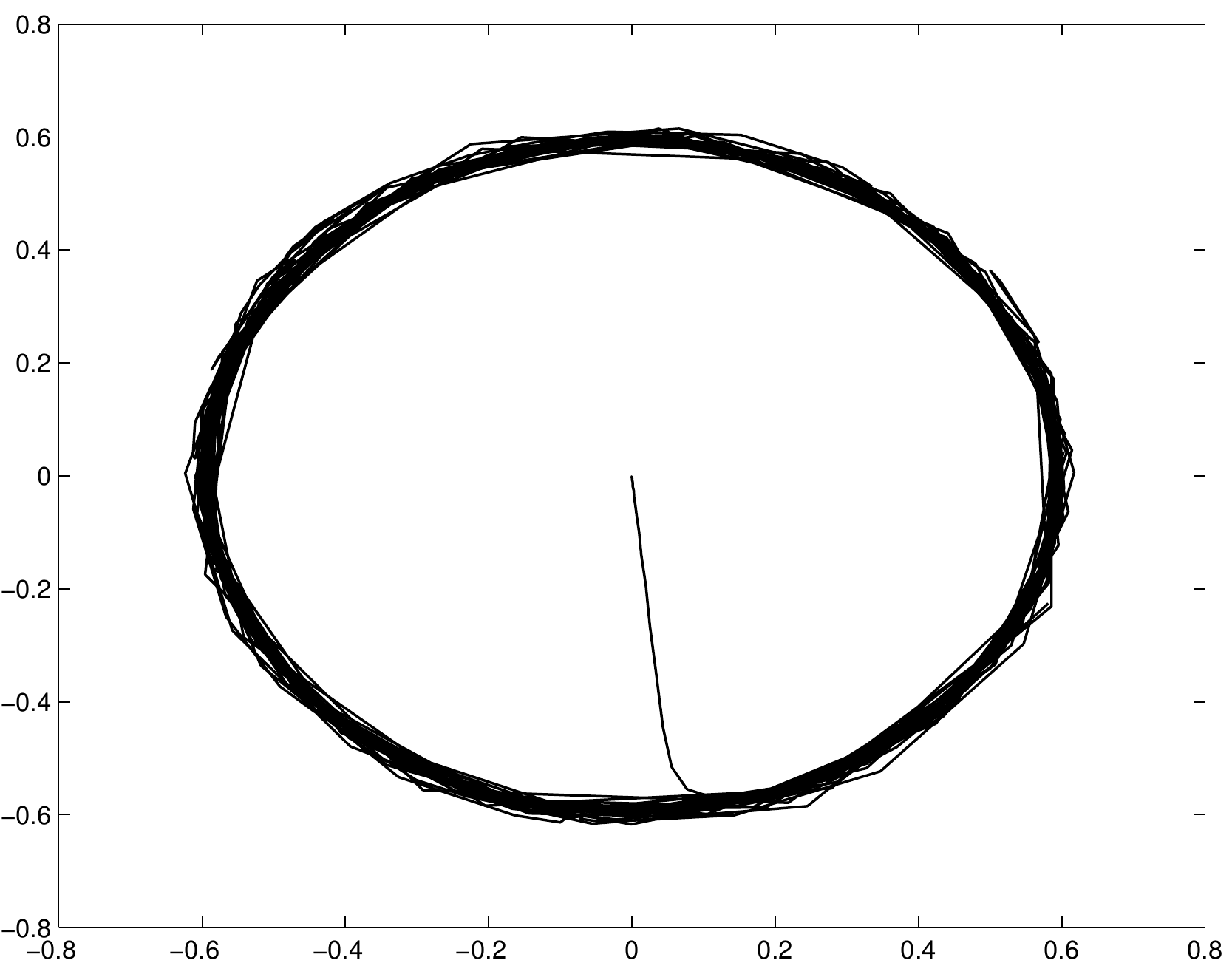} \\
\hbox{\hspace{0.8in}   (a)  $\alpha=0.0002$ \hspace{0.3in} (b)
$\alpha=0.0002$ \hspace{0.3in} (c)  $\alpha=0.00025$  \hspace{0.3in}
} \caption{\footnotesize  {\bf xy-position of magnetic ball
(non-uniform temperature).}  An aerial view of the path of the
center of mass of the ball in the xy-plane. }
\label{fig:nctmagmotorxyplot}
\end{center}
\end{figure}

\begin{figure}[htbp]
\begin{center}
\includegraphics[scale=0.2,angle=0]{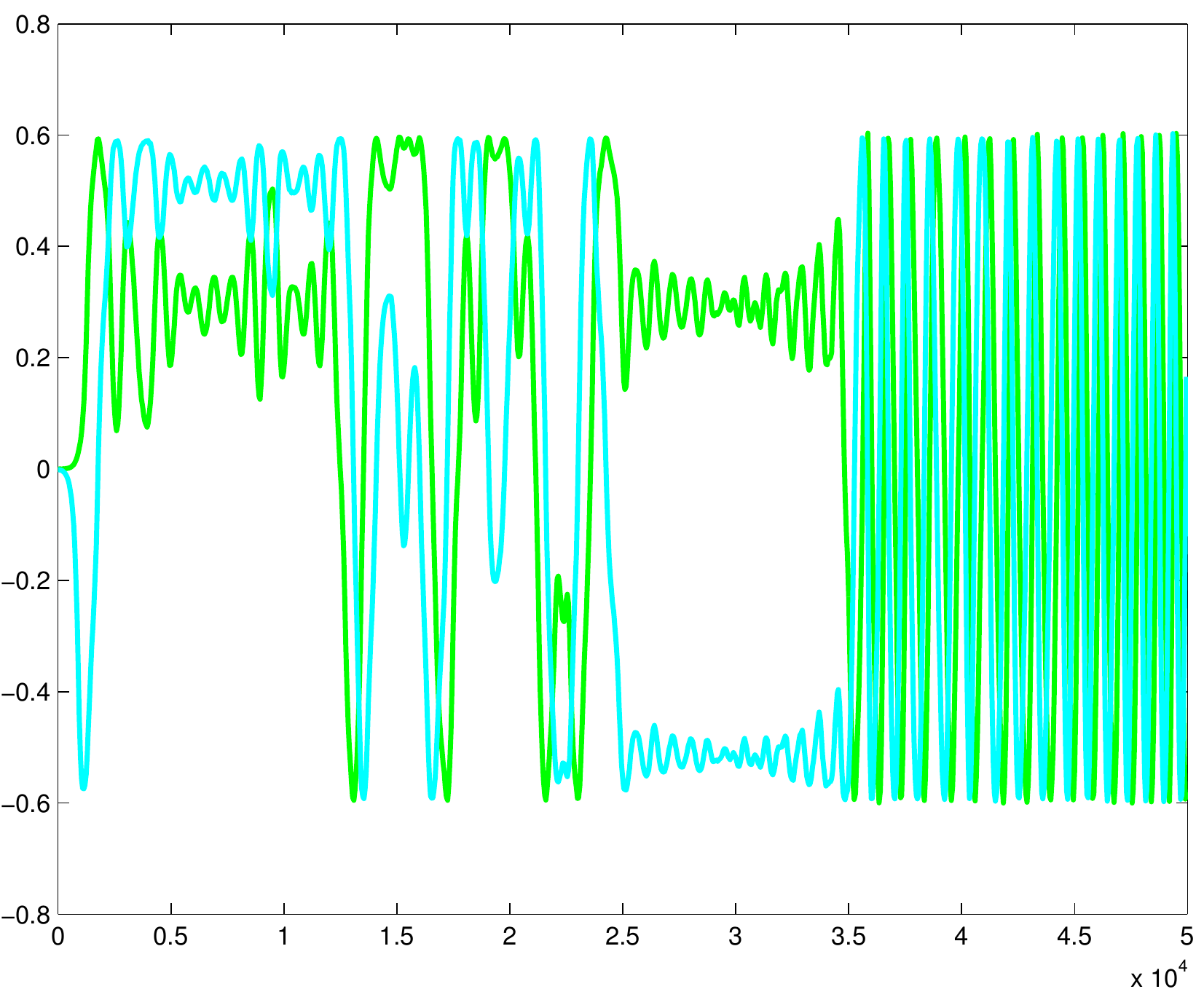}
\includegraphics[scale=0.2,angle=0]{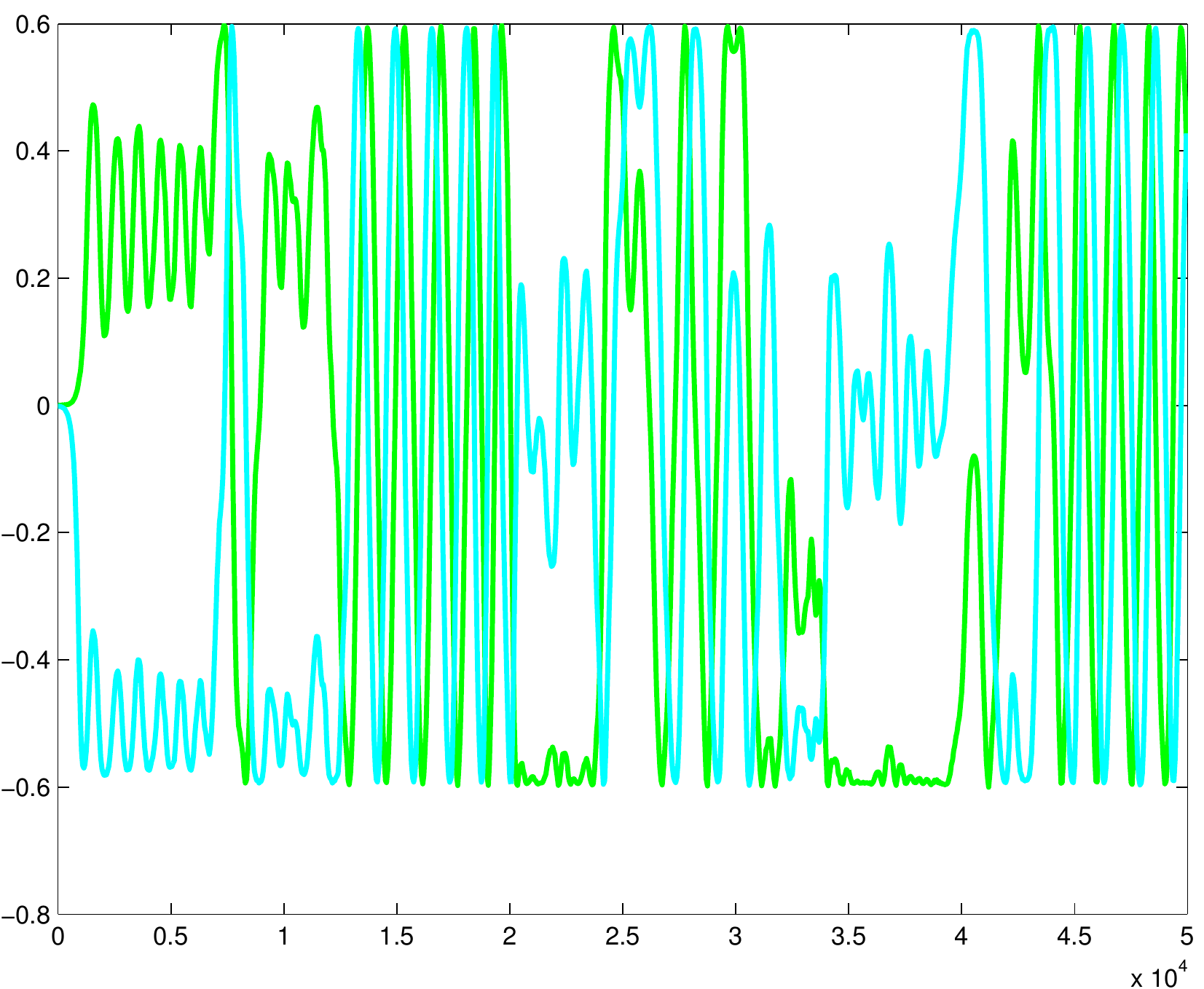}
\includegraphics[scale=0.2,angle=0]{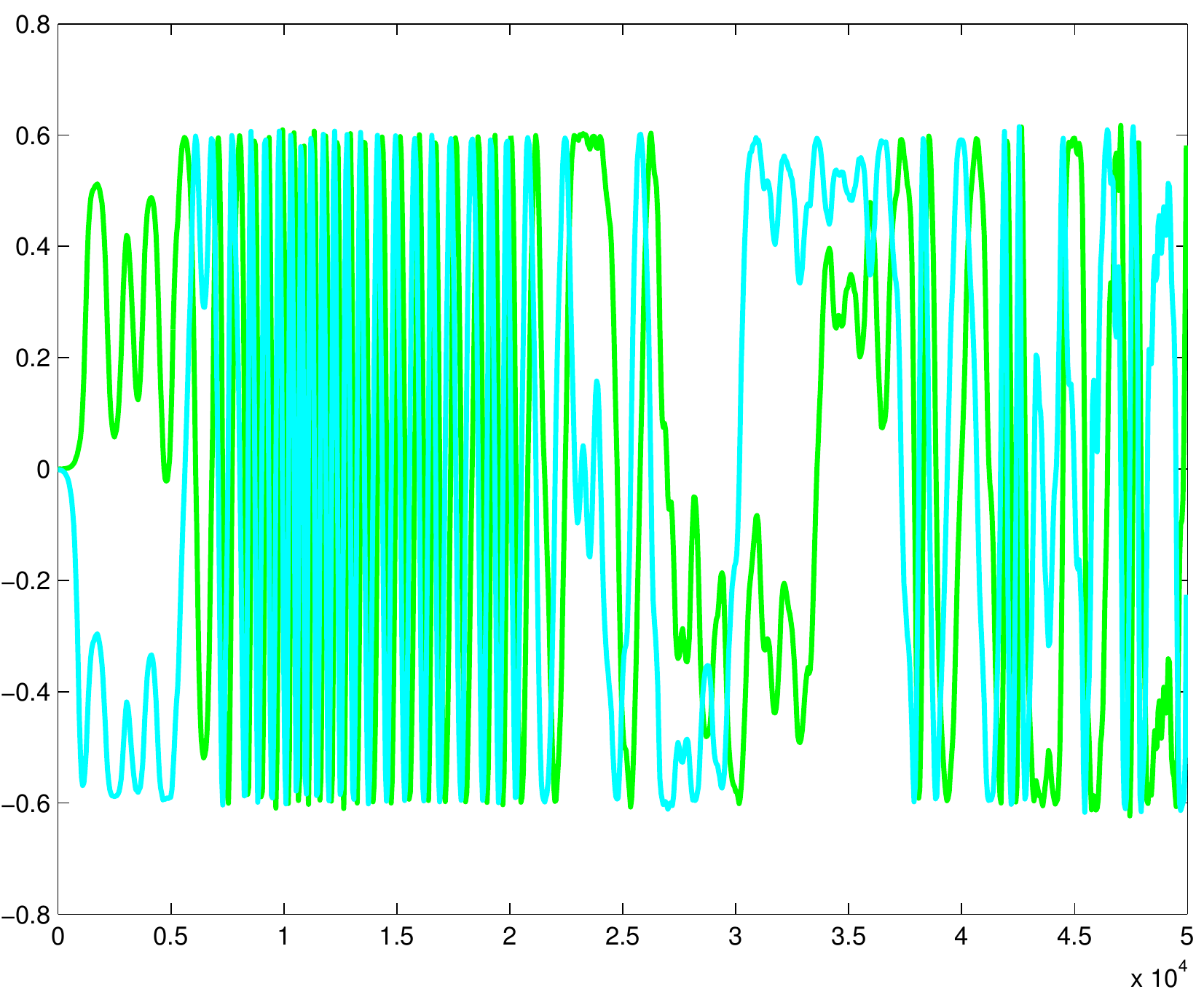}  \\
\hbox{\hspace{0.8in}   (a)  $\alpha=0.0002$ \hspace{0.3in} (b)
$\alpha=0.0002$ \hspace{0.3in} (c)  $\alpha=0.00025$  \hspace{0.3in}
} \caption{\footnotesize  {\bf xy-position of magnetic ball planar
view (non-uniform temperature).} The  $x$ and $y$ components of the
center of mass are plotted in green and blue respectively.   It is
very clear from this plot that the motion transitions from noise
driven to being inertia driven. } \label{fig:nctmagmotorxyplot}
\end{center}
\end{figure}

\begin{figure}[htbp]
\begin{center}
\includegraphics[scale=0.2,angle=0]{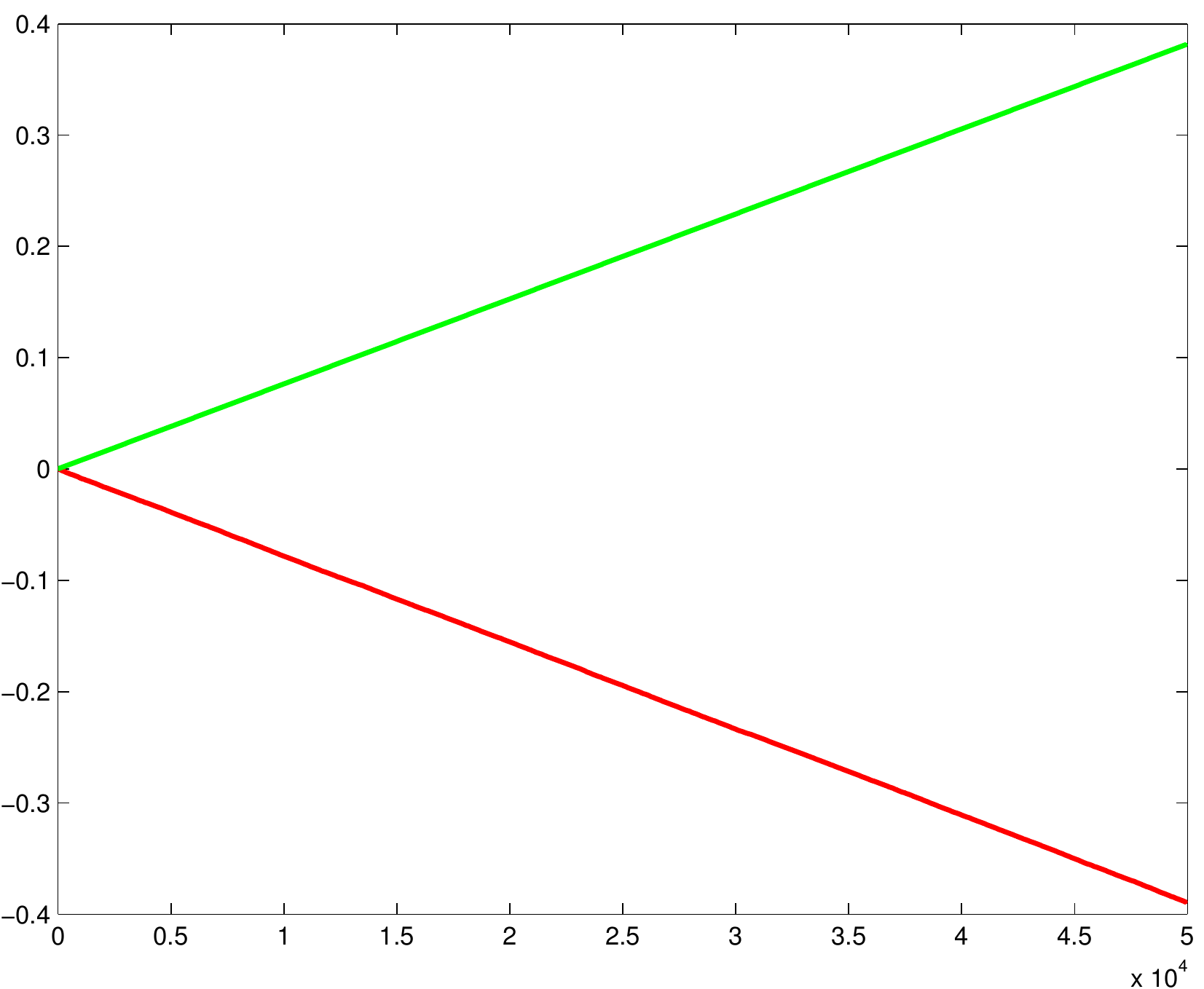}
\includegraphics[scale=0.2,angle=0]{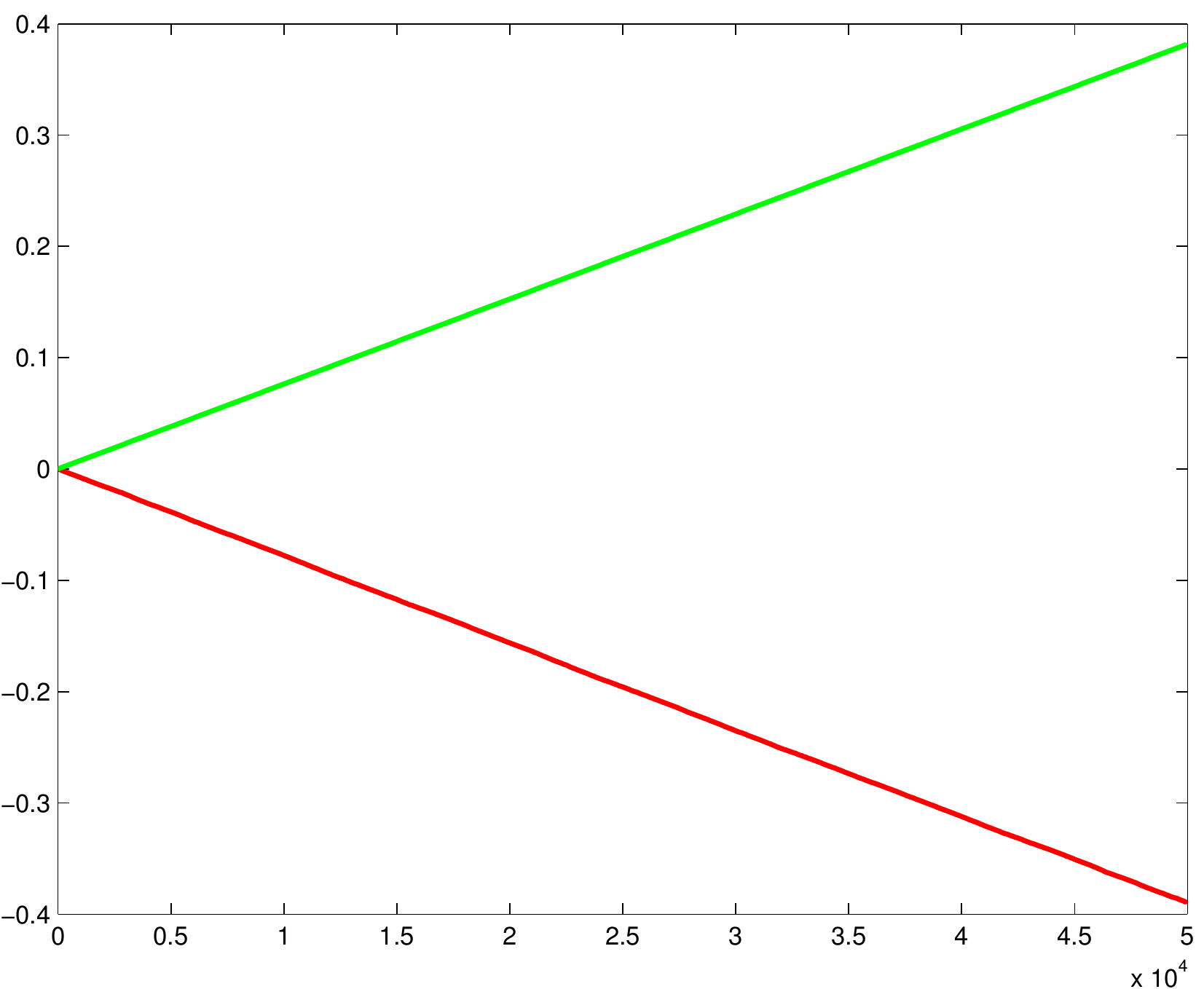}
\includegraphics[scale=0.2,angle=0]{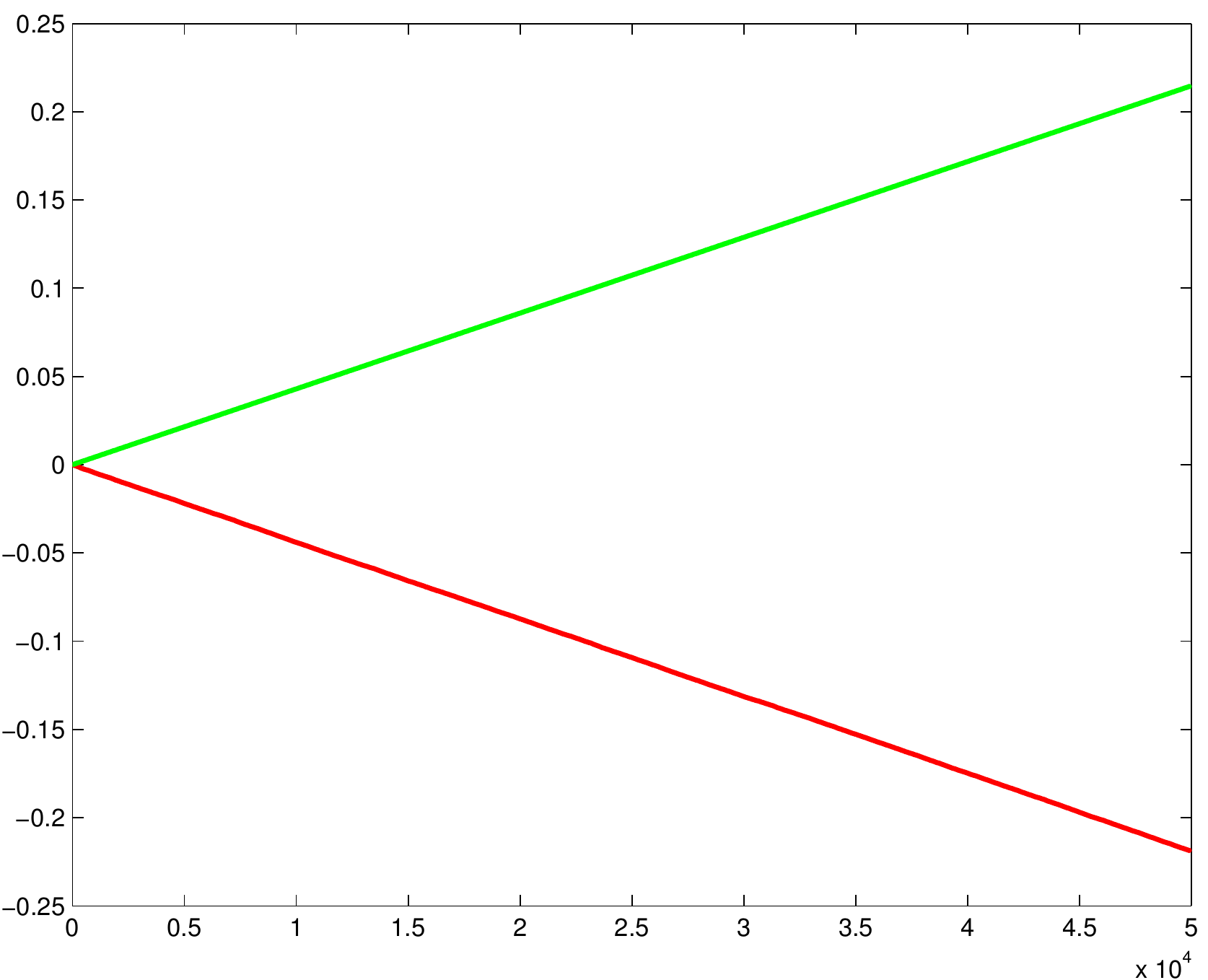}
\includegraphics[scale=0.2,angle=0]{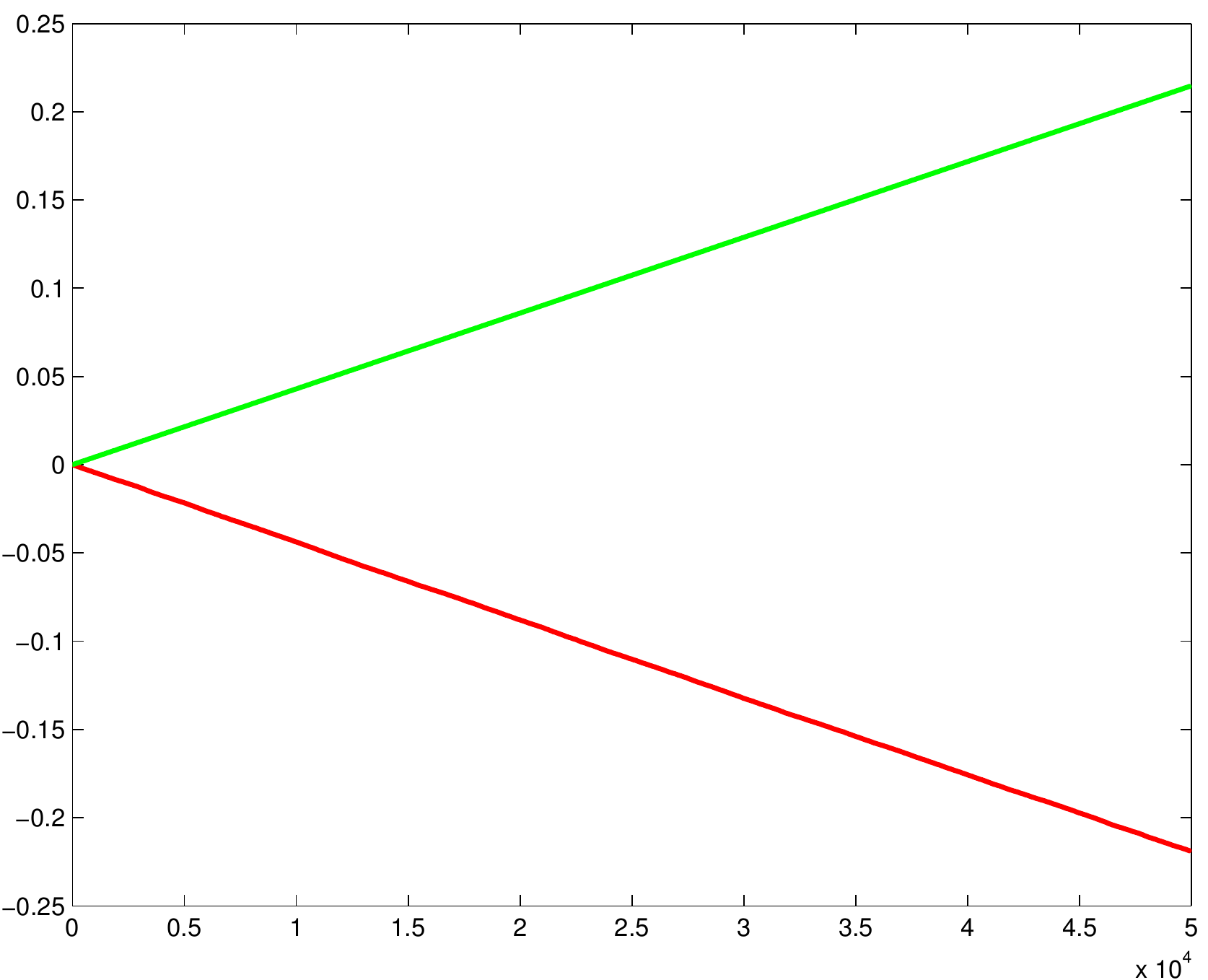} \\
\hbox{\hspace{0.15in}   (a)  $\alpha=0.001$ \hspace{0.3in} (b)
$\alpha=0.001$ \hspace{0.3in} (c)  $\alpha=0.00075$ \hspace{0.3in}
(d) $\alpha=0.00075$  \hspace{0.3in}    } \caption{\footnotesize
{\bf Energy injected and dissipated (uniform temperature).} A plot
of the energy injected by the white noise as computed analytically
and dissipated by friction computed numerically.  In the constant
temperature case we see an approximate balance between these
energies. } \label{fig:ctmagmotorenergyloss}
\end{center}
\end{figure}

\begin{figure}[htbp]
\begin{center}
\includegraphics[scale=0.2,angle=0]{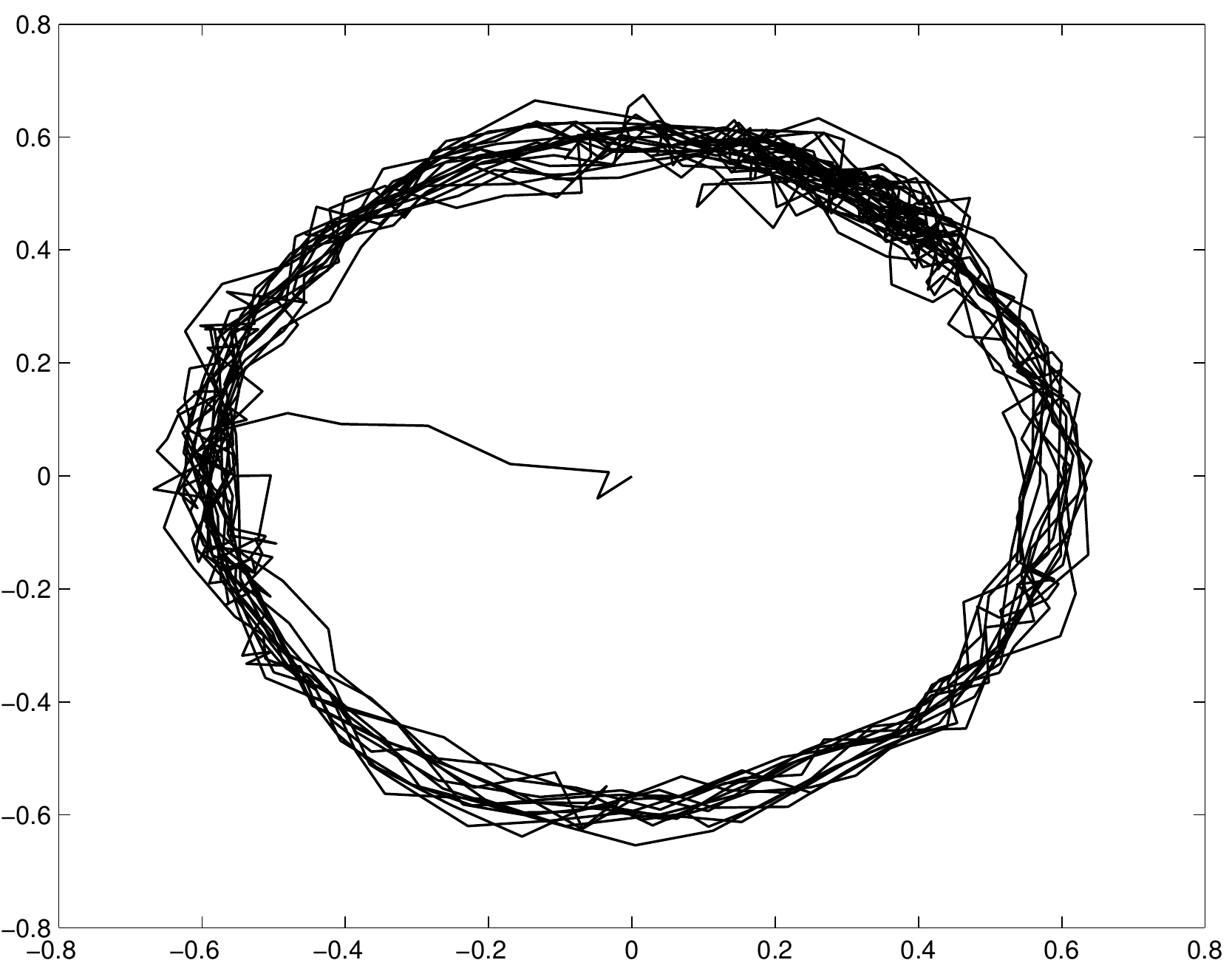}
\includegraphics[scale=0.2,angle=0]{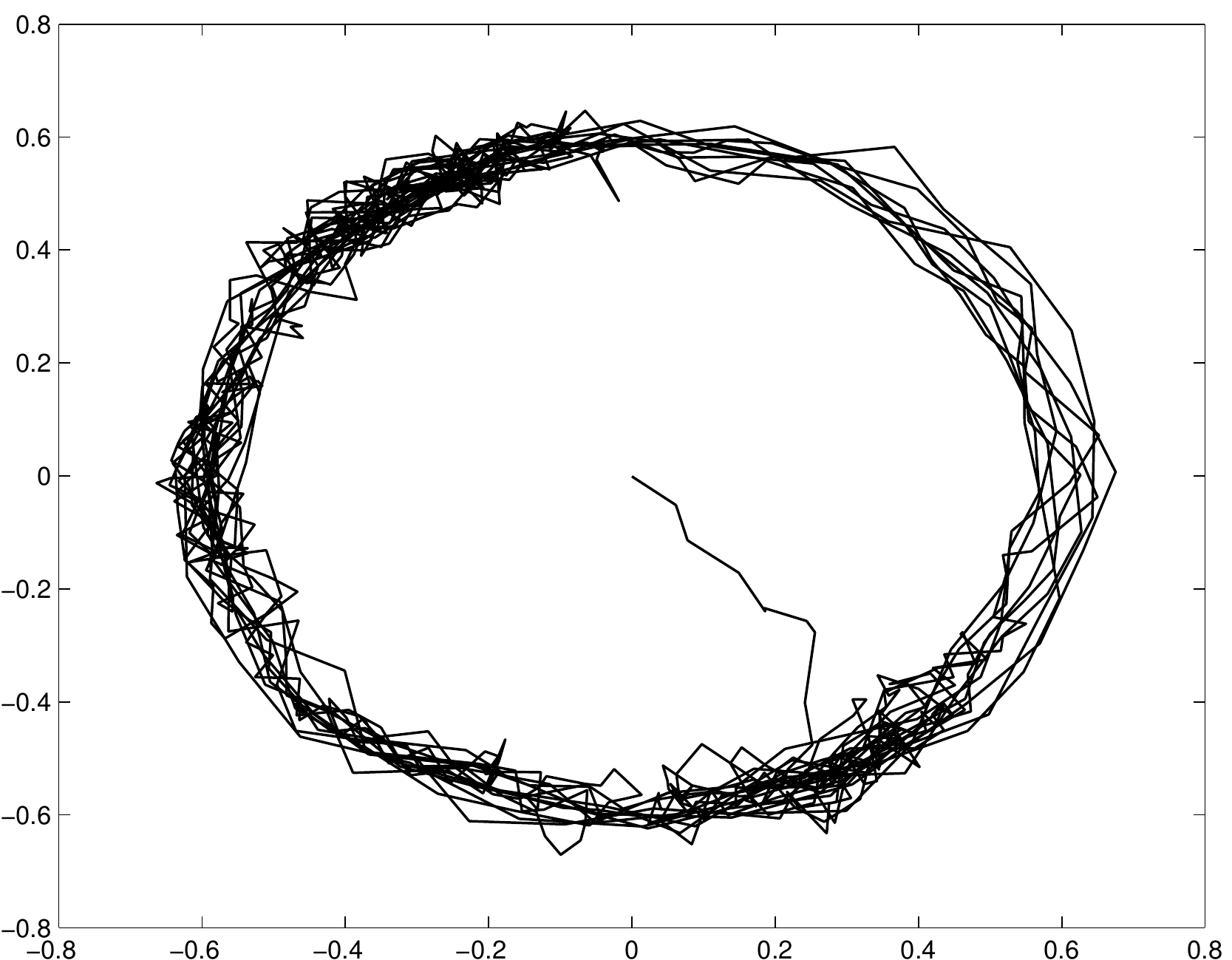}
\includegraphics[scale=0.2,angle=0]{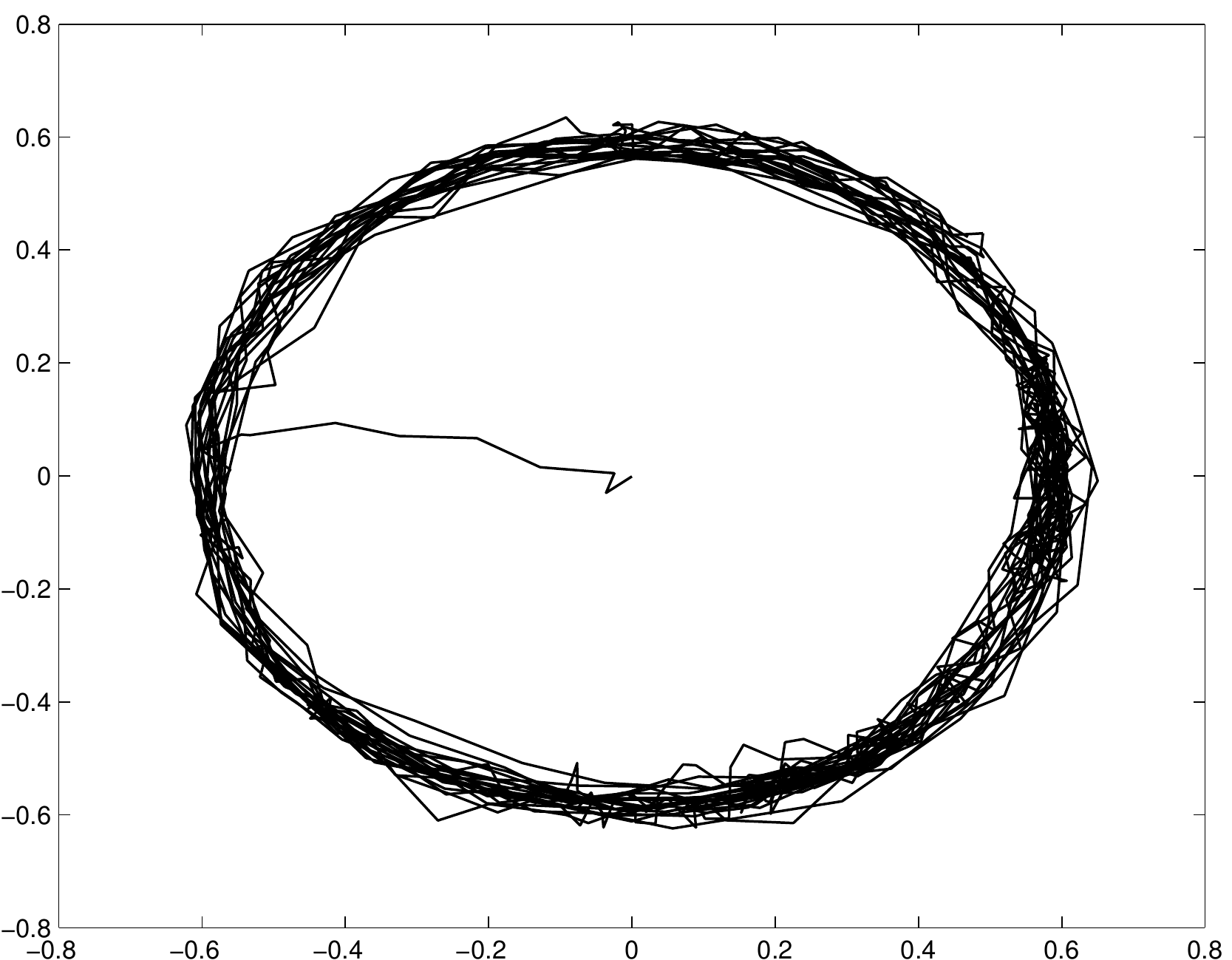}
\includegraphics[scale=0.2,angle=0]{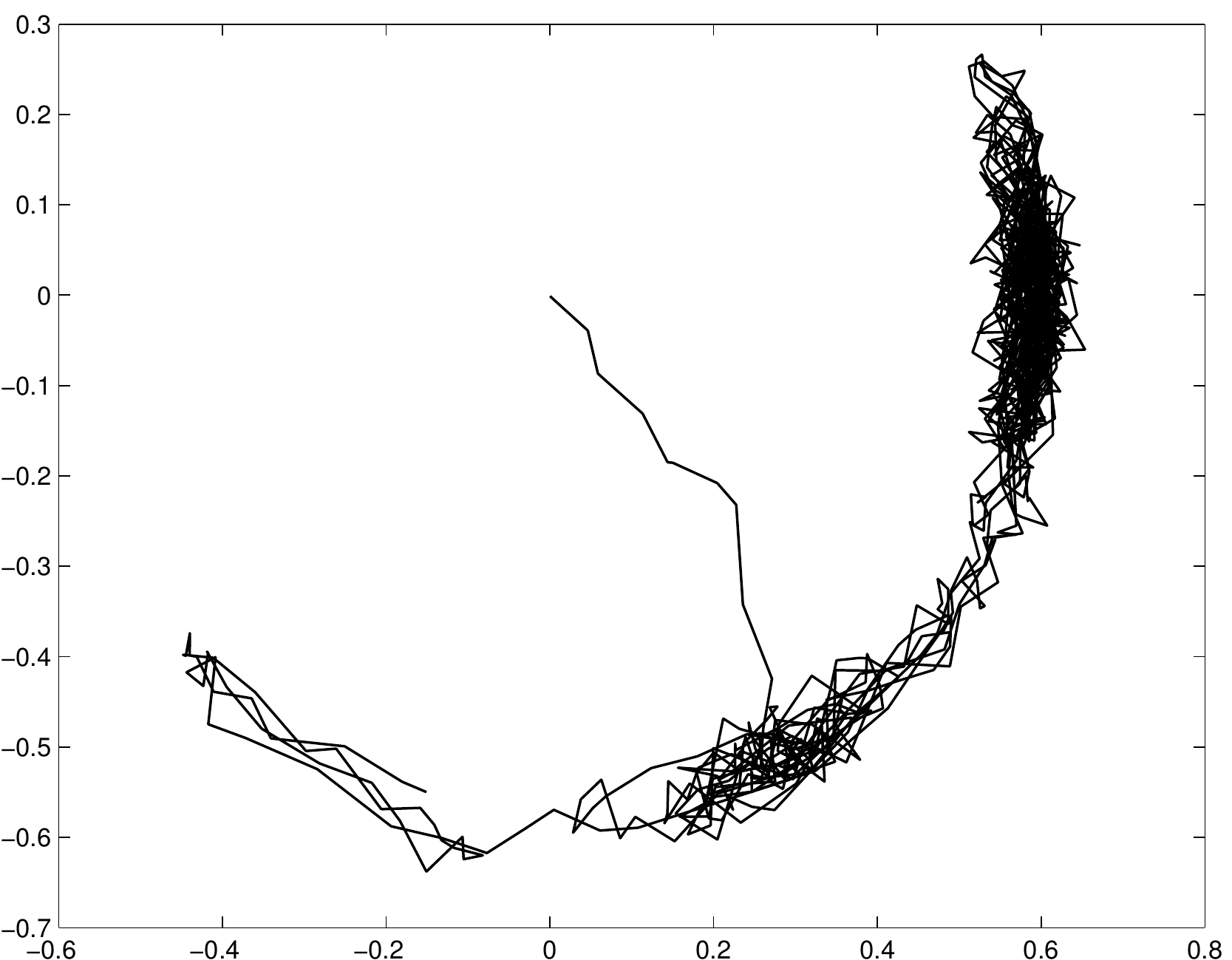} \\
\hbox{\hspace{0.15in}   (a)  $\alpha=0.001$ \hspace{0.3in} (b)
$\alpha=0.001$ \hspace{0.3in} (c)  $\alpha=0.00075$ \hspace{0.3in}
(d) $\alpha=0.00075$  \hspace{0.3in}    } \caption{\footnotesize
{\bf xy-position of magnetic ball (uniform temperature).}  An aerial
view of the path of the center of mass of the ball in the xy-plane.
} \label{fig:ctmagmotorxyplot}
\end{center}
\end{figure}

\newpage

\section{Conclusion}

This paper introduces a novel mechanism to rectify random perturbations to achieve
directed motion as a meta-stable state and ballistic mean-squared displacement with respect to the invariant law.   The basic idea behind the mechanism is explained using the simple sliding disk model.  With this model one can prove that if the potential energy is non-constant, then the invariant Gibbs measure of the system, $\mu$, is ergodic and strongly mixing.  As a consequence $\mu$ a.s.~it is shown that the translational displacement of the sliding disk is not ballistic. However, it is shown that the mean-squared displacement with respect to the invariant law is ballistic.

This anomalous diffusion manifests itself numerically.  That is, starting from a position at rest and when $U$ is non-constant, one observes flights in the $x$-displacement over time-scales which are random and the probability distribution of the flight-durations have a heavy tail.  Moreover, from a numerical statistical analysis, it is observed that the mean-squared translational displacement of the sliding disk is ballistic and its translational displacement is characterized by long time memory/correlation.  The paper proceeds to show that this basic phenomenon arises in a more complex, rigid-body system consisting of two rigid bodies interacting via magnetostatic effects.  Along the way we explain the observed dynamics of Hamel's magnetic device.

\section{Appendix}

\paragraph{$\operatorname{SO}(3)$ Preliminaries}

In the paper we use the following isomorphism between an element of
the Lie algebra of $\operatorname{SO}(3)$,
$T_e\operatorname{SO(3)}=\mathfrak{so}(3)$, and $\mathbb{R}^3$.
Recall that elements of $\mathfrak{so}(3)$ are skew-symmetric
matrices with Lie bracket given by the matrix commutator.  Let
$\boldsymbol{\omega} = (\omega_1, \omega_2, \omega_3)$.  Then one
can relate $\mathbb{R}^3$ with a skew-symmetric matrix via the hat
map   $~~\widehat{}:\mathbb{R}^3 \rightarrow so(3)$,
\[
\widehat{\boldsymbol{\omega}} = \begin{bmatrix} 0 & - \omega_3 & \omega_2  \\
                                  \omega_3 & 0  & -\omega_1 \\
                                -\omega_2  & \omega_1 & 0 \end{bmatrix} \text{.}
\]
Let $g(t)$ be a curve in  $\operatorname{SO}(3)$.  With this
identification of $\mathfrak{so}(3)$ to $\mathbb{R}^3$, the
right-trivialization of a tangent vector $\dot{g}$ to this curve,
given by  $\xi = T R_{g^{-1}} \cdot \dot{g} \in \mathfrak{so}(3)$,
can be written in terms of the spatial angular velocity vector
$\boldsymbol{\omega} \in \mathbb{R}^3$, i.e.,  as $\xi =
\widehat{\boldsymbol{\omega}} \in so(3)$.

\paragraph{Magnetic motor at non-uniform temperature governing equations}

Below the governing equations of the ring-ball system with
nonconservative effects due to surface friction and white noise are
written in Ito form.   We introduce the potential energy function
$U:\mathbb{R}^3 \times \operatorname{SO}(3) \times
\operatorname{SO}(3) \to \mathbb{R}$ which represents the total
potential energy of the ring-ball system.  In terms of $U$ the
governing equations can be written as,
\begin{align*}
 d \mathbf{x} =& \mathbf{v} dt \\
 d \mathbf{v} =& \frac{1}{m} \left( -  \frac{ \partial U}{\partial \mathbf{x} } + ( \lambda  - mg ) \mathbf{e}_3 - c \mathbf{V}_Q  \right) dt \\
\mathbf{x}^\mathrm{T} \mathbf{e}_3 =& r \\
     d  R_B =&   \widehat{\boldsymbol{\omega}_B} R_B dt \\
       d    R_R =&   \widehat{\boldsymbol{\omega}_R} R_R dt \\
      \boldsymbol{\pi}_B =& J_B \boldsymbol{\omega}_B  \\
      \boldsymbol{\pi}_R =& R_R \mathbb{I}_R R_R^{\mathrm{T}} \boldsymbol{\omega}_R  \\
d \widehat{\boldsymbol{\pi}_B}  =&    \left( -  \frac{\partial U}{\partial R_B} R_B^{\mathrm{T}} - c \widehat{ \widehat{\mathbf{q}} \mathbf{V}_Q } \right) dt   \\
d \widehat{\boldsymbol{\pi}_R}  =&   \left(  - \frac{\partial
U}{\partial R_R} R_R^{\mathrm{T}} \right) dt + \alpha \widehat{d
\mathbf{W}}
 \end{align*}
Evaluating these governing equations at $\ell$ defined in
(\ref{eq:ellringnball}) yields
 \begin{align}
 d \mathbf{x} =& \mathbf{v} dt \\
d \mathbf{v} =& \frac{1}{m} \left( 2 \sum_{i=1}^{N}  \mathbf{D} \mathbf{B}_i^{\text{inner}}( \mathbf{r}_{i0} )^\mathrm{T}  \boldsymbol{\xi}_3  + \sum_{j=1}^{M}   \mathbf{D} \mathbf{B}_j^{\text{outer}}( \mathbf{r}_{0j} )^\mathrm{T}   \boldsymbol{\xi}_3  \right) dt \nonumber \\
& + \frac{1}{m} \left(- c \mathbf{V}_Q  + (\lambda - mg ) \mathbf{e}_3 \right) dt \label{eq:eulerlagrangeb} \\
\mathbf{x}^\mathrm{T}  \mathbf{e}_3 &= r \label{eq:constraintballb} \\
      d R_B =&   \widehat{\boldsymbol{\omega}_B} R_B dt \label{eq:reconstructionballb} \\
            d R_R =&   \widehat{\boldsymbol{\omega}_R} R_R dt \label{eq:reconstructionringb} \\
      \boldsymbol{\pi}_B =& J_B \boldsymbol{\omega}_B \label{eq:legendreballb} \\
      \boldsymbol{\pi}_R =& R_R \mathbb{I}_R R_R^{\mathrm{T}} \boldsymbol{\omega}_R \label{eq:legendreringb} \\
d \boldsymbol{\pi}_B  =&  \widehat{\boldsymbol{\xi}_3} \left( 2
\sum_{i=1}^{N} \mathbf{B}_i^{\text{inner}}(  \mathbf{r}_{i0}) +
\sum_{j=1}^{M}   \mathbf{B}_j^{\text{outer}}(  \mathbf{r}_{0j}) \right) dt - c \widehat{\mathbf{q}} \mathbf{V}_Q  dt \label{eq:lpballb} \text{,} \\
d \boldsymbol{\pi}_R  =&  2 \sum_{i=1}^N \left[ \widehat{ \mathbf{m}_i^{\text{inner}}} \mathbf{B}_0 (  \mathbf{r}_{i0} )- \widehat{ \mathbf{d}_i^{\text{inner}}}   \left( \mathbf{D} \mathbf{B}_i^{\text{inner}}( \mathbf{r}_{i0} ) ^\mathrm{T}  \boldsymbol{\xi}_3 \right) \right] dt \nonumber \\
&+ \sum_{i=1}^N \sum_{j=1}^M \left[ \widehat{ \mathbf{m}_i^{\text{inner}} }  \mathbf{B}_j^{\text{outer}} \left(  \mathbf{r}_{ij} \right) + \widehat{ \mathbf{d}_i^{\text{inner}}}  \left( \mathbf{D}  \mathbf{B}_j^{\text{outer}} (   \mathbf{r}_{ij} ) ^\mathrm{T}  \mathbf{m}_i^{\text{inner}} \right)  \right] dt \nonumber \\
&+ \alpha d \mathbf{W} \label{eq:lpringb} \text{.}
 \end{align}
 Since the ring is axisymmetric its Legendre transform simplifies:
 \begin{align}
\boldsymbol{\pi}_R &=  J \boldsymbol{\omega}_R + (J_3 - J)  (
\boldsymbol{\omega}_R^{\mathrm{T}}  \boldsymbol{\zeta}_3 )
\boldsymbol{\zeta}_3  \implies \boldsymbol{\omega}_R  = \frac{1}{J}
\left( \boldsymbol{\pi}_R + \frac{J - J_3}{J_3} (
\boldsymbol{\pi}_R^{\mathrm{T}}  \boldsymbol{\zeta}_3 )
\boldsymbol{\zeta}_3 \right) \label{eq:spatialltring} \text{.}
 \end{align}

\begin{Remark}
The energy of the magnetic motor is given by
\[
E(\mathbf{x}, \mathbf{v}, R_B, \boldsymbol{\pi}_B, R_R,
\boldsymbol{\pi}_R) =  \frac{m}{2} \mathbf{v}^{\mathrm{T}}
\mathbf{v} + \frac{1}{2} J_B^{-1} \boldsymbol{\pi}_B^{\mathrm{T}}
\boldsymbol{\pi}_B + \frac{1}{2} \boldsymbol{\pi}_R^{\mathrm{T}} R_R
\mathbb{I}_R^{-1} R_R^{\mathrm{T}}  \boldsymbol{\pi}_R +
U(\mathbf{x}, R_B, R_R)    \text{.}
\]
Using Ito's formula one can calculate the stochastic differential of the
energy in the case when the magnetic motor is at non-uniform temperature
\[
d E = \left( - c \mathbf{V}_Q^{\mathrm{T}}  \mathbf{V}_Q +
\frac{1}{2} \alpha^2 \operatorname{trace} \left[ \mathbb{I}_R^{-1}
\right]  \right) dt + \alpha \boldsymbol{\omega}_R^{\mathrm{T}}  d
\mathbf{W} \text{.}
\]
Integrating yields,
\[
E(t) - E(0) = - c \int_0^t \mathbf{V}_Q^{\mathrm{T}}  \mathbf{V}_Q
ds  + \frac{1}{2} \alpha^2 \operatorname{trace} \left[
\mathbb{I}_R^{-1} \right] t + \text{martingales.}
\]
The first term represents the energy dissipated by friction.  The
next terms represent the work done by the white noise torque as
computed by Ito's integral.
\end{Remark}

\paragraph{Isothermal, magnetic motor governing equations}

Similar to the sliding disk, to put the magnetic motor at constant
temperature we define the following dissipation matrix,
\[
\mathbf{C} = \begin{bmatrix} 1/m^2 & 0 & 0  & - r/(m J) \\
                               0  & 1/m^2 & r/(m J) & 0 \\
                                 0 & r/(m J) & r^2/J^2 & 0 \\
                                 -r/(m J) & 0 & 0 & r^2/J^2 \end{bmatrix}
\]
The translational position of the ball is written in coordinates as
$\mathbf{x} = (x_1, x_2, r)$.  The translational and angular
velocities of the ball are given by $\mathbf{v} = (v_1, v_2, 0)$ and
$\boldsymbol{\omega}_B = ( \omega_B^{(1)} ,  \omega_B^{(2)} ,
\omega_B^{(3)} )$.  The dynamical equations for the constant
temperature magnetic motor are given by:
 \begin{align}
 d x_1 &= v_1 dt \label{eq:sdea} \\
 d x_2 &= v_2 dt \label{eq:sdeb} \\
d R_B &=  \widehat{\boldsymbol{\omega}_B}  R_B dt \label{eq:sdec}  \\
d R_R &=  \widehat{\boldsymbol{\omega}_R}  R_R dt \label{eq:sded}  \\
\begin{bmatrix} d v_1 \\ d v_2 \\   d  \omega_B^{(1)} \\ d \omega_B^{(2)} \end{bmatrix} &=
-\begin{bmatrix} U_x^{\mathrm{T}} \mathbf{e}_1/m \\ U_x^{\mathrm{T}}
\mathbf{e}_2/m  \\ U_B^{\mathrm{T}} \mathbf{e}_1/J \\
U_B^{\mathrm{T}} \mathbf{e}_2/J \end{bmatrix} dt - c_B \mathbf{C}
\begin{bmatrix}  m v_1 \\ m v_2 \\    J \omega_B^{(1)} \\ J
\omega_B^{(2)}  \end{bmatrix} +
\alpha_B \mathbf{C}^{1/2} \begin{bmatrix} d B_{v_1} \\ d B_{v_2}  \\ d B_{\omega_B^{(1)}} \\ d B_{\omega_B^{(2)}} \end{bmatrix} \label{eq:sdee} \\
d \omega_B^{(3)} &= - U_B^{\mathrm{T}} \mathbf{e}_3/J \label{eq:sdef} \\
\boldsymbol{\pi}_R &= R_R \mathbb{I}_R R_R^{\mathrm{T}} \boldsymbol{\omega}_R \label{eq:sdeg} \\
d \boldsymbol{\pi}_R &= - U_R - c_R \boldsymbol{\omega}_R + \alpha_R
d \mathbf{B}_R \label{eq:sdeh}
 \end{align}

 The terms $U_x$, $U_B$, and $U_R$  are defined in terms of the inner product on $\mathbb{R}^3$ as,
 \begin{align*}
 U_x^{\mathrm{T}} y &= \left\langle \frac{\partial U}{\partial x}, y \right\rangle  = \partial_x U(x, R_B, R_R) \cdot y \\
 U_B^{\mathrm{T}} y &= \left\langle \frac{\partial U}{\partial R_B} R_B^{\mathrm{T}}, \widehat{y} \right\rangle =   \partial_{R_B} U(x, R_B, R_R) \cdot \widehat{y} R_B \\
U_R^{\mathrm{T}} y &= \left\langle \frac{\partial U}{\partial R_R}
R_R^{\mathrm{T}}, \widehat{y} \right\rangle =   \partial_{R_R} U(x,
R_B, R_R) \cdot \widehat{y} R_R
 \end{align*}
 and $\partial_{R_B} U, \partial_{R_R} U :  \operatorname{SO}(3) \to T_R^*\operatorname{SO}(3)$, and  $\partial_x U: \mathbb{R}^3 \to T_x^* \mathbb{R}^3$.  Similar to the sliding disk at uniform temperature, one can prove the following.

\begin{theorem}
Let $\chi$ denote the phase space of the magnetic motor.  Set $\xi$ to be the solution of  (\ref{eq:sdea})-(\ref{eq:sdeh}).  Suppose that $ c_R / \alpha_R^2 = c_B / \alpha_B^2$.  Let $E$ denote the energy of the magnetic motor given by,
\[
E(\mathbf{x}, \mathbf{v}, R_B, \boldsymbol{\pi}_B, R_R,
\boldsymbol{\pi}_R) =  \frac{m}{2}  \mathbf{v}^{\mathrm{T}}
\mathbf{v} + \frac{1}{2} J_B^{-1} \boldsymbol{\pi}_B^{\mathrm{T}}
\boldsymbol{\pi}_B + \frac{1}{2} \boldsymbol{\pi}_R^{\mathrm{T}} R_R
\mathbb{I}_R^{-1} R_R^{\mathrm{T}}  \boldsymbol{\pi}_R +
U(\mathbf{x}, R_B, R_R)    \text{.}
\]
Let $\beta = 2 c_B/ \alpha_B^2$.  The following measure,
\begin{equation}\label{kjhss}
\mu(d\xi):=\frac{e^{-\beta E}}{Z} d\xi
\end{equation}
is the unique Gibbs invariant measure of the stochastic process $\xi$.
\end{theorem}

\begin{Remark}
By using Ito's formula one can show that:
\[
d E = \left( - c_B \mathbf{V}_Q^{\mathrm{T}}  \mathbf{V}_Q  -c_R
\boldsymbol{\omega}_R^{\mathrm{T}} \boldsymbol{\omega}_R +
\alpha_B^2 \left( \frac{1}{m} + \frac{r^2}{J} \right) dt +
\frac{1}{2} \alpha_R^2 \operatorname{trace} \left[ \mathbb{I}_R^{-1}
\right]  \right) dt + \text{martingales} \text{.}
\]
Integrating yields,
\begin{align*}
E(t) - E(0) =& - c_B \int_0^t \mathbf{V}_Q^{\mathrm{T}} \mathbf{V}_Q
ds
- c_R \int_0^t \boldsymbol{\omega}_R^{\mathrm{T}} \boldsymbol{\omega}_R  ds \\
&+ \alpha_B^2 \left( \frac{1}{m} + \frac{r^2}{J} \right) t +
\frac{1}{2} \alpha_R^2 \operatorname{trace} \left[ \mathbb{I}_R^{-1}
\right] t  + \text{martingales} \text{.}
\end{align*}
\end{Remark}

\newpage

\bibliographystyle{plain}
\bibliography{nawaf}

 \end{document}